\documentclass[11pt, reqno]{amsart}

\usepackage{amsmath, amsthm, amssymb}
\usepackage{enumitem}
\usepackage{pdflscape}
\usepackage{caption}
\usepackage{bm}

\usepackage{ifpdf}
\ifpdf
\usepackage[pdftex]{graphicx}
\else
\usepackage[dvips]{graphicx}
\fi
\usepackage[all]{xy}

\usepackage{multicol}
\usepackage{mathtools}

\usepackage{tocvsec2}

\usepackage{bbm}

\input xy
\xyoption{all}

\usepackage[pdftex,plainpages=false,hypertexnames=false,pdfpagelabels]{hyperref}
\newcommand{\arxiv}[1]{\href{http://arxiv.org/abs/#1}{\tt arXiv:\nolinkurl{#1}}}
\newcommand{\arXiv}[1]{\href{http://arxiv.org/abs/#1}{\tt arXiv:\nolinkurl{#1}}}
\newcommand{\doi}[1]{\href{http://dx.doi.org/#1}{{\tt DOI:#1}}}

\newcommand{\googlebooks}[1]{(preview at \href{http://books.google.com/books?id=#1}{google books})}

\usepackage{xcolor}
\definecolor{dark-red}{rgb}{0.7,0.25,0.25}
\definecolor{dark-blue}{rgb}{0.15,0.15,0.55}
\definecolor{medium-blue}{rgb}{0,0,.8}
\definecolor{DarkGreen}{RGB}{0,150,0}
\definecolor{rho}{named}{red}
\hypersetup{
   colorlinks, linkcolor={purple},
   citecolor={medium-blue}, urlcolor={medium-blue}
}

\usepackage{longtable}
\usepackage{fullpage}

\theoremstyle{plain}
\newtheorem{thm}{Theorem}[section]
\newtheorem*{thm*}{Theorem}
\newtheorem{thmalpha}{Theorem}

\newtheorem*{cor*}{Corollary}

\newtheorem*{conj*}{Conjecture}
\newtheorem{lem}[thm]{Lemma}
\newtheorem{prop}[thm]{Proposition}

\newtheorem*{quest*}{Question}
\newtheorem*{claim*}{Claim}

\theoremstyle{definition}
\newtheorem{defn}[thm]{Definition}

\newtheorem{construction}[thm]{Construction}

\newtheorem{nota}[thm]{Notation}

\newtheorem{sub-ex}[thm]{Sub-Example}
\newtheorem{counter-ex}[thm]{Counter-Example}
\newtheorem{rem}[thm]{Remark}
\newtheorem*{rem*}{Remark}

\DeclareMathOperator{\core}{core}

\DeclareMathOperator{\End}{End}

\DeclareMathOperator{\id}{id}

\newcommand{\comment}[1]{}

\newcommand{\noshow}[1]{}
\renewcommand{\MR}[1]{}
\newcommand{\condense}{\mathrel{\,\hspace{.75ex}\joinrel\rhook\joinrel\hspace{-.75ex}\joinrel\rightarrow}}

\def\semicolon{;}
\def\applytolist#1{
    \expandafter\def\csname multi#1\endcsname##1{
        \def\multiack{##1}\ifx\multiack\semicolon
            \def\next{\relax}
        \else
            \csname #1\endcsname{##1}
            \def\next{\csname multi#1\endcsname}
        \fi
        \next}
    \csname multi#1\endcsname}

\def\calc#1{\expandafter\def\csname c#1\endcsname{{\mathcal #1}}}
\applytolist{calc}QWERTYUIOPLKJHGFDSAZXCVBNM;
\def\bbc#1{\expandafter\def\csname bb#1\endcsname{{\mathbb #1}}}
\applytolist{bbc}QWERTYUIOPLKJHGFDSAZXCVBNM;
\def\bfc#1{\expandafter\def\csname bf#1\endcsname{{\mathbf #1}}}
\applytolist{bfc}QWERTYUIOPLKJHGFDSAZXCVBNM;
\def\sfc#1{\expandafter\def\csname s#1\endcsname{{\sf #1}}}
\applytolist{sfc}QWERTYUIOPLKJHGFDSAZXCVBNM;
\def\fc#1{\expandafter\def\csname f#1\endcsname{{\mathfrak #1}}}
\applytolist{fc}QWERTYUIOPLKJHGFDSAZXCVBNM;

\newcommand{\Fun}{{\sf Fun}}

\newcommand{\Gray}{{\mathsf{Gray}}}

\newcommand{\Cat}{{\mathsf{Cat}}}

\newcommand{\QSys}{{\sf QSys}}

\usepackage{tikz}
\usepackage{tikz-cd}
\usetikzlibrary{arrows,backgrounds,patterns.meta}
\usetikzlibrary{positioning,shadings,cd}
\usetikzlibrary{shapes}
\usetikzlibrary{backgrounds}
\usetikzlibrary{decorations,decorations.pathreplacing,decorations.markings,decorations.pathmorphing}
\usetikzlibrary{fit,calc,through}
\usetikzlibrary{external}
\usetikzlibrary{arrows}
\tikzset{vertex/.style = {shape=circle,draw,fill=black,inner sep=0pt,minimum size=5pt}}
\tikzset{edge/.style = {->,> = latex', bend right}}
\tikzset{
	super thick/.style={line width=3pt}
}
\tikzset{
    quadruple/.style args={[#1] in [#2] in [#3] in [#4]}{
        #1,preaction={preaction={preaction={draw,#4},draw,#3}, draw,#2}
    }
}
\tikzstyle{shaded}=[fill=red!10!blue!20!gray!30!white]
\tikzstyle{unshaded}=[fill=white]
\tikzstyle{empty box}=[circle, draw, thick, fill=white, opaque, inner sep=2mm]
\tikzstyle{annular}=[scale=.7, inner sep=1mm, baseline]
\tikzstyle{rectangular}=[scale=.75, inner sep=1mm, baseline=-.1cm]
\tikzstyle{mid>}=[decoration={markings, mark=at position 0.5 with {\arrow{>}}}, postaction={decorate}]
\tikzstyle{mid<}=[decoration={markings, mark=at position 0.5 with {\arrow{<}}}, postaction={decorate}]
\tikzstyle{over}=[double, draw=white, super thick, double=]
\tikzstyle{snake}=[decorate, decoration={snake, segment length=1mm, amplitude=.3mm}]
\tikzstyle{saw}=[decorate, decoration={saw, segment length=.7mm, amplitude=.25mm}]

\tikzstyle{coupon}=[draw, very thick, rectangle, rounded corners=5pt]
\tikzset{Rightarrow/.style={double equal sign distance,>={Implies},->},
triplecd/.style={-,preaction={draw,Rightarrow}},
quadruplecd/.style={preaction={draw,Rightarrow,
shorten >=0pt
},
shorten >=1pt,
-,double,double
distance=0.2pt}}
\tikzset{
    tripleline/.style args={[#1] in [#2] in [#3]}{
        #1,preaction={preaction={draw,#3},draw,#2}
    }
}
\tikzstyle{triple}=[tripleline={[line width=.15mm,black] in
      [line width=.7mm,white] in
      [line width=1mm,black]}] 
\tikzset{
    quadrupleline/.style args={[#1] in [#2] in [#3] in [#4]}{
        #1,preaction={preaction={preaction={draw,#4},draw,#3}, draw,#2}
    }
}
\tikzstyle{quadruple}=[quadrupleline={[line width=.3mm,white] in
      [line width=.6mm,black] in
      [line width=1.2mm,white] in
      [line width=1.5mm,black]}]

\tikzdeclarepattern{
  name=primeddots,
  type=uncolored,
  bounding box={(-.6pt,-.6pt) and (.6pt,.6pt)},
  tile size={(5pt,5pt)},
  tile transformation={rotate=60},
  parameters={none}, 
  code={
    \fill(0pt,0pt) circle (.35pt);
  }
}

\tikzdeclarepattern{
  name=primedbox,
  type=uncolored,
  bounding box={(-1pt,-1pt) and (1pt,1pt)},
  tile size={(5pt,5pt)},
  tile transformation={rotate=60},
  parameters={none}, 
  code={
    \draw[very thin] (-.9pt,-.9pt) -- (-.9pt,.9pt) -- (.9pt,.9pt) -- (.9pt,-.9pt) -- (-.9pt,-.9pt);
  }
}

\tikzdeclarepattern{
  name=primedplus,
  type=uncolored,
  bounding box={(-1pt,-1pt) and (1pt,1pt)},
  tile size={(5pt,5pt)},
  tile transformation={rotate=30},
  parameters={none}, 
  code={
    \draw[very thin] (-.9pt,-.9pt) -- (.9pt,.9pt);
    \draw[very thin] (.9pt,-.9pt) -- (-.9pt,.9pt);
  }
}

\tikzdeclarepattern{
  name=primedstar,
  type=uncolored,
  bounding box={(-1.2pt,-1.2pt) and (1.2pt,1.2pt)},
  tile size={(5pt,5pt)},
  tile transformation={rotate=30},
  parameters={none}, 
  code={
    \draw[very thin] (-1.2pt,0pt) -- (1.2pt,0pt);
    \draw[very thin] (-.6pt,-1.04pt) -- (.6pt,1.04pt);
    \draw[very thin] (-.6pt,1.04pt) -- (.6pt,-1.04pt);
  }
}

\tikzdeclarepattern{
  name=primeddots2,
  type=uncolored,
  bounding box={(-.3pt,-.3pt) and (.3pt,.3pt)},
  tile size={(2.5pt,2.5pt)},
  tile transformation={rotate=60},
  parameters={none}, 
  code={
    \fill(0pt,0pt) circle (.35pt);
  }
}

\tikzstyle{primedregion}[none]=[
	preaction={fill=#1},
	pattern=primeddots,
  draw=#1,
]

\tikzstyle{boxregion}[none]=[
	preaction={fill=#1},
	pattern=primedbox,
  draw=#1,
]

\tikzstyle{plusregion}[none]=[
	preaction={fill=#1},
	pattern=primedplus,
  draw=#1,
]

\tikzstyle{starregion}[none]=[
	preaction={fill=#1},
	pattern=primedstar,
  draw=#1,
]

\tikzstyle{primedregion2}[none]=[
	preaction={fill=#1},
	pattern=primeddots2,
  draw=#1,
]

\newcommand{\roundNbox}[6]{
	\draw[rounded corners=5pt, very thick, #1] ($#2+(-#3,-#3)+(-#4,0)$) rectangle ($#2+(#3,#3)+(#5,0)$);
	\coordinate (ZZa) at ($#2+(-#4,0)$);
	\coordinate (ZZb) at ($#2+(#5,0)$);
	\node at ($1/2*(ZZa)+1/2*(ZZb)$) {#6};
}

\newcommand{\tikzmath}[2][]
     {\vcenter{\hbox{\begin{tikzpicture}[#1]#2
                     \end{tikzpicture}}}
     }

\newcommand{\XColor}{red} 
\newcommand{\YColor}{orange}
\newcommand{\ZColor}{violet}
\newcommand{\WColor}{blue}
\newcommand{\PsColor}{brown!80} 
\newcommand{\QsColor}{DarkGreen}
\newcommand{\RsColor}{cyan!90}

\newcommand{\AColor}{gray!30}
\newcommand{\BColor}{gray!55}
\newcommand{\CColor}{gray!75}
\newcommand{\DColor}{gray!95}
\newcommand{\PrColor}{yellow!50} 
\newcommand{\QrColor}{green!30}
\newcommand{\RrColor}{cyan!30}

\newcommand{\phiColor}{black}
\newcommand{\psiColor}{snake}
\newcommand{\gammaColor}{saw}

\newcommand{\xz}{\circ}
\newcommand{\xo}{\otimes}
\newcommand{\xt}{\star}

\newcommand{\xxo}{\otimes}
\newcommand{\xxt}{\star}

\tikzcdset{scale cd/.style={every label/.append style={scale=#1},
    cells={nodes={scale=#1}}}}

\begin{document}
\title{Q-system completion is a 3-functor}
\author{Quan Chen and David Penneys}
\date{\today}
\maketitle
\begin{abstract}
Q-systems are unitary versions of Frobenius algebra objects which appeared in the theory of subfactors.
In recent joint work with R.~Hern\'andez Palomares and C.~Jones, the authors defined a notion of Q-system completion for $\rm C^*/W^*$ 2-categories, which is a unitary version of a higher idempotent completion in the spirit of Douglas--Reutter and Gaiotto--Johnson-Freyd.
In this article, we prove that Q-system completion is a $\dag$ 3-functor on the $\dag$ 3-category of $\rm C^*/W^*$ 2-categories.
We also prove that Q-system completion satisfies a universal property analogous to the universal property satisfied by idempotent completion for 1-categories.
\end{abstract}

\section{Introduction}

Idempotent completions for higher categories has seen tremendous recent progress.
For 2-categories (which we always assume are locally idempotent complete), completing with respect to the two notions of condensation monads \cite{1905.09566} and separable monads \cite{1812.11933} produces equivalent 2-categories by \cite[Thm.~3.3.3]{1905.09566}.
The major difference is that condensation monads are \emph{non-unital} and include the \emph{data} of the separating structure, while separable monads are \emph{unital} and include only the \emph{existence} of separating structure, the choice of which is contractible.

In the setting of $\rm C^*/W^*$ 2-categories (which we always assume are locally orthogonal projection complete), the analogous notion of separable monad is Longo's \emph{Q-system} \cite{MR1257245,MR1444286}, which was originally studied for its role in subfactor theory.
In our recent joint article \cite{2105.12010}, we introduced the notion of \emph{Q-system completion} $\QSys(\cC)$ for a $\rm C^*/W^*$ 2-category $\cC$, which comes equipped with a canonical $\dag$ 2-functor $\iota_\cC : \cC \hookrightarrow \QSys(\cC)$.
While we analyzed some of the general theory of Q-system completion in that article, we focused more on applications to $\rm C^*$-algebra theory, showing the $\rm C^*$ 2-category of $\rm C^*$-algebras is Q-system complete.
As an application, we used Q-system completion to induce actions of unitary fusion categories on $\rm C^*$-algebras, similar to the spirit of \cite{2010.01072}.

In this article, we study some basic formal properties of Q-system completion, and our proofs can easily be adapted to the separable monad setting.
Our main results extend the treatment of idempotent completion for 2-categories in \cite[Appendix A]{1812.11933}.
Here is our first main theorem:

\begin{thmalpha}
\label{thm:QSys3Functor}
Q-system completion is a $\dag$ 3-endofunctor on the $\dag$ 3-category of $\rm C^*/W^*$ 2-categories.
\end{thmalpha}

In \cite[Prop.~A.6.3]{1812.11933}, Douglas and Reutter provided strong evidence towards this theorem, and they mentioned they expect such a result to be true. 
To prove this theorem, we introduce an \emph{overlay} compatibility between the 2D graphical calculi for a $\rm C^*/W^*$ 2-category $\cC$ and the $\rm C^*/W^*$ 2-category $\Fun^\dag(\cC\to \cD)$ for another $\cD$.
(We show in Proposition \ref{prop:FunDagC*W*} below that $\Fun^\dag(\cC\to \cD)$ is $\rm C^*/W^*$ whenever $\cC,\cD$ are.)
See \S\ref{sec:Fun(A,B)} below for more details.

Our second main theorem regards the universal property for idempotent completion for 2-categories discussed in \cite[\S1.2]{2012.15774}, proving the best possible uniqueness statement.
Given 2-categories $\cC,\cD,\cE$ and 2-functors $F: \cC \to \cD$ and $G: \cC\to \cE$,
the 2-category of \emph{lifts} of $F$ to $\cE$ along $G$ is the homotopy fiber at $F$ of the functor
$$
-\xz G : \Fun(\cE \to \cD) \to \Fun(\cC\to \cD).
$$
Objects in this lift 2-category are pairs $(\widetilde{F}, \theta)$ where $\widetilde{F}:\cE\to \cD$ is a 2-functor and $\theta: F \Rightarrow \widetilde{F}\xz G$ is an invertible 2-transformation. 
We refer the reader to \S\ref{sec:UniversalProperty} for the rest of the unpacked definition.

\begin{thmalpha}
\label{thm:UniqueLift}
Suppose $\cC$ is a $\rm C^*/W^*$ 2-category.
The Q-system completion $\QSys(\cC)$ satisfies the following universal property.
For any $\dag$ 2-functor $F:\cC\to \cD$ where $\cD$ is Q-system complete, 
the 2-category of lifts of $F$ along $\iota_\cC$ is \emph{$(-2)$-truncated}, i.e., equivalent to a point.
That is, $-\circ \iota_\cC: \Fun^\dag(\QSys(\cC) \to \cD) \to \Fun^\dag(\cC\to \cD)$ is a $\dag$ 2-equivalence.
\end{thmalpha}

The main idea of the proof of this theorem comes from \cite[\S3.1]{1910.03178}.
By a version of Grothendieck's \emph{Homotopy Hypothesis} for 2-categories \cite{MR1239560}, the homotopy category of strict 2-groupoids and strict 2-functors localized at the strict equivalences is equivalent to the 1-category of homotopy 2-types.
Hence the homotopy fiber of $-\xz G$ restricted to the \emph{core} 2-groupoids
$$
-\xz G :
\core(\Fun(\cE \to \cD))
\to
\core(\Fun(\cC \to \cD))
$$
is $k$-truncated for $-2\leq k\leq 1$ if and only if various (essential) surjectivity properties hold for $-\xz G$.
In turn, these surjectivity properties for $-\xz G$ are ensured by various levels of \emph{dominance} for the 2-functor $G$.
We make these notions precise in \S\ref{sec:UniversalProperty}.

While we work in the $\rm C^*/W^*$ setting both for novelty and for applications to the world of operator algebras, we re-emphasize that these results do not depend on the dagger structure.

\subsection*{Acknowledgements}
The authors would like to thank Thibault D\'ecoppet, Brent Nelson, and David Reutter for helpful discussions.
The authors were supported by 
NSF grants DMS 1654159, 1927098, and 2051170.

\section{Preliminaries}

In this article, \emph{2-category} will always mean a weak 2-category/bicategory which is locally idempotent complete,
and a $\rm C^*/W^*$ 2-category will always mean a weak $\rm C^*/W^*$ 2-category which is locally orthogonal projection complete.
We refer the reader to \cite{2002.06055} for background on 2-categories and to \cite{2105.12010} for background on $\rm C^*/W^*$ 2-categories.
We refer the reader to \cite{MR3971584} or \cite{2105.12010} for a detailed discussion of the graphical calculus of string diagrams for 2-categories.
The only 3-categories in this article are the 3-category $2\Cat$ of 2-categories \cite[\S5.1]{MR3076451} and its 3-subcategories $\rm C^*2\Cat$ and $\rm W^*2\Cat$ of $\rm C^*/W^*$ 2-categories respectively.

\begin{nota}
In a 2-category $\cC$, we refer to its objects, 1-morphisms, and 2-morphisms as \emph{0-cells}, \emph{1-cells}, and \emph{2-cells} respectively.
We denote 0-cells in a 2-category $\cC$ by lowercase Roman letters $a,b,c$, 1-cells by uppercase Roman letters ${}_aX_b,{}_bY_c$ using bimodule notation for source (left) and target (right), and 2-cells by lowercase Roman letters later in the alphabet $f,m,n,t$.
We write 1-composition as $\xxo$ read \emph{left to right}, and we write 2-composition as $\xxt$, 
which is read \emph{right to left}.
In the graphical calculus of string diagrams in 2-categories, which is formally dual to the manipulation of pasting diagrams, we read 1-composition \emph{left to right} and 2-composition \emph{bottom to top}.
$$
f: {}_aX\xxo_b Y_c \Rightarrow {}_aZ_c
\qquad\rightsquigarrow\qquad
\begin{tikzcd}[row sep=0]
&\mbox{}
\\
a
\arrow[rr, bend left = 30, "Z"]
\arrow[dr, swap, "X"]
&
\mbox{}
&
c
\\
&b
\arrow[uu, Rightarrow, swap, "f"]
\arrow[ur, swap, "Y"]
\end{tikzcd}
\qquad\rightsquigarrow\qquad
\tikzmath{
\begin{scope}
\clip[rounded corners=5pt] (-.7,-.7) rectangle (.7,.7);
\filldraw[\AColor] (-.2,-.7) -- (-.2,0) -- (0,0) -- (0,.7) -- (-.7,.7) -- (-.7,-.7);
\filldraw[\BColor] (-.2,-.7) rectangle (.2,0);
\filldraw[\CColor] (.2,-.7) -- (.2,0) -- (0,0) -- (0,.7) -- (.7,.7) -- (.7,-.7);
\end{scope}
\node at (0,-.9) {$\to$};
\node at (-.9,0) {$\Uparrow$};
\draw[\ZColor,thick] (0,.3) -- node[left]{$\scriptstyle Z$}(0,.7);
\draw[\XColor,thick] (-.2,-.3) -- node[left]{$\scriptstyle X$} (-.2,.-.7);
\draw[\YColor,thick] (.2,-.3) -- node[right]{$\scriptstyle Y$} (.2,-.7);
\roundNbox{unshaded}{(0,0)}{.3}{.1}{.1}{\scriptsize{$f$}};
}
$$

In the 3-category $2\Cat$ of 2-categories, the object 2-categories are denoted by math calligraphic letters $\cC,\cD,\cE$,
the 2-functor 1-morphisms are denoted by capital Roman letters $F,G,H$,
the 2-transformation 2-morphisms are denoted by lowercase Greek letters $\varphi,\psi$,
and 2-modification 3-morphisms are denoted by lowercase Roman letters $m,n$.
We write 1-composition of 2-functors as $\xz$, which we read \emph{right to left}, i.e., if $F: \cA \to \cB$ and $G: \cB \to \cC$, then $G\xz F : \cA \to \cC$.
We write 2-composition of 2-transformations as $\xo$, and we write 3-composition of 2-modification as $\xt$.
\end{nota}

\begin{rem}
While we will not rely on any 3D string diagram graphical calculus in this article, 
its use for weak 3-categories can be justified using the article \cite{1903.05777}.
In several locations, we provide 3D diagrams for conceptual clarity.
Our conventions for 1-, 2-, and 3-composition in these 3D diagrams are indicated in the figure below.
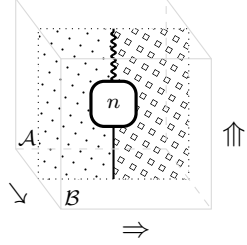
\[
\tikzmath{
\begin{scope}
\filldraw[primedregion=white] (-.3,.4) rectangle (.7,2.4);
\filldraw[boxregion=white] (.7,.4) rectangle (1.7,2.4);
\end{scope}
\draw[\phiColor,thick] (.7,.4) -- (.7,1.4);
\draw[\psiColor,thick] (.7,1.4) -- (.7,2.4);
\roundNbox{unshaded}{(.7,1.4)}{.3}{0}{0}{\scriptsize{$n$}};
\draw[dotted] (-.3,.4) rectangle (1.7,2.4); 
\node at (-.6,.2) [rotate=-53] {$\to$};
\node at (1,-.3) {$\Rightarrow$};
\node at (2.3,1) [rotate=90] {$\Rrightarrow$};
\draw[\AColor] (0,0) rectangle (2,2);
\draw[\AColor] (0,0) -- (-.6,.8) -- (-.6,2.8) -- (1.4,2.8) -- (2,2);
\draw[\AColor] (-.6,2.8) -- (0,2);
\draw[\AColor,dashed] (-.6,.8) -- (1.4,.8) -- (1.4,2.8);
\draw[\AColor,dashed] (1.4,.8) -- (2,0);
\node at (-.45,.95) {\scriptsize{$\cA$}};
\node at (.15,.15) {\scriptsize{$\cB$}};
}
\]
\end{rem}

\subsection{The 2-category \texorpdfstring{$\Fun(\cA\to \cB)$}{Fun(C->D)} of 2-functors, 2-transformations, and 2-modifications}
\label{sec:Fun(A,B)}

In this section, we first describe our graphical conventions for working with 2-functors, 2-transformations, and 2-modifications.
We then use our graphical notation to unpack their definitions.

\begin{nota}
To define 2-transformations between 2-functors and 2-modifications between 2-transformations in a diagrammatic language, we overlay the 2D diagrammatic calculus for the hom 2-category $\Fun(\cA\to \cB)$ between 2-categories $\cA,\cB$ with the 2D diagrammatic calculus for $\cB$.

For our 2D diagrammatic calculus for the hom 2-category $\Fun(\cA\to \cB)$, we represent the object functors by \emph{unshaded} regions with \emph{textured} decorations, e.g.,
$$
\tikzmath{
\filldraw[primedregion=white, rounded corners = 5pt] (0,0) rectangle (.6,.6);
\draw[thin, dotted, rounded corners = 5pt] (0,0) rectangle (.6,.6);
}
=
F
\qquad\qquad
\tikzmath{
\filldraw[boxregion=white, rounded corners = 5pt] (0,0) rectangle (.6,.6);
\draw[thin, dotted, rounded corners = 5pt] (0,0) rectangle (.6,.6);
}
=
F'
\qquad\qquad
\tikzmath{
\filldraw[plusregion=white, rounded corners = 5pt] (0,0) rectangle (.6,.6);
\draw[thin, dotted, rounded corners = 5pt] (0,0) rectangle (.6,.6);
}
=
F''
\qquad\qquad
\tikzmath{
\filldraw[starregion=white, rounded corners = 5pt] (0,0) rectangle (.6,.6);
\draw[thin, dotted, rounded corners = 5pt] (0,0) rectangle (.6,.6);
}
=
F'''
$$
We represent 2-transformations (see Definition \ref{defn:2Transformation} below) by \emph{textured} strings between these textured regions, e.g.,
$$
\tikzmath{
\begin{scope}
\filldraw[primedregion=white, rounded corners = 5pt] (0,0) rectangle (.3,.6);
\filldraw[boxregion=white, rounded corners = 5pt] (.3,0) rectangle (.6,.6);
\end{scope}
\draw[\phiColor,thick] (.3,0) -- (.3,.6);
\draw[thin, dotted, rounded corners = 5pt] (0,0) rectangle (.6,.6);
}
=
\varphi: F\Rightarrow F'
\qquad\qquad
\tikzmath{
\begin{scope}
\filldraw[boxregion=white, rounded corners = 5pt] (0,0) rectangle (.3,.6);
\filldraw[plusregion=white, rounded corners = 5pt] (.3,0) rectangle (.6,.6);
\end{scope}
\draw[\psiColor,thick] (.3,0) -- (.3,.6);
\draw[thin, dotted, rounded corners = 5pt] (0,0) rectangle (.6,.6);
}
=
\psi: F'\Rightarrow F''
\qquad\qquad
\tikzmath{
\begin{scope}
\filldraw[plusregion=white, rounded corners = 5pt] (0,0) rectangle (.3,.6);
\filldraw[starregion=white, rounded corners = 5pt] (.3,0) rectangle (.6,.6);
\end{scope}
\draw[\gammaColor,thick] (.3,0) -- (.3,.6);
\draw[thin, dotted, rounded corners = 5pt] (0,0) rectangle (.6,.6);
}
=
\gamma: F''\Rightarrow F'''
$$
We represent 2-modifications (see Definition \ref{Defn:2Modification} below) by coupons as usual.

To depict a 2-morphism in $\cB$ in the image of $F$, we \emph{overlay} the 2D string diagrammatic calculus for $\Fun(\cA\to \cB)$ on top of the 2D string diagrammatic calculus for $\cA$.
For example, given $F,F': \cA \to \cB$, $\varphi,\varphi' : F\Rightarrow F'$, and $m: \varphi\Rrightarrow \varphi'$, we can `overlay' the coupon for $m$ over the shaded region for $a\in \cA$ to obtain the 2-cell $m_a:\varphi_a\Rightarrow \varphi'_a$:
$$
\left(\,
\tikzmath[scale=.75, transform shape]{
\begin{scope}
\clip[rounded corners = 5] (-.7,-.7) rectangle (.7,.7);
\filldraw[primedregion=white] (-.8,-.9) rectangle (0,.9);
\filldraw[boxregion=white] (0,-.9) rectangle (.8,.9);
\end{scope}
\draw[\phiColor,thick] (0,-.7) -- (0,0);
\draw[\phiColor,thick] (0,.7) -- (0,0);
\roundNbox{unshaded}{(0,0)}{.3}{0}{0}{\scriptsize{$m$}}; 
\node at (0,-.9) {\scriptsize{$\varphi$}};
\node at (0,.9) {\scriptsize{$\varphi'$}};
\draw[thin, dotted, rounded corners = 5pt] (-.7,-.7) rectangle (.7,.7);
}
\,\right)
\left(
\tikzmath{
\fill[\AColor,rounded corners=5pt] (0,0) rectangle (.6,.6);
\node at (.3,.8) {\phantom{\scriptsize{$a$}}};
\node at (.3,-.2) {\scriptsize{$a$}};
}
\right)
=
\tikzmath[scale=.75, transform shape]{
\begin{scope}
\clip[rounded corners = 5] (-.7,-.7) rectangle (.7,.7);
\filldraw[primedregion=\AColor] (-.8,-.9) rectangle (0,.9);
\filldraw[boxregion=\AColor] (0,-.9) rectangle (.8,.9);
\end{scope}
\draw[\phiColor,thick] (0,-.7) -- (0,0);
\draw[\phiColor,thick] (0,.7) -- (0,0);
\roundNbox{unshaded}{(0,0)}{.3}{0}{0}{\scriptsize{$m_a$}}; 
\node at (0,-.9) {\scriptsize{$\varphi_{a}$}};
\node at (0,.9) {\scriptsize{$\varphi'_a$}};
}
\qquad\qquad\qquad
\tikzmath{
\filldraw[primedregion=\AColor, rounded corners = 5pt] (0,0) rectangle (.6,.6);
}
=F(a),
\quad
\tikzmath{
\filldraw[boxregion=\AColor, rounded corners = 5pt] (0,0) rectangle (.6,.6);
}
=F'(a)
\,.
$$
We do not attempt to formalize this `overlay' operation, as all string diagrams can be interpreted uniquely as 2-cells in $\cB$; see Remark \ref{rem:Overlay} below for further discussion.
\end{nota}

\begin{defn}
\label{defn:2Functor}
Suppose $\cA,\cB$ are 2-categories.
We use the following conventions for the coheretors of a 2-functor
$F=(F,F^2,F^1):\cA \to \cB$:
$$
F^2_{X,Y}\in \cB( F(X)\xxo_{F(b)} F(Y) \Rightarrow F(X\xxo_b Y))
\qquad\text{and}\qquad
F^1_a \in \cB(1_{F(a)} \Rightarrow F(1_a)),
$$
which satisfy the hexagon associativity equation and triangle unit equations.
We depict these axioms below in the graphical calculus for $\cB$.
Denoting objects in $\cB$ by the shaded regions
$$
\tikzmath{
\filldraw[primedregion=\AColor, rounded corners = 5pt] (0,0) rectangle (.6,.6);
\draw[rounded corners=5, thin, dotted] (0,0) rectangle (.6,.6);
}
=F(a)
\qquad\qquad
\tikzmath{
\filldraw[primedregion=\BColor, rounded corners = 5pt] (0,0) rectangle (.6,.6);
\draw[rounded corners=5, thin, dotted] (0,0) rectangle (.6,.6);
}
=F(b)
\qquad\qquad
\tikzmath{
\filldraw[primedregion=\CColor, rounded corners = 5pt] (0,0) rectangle (.6,.6);
\draw[rounded corners=5, thin, dotted] (0,0) rectangle (.6,.6);
}
=F(c)
\qquad\qquad
\tikzmath{
\filldraw[primedregion=\DColor, rounded corners = 5pt] (0,0) rectangle (.6,.6);
\draw[rounded corners=5, thin, dotted] (0,0) rectangle (.6,.6);
}
=F(d),
$$
and 1-cells in $\cB$ by shaded strands, e.g.
\[
\tikzmath{
\begin{scope}
\clip[rounded corners=5pt] (-.3,0) rectangle (.3,.6);
\filldraw[\AColor] (0,0) rectangle (-.3,.6);
\filldraw[\BColor] (0,0) rectangle (.3,.6);
\end{scope}
\draw[thick, \XColor] (0,0) -- (0,.6);
}={}_aX_b
\qquad
\tikzmath{
\begin{scope}
\clip[rounded corners=5pt] (-.3,0) rectangle (.3,.6);
\filldraw[\BColor] (0,0) rectangle (-.3,.6);
\filldraw[\CColor] (0,0) rectangle (.3,.6);
\end{scope}
\draw[thick, \YColor] (0,0) -- (0,.6);
}={}_bY_c
\qquad
\tikzmath{
\begin{scope}
\clip[rounded corners=5pt] (-.3,0) rectangle (.3,.6);
\filldraw[\CColor] (0,0) rectangle (-.3,.6);
\filldraw[\DColor] (0,0) rectangle (.3,.6);
\end{scope}
\draw[thick, \ZColor] (0,0) -- (0,.6);
}={}_cZ_d
\qquad
\tikzmath{
\begin{scope}
\clip[rounded corners=5pt] (-.5,-.3) rectangle (.5,.3);
\filldraw[primedregion=\AColor] (-.5,-.3) rectangle (-.2,.3);
\filldraw[primedregion=\BColor] (-.2,-.3) rectangle (.2,.3);
\filldraw[primedregion=\CColor] (.5,-.3) rectangle (.2,.3);
\end{scope}
\draw[\XColor,thick] (-.2,-.3) -- (-.2,.3);
\draw[\YColor,thick] (.2,-.3) -- (.2,.3);
}
=
F(X)\xxo_{F(b)} F(Y)
\qquad
\tikzmath{
\begin{scope}
\clip[rounded corners=5pt] (-.35,-.3) rectangle (.35,.3);
\filldraw[primedregion=\AColor] (-.05,-.3) rectangle (-.35,.3);
\filldraw[primedregion=\BColor] (-.05,-.3) rectangle (.05,.3);
\filldraw[primedregion=\CColor] (.05,-.3) rectangle (.35,.3);
\end{scope}
\draw[\XColor,thick] (-.05,-.3) -- (-.05,.3);
\draw[\YColor,thick] (.05,-.3) -- (.05,.3);
}
=
F(X\xxo_b Y),
\]
the hexagon and triangle equations are given by
\[
\underset{
\tiny (F(X)\xxo_{F(b)} F(Y))\xxo_{F(c)}F(Z)\Rightarrow F(X\xxo_b(Y\xxo_c Z))
}
{
\tikzmath{
\begin{scope}
\clip[rounded corners=5pt] (-.8,.3) rectangle (.8,3.7);
\filldraw[primedregion=\AColor] (-.4,0) -- (-.4,1) -- (-.25,1) -- (-.25,2) -- (-.1,2) -- (-.1,4) -- (-.8,4) -- (-.8,0);
\filldraw[primedregion=\BColor] (0,0) -- (0,1) -- (-.15,1) -- (-.15,2) -- (0,2) -- (0,4) -- (-.1,4) -- (-.1,2) -- (-.25,2) -- (-.25,1) -- (-.4,1) -- (-.4,0);
\filldraw[primedregion=\CColor] (.4,0) -- (.4,2) -- (.1,2) -- (.1,4) -- (0,4) -- (0,2) -- (-.15,2) -- (-.15,1) -- (0,1) -- (0,0);
\filldraw[primedregion=\DColor] (.4,0) -- (.4,2) -- (.1,2) -- (.1,4) -- (.8,4) -- (.8,0);
\end{scope}
\draw[\XColor,thick] (-.4,.3) -- (-.4,1);
\draw[\YColor,thick] (0,.3) -- (0,1);
\draw[\ZColor,thick] (.4,.3) -- (.4,2);
\draw[\XColor,thick] (-.25,1) -- (-.25,2);
\draw[\YColor,thick] (-.15,1) -- (-.15,2);
\draw[\XColor,thick] (-.1,2) -- (-.1,3.7);
\draw[\YColor,thick] (0,2) -- (0,3.7);
\draw[\ZColor,thick] (.1,2) -- (.1,3.7);
\roundNbox{unshaded}{(-.2,1)}{.3}{.05}{.05}{\scriptsize{$F^2_{X,Y}$}};
\roundNbox{unshaded}{(0,2)}{.3}{.25}{.25}{\scriptsize{$F^2_{X\xxo Y,Z}$}};
\roundNbox{unshaded}{(0,3)}{.3}{.25}{.25}{\scriptsize{$F(\alpha^{\cC})$}};
}
=
\tikzmath{
\begin{scope}
\clip[rounded corners=5pt] (-.8,.3) rectangle (.8,3.7);
\filldraw[primedregion=\AColor] (-.4,0) -- (-.4,3) -- (-.1,3) -- (-.1,4) -- (-.8,4) -- (-.8,0);
\filldraw[primedregion=\BColor] (-.4,0) -- (-.4,3) -- (-.1,3) -- (-.1,4) -- (0,4) -- (0,3) -- (.15,3) -- (.15,2) -- (0,2) -- (0,0);
\filldraw[primedregion=\CColor] (.4,0) -- (.4,2) -- (.25,2) -- (.25,3) -- (.1,3) -- (.1,4) -- (0,4) -- (0,3) -- (.15,3) -- (.15,2) -- (0,2) -- (0,0);
\filldraw[primedregion=\DColor] (.4,0) -- (.4,2) -- (.25,2) -- (.25,3) -- (.1,3) -- (.1,4) -- (.8,4) -- (.8,0);
\end{scope}
\draw[\XColor,thick] (-.4,.3) -- (-.4,3);
\draw[\YColor,thick] (0,.3) -- (0,2);
\draw[\ZColor,thick] (.4,.3) -- (.4,2);
\draw[\YColor,thick] (.15,2) -- (.15,3);
\draw[\ZColor,thick] (.25,2) -- (.25,3);
\draw[\XColor,thick] (-.1,3) -- (-.1,3.7);
\draw[\YColor,thick] (0,3) -- (0,3.7);
\draw[\ZColor,thick] (.1,3) -- (.1,3.7);
\roundNbox{unshaded}{(0,1)}{.3}{.25}{.25}{\scriptsize{$\alpha^\cB$}};
\roundNbox{unshaded}{(.2,2)}{.3}{.05}{.05}{\scriptsize{$F^2_{Y,Z}$}};
\roundNbox{unshaded}{(-.1,3)}{.3}{.15}{.35}{\scriptsize{$F^2_{X,Y\xxo Z}$}};
}
}
\qquad\qquad
\underset{
\tiny F(X)\xxo_{F(b)}1_{F(b)}\Rightarrow F(X)
}{
\tikzmath{
\begin{scope}
\clip[rounded corners=5pt] (-.7,.3) rectangle (.7,3.7);
\filldraw[primedregion=\AColor] (-.2,0) -- (-.2,2) -- (-.05,2) -- (-.05,3) -- (0,3) -- (0,4) -- (-.7,4) -- (-.7,0);
\filldraw[primedregion=\BColor] (-.2,0) -- (-.2,2) -- (-.05,2) -- (-.05,3) -- (0,3) -- (0,4) -- (.7,4) -- (.7,0);
\end{scope}
\draw[\XColor,thick] (-.2,.3) -- (-.2,2);
\draw[dotted,thick] (.2,.3) -- (.2,2);
\draw[\XColor,thick] (-.05,2) -- (-.05,3);
\draw[dotted,thick] (.05,2) -- (.05,3);
\draw[\XColor,thick] (0,3) -- (0,3.7);
\roundNbox{unshaded}{(.2,1)}{.25}{0}{0}{\scriptsize{$F^1_b$}};
\roundNbox{unshaded}{(0,2)}{.3}{.15}{.15}{\scriptsize{$F^2_{X,1_b}$}};
\roundNbox{unshaded}{(0,3)}{.3}{.15}{.15}{\scriptsize{$F(\rho_X^b)$}};
}
=
\tikzmath{
\begin{scope}
\clip[rounded corners=5pt] (-.7,0) rectangle (.7,2);
\filldraw[primedregion=\AColor] (-.2,0) -- (-.2,1) -- (0,1) -- (0,2) -- (-.7,2) -- (-.7,0);
\filldraw[primedregion=\BColor] (-.2,0) -- (-.2,1) -- (0,1) -- (0,2) -- (.7,2) -- (.7,0);
\end{scope}
\draw[\XColor,thick] (-.2,0) -- (-.2,1);
\draw[dotted,thick] (.2,0) -- (.2,1);
\draw[\XColor,thick] (0,1) -- (0,2);
\roundNbox{unshaded}{(0,1)}{.3}{.15}{.15}{\scriptsize{$\rho_{F(X)}^{F(b)}$}};
}
}
\qquad\qquad
\underset{
\tiny 1_{F(a)}\xxo_{F(a)}F(X)\Rightarrow F(X)
}{
\tikzmath{
\begin{scope}
\clip[rounded corners=5pt] (-.7,.3) rectangle (.7,3.7);
\filldraw[primedregion=\AColor] (.2,0) -- (.2,2) -- (.05,2) -- (.05,3) -- (0,3) -- (0,4) -- (-.7,4) -- (-.7,0);
\filldraw[primedregion=\BColor] (.2,0) -- (.2,2) -- (.05,2) -- (.05,3) -- (0,3) -- (0,4) -- (.7,4) -- (.7,0);
\end{scope}
\draw[\XColor,thick] (.2,.3) -- (.2,2);
\draw[dotted,thick] (-.2,.3) -- (-.2,2);
\draw[\XColor,thick] (.05,2) -- (.05,3);
\draw[dotted,thick] (-.05,2) -- (-.05,3);
\draw[\XColor,thick] (0,3) -- (0,3.7);
\roundNbox{unshaded}{(-.2,1)}{.25}{0}{0}{\scriptsize{$F^1_a$}};
\roundNbox{unshaded}{(0,2)}{.3}{.15}{.15}{\scriptsize{$F^2_{1_a,X}$}};
\roundNbox{unshaded}{(0,3)}{.3}{.15}{.15}{\scriptsize{$F(\lambda_X^a)$}};
}
=
\tikzmath{
\begin{scope}
\clip[rounded corners=5pt] (-.7,0) rectangle (.7,2);
\filldraw[primedregion=\AColor] (.2,0) -- (.2,1) -- (0,1) -- (0,2) -- (-.7,2) -- (-.7,0);
\filldraw[primedregion=\BColor] (.2,0) -- (.2,1) -- (0,1) -- (0,2) -- (.7,2) -- (.7,0);
\end{scope}
\draw[\XColor,thick] (.2,0) -- (.2,1);
\draw[dotted,thick] (-.2,0) -- (-.2,1);
\draw[\XColor,thick] (0,1) -- (0,2);
\roundNbox{unshaded}{(0,1)}{.3}{.15}{.15}{\scriptsize{$\lambda_{F(X)}^{F(a)}$}};
}
}\,.
\]
Whenever possible, we will suppress the associator and unitor coheretors in our 2-categories.
\end{defn}

\begin{defn}
\label{defn:2Transformation}
Suppose $\cA,\cB$ are 2-categories, $F,F': \cA \to \cB$ are 2-functors.
A \emph{2-transformation} $\varphi:F\Rightarrow F'$ consists of:
\begin{itemize}
\item 
for every 0-cell $c\in \cA$,
a 1-cell $\varphi_c \in \cB(F(c)\to F'(c))$, and
\item
for every 1-cell ${}_aX_b\in \cA(a\to b)$,
an invertible $F(a)-F'(b)$ bimodular 2-cell 
$$
\tikzmath[scale=.7, transform shape]{
\begin{scope}
\clip[rounded corners = 5] (-.6,0) rectangle (1.8,2.4);
\filldraw[primedregion=\AColor] (0,0) -- (0,.6) .. controls ++(90:.4cm) and ++(-135:.2cm) .. (.6,1.2) .. controls ++(135:.2cm) and ++(270:.4cm) .. (0,1.8) -- (0,3) -- (-.6,3) -- (-.6,0); 
\filldraw[primedregion=\BColor] (1.2,0) -- (1.2,.6) .. controls ++(90:.4cm) and ++(-45:.2cm) .. (.6,1.2) .. controls ++(-135:.2cm) and ++(90:.4cm) .. (0,.6) -- (0,0);
\filldraw[boxregion=\AColor] (0,3) -- (0,1.8) .. controls ++(270:.4cm) and ++(135:.2cm) .. (.6,1.2) .. controls ++(45:.2cm) and ++(270:.4cm) .. (1.2,1.8) -- (1.2,3);
\filldraw[boxregion=\BColor] (1.2,0) -- (1.2,.6) .. controls ++(90:.4cm) and ++(-45:.2cm) .. (.6,1.2) .. controls ++(45:.2cm) and ++(270:.4cm) .. (1.2,1.8) -- (1.2,3) -- (1.8,3) -- (1.8,0);
\end{scope}
\draw[\XColor,thick] (0,0) -- (0,.6) .. controls ++(90:.6cm) and ++(270:.6cm) .. (1.2,1.8) -- (1.2,2.4);
\draw[\phiColor,thick] (1.2,0) -- (1.2,.6) .. controls ++(90:.6cm) and ++(270:.6cm) .. (0,1.8) -- (0,2.4);
\filldraw[white] (.6,1.2) circle (.1cm);
\draw[thick] (.6,1.2) circle (.1cm); 
\node at (0,-.2) {\scriptsize{$F(X)$}};
\node at (1.2,2.6) {\scriptsize{$F'(X)$}};
%
\node at (1.2,-.2) {\scriptsize{$\varphi_b$}};
\node at (0,2.6) {\scriptsize{$\varphi_a$}};
} 
=\varphi_X 
\in 
\cB(F(X)\xxo_{F(b)} \varphi_b \Rightarrow \varphi_a \xxo_{F'(a)} F'(X)).
$$
\end{itemize}
This data satisfies the following coherence properties:
$$
\tikzmath[scale=.75, transform shape]{
\begin{scope}
\clip[rounded corners = 5] (-.6,0) rectangle (3,4.2);
\filldraw[primedregion=\AColor] (0,0) -- (0,2) .. controls ++(90:.4cm) and ++(-135:.2cm) .. (.6,2.6) .. controls ++(135:.2cm) and ++(270:.4cm) .. (0,3.2) -- (0,4.2) -- (-.6,4.2) -- (-.6,0); 
\filldraw[primedregion=\BColor] (1.2,0) -- (1.2,.4) .. controls ++(90:.4cm) and ++(-135:.2cm) .. (1.8,1) .. controls ++(135:.2cm) and ++(270:.4cm) .. (1.2,1.6) -- (1.2,2) .. controls ++(90:.4cm) and ++(-45:.2cm) .. (.6,2.6) .. controls ++(-135:.2cm) and ++(90:.4cm) .. (0,2) -- (0,0);
\filldraw[primedregion=\CColor] (2.4,0) -- (2.4,.4) .. controls ++(90:.4cm) and ++(-45:.2cm) .. (1.8,1) .. controls ++(-135:.2cm) and ++(90:.4cm) .. (1.2,.4) -- (1.2,0);
\filldraw[boxregion=\AColor] (0,4.2) -- (0,3.2) .. controls ++(270:.4cm) and ++(135:.2cm) .. (.6,2.6) .. controls ++(45:.2cm) and ++(270:.4cm) .. (1.2,3.2) -- (1.75,3.6) -- (1.75,4.2);
\filldraw[boxregion=\BColor] (1.85,4.2) -- (1.85,3.6) -- (2.4,3.2) -- (2.4,1.6) .. controls ++(270:.4cm) and ++(45:.2cm) .. (1.8,1) .. controls ++(135:.2cm) and ++(270:.4cm) .. (1.2,1.6) -- (1.2,2) .. controls ++(90:.4cm) and ++(-45:.2cm) .. (.6,2.6) .. controls ++(45:.2cm) and ++(270:.4cm) .. (1.2,3.2) -- (1.75,3.6) -- (1.75,4.2);
\filldraw[boxregion=\CColor] (2.4,0) -- (2.4,.4) .. controls ++(90:.4cm) and ++(-45:.2cm) .. (1.8,1) .. controls ++(45:.2cm) and ++(270:.4cm) .. (2.4,1.6) -- (2.4,3.2) -- (1.85,3.6) -- (1.85,4.2) -- (3,4.2) -- (3,0);
\filldraw[\BColor] (1,1.8) circle (.2cm);
\end{scope}
\draw[\XColor,thick] (0,0) -- (0,.4) -- (0,2) .. controls ++(90:.6cm) and ++(270:.6cm) .. (1.2,3.2) -- (1.2,3.6);
\draw[\YColor,thick] (1.2,0) -- (1.2,.4) .. controls ++(90:.6cm) and ++(270:.6cm) .. (2.4,1.6) -- (2.4,3.6);
\draw[\phiColor,thick] (2.4,0) -- (2.4,.4) .. controls ++(90:.6cm) and ++(270:.6cm) .. (1.2,1.6) -- (1.2,2) .. controls ++(90:.6cm) and ++(270:.6cm) .. (0,3.2) -- (0,4.2);
\draw[\XColor,thick] (1.733,3.6) -- (1.733,4.2);
\draw[\YColor,thick] (1.867,3.6) -- (1.867,4.2);
\roundNbox{unshaded}{(1.8,3.4)}{.3}{.6}{.6}{\scriptsize{${F'}^2_{X,Y}$}}; 
\filldraw[white] (1.8,1) circle (.1cm);
\draw[thick] (1.8,1) circle (.1cm); 
\filldraw[white] (.6,2.6) circle (.1cm);
\draw[thick] (.6,2.6) circle (.1cm); 
\node at (0,-.2) {\scriptsize{$F(X)$}};
\node at (1.2,-.2) {\scriptsize{$F(Y)$}};
\node at (1.8,4.4) {\scriptsize{$F'(X\xxo_b Y)$}};
\node at (2.4,-.2) {\scriptsize{$\varphi_c$}};
\node at (1,1.8) {\scriptsize{$\varphi_b$}};
\node at (0,4.4) {\scriptsize{$\varphi_a$}};
}
=
\tikzmath[scale=.75, transform shape]{
\begin{scope}
\clip[rounded corners = 5] (-.6,0) rectangle (2.4,3);
\filldraw[primedregion=\AColor] (0,0) -- (0,.8) -- (.533,.8) -- (.533,1.2) .. controls ++(90:.4cm) and ++(-135:.2cm) .. (1.2,1.8) .. controls ++(135:.2cm) and ++(270:.4cm) .. (.6,2.4) -- (.6,3) -- (-.6,3) -- (-.6,0); 
\filldraw[primedregion=\CColor] (1.8,0) -- (1.8,1.2) .. controls ++(90:.4cm) and ++(-45:.2cm) .. (1.2,1.8) .. controls ++(-135:.2cm) and ++(90:.4cm) .. (.65,1.2) -- (.65,.8) -- (1.2,.8) -- (1.2,0);
\filldraw[primedregion=\BColor] (0,0) -- (0,.8) -- (.533,.8) -- (.533,1.2) .. controls ++(90:.4cm) and ++(-135:.2cm) .. (1.133,1.8) -- (1.267,1.8) .. controls ++(-135:.2cm) and ++(90:.4cm) .. (.667,1.2) -- (.667,.8) -- (1.2,.8) -- (1.2,0);
\filldraw[boxregion=\AColor] (.6,3) -- (.6,2.4) .. controls ++(270:.4cm) and ++(135:.2cm) .. (1.2,1.8) .. controls ++(45:.2cm) and ++(270:.4cm) .. (1.75,2.4) -- (1.75,3);
\filldraw[boxregion=\CColor] (1.8,0) -- (1.8,1.2) .. controls ++(90:.4cm) and ++(-45:.2cm) .. (1.2,1.8) .. controls ++(45:.2cm) and ++(270:.4cm) .. (1.85,2.4) -- (1.85,3) -- (2.4,3) -- (2.4,0);
\filldraw[boxregion=\BColor] (1.733,3) -- (1.733,2.4) .. controls ++(270:.4cm) and ++(45:.2cm) .. (1.133,1.8) -- (1.267,1.8) .. controls ++(45:.2cm) and ++(270:.4cm) .. (1.867,2.4) -- (1.867,3); 
\end{scope}
\draw[\XColor,thick] (0,0) -- (0,.6);
\draw[\YColor,thick] (1.2,0) -- (1.2,.6);
\draw[\XColor,thick] (.533,.6) -- (.533,1.2) .. controls ++(90:.6cm) and ++(270:.6cm) .. (1.733,2.4) -- (1.733,3);
\draw[\YColor,thick] (.667,.6) -- (.667,1.2) .. controls ++(90:.6cm) and ++(270:.6cm) .. (1.867,2.4) -- (1.867,3);
\draw[\phiColor,thick] (1.8,0) -- (1.8,1.2) .. controls ++(90:.6cm) and ++(270:.6cm) .. (.6,2.4) -- (.6,3);
\roundNbox{unshaded}{(.6,.8)}{.3}{.6}{.6}{\scriptsize{$F^2_{X,Y}$}}; 
\filldraw[white] (1.2,1.8) circle (.1cm);
\draw[thick] (1.2,1.8) circle (.1cm); 
\node at (0,-.2) {\scriptsize{$F(X)$}};
\node at (1.2,-.2) {\scriptsize{$F(Y)$}};
\node at (1.8,3.2) {\scriptsize{$F'(X\xxo_b Y)$}};
\node at (1.8,-.2) {\scriptsize{$\varphi_c$}};
\node at (.6,3.2) {\scriptsize{$\varphi_a$}};
}
\qquad
\tikzmath[scale=.75, transform shape]{
\begin{scope}
\clip[rounded corners = 5] (-.6,0) rectangle (1.8,3);
\filldraw[primedregion=\BColor] (1.2,0) -- (1.2,.6) .. controls ++(90:.6cm) and ++(270:.6cm) .. (0,1.8) -- (0,3) -- (-.6,3) -- (-.6,0);
\filldraw[boxregion=\BColor] (1.2,0) -- (1.2,.6) .. controls ++(90:.6cm) and ++(270:.6cm) .. (0,1.8) -- (0,3) -- (1.8,3) -- (1.8,0);
\end{scope}
\draw[\XColor,thick,dotted] (0,0) -- (0,.6) .. controls ++(90:.6cm) and ++(270:.6cm) .. (1.2,1.8) -- (1.2,2.2);
\draw[\XColor,thick] (1.2,2.2) -- (1.2,3);
\draw[\phiColor,thick] (1.2,0) -- (1.2,.6) .. controls ++(90:.6cm) and ++(270:.6cm) .. (0,1.8) -- (0,3);
\roundNbox{unshaded}{(1.2,2.2)}{.3}{0}{0}{\scriptsize{${F'}^1_b$}}; 
\filldraw[white] (.6,1.2) circle (.1cm);
\draw[thick] (.6,1.2) circle (.1cm); 
\node at (0,-.2) {\scriptsize{$1_{F(b)}$}};
\node at (1.2,3.2) {\scriptsize{$F'(1_b)$}};
%
\node at (1.2,-.2) {\scriptsize{$\varphi_b$}};
\node at (0,3.2) {\scriptsize{$\varphi_b$}};
\node at (.2,1.2) {\scriptsize{$\id$}};
} 
=
\tikzmath[scale=.75, transform shape]{
\begin{scope}
\clip[rounded corners = 5] (-.6,0) rectangle (1.8,3);
\filldraw[primedregion=\BColor] (1.2,0) -- (1.2,1.2) .. controls ++(90:.6cm) and ++(270:.6cm) .. (0,2.4) -- (0,3) -- (-.6,3) -- (-.6,0);
\filldraw[boxregion=\BColor] (1.2,0) -- (1.2,1.2) .. controls ++(90:.6cm) and ++(270:.6cm) .. (0,2.4) -- (0,3) -- (1.8,3) -- (1.8,0);
\end{scope}
\draw[\XColor,thick,dotted] (0,0) -- (0,.8);
\draw[\XColor,thick] (0,.8) -- (0,1.2) .. controls ++(90:.6cm) and ++(270:.6cm) .. (1.2,2.4) -- (1.2,3);
\draw[\phiColor,thick] (1.2,0) -- (1.2,1.2) .. controls ++(90:.6cm) and ++(270:.6cm) .. (0,2.4) -- (0,3);
\roundNbox{unshaded}{(0,.8)}{.3}{0}{0}{\scriptsize{$F^1_b$}}; 
\filldraw[white] (.6,1.8) circle (.1cm);
\draw[thick] (.6,1.8) circle (.1cm); 
\node at (0,-.2) {\scriptsize{$1_{F(b)}$}};
\node at (1.2,3.2) {\scriptsize{$F'(1_b)$}};
\node at (1.2,-.2) {\scriptsize{$\varphi_b$}};
\node at (0,3.2) {\scriptsize{$\varphi_b$}};
}
\qquad
\tikzmath[scale=.75, transform shape]{
\begin{scope}
\clip[rounded corners = 5] (-.6,0) rectangle (1.8,3);
\filldraw[primedregion=\AColor] (0,0) -- (0,.6) .. controls ++(90:.4cm) and ++(-135:.2cm) .. (.6,1.2) .. controls ++(135:.2cm) and ++(270:.4cm) .. (0,1.8) -- (0,3) -- (-.6,3) -- (-.6,0); 
\filldraw[primedregion=\BColor] (1.2,0) -- (1.2,.6) .. controls ++(90:.4cm) and ++(-45:.2cm) .. (.6,1.2) .. controls ++(-135:.2cm) and ++(90:.4cm) .. (0,.6) -- (0,0);
\filldraw[boxregion=\AColor] (0,3) -- (0,1.8) .. controls ++(270:.4cm) and ++(135:.2cm) .. (.6,1.2) .. controls ++(45:.2cm) and ++(270:.4cm) .. (1.2,1.8) -- (1.2,3);
\filldraw[boxregion=\BColor] (1.2,0) -- (1.2,.6) .. controls ++(90:.4cm) and ++(-45:.2cm) .. (.6,1.2) .. controls ++(45:.2cm) and ++(270:.4cm) .. (1.2,1.8) -- (1.2,3) -- (1.8,3) -- (1.8,0);
\end{scope}
\draw[\XColor,thick] (0,0) -- (0,.6) .. controls ++(90:.6cm) and ++(270:.6cm) .. (1.2,1.8) -- (1.2,2.2);
\draw[\ZColor,thick] (1.2,2.2) -- (1.2,3);
\draw[\phiColor,thick] (1.2,0) -- (1.2,.6) .. controls ++(90:.6cm) and ++(270:.6cm) .. (0,1.8) -- (0,3);
\roundNbox{unshaded}{(1.2,2.2)}{.3}{.1}{.1}{\scriptsize{$F'(f)$}}; 
\filldraw[white] (.6,1.2) circle (.1cm);
\draw[thick] (.6,1.2) circle (.1cm); 
\node at (0,-.2) {\scriptsize{$F(X)$}};
\node at (1.2,3.2) {\scriptsize{$F'(Z)$}};
%
\node at (1.2,-.2) {\scriptsize{$\varphi_b$}};
\node at (0,3.2) {\scriptsize{$\varphi_a$}};
} 
=
\tikzmath[scale=.75, transform shape]{
\begin{scope}
\clip[rounded corners = 5] (-.6,0) rectangle (1.8,3);
\filldraw[primedregion=\AColor] (0,0) -- (0,1.2) .. controls ++(90:.4cm) and ++(-135:.2cm) .. (.6,1.8) .. controls ++(135:.2cm) and ++(270:.4cm) .. (0,2.4) -- (0,3) -- (-.6,3) -- (-.6,0); 
\filldraw[primedregion=\BColor] (1.2,0) -- (1.2,1.2) .. controls ++(90:.4cm) and ++(-45:.2cm) .. (.6,1.8) .. controls ++(-135:.2cm) and ++(90:.4cm) .. (0,1.2) -- (0,0);
\filldraw[boxregion=\AColor] (0,3) -- (0,2.4) .. controls ++(270:.4cm) and ++(135:.2cm) .. (.6,1.8) .. controls ++(45:.2cm) and ++(270:.4cm) .. (1.2,2.4) -- (1.2,3);
\filldraw[boxregion=\BColor] (1.2,0) -- (1.2,1.2) .. controls ++(90:.4cm) and ++(-45:.2cm) .. (.6,1.8) .. controls ++(45:.2cm) and ++(270:.4cm) .. (1.2,2.4) -- (1.2,3) -- (1.8,3) -- (1.8,0);
\end{scope}
\draw[\XColor,thick] (0,0) -- (0,.8);
\draw[\ZColor,thick] (0,.8) -- (0,1.2) .. controls ++(90:.6cm) and ++(270:.6cm) .. (1.2,2.4) -- (1.2,3);
\draw[\phiColor,thick] (1.2,0) -- (1.2,1.2) .. controls ++(90:.6cm) and ++(270:.6cm) .. (0,2.4) -- (0,3);
\roundNbox{unshaded}{(0,.8)}{.3}{.1}{.1}{\scriptsize{$F(f)$}}; 
\filldraw[white] (.6,1.8) circle (.1cm);
\draw[thick] (.6,1.8) circle (.1cm); 
\node at (0,-.2) {\scriptsize{$F(X)$}};
\node at (1.2,3.2) {\scriptsize{$F'(Z)$}};
\node at (1.2,-.2) {\scriptsize{$\varphi_b$}};
\node at (0,3.2) {\scriptsize{$\varphi_a$}};
}\,.
$$
\end{defn}

\begin{defn}
\label{Defn:2Modification}
Suppose $\cA,\cB$ are 2-categories, $F,F': \cA \to \cB$ are 2-functors, and $\varphi, \psi : F\Rightarrow F'$ are 2-transformations.
A \emph{2-modification} $n: \varphi\Rrightarrow \psi$
consists of a 2-cell $n_a\in \cB(\varphi_a \Rightarrow \psi_a)$ for all $a\in \cA$ such that
$$
\tikzmath[scale=.75, transform shape]{
\begin{scope}
\clip[rounded corners = 5] (-.6,0) rectangle (1.8,3);
\filldraw[primedregion=\AColor] (0,0) -- (0,.6) .. controls ++(90:.4cm) and ++(-135:.2cm) .. (.6,1.2) .. controls ++(135:.2cm) and ++(270:.4cm) .. (0,1.8) -- (0,3) -- (-.6,3) -- (-.6,0); 
\filldraw[primedregion=\BColor] (1.2,0) -- (1.2,.6) .. controls ++(90:.4cm) and ++(-45:.2cm) .. (.6,1.2) .. controls ++(-135:.2cm) and ++(90:.4cm) .. (0,.6) -- (0,0);
\filldraw[boxregion=\AColor] (0,3) -- (0,1.8) .. controls ++(270:.4cm) and ++(135:.2cm) .. (.6,1.2) .. controls ++(45:.2cm) and ++(270:.4cm) .. (1.2,1.8) -- (1.2,3);
\filldraw[boxregion=\BColor] (1.2,0) -- (1.2,.6) .. controls ++(90:.4cm) and ++(-45:.2cm) .. (.6,1.2) .. controls ++(45:.2cm) and ++(270:.4cm) .. (1.2,1.8) -- (1.2,3) -- (1.8,3) -- (1.8,0);
\end{scope}
\draw[\XColor,thick] (0,0) -- (0,.6) .. controls ++(90:.6cm) and ++(270:.6cm) .. (1.2,1.8) -- (1.2,3);
\draw[\phiColor,thick] (1.2,0) -- (1.2,.6) .. controls ++(90:.6cm) and ++(270:.6cm) .. (0,1.8) -- (0,2.2);
\draw[\psiColor,thick] (0,2.2) -- (0,3);
\roundNbox{unshaded}{(0,2.2)}{.3}{0}{0}{\scriptsize{$n_a$}}; 
\filldraw[white] (.6,1.2) circle (.1cm);
\draw[thick] (.6,1.2) circle (.1cm); 
\node at (0,-.2) {\scriptsize{$F(X)$}};
\node at (1.2,3.2) {\scriptsize{$F'(X)$}};
\node at (1.2,-.2) {\scriptsize{$\varphi_b$}};
\node at (0,3.2) {\scriptsize{$\psi_a$}};
} 
=
\tikzmath[scale=.75, transform shape]{
\begin{scope}
\clip[rounded corners = 5] (-.6,0) rectangle (1.8,3);
\filldraw[primedregion=\AColor] (0,0) -- (0,1.2) .. controls ++(90:.4cm) and ++(-135:.2cm) .. (.6,1.8) .. controls ++(135:.2cm) and ++(270:.4cm) .. (0,2.4) -- (0,3) -- (-.6,3) -- (-.6,0); 
\filldraw[primedregion=\BColor] (1.2,0) -- (1.2,1.2) .. controls ++(90:.4cm) and ++(-45:.2cm) .. (.6,1.8) .. controls ++(-135:.2cm) and ++(90:.4cm) .. (0,1.2) -- (0,0);
\filldraw[boxregion=\AColor] (0,3) -- (0,2.4) .. controls ++(270:.4cm) and ++(135:.2cm) .. (.6,1.8) .. controls ++(45:.2cm) and ++(270:.4cm) .. (1.2,2.4) -- (1.2,3);
\filldraw[boxregion=\BColor] (1.2,0) -- (1.2,1.2) .. controls ++(90:.4cm) and ++(-45:.2cm) .. (.6,1.8) .. controls ++(45:.2cm) and ++(270:.4cm) .. (1.2,2.4) -- (1.2,3) -- (1.8,3) -- (1.8,0);
\end{scope}
\draw[\XColor,thick] (0,0) -- (0,1.2) .. controls ++(90:.6cm) and ++(270:.6cm) .. (1.2,2.4) -- (1.2,3);
\draw[\phiColor,thick] (1.2,0) -- (1.2,.8);
\draw[\psiColor,thick] (1.2,.8) -- (1.2,1.2) .. controls ++(90:.6cm) and ++(270:.6cm) .. (0,2.4) -- (0,3);
\roundNbox{unshaded}{(1.2,.8)}{.3}{0}{0}{\scriptsize{$n_b$}}; 
\filldraw[white] (.6,1.8) circle (.1cm);
\draw[thick] (.6,1.8) circle (.1cm); 
\node at (0,-.2) {\scriptsize{$F(X)$}};
\node at (1.2,3.2) {\scriptsize{$F'(X)$}};
\node at (1.2,-.2) {\scriptsize{$\varphi_b$}};
\node at (0,3.2) {\scriptsize{$\psi_a$}};
}
\qquad\qquad
\forall\,X\in\cA(a\to b)
\qquad\qquad
\begin{aligned}
\tikzmath{
\filldraw[primedregion=white, rounded corners = 5pt] (0,0) rectangle (.6,.6);
\draw[thin, dotted, rounded corners = 5pt] (0,0) rectangle (.6,.6);
}
&=
F
\\
\tikzmath{
\filldraw[boxregion=white, rounded corners = 5pt] (0,0) rectangle (.6,.6);
\draw[thin, dotted, rounded corners = 5pt] (0,0) rectangle (.6,.6);
}
&=
F'
\end{aligned}
$$
The 2-composition of 2-modifications in $\Fun(\cA\to \cB)$ is defined as follows.
Suppose $F,F'\in \Fun(\cA\to \cB)$
and $\varphi,\varphi',\varphi''$ are 2-transformations $F\Rightarrow G$.
Let $n:\varphi\Rrightarrow\varphi'$ and $n':\varphi'\Rrightarrow\varphi''$ be 2-modifications.
The 2-composition in $\Fun(\cA\to \cB)$, denoted by $n'\xt n:\varphi\Rrightarrow\varphi''$ is defined by 
$(n'\xt n)_a:=n'_a\xxt n_a$ for $a\in \cA$ as composition of 2-cells in $\cB$.
\end{defn}

\begin{defn}[1-composition in $\Fun(\cA\to \cB)$]
Suppose $F,F',F''\in \Fun(\cA\to \cB)$ are 2-functors, 
and
let $\varphi:F\Rightarrow F'$ and $\psi:F'\Rightarrow F''$ be 2-transformations.
The 1-composite $\varphi\xo\psi :F\Rightarrow F''$ of 2-transformations is defined as follows.
Let $X\in\cA(a\to b)$, we define $(\varphi\xo\psi)_a:=\varphi_a\xxo\psi_a$ as 1-composition of 1-cells in $\cB$, and $(\varphi\xo\psi)_X$ by
\[
(\varphi\xo\psi)_X
:=
\tikzmath[scale=.7, transform shape]{
\begin{scope}
\clip[rounded corners = 5] (-.6,0) rectangle (1.8,2.4);
\filldraw[primedregion=\AColor] (0,0) -- (0,.6) .. controls ++(90:.4cm) and ++(-135:.2cm) .. (.6,1.2) .. controls ++(135:.2cm) and ++(270:.4cm) .. (0,1.8) -- (0,3) -- (-.6,3) -- (-.6,0); 
\filldraw[primedregion=\BColor] (1.2,0) -- (1.2,.6) .. controls ++(90:.4cm) and ++(-45:.2cm) .. (.6,1.2) .. controls ++(-135:.2cm) and ++(90:.4cm) .. (0,.6) -- (0,0);
\filldraw[plusregion=\AColor] (0,3) -- (0,1.8) .. controls ++(270:.4cm) and ++(135:.2cm) .. (.6,1.2) .. controls ++(45:.2cm) and ++(270:.4cm) .. (1.2,1.8) -- (1.2,3);
\filldraw[plusregion=\BColor] (1.2,0) -- (1.2,.6) .. controls ++(90:.4cm) and ++(-45:.2cm) .. (.6,1.2) .. controls ++(45:.2cm) and ++(270:.4cm) .. (1.2,1.8) -- (1.2,3) -- (1.8,3) -- (1.8,0);
\end{scope}
\draw[\XColor,thick] (0,0) -- (0,.6) .. controls ++(90:.6cm) and ++(270:.6cm) .. (1.2,1.8) -- (1.2,2.4);
\draw[\phiColor,thick] (1.15,0) -- (1.15,.6) .. controls ++(90:.6cm) and ++(270:.6cm) .. (-.05,1.8) -- (-.05,2.4);
\draw[\psiColor,thick] (1.25,0) -- (1.25,.6) .. controls ++(90:.6cm) and ++(270:.6cm) .. (.05,1.8) -- (.05,2.4);
\filldraw[white] (.6,1.2) circle (.1cm);
\draw[thick] (.6,1.2) circle (.1cm); 
\node at (0,-.2) {\scriptsize{$F(X)$}};
\node at (1.2,2.6) {\scriptsize{$F''(X)$}};
\node at (1.2,-.2) {\scriptsize{$(\varphi\xo\psi)_b$}};
\node at (0,2.6) {\scriptsize{$(\varphi\xo\psi)_a$}};
} 
:=
\tikzmath[scale=.7, transform shape]{
\begin{scope}
\clip[rounded corners = 5] (-1.8,-1.8) rectangle (1.8,1.8);
\filldraw[primedregion=\AColor] (-1.2,-1.8) -- (-1.2,-1.4) .. controls ++(90:.4cm) and ++(-135:.2cm) .. (-.6,-.8) .. controls ++(135:.2cm) and ++(270:.4cm) .. (-1.2,-.2) -- (-1.2,1.8) -- (-.55,1.8) -- (-.55,2.7) -- (-1.8,2.7) -- (-1.8,-1.8);
\filldraw[primedregion=\BColor](0,-1.8) -- (0,-1.4) .. controls ++(90:.4cm) and ++(-45:.2cm) .. (-.6,-.8) .. controls ++(-135:.2cm) and ++(90:.4cm) .. (-1.2,-1.4) -- (-1.2,-1.8);
\filldraw[boxregion=\AColor] (-1.2,1.8) -- (-1.2,-.2) .. controls ++(270:.4cm) and ++(135:.2cm) .. (-.6,-.8) .. controls ++(45:.2cm) and ++(270:.4cm) .. (0,-.2) -- (0,.2) .. controls ++(90:.4cm) and ++(-135:.2cm) .. (.6,.8) .. controls ++(135:.2cm) and ++(270:.4cm) .. (0,1.4) -- (0,1.8);
\filldraw[boxregion=\BColor] (0,-1.8) -- (0,-1.4) .. controls ++(90:.4cm) and ++(-45:.2cm) .. (-.6,-.8) .. controls ++(45:.2cm) and ++(270:.4cm) .. (0,-.2) -- (0,.2) .. controls ++(90:.4cm) and ++(-135:.2cm) .. (.6,.8) .. controls ++(-45:.2cm) and ++(90:.4cm) .. (1.2,.2) -- (1.2,-1.8);
\filldraw[plusregion=\AColor] (-.55,2.7) -- (-.55,1.8) -- (0,1.8) -- (0,1.4) .. controls ++(270:.4cm) and ++(135:.2cm) .. (.6,.8) .. controls ++(45:.2cm) and ++(270:.4cm) .. (1.2,1.4) -- (1.2,2.7);
\filldraw[plusregion=\BColor] (1.2,-1.8) -- (1.2,.2) .. controls ++(90:.4cm) and ++(-45:.2cm) .. (.6,.8) .. controls ++(45:.2cm) and ++(270:.4cm) .. (1.2,1.4) -- (1.2,2.7) -- (1.8,2.7) -- (1.8,-1.8);
\end{scope}
\draw[\XColor,thick] (-1.2,-1.8) -- (-1.2,-1.4) .. controls ++(90:.6cm) and ++(270:.6cm) .. (0,-.2) -- (0,.2) .. controls ++(90:.6cm) and ++(270:.6cm) .. (1.2,1.4) -- (1.2,1.8);
\draw[\phiColor,thick] (0,-1.8) -- (0,-1.4) .. controls ++(90:.6cm) and ++(270:.6cm) .. (-1.2,-.2) -- (-1.2,1.8);
\draw[\psiColor,thick] (1.2,-1.8) -- (1.2,.2) .. controls ++(90:.6cm) and ++(270:.6cm) .. (0,1.4) -- (0,1.8);
\filldraw[white] (-.6,-.8) circle (.1cm);
\draw[thick] (-.6,-.8) circle (.1cm); 
\filldraw[white] (.6,.8) circle (.1cm);
\draw[thick] (.6,.8) circle (.1cm); 
\node at (-1.2,-2) {\scriptsize{$F(X)$}};
\node at (0,-2) {\scriptsize{$\varphi_b$}};
\node at (1.2,-2) {\scriptsize{$\psi_b$}};
\node at (-1.2,2) {\scriptsize{$\varphi_a$}};
\node at (0,2) {\scriptsize{$\psi_a$}};
\node at (1.2,2) {\scriptsize{$F''(X)$}};
}
\qquad\qquad
\forall\,X\in\cA(a\to b)
\qquad\qquad
\begin{aligned}
\tikzmath{
\filldraw[primedregion=white, rounded corners = 5pt] (0,0) rectangle (.6,.6);
\draw[thin, dotted, rounded corners = 5pt] (0,0) rectangle (.6,.6);
}
&=
F
\\
\tikzmath{
\filldraw[boxregion=white, rounded corners = 5pt] (0,0) rectangle (.6,.6);
\draw[thin, dotted, rounded corners = 5pt] (0,0) rectangle (.6,.6);
}
&=
F'
\\
\tikzmath{
\filldraw[plusregion=white, rounded corners = 5pt] (0,0) rectangle (.6,.6);
\draw[thin, dotted, rounded corners = 5pt] (0,0) rectangle (.6,.6);
}
&=
F''
\end{aligned}
\]

Suppose $\varphi,\varphi': F\Rightarrow F'$ and $\psi,\psi':F'\Rightarrow F''$ are 2-transformations, 
and
let $n:\varphi\Rrightarrow\varphi'$ and 
$t:\psi\Rrightarrow\psi'$ be 2-modifications.
The 1-composite $n\xo t:\varphi\xo\psi \Rrightarrow \varphi'\xo\psi'$ 
of 2-modifications is defined component-wise 
as 1-composition of 2-cells in $\cB$ by
$(n\xo t)_a:=n_a\xxo t_a$ for $a\in\cA$.

Finally, we define the associator for 1-composition in $\Fun(\cA\to \cB)$ as follows.
Suppose $\varphi:F\Rightarrow F'$, $\psi:F'\Rightarrow F'':$ and $\gamma:F''\Rightarrow F''':$ are 2-transformations.
The associator $\alpha^\xo_{\varphi,\psi,\gamma}$ is an invertible modification 
$(\varphi\xo\psi)\xo\gamma\Rrightarrow \varphi\xo(\psi\xo\gamma)$
which is given component-wise by 
\begin{equation}
\label{eq:ComponentwiseAssociator}
\left(\alpha^\xo_{\varphi,\psi,\gamma}\right)_a:= \alpha^\cB_{\varphi(a),\psi(a),\gamma(a)},
\end{equation}
which is the associator in $\cB$ between 1-cells $\varphi(a),\psi(a),\gamma(a)$.
One checks that $\alpha^\otimes_{\varphi,\psi,\gamma}$ is a modification, and that $\alpha^\otimes$ satisfies the pentagon axiom. 

The left/right unitors $\lambda_\varphi^F:1_F\xo\varphi\Rrightarrow \varphi$ and $\rho_\varphi^{F'}:\varphi\xo 1_{F'}\Rrightarrow \varphi$ are an invertible 2-modifications which are given component-wise by
\begin{equation}
\label{eq:ComponentwiseUnitor}
\left(\lambda_\varphi^F\right)_a:=\lambda_{\varphi(a)}^{{F(a)}} 
\qquad\qquad 
\left(\rho_\varphi^{F'}\right)_a:=\rho_{\varphi(a)}^{{F'(a)}},
\end{equation}
which are the unitors in $\cB$ for 1-cell $\varphi(a)$.
\end{defn}

\begin{rem}
\label{rem:Overlay}
We do not attempt to formalize this overlay operation in this article, as all such string diagrams can be interpreted uniquely as a 2-cell in $\cB$ without confusion.
However, we sketch the following strategy to formalize this graphical calculus, which was communicated to us by David Reutter.

First, by \cite{1903.05777}, the 3D graphical calculus for $\Gray$-categories \cite{1211.0529,1409.2148} may be applied in any 3-category, in particular, to $2\Cat$.
Second, given a 2-category $\cA\in 2\Cat$, we may identify
$\cA = \Fun(* \to \cA)$ where $*$ is the trivial 2-category.
This identification allows us to identify the \emph{internal} 2D string diagrammatic calculus for $\cA$ with the \emph{external} 2D string diagrammatic calculus for $\Fun(* \to \cA)$ as a hom 2-category of $2\Cat$.
Finally, identifying a 2-functor $F: \cA\to \cB$ with the 2-functor $\Fun(*\to \cA) \to \Fun(*\to \cB)$ given by post-composition with $F$, and similarly for transformations and modifications, we see that our overlay graphical calculus is exactly stacking of 2D sheets in the 3D graphical calculus for $2\Cat$.
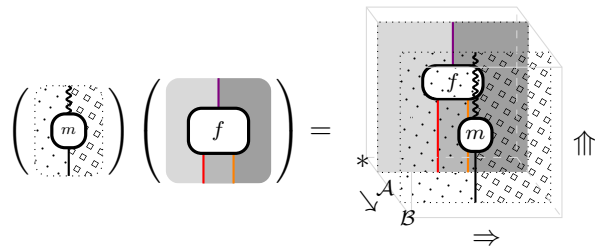
\[
\left(\,
\tikzmath[scale=.75, transform shape]{
\begin{scope}
\clip[rounded corners = 5] (-.6,-.8) rectangle (.6,.8);
\filldraw[primedregion=white] (-.6,-.9) rectangle (0,.9);
\filldraw[boxregion=white] (0,-.9) rectangle (.6,.9);
\end{scope}
\draw[\phiColor,thick] (0,-.8) -- (0,0);
\draw[\psiColor,thick] (0,.8) -- (0,0);
\roundNbox{unshaded}{(0,0)}{.3}{0}{0}{\scriptsize{$m$}}; 
\draw[thin, dotted, rounded corners = 5pt] (-.6,-.8) rectangle (.6,.8);
}
\,\right)
\left(\,
\tikzmath{
\begin{scope}
\clip[rounded corners=5pt] (-.7,-.7) rectangle (.7,.7);
\filldraw[\AColor] (-.2,-.7) -- (-.2,0) -- (0,0) -- (0,.7) -- (-.7,.7) -- (-.7,-.7);
\filldraw[\BColor] (-.2,-.7) rectangle (.2,0);
\filldraw[\CColor] (.2,-.7) -- (.2,0) -- (0,0) -- (0,.7) -- (.7,.7) -- (.7,-.7);
\end{scope}
\draw[\ZColor,thick] (0,.3) -- (0,.7);
\draw[\XColor,thick] (-.2,-.3) -- (-.2,.-.7);
\draw[\YColor,thick] (.2,-.3) -- (.2,-.7);
\roundNbox{unshaded}{(0,0)}{.3}{.1}{.1}{\scriptsize{$f$}};
}
\,\right)
=
\tikzmath{
\begin{scope}
\filldraw[\AColor] (-.45,.6) -- (-.45,2.6) -- (.55,2.6) -- (.55,2.2) -- (-.15,2.2) -- (-.15,.6);
\filldraw[\CColor] (.55,2.2) rectangle (1.55,2.6);
\filldraw[primedregion=\AColor] (-.15,.6) -- (-.15,2.2) -- (.55,2.2) -- (.55,1.8) -- (.35,1.8) -- (.35,.6);
\filldraw[primedregion=\CColor] (.55,2.2) rectangle (.85,1.2);
\filldraw[primedregion=\BColor] (.35,.6) rectangle (.75,1.8);
\filldraw[primedregion=\CColor] (.75,.6) rectangle (.85,1.2);
\filldraw[boxregion=\CColor] (.85,.6) rectangle (1.55,2.2);
\filldraw[primedregion=white] (-.15,.2) rectangle (.85,.6);
\filldraw[boxregion=white] (.85,.2) -- (.85,.6) -- (1.55,.6) -- (1.55,2.2) -- (1.85,2.2) -- (1.85,.2);
\end{scope}
\draw[\XColor,thick] (.35,1.8) -- (.35,.6);
\draw[\YColor,thick] (.75,1.8) -- (.75,.6);
\draw[\ZColor,thick] (.55,1.8) -- (.55,2.6);
\roundNbox{unshaded}{(.55,1.8)}{.22}{.18}{.18}{\scriptsize{$f$}}; 
\filldraw[primedregion=white,rounded corners = 5] (.18,1.6) rectangle (.92,2.0);
\node at (.55,1.8) {\scriptsize{$f$}};
\draw[\psiColor,thick] (.85,1.2) -- (.85,2.2);
\draw[\phiColor,thick] (.85,1.2) -- (.85,.2);
\roundNbox{unshaded}{(.85,1.1)}{.22}{0}{0}{\scriptsize{$m$}};
\draw[dotted] (-.15,.2) rectangle (1.85,2.2); 
\draw[dotted] (-.45,.6) rectangle (1.55,2.6); 
\node at (-.6,.2) [rotate=-53] {$\to$};
\node at (1,-.3) {$\Rightarrow$};
\node at (2.3,1) [rotate=90] {$\Rrightarrow$};
\draw[\AColor] (0,0) rectangle (2,2);
\draw[\AColor] (0,0) -- (-.6,.8) -- (-.6,2.8) -- (1.4,2.8) -- (2,2);
\draw[\AColor] (-.6,2.8) -- (0,2);
\draw[\AColor,dashed] (-.6,.8) -- (1.4,.8) -- (1.4,2.8);
\draw[\AColor,dashed] (1.4,.8) -- (2,0);
\node at (-.65,.75) {$*$};
\node at (-.35,.4) {\scriptsize{$\cA$}};
\node at (-.05,0) {\scriptsize{$\cB$}};
}
\]

Now in order to interpret each diagram as a unique 2-morphism in $\cB$, one should require the strings and coupons of our $\cA$-diagram and our $\Fun(\cA\to \cB)$ diagram not overlap, except at finitely many points where strings can cross transversely.
The axioms of 2-functor, 2-transformation, and 2-modification will then ensure than any two ways of resolving non-generic intersections agree.
For example, we may overlay the 2-transformation $\varphi: F\Rightarrow F'$
on the identity 2-morphism $\id_X \otimes_b \id_Y$ in $\cA$ in several ways.
The equality of two such ways below produces the monoidal coherence axiom:
\[
\left(
\underset{\varphi: F\Rightarrow F'}{
\tikzmath{
\begin{scope}
\filldraw[primedregion=white, rounded corners = 5pt] (0,0) rectangle (.3,.6);
\filldraw[boxregion=white, rounded corners = 5pt] (.3,0) rectangle (.6,.6);
\end{scope}
\draw[\phiColor,thick] (.3,0) -- (.3,.6);
\draw[thin, dotted, rounded corners = 5pt] (0,0) rectangle (.6,.6);
}
}
\right)
\left(
\tikzmath{
\begin{scope}
\clip[rounded corners = 5] (0,0) rectangle (.9,.6);
\fill[\AColor] (0,0) rectangle (.3,.6);
\fill[\BColor] (.3,0) rectangle (.6,.6);
\fill[\CColor] (.6,0) rectangle (.9,.6);
\end{scope}
\draw[\XColor,thick] (.3,0) -- (.3,.6);
\draw[\YColor,thick] (.6,0) -- (.6,.6);
\node at (.3,-.2) {\scriptsize{$X$}};
\node at (.6,-.2) {\scriptsize{$Y$}};
}
\right)
=
\tikzmath[scale=.75, transform shape]{
\begin{scope}
\clip[rounded corners = 5] (-.6,0) rectangle (3,4.2);
\filldraw[primedregion=\AColor] (0,0) -- (0,2) .. controls ++(90:.4cm) and ++(-135:.2cm) .. (.6,2.6) .. controls ++(135:.2cm) and ++(270:.4cm) .. (0,3.2) -- (0,4.2) -- (-.6,4.2) -- (-.6,0); 
\filldraw[primedregion=\BColor] (1.2,0) -- (1.2,.4) .. controls ++(90:.4cm) and ++(-135:.2cm) .. (1.8,1) .. controls ++(135:.2cm) and ++(270:.4cm) .. (1.2,1.6) -- (1.2,2) .. controls ++(90:.4cm) and ++(-45:.2cm) .. (.6,2.6) .. controls ++(-135:.2cm) and ++(90:.4cm) .. (0,2) -- (0,0);
\filldraw[primedregion=\CColor] (2.4,0) -- (2.4,.4) .. controls ++(90:.4cm) and ++(-45:.2cm) .. (1.8,1) .. controls ++(-135:.2cm) and ++(90:.4cm) .. (1.2,.4) -- (1.2,0);
\filldraw[boxregion=\AColor] (0,4.2) -- (0,3.2) .. controls ++(270:.4cm) and ++(135:.2cm) .. (.6,2.6) .. controls ++(45:.2cm) and ++(270:.4cm) .. (1.2,3.2) -- (1.75,3.6) -- (1.75,4.2);
\filldraw[boxregion=\BColor] (1.85,4.2) -- (1.85,3.6) -- (2.4,3.2) -- (2.4,1.6) .. controls ++(270:.4cm) and ++(45:.2cm) .. (1.8,1) .. controls ++(135:.2cm) and ++(270:.4cm) .. (1.2,1.6) -- (1.2,2) .. controls ++(90:.4cm) and ++(-45:.2cm) .. (.6,2.6) .. controls ++(45:.2cm) and ++(270:.4cm) .. (1.2,3.2) -- (1.75,3.6) -- (1.75,4.2);
\filldraw[boxregion=\CColor] (2.4,0) -- (2.4,.4) .. controls ++(90:.4cm) and ++(-45:.2cm) .. (1.8,1) .. controls ++(45:.2cm) and ++(270:.4cm) .. (2.4,1.6) -- (2.4,3.2) -- (1.85,3.6) -- (1.85,4.2) -- (3,4.2) -- (3,0);
\filldraw[\BColor] (1,1.8) circle (.2cm);
\end{scope}
\draw[\XColor,thick] (0,0) -- (0,.4) -- (0,2) .. controls ++(90:.6cm) and ++(270:.6cm) .. (1.2,3.2) -- (1.2,3.6);
\draw[\YColor,thick] (1.2,0) -- (1.2,.4) .. controls ++(90:.6cm) and ++(270:.6cm) .. (2.4,1.6) -- (2.4,3.6);
\draw[\phiColor,thick] (2.4,0) -- (2.4,.4) .. controls ++(90:.6cm) and ++(270:.6cm) .. (1.2,1.6) -- (1.2,2) .. controls ++(90:.6cm) and ++(270:.6cm) .. (0,3.2) -- (0,4.2);
\draw[\XColor,thick] (1.733,3.6) -- (1.733,4.2);
\draw[\YColor,thick] (1.867,3.6) -- (1.867,4.2);
\roundNbox{unshaded}{(1.8,3.4)}{.3}{.6}{.6}{\scriptsize{${F'}^2_{X,Y}$}}; 
\filldraw[white] (1.8,1) circle (.1cm);
\draw[thick] (1.8,1) circle (.1cm); 
\filldraw[white] (.6,2.6) circle (.1cm);
\draw[thick] (.6,2.6) circle (.1cm); 
\node at (0,-.2) {\scriptsize{$F(X)$}};
\node at (1.2,-.2) {\scriptsize{$F(Y)$}};
\node at (1.8,4.4) {\scriptsize{$F'(X\xxo_b Y)$}};
\node at (2.4,-.2) {\scriptsize{$\varphi_c$}};
\node at (1,1.8) {\scriptsize{$\varphi_b$}};
\node at (0,4.4) {\scriptsize{$\varphi_a$}};
}
\quad
\text{ or }
\quad
\tikzmath[scale=.75, transform shape]{
\begin{scope}
\clip[rounded corners = 5] (-.6,0) rectangle (2.4,3);
\filldraw[primedregion=\AColor] (0,0) -- (0,.8) -- (.533,.8) -- (.533,1.2) .. controls ++(90:.4cm) and ++(-135:.2cm) .. (1.2,1.8) .. controls ++(135:.2cm) and ++(270:.4cm) .. (.6,2.4) -- (.6,3) -- (-.6,3) -- (-.6,0); 
\filldraw[primedregion=\CColor] (1.8,0) -- (1.8,1.2) .. controls ++(90:.4cm) and ++(-45:.2cm) .. (1.2,1.8) .. controls ++(-135:.2cm) and ++(90:.4cm) .. (.65,1.2) -- (.65,.8) -- (1.2,.8) -- (1.2,0);
\filldraw[primedregion=\BColor] (0,0) -- (0,.8) -- (.533,.8) -- (.533,1.2) .. controls ++(90:.4cm) and ++(-135:.2cm) .. (1.133,1.8) -- (1.267,1.8) .. controls ++(-135:.2cm) and ++(90:.4cm) .. (.667,1.2) -- (.667,.8) -- (1.2,.8) -- (1.2,0);
\filldraw[boxregion=\AColor] (.6,3) -- (.6,2.4) .. controls ++(270:.4cm) and ++(135:.2cm) .. (1.2,1.8) .. controls ++(45:.2cm) and ++(270:.4cm) .. (1.75,2.4) -- (1.75,3);
\filldraw[boxregion=\CColor] (1.8,0) -- (1.8,1.2) .. controls ++(90:.4cm) and ++(-45:.2cm) .. (1.2,1.8) .. controls ++(45:.2cm) and ++(270:.4cm) .. (1.85,2.4) -- (1.85,3) -- (2.4,3) -- (2.4,0);
\filldraw[boxregion=\BColor] (1.733,3) -- (1.733,2.4) .. controls ++(270:.4cm) and ++(45:.2cm) .. (1.133,1.8) -- (1.267,1.8) .. controls ++(45:.2cm) and ++(270:.4cm) .. (1.867,2.4) -- (1.867,3); 
\end{scope}
\draw[\XColor,thick] (0,0) -- (0,.6);
\draw[\YColor,thick] (1.2,0) -- (1.2,.6);
\draw[\XColor,thick] (.533,.6) -- (.533,1.2) .. controls ++(90:.6cm) and ++(270:.6cm) .. (1.733,2.4) -- (1.733,3);
\draw[\YColor,thick] (.667,.6) -- (.667,1.2) .. controls ++(90:.6cm) and ++(270:.6cm) .. (1.867,2.4) -- (1.867,3);
\draw[\phiColor,thick] (1.8,0) -- (1.8,1.2) .. controls ++(90:.6cm) and ++(270:.6cm) .. (.6,2.4) -- (.6,3);
\roundNbox{unshaded}{(.6,.8)}{.3}{.6}{.6}{\scriptsize{$F^2_{X,Y}$}}; 
\filldraw[white] (1.2,1.8) circle (.1cm);
\draw[thick] (1.2,1.8) circle (.1cm); 
\node at (0,-.2) {\scriptsize{$F(X)$}};
\node at (1.2,-.2) {\scriptsize{$F(Y)$}};
\node at (1.8,3.2) {\scriptsize{$F'(X\xxo_b Y)$}};
\node at (1.8,-.2) {\scriptsize{$\varphi_c$}};
\node at (.6,3.2) {\scriptsize{$\varphi_a$}};
}\,.
\]
For another example, when we have a 2-modification between 2-transformations, we may overlay it on an identity 2-morphism $\id_X$ in many ways.
The equality of two such ways below produces the modification coherence axiom:
\[
\left(
\tikzmath[scale=.75, transform shape]{
\begin{scope}
\clip[rounded corners = 5] (-.6,-.8) rectangle (.6,.8);
\filldraw[primedregion=white] (-.6,-.9) rectangle (0,.9);
\filldraw[boxregion=white] (0,-.9) rectangle (.6,.9);
\end{scope}
\draw[\phiColor,thick] (0,-.8) -- (0,0);
\draw[\psiColor,thick] (0,.8) -- (0,0);
\roundNbox{unshaded}{(0,0)}{.3}{0}{0}{\scriptsize{$n$}}; 
\node at (0,1) {\scriptsize{$\varphi'$}};
\node at (0,-1) {\scriptsize{$\varphi$}};
\draw[thin, dotted, rounded corners = 5pt] (-.6,-.8) rectangle (.6,.8);
}
\right)
\left(
\tikzmath{
\begin{scope}
\clip[rounded corners=5pt] (-.3,0) rectangle (.3,.6);
\fill[\AColor] (0,0) rectangle (-.3,.6);
\fill[\BColor] (0,0) rectangle (.3,.6);
\end{scope}
\draw[thick, \XColor] (0,0) -- (0,.6);
\node at (0,-.2) {\scriptsize{$X$}};
}
\right)
=
\tikzmath[scale=.75, transform shape]{
\begin{scope}
\clip[rounded corners = 5] (-.6,0) rectangle (1.8,3);
\filldraw[primedregion=\AColor] (0,0) -- (0,.6) .. controls ++(90:.4cm) and ++(-135:.2cm) .. (.6,1.2) .. controls ++(135:.2cm) and ++(270:.4cm) .. (0,1.8) -- (0,3) -- (-.6,3) -- (-.6,0); 
\filldraw[primedregion=\BColor] (1.2,0) -- (1.2,.6) .. controls ++(90:.4cm) and ++(-45:.2cm) .. (.6,1.2) .. controls ++(-135:.2cm) and ++(90:.4cm) .. (0,.6) -- (0,0);
\filldraw[boxregion=\AColor] (0,3) -- (0,1.8) .. controls ++(270:.4cm) and ++(135:.2cm) .. (.6,1.2) .. controls ++(45:.2cm) and ++(270:.4cm) .. (1.2,1.8) -- (1.2,3);
\filldraw[boxregion=\BColor] (1.2,0) -- (1.2,.6) .. controls ++(90:.4cm) and ++(-45:.2cm) .. (.6,1.2) .. controls ++(45:.2cm) and ++(270:.4cm) .. (1.2,1.8) -- (1.2,3) -- (1.8,3) -- (1.8,0);
\end{scope}
\draw[\XColor,thick] (0,0) -- (0,.6) .. controls ++(90:.6cm) and ++(270:.6cm) .. (1.2,1.8) -- (1.2,3);
\draw[\phiColor,thick] (1.2,0) -- (1.2,.6) .. controls ++(90:.6cm) and ++(270:.6cm) .. (0,1.8) -- (0,2.2);
\draw[\psiColor,thick] (0,2.2) -- (0,3);
\roundNbox{unshaded}{(0,2.2)}{.3}{0}{0}{\scriptsize{$n_a$}}; 
\filldraw[white] (.6,1.2) circle (.1cm);
\draw[thick] (.6,1.2) circle (.1cm); 
\node at (0,-.2) {\scriptsize{$F(X)$}};
\node at (1.2,3.2) {\scriptsize{$F'(X)$}};
\node at (1.2,-.2) {\scriptsize{$\varphi_b$}};
\node at (0,3.2) {\scriptsize{$\psi_a$}};
} 
\quad
\text{ or }
\quad
\tikzmath[scale=.75, transform shape]{
\begin{scope}
\clip[rounded corners = 5] (-.6,0) rectangle (1.8,3);
\filldraw[primedregion=\AColor] (0,0) -- (0,1.2) .. controls ++(90:.4cm) and ++(-135:.2cm) .. (.6,1.8) .. controls ++(135:.2cm) and ++(270:.4cm) .. (0,2.4) -- (0,3) -- (-.6,3) -- (-.6,0); 
\filldraw[primedregion=\BColor] (1.2,0) -- (1.2,1.2) .. controls ++(90:.4cm) and ++(-45:.2cm) .. (.6,1.8) .. controls ++(-135:.2cm) and ++(90:.4cm) .. (0,1.2) -- (0,0);
\filldraw[boxregion=\AColor] (0,3) -- (0,2.4) .. controls ++(270:.4cm) and ++(135:.2cm) .. (.6,1.8) .. controls ++(45:.2cm) and ++(270:.4cm) .. (1.2,2.4) -- (1.2,3);
\filldraw[boxregion=\BColor] (1.2,0) -- (1.2,1.2) .. controls ++(90:.4cm) and ++(-45:.2cm) .. (.6,1.8) .. controls ++(45:.2cm) and ++(270:.4cm) .. (1.2,2.4) -- (1.2,3) -- (1.8,3) -- (1.8,0);
\end{scope}
\draw[\XColor,thick] (0,0) -- (0,1.2) .. controls ++(90:.6cm) and ++(270:.6cm) .. (1.2,2.4) -- (1.2,3);
\draw[\phiColor,thick] (1.2,0) -- (1.2,.8);
\draw[\psiColor,thick] (1.2,.8) -- (1.2,1.2) .. controls ++(90:.6cm) and ++(270:.6cm) .. (0,2.4) -- (0,3);
\roundNbox{unshaded}{(1.2,.8)}{.3}{0}{0}{\scriptsize{$n_b$}}; 
\filldraw[white] (.6,1.8) circle (.1cm);
\draw[thick] (.6,1.8) circle (.1cm); 
\node at (0,-.2) {\scriptsize{$F(X)$}};
\node at (1.2,3.2) {\scriptsize{$F'(X)$}};
\node at (1.2,-.2) {\scriptsize{$\varphi_b$}};
\node at (0,3.2) {\scriptsize{$\psi_a$}};
}\,.
\]
Here, the white dots which appear may be interpreted as interchangers in $2\Cat$ (see Construction \ref{const:1CompositionIn2Cat} below) which arise from resolving the two stacked 2D diagrams in $2\Cat$.
(Recall that ${}_aX_b\in \cA$ is a transformation when viewed as a 1-morphism in $\Fun(*\to A)$.)

We leave a rigorous proof of our formalization strategy of this `overlay' graphical calculus to the interested reader.
\end{rem}

\subsection{The \texorpdfstring{$\rm C^*/W^*$}{C*/W*} 2-category \texorpdfstring{$\Fun^\dag(\cA\to \cB)$}{Fun(C->D)} between \texorpdfstring{$\rm C^*/W^*$}{C*/W*} 2-categories}

To the best of our knowledge, the notion of $\rm C^*$ 2-category first appeared in \cite{MR1444286}, and the notion of $\rm W^*$ 2-category first appeared in \cite{MR2325696}.
The notion of $\rm W^*$-category was studied in detail in \cite{MR808930}.
We refer the reader to \cite[\S2.1]{2105.12010} for an introduction to $\rm C^*/W^*$ 2-categories.

\begin{defn}
Suppose $\cA,\cB$ are $\rm C^*/W^*$ 2-categories.
A $\dag$ 2-functor $F:\cA\to\cB$ is a 2-functor
$F=(F,F^2,F^1):\cA \to \cB$ such that 
$F^2_{X,Y}$ and $F^1_a$ are unitary for all composable 1-cells $X,Y$ in $\cA$ and all objects $a\in \cA$.
When $\cA,\cB$ are $\rm W^*$, we call a $\dag$ 2-functor \emph{normal} when each hom functor $F_{a\to b}: \cA(a\to b) \to \cB(F(a)\to F(b))$ is a normal $\dag$ functor.

Suppose now $F,G: \cA \to \cB$ are $\dag$-2-functors.
A $\dag$-2-transformation $\varphi:F\Rightarrow G$ consists of
a 2-transformation $\varphi=(\varphi_c, \varphi_X): F\Rightarrow G$ such that every invertible 2-cell $\varphi_X \in \cB(F(X)\xxo_{F(b)} \varphi_b \Rightarrow \varphi_a \xxo_{G(a)} G(X))$ is unitary.

Given two $\dag$-2-transformations $\varphi, \psi : F\Rightarrow G$, a 2-modification $n: \varphi\Rrightarrow \psi$ is \emph{(uniformly) bounded} if 
the 2-cells $n_a\in \cB(\varphi_a \Rightarrow \psi_a)$ for all $a\in \cA$ are uniformly bounded.

Now consider the 2-subcategory $\Fun^\dag(\cA\to \cB)$ of $\Fun(\cA\to \cB)$ consisting of $\dag$ 2-functors, $\dag$ 2-transformations, and uniformly bounded modifications.
When $\cA,\cB$ are $\rm W^*$, we further require all $\dag$ 2-functors to be normal.
\end{defn}

\begin{rem}
\label{rem:Underlying2FunctorEquivalence}
It is well known (e.g., see \cite[Thm.~7.4.1]{2002.06055})
that a 2-functor is an equivalence if and only if it is an equivalence on hom 1-categories (fully faithful on 2-morphisms and essentially surjective on 1-morphisms) and essentially surjective on objects.
Similarly, a $\dag$ 2-functor is an equivalence if and only if it is a $\dag$-equivalence on hom categories (fully faithful on 2-morphisms and unitarily essentially surjective on 1-morphisms) and unitarily essentially surjective on objects.

When $F: \cC\to \cD$ is a $\dag$ 2-functor between $\rm C^*$ 2-categories, observe that $F$ is a dagger equivalence if and only if the underlying 2-fucntor is an equivalence.
Indeed, $F$ is unitarily essentially surjective on 1-morphisms and objects if and only if it is essentially surjective on 1-morphisms and objects by the existence of polar decomposition for invertible 2-morphisms in $\cD$.

Finally, observe that when $\cC,\cD$ are $\rm W^*$, any inverse $\dag$ 2-functor will automatically be normal.
This is an immediate consequence of the fact that every unital $*$-isomorphism between von Neumann algebras is automatically normal using Roberts' $2\times 2$ trick \cite[Lem.~2.6]{MR808930} on linking algebras of hom 1-categories.
\end{rem}

In Proposition \ref{prop:FunDagC*W*} below, we prove that whenever $\cA,\cB$ are $\rm C^*/W^*$, then so is $\Fun^\dag(\cA\to \cB)$ respectively.
In order to prove this result, we prove Lemma \ref{Lem:ProdvNa} on weak* convergence in a product von Neumann algebra, which is certainly known to experts.

Suppose that $(M_i)_{i\in I}$ is a family of von Neumann algebras, and consider the product von Neumann algebra $\prod_{i\in I} M_i$, which is defined as the double commutant of the unital $*$-algebra of uniformly bounded elements $(m_i)$ in the algebraic product of the $M_i$
acting on the Hilbert space $\prod_{i\in I} H_i$, which consists of $L^2$-summable sequences of vectors.
For $j\in I$, there are mutually orthogonal projections $p_j : \prod_i H_i \to H_j$ such that $\sum p_j = 1$ SOT, so every element $m\in \prod_i M_i$ is diagonal, i.e., $m$ be written as $m=(m_i:= p_imp_i)_{i\in I}$.
\begin{lem}
\label{Lem:ProdvNa}
A norm-bounded net $(m_i)^j \to (m_i)$ in the weak* topology on $\prod M_i$ if and only if every component net $m_i^j \to m_i$ in the weak* topology on $M_i$.
\end{lem}
\begin{proof}
On norm-bounded sets in a von Neumann algebra, the weak* topology agrees with the weak operator topology.
Suppose $\eta, \xi \in \prod_i H_i$.
It is clear that
$\langle (m_i)^j\eta, \xi \rangle \to \langle (m_i)\eta, \xi\rangle$
for all $\eta, \xi$ implies 
$\langle m_i^j\eta_i, \xi_i \rangle \to \langle m_i\eta_i, \xi_i\rangle$
for all $i$.

For the converse, let $\varepsilon>0$.
Suppose $M$ is the norm bound for $(m_i)^j$ and $(m_i)$. 
It suffices to show $\langle (m_i)^j\eta,\xi\rangle\to (m_i)\eta,\xi\rangle$ for all given $\eta,\xi\in \prod_i H_i$ with $\|\eta\|,\|\xi\|<1$.
Now choose $\eta',\xi'$ in a finite product with $\|\eta'\|<1$ and $\|\xi'\|<1$ such that 
$$
\|\eta-\eta'\|<\frac{\varepsilon}{5M} 
\qquad
\text{and}
\qquad
\|\xi-\xi'\|<\frac{\varepsilon}{5M}.
$$
Since $\eta',\xi'$ are finitely supported and $m_i^j\to m_i$ weak* for all components $i\in I$ by assumption, we can choose
$j_0$ such that for all $j\geq j_0$,
$$
|\langle [(m_i)^j-(m_i)]\eta',\xi'\rangle|<\frac{\varepsilon}{5}.
$$
Then for all $j\geq j_0$, we have
\begin{align*}
|\langle [(m_i)^j - (m_i)]\eta, \xi \rangle |
&\leq
|\langle (m_i)^j(\eta-\eta'), \xi \rangle |
+
|\langle (m_i)^j\eta', (\xi-\xi') \rangle |
\\
&\qquad+
|\langle [(m_i)^j - (m_i)]\eta', \xi' \rangle |
+
|\langle (m_i)(\eta-\eta'), \xi \rangle |
+
|\langle (m_i)\eta', (\xi-\xi') \rangle | 
\\
& \le
\|(m_i)^j\| \|\eta-\eta'\| \|\xi \|
+
\|(m_i)^j\| \|\eta'\| \|(\xi-\xi')\|
\\
&\qquad+
|\langle [(m_i)^j - (m_i)]\eta', \xi' \rangle |
+
\|(m_i)\| \|\eta-\eta'\| \|\xi\|
+
\|(m_i)\| \|\eta'\| \|\xi-\xi'\| 
\\
& 
< \varepsilon. \qedhere
\end{align*}
\end{proof}

\begin{construction}
\label{construction:DaggerStructureOnFunAB}
We construct a $\dag$-structure on $\Fun^\dag(\cA\to \cB)$ (c.f.~\cite{2004.12760}).
Suppose $F,F'\in \Fun^\dag(\cA\to\cB)$,
$\varphi,\psi: F\Rightarrow F'$,
and $n:\varphi\Rrightarrow \psi$ is a uniformly bounded modification.
For each 0-cell $b\in\cB$,
we define $(n^\dag)_b:=(n_b)^\dag$, where $(n_b)^\dag$ is the dagger in $\cB$.

We now verify that $n^\dag$ is a modification $\psi\Rightarrow\varphi$ with $\|n^\dag\|=\|n\|$.
First, note that $\varphi_X,\psi_X$ are unitaries for all $X\in \cA(a\to b)$. 
We compose $\psi_X^\dag$ on the top and $\varphi_X^\dag$ on the bottom, and apply the dagger in $\cB$, to obtain
\[
\tikzmath[scale=.75, transform shape]{
\begin{scope}
\clip[rounded corners = 5] (-.6,0) rectangle (1.8,3);
\filldraw[primedregion=\AColor] (0,0) -- (0,1.2) .. controls ++(90:.4cm) and ++(-135:.2cm) .. (.6,1.8) .. controls ++(135:.2cm) and ++(270:.4cm) .. (0,2.4) -- (0,3) -- (-.6,3) -- (-.6,0); 
\filldraw[boxregion=\AColor] (1.2,0) -- (1.2,1.2) .. controls ++(90:.4cm) and ++(-45:.2cm) .. (.6,1.8) .. controls ++(-135:.2cm) and ++(90:.4cm) .. (0,1.2) -- (0,0);
\filldraw[primedregion=\BColor] (0,3) -- (0,2.4) .. controls ++(270:.4cm) and ++(135:.2cm) .. (.6,1.8) .. controls ++(45:.2cm) and ++(270:.4cm) .. (1.2,2.4) -- (1.2,3);
\filldraw[boxregion=\BColor] (1.2,0) -- (1.2,1.2) .. controls ++(90:.4cm) and ++(-45:.2cm) .. (.6,1.8) .. controls ++(45:.2cm) and ++(270:.4cm) .. (1.2,2.4) -- (1.2,3) -- (1.8,3) -- (1.8,0);
\filldraw[\AColor] (0,1.8) circle (.3cm);
\end{scope}
\draw[\phiColor,thick] (0,0) -- (0,.8);
\draw[\psiColor,thick] (0,.8) -- (0,1.2) .. controls ++(90:.6cm) and ++(270:.6cm) .. (1.2,2.4) -- (1.2,3);
\draw[\XColor,thick] (1.2,0) -- (1.2,.8);
\draw[\XColor,thick] (1.2,.8) -- (1.2,1.2) .. controls ++(90:.6cm) and ++(270:.6cm) .. (0,2.4) -- (0,3);
\roundNbox{unshaded}{(0,.8)}{.3}{0}{0}{\scriptsize{$n_a$}}; 
\filldraw[white] (.6,1.8) circle (.1cm);
\draw[thick] (.6,1.8) circle (.1cm); 
\node at (0,-.2) {\scriptsize{$\varphi_a$}};
\node at (1.2,3.2) {\scriptsize{$\psi_b$}};
\node at (1.2,-.2) {\scriptsize{$F'(X)$}};
\node at (0,3.2) {\scriptsize{$F(X)$}};
\node at (0,1.8) {\scriptsize{$\psi_X^\dag$}};
}
=
\tikzmath[scale=.75, transform shape]{
\begin{scope}
\clip[rounded corners = 5] (-.6,0) rectangle (1.8,3);
\filldraw[primedregion=\AColor] (0,0) -- (0,.6) .. controls ++(90:.4cm) and ++(-135:.2cm) .. (.6,1.2) .. controls ++(135:.2cm) and ++(270:.4cm) .. (0,1.8) -- (0,3) -- (-.6,3) -- (-.6,0); 
\filldraw[boxregion=\AColor] (1.2,0) -- (1.2,.6) .. controls ++(90:.4cm) and ++(-45:.2cm) .. (.6,1.2) .. controls ++(-135:.2cm) and ++(90:.4cm) .. (0,.6) -- (0,0);
\filldraw[primedregion=\BColor] (0,3) -- (0,1.8) .. controls ++(270:.4cm) and ++(135:.2cm) .. (.6,1.2) .. controls ++(45:.2cm) and ++(270:.4cm) .. (1.2,1.8) -- (1.2,3);
\filldraw[boxregion=\BColor] (1.2,0) -- (1.2,.6) .. controls ++(90:.4cm) and ++(-45:.2cm) .. (.6,1.2) .. controls ++(45:.2cm) and ++(270:.4cm) .. (1.2,1.8) -- (1.2,3) -- (1.8,3) -- (1.8,0);
\filldraw[\AColor] (0,1.2) circle (.3cm);
\end{scope}
\draw[\phiColor,thick] (0,0) -- (0,.6) .. controls ++(90:.6cm) and ++(270:.6cm) .. (1.2,1.8) -- (1.2,2.2);
\draw[\psiColor,thick] (1.2,2.2) -- (1.2,3);
\draw[\XColor,thick] (1.2,0) -- (1.2,.6) .. controls ++(90:.6cm) and ++(270:.6cm) .. (0,1.8) -- (0,2.2);
\draw[\XColor,thick] (0,2.2) -- (0,3);
\roundNbox{unshaded}{(1.2,2.2)}{.3}{0}{0}{\scriptsize{$n_b$}}; 
\filldraw[white] (.6,1.2) circle (.1cm);
\draw[thick] (.6,1.2) circle (.1cm); 
\node at (0,-.2) {\scriptsize{$\varphi_a$}};
\node at (1.2,3.2) {\scriptsize{$\psi_b$}};
\node at (1.2,-.2) {\scriptsize{$F'(X)$}};
\node at (0,3.2) {\scriptsize{$F(X)$}};
\node at (0,1.2) {\scriptsize{$\varphi_{X}^\dag$}};
}
\qquad\overset{\dag}{\Longrightarrow}\qquad
\tikzmath[scale=.75, transform shape]{
\begin{scope}
\clip[rounded corners = 5] (-.6,0) rectangle (1.8,3);
\filldraw[primedregion=\AColor] (0,0) -- (0,.6) .. controls ++(90:.4cm) and ++(-135:.2cm) .. (.6,1.2) .. controls ++(135:.2cm) and ++(270:.4cm) .. (0,1.8) -- (0,3) -- (-.6,3) -- (-.6,0); 
\filldraw[primedregion=\BColor] (1.2,0) -- (1.2,.6) .. controls ++(90:.4cm) and ++(-45:.2cm) .. (.6,1.2) .. controls ++(-135:.2cm) and ++(90:.4cm) .. (0,.6) -- (0,0);
\filldraw[boxregion=\AColor] (0,3) -- (0,1.8) .. controls ++(270:.4cm) and ++(135:.2cm) .. (.6,1.2) .. controls ++(45:.2cm) and ++(270:.4cm) .. (1.2,1.8) -- (1.2,3);
\filldraw[boxregion=\BColor] (1.2,0) -- (1.2,.6) .. controls ++(90:.4cm) and ++(-45:.2cm) .. (.6,1.2) .. controls ++(45:.2cm) and ++(270:.4cm) .. (1.2,1.8) -- (1.2,3) -- (1.8,3) -- (1.8,0);
\filldraw[\AColor] (0,1.2) circle (.3cm);
\end{scope}
\draw[\XColor,thick] (0,0) -- (0,.6) .. controls ++(90:.6cm) and ++(270:.6cm) .. (1.2,1.8) -- (1.2,3);
\draw[\psiColor,thick] (1.2,0) -- (1.2,.6) .. controls ++(90:.6cm) and ++(270:.6cm) .. (0,1.8) -- (0,2.2);
\draw[\phiColor,thick] (0,2.2) -- (0,3);
\roundNbox{unshaded}{(0,2.2)}{.3}{0}{0}{\scriptsize{$n_a^\dag$}}; 
\filldraw[white] (.6,1.2) circle (.1cm);
\draw[thick] (.6,1.2) circle (.1cm); 
\node at (0,-.2) {\scriptsize{$F(X)$}};
\node at (1.2,3.2) {\scriptsize{$F'(X)$}};
\node at (1.2,-.2) {\scriptsize{$\psi_b$}};
\node at (0,3.2) {\scriptsize{$\varphi_a$}};
\node at (0,1.2) {\scriptsize{$\psi_{X}$}};
} 
=
\tikzmath[scale=.75, transform shape]{
\begin{scope}
\clip[rounded corners = 5] (-.6,0) rectangle (1.8,3);
\filldraw[primedregion=\AColor] (0,0) -- (0,1.2) .. controls ++(90:.4cm) and ++(-135:.2cm) .. (.6,1.8) .. controls ++(135:.2cm) and ++(270:.4cm) .. (0,2.4) -- (0,3) -- (-.6,3) -- (-.6,0); 
\filldraw[primedregion=\BColor] (1.2,0) -- (1.2,1.2) .. controls ++(90:.4cm) and ++(-45:.2cm) .. (.6,1.8) .. controls ++(-135:.2cm) and ++(90:.4cm) .. (0,1.2) -- (0,0);
\filldraw[boxregion=\AColor] (0,3) -- (0,2.4) .. controls ++(270:.4cm) and ++(135:.2cm) .. (.6,1.8) .. controls ++(45:.2cm) and ++(270:.4cm) .. (1.2,2.4) -- (1.2,3);
\filldraw[boxregion=\BColor] (1.2,0) -- (1.2,1.2) .. controls ++(90:.4cm) and ++(-45:.2cm) .. (.6,1.8) .. controls ++(45:.2cm) and ++(270:.4cm) .. (1.2,2.4) -- (1.2,3) -- (1.8,3) -- (1.8,0);
\filldraw[\AColor] (0,1.8) circle (.3cm);
\end{scope}
\draw[\XColor,thick] (0,0) -- (0,1.2) .. controls ++(90:.6cm) and ++(270:.6cm) .. (1.2,2.4) -- (1.2,3);
\draw[\psiColor,thick] (1.2,0) -- (1.2,.8);
\draw[\phiColor,thick] (1.2,.8) -- (1.2,1.2) .. controls ++(90:.6cm) and ++(270:.6cm) .. (0,2.4) -- (0,3);
\roundNbox{unshaded}{(1.2,.8)}{.3}{0}{0}{\scriptsize{$n_b^\dag$}}; 
\filldraw[white] (.6,1.8) circle (.1cm);
\draw[thick] (.6,1.8) circle (.1cm); 
\node at (0,-.2) {\scriptsize{$F(X)$}};
\node at (1.2,3.2) {\scriptsize{$F'(X)$}};
\node at (1.2,-.2) {\scriptsize{$\psi_b$}};
\node at (0,3.2) {\scriptsize{$\varphi_a$}};
\node at (0,1.8) {\scriptsize{$\varphi_X$}};
}
\]
Thus, $n^\dag$ is a 2-modification $\psi\Rightarrow\varphi$.
Since $\dag$ preserves the norm on all 2-cells of $\cB$, we have $\|n_b\|=\|n_b^\dag\|$ for all $b\in \cB$, and thus $n^\dag$ is uniformly bounded with $\|n^\dag\|=\|n\|$.

We show $(n\xo k)^\dag = n^\dag\xo k^\dag$ and $(n\xt t)^\dag = t^\dag\xt n^\dag$, and clearly $n^{\dag\dag}=n$ by construction.
For $a\in\cA$,
\begin{align*}
(n\xo k)^\dag_a 
&= 
((n\xo k)_a)^\dag 
= 
(n_a \xxo k_a)^\dag 
= 
n_a^\dag\xxo k_a^\dag = (n^\dag)_a \xxo(k^\dag)_a = (n^\dag\xo k^\dag)_a \\
(n\xt t)^\dag_a &= ((n\xt t)_a)^\dag = (n_a \xxt t_a)^\dag = t_a^\dag \xxt n_a^\dag = (t^\dag)_a \xxt (n^\dag)_a = (t^\dag\xt n^\dag)_a.
\end{align*}
Finally, we observe that since all associators and unitors in $\cB$ are unitary, so are the associators and unitors in $\Fun^\dag(\cA\to \cB)$, as all their components are unitary by (\ref{eq:ComponentwiseAssociator},\ref{eq:ComponentwiseUnitor}).
\end{construction}

\begin{prop}
\label{prop:FunDagC*W*}
When $\cA,\cB$ are $\rm C^*/W^*$ 2-categories, so is $\Fun^\dag(\cA\to \cB)$.
\end{prop}
\begin{proof}
By Construction \ref{construction:DaggerStructureOnFunAB}, $\Fun^\dag(\cA\to \cB)$ is a $\dag$ 2-category.
Since $\Fun^\dag(\cA\to\cB)$ admits direct sums of 1-morphisms, to show $\Fun^\dag(\cA\to \cB)$ is $\rm C^*$, by Roberts' $2\times 2$ trick \cite[Lem.~2.6]{MR808930}, 
it suffices to show that for each 1-morphism/2-transformation $\varphi : F\Rightarrow G$, 
$\End_{\Fun^\dag(\cA\to\cB)}(\varphi)$
is a $\rm C^*$ algebra.
Indeed, the uniformly bounded modifications $n: \varphi\Rrightarrow \varphi$ do form a $\rm C^*$-algebra under the supreme norm:
$$
\|n^\dag\cdot n\|
=
\sup_{a\in\cA}\|(n^\dag\cdot n)_a\| 
= 
\sup_{a\in\cA}\|(n^\dag)_a\xxt n_a\| 
= 
\sup_{a\in\cA}\|(n_a)^\dag\xxt n_a\| 
= 
\sup_{a\in\cA}\|n_a\|^2 
= 
\|n\|^2.
$$

Now suppose $\cA,\cB$ are $\rm W^*$ 2-categories.
It remains to prove $\End_{\Fun^\dag(\cA\to \cB)}(\varphi)$ is a $\rm W^*$-algebra
and that 1-compositions with identity 2-transformations is a normal $\dag$ functor on hom categories.
Note that 
$$
n
=
(n_a)_{a\in\cA}
\in
\End(\varphi: F\to G)
\subset
\prod_{a\in\cA} \End(\varphi_a),
$$
where $n$ satisfies
$\varphi_X\xxt(1_{F(X)}\xxo_{F(b)} n_b)=(n_a\xxo_{G(a)}1_{G(X)})\xxt\varphi_X$, for all $X\in\cA(a\to b)$.

By either the Krein-Smulian or Kaplansky Density Theorems, 
to prove $\End_{\Fun^\dag(\cA\to \cB)}(\varphi)$ is a $\rm W^*$-algebra,
it suffices to show the unit ball in $\End(\varphi)$ is weak* closed. 
Let $(n_j=(n_a^j))$ be a weak* convergent net in the unit ball of $\End(\varphi)\subset\prod_a \End(\varphi_a)$, a $\rm W^*$-algebra.
By Lemma \ref{Lem:ProdvNa}, each component net $(n_a^j)$ converges weak* to an element $n_a$ in the unit ball of $\prod_a \End(\varphi_a)$.
We verify that $n:=(n_a)$ is a 2-modification in $\End(\varphi)$.
By the axioms of a $\rm W^*$ 2-category (see ($\rm W^*$2') in \cite[Prop.~2.4]{2105.12010}), $1_{F(X)}\xxo_{F(b)}-$, $-\xxo_{G(a)}1_{G(X)}$, $\varphi_X\xxt -$, and $-\xxt\varphi_X$ are normal operations on 2-cells in $\cB$.
We thus have
$$
\varphi_X\xxt(1_{F(X)}\xxo_{F(b)} n_b) = \lim_k\varphi_X\xxt(1_{F(X)}\xxo_{F(b)} (n_k)_b) = \lim_k((n_k)_a\xxo_{G(a)}1_{G(X)})\xxt\varphi_X= (n_a\xxo_{G(a)}1_{G(X)})\xxt\varphi_X,
$$
which implies that $n$ is a 2-modification $\varphi\Rrightarrow\varphi$.

We now show that 1-composition with an identity 2-transformation is normal.
Let $\varphi:F\Rightarrow G$; we show $1_\varphi\xo-$ is normal.
Suppose $n^j,n:\psi\Rrightarrow\gamma$ are modifications with $n^j \to n$ weak*.
Again by Lemma \ref{Lem:ProdvNa}, $n^j_a \to n_a$ weak* for all $a\in \cA$.
Since $1_{\varphi(a)}\xxo-$ is normal,
$$
(1_\varphi\xo n^j)_a=1_{\varphi(a)}\xxo n^j_a\to 1_{\varphi(a)}\xxo n_a=(1_\varphi\xo n)_a,
$$
for each $a\in\cA$, 
which implies $1_\varphi\xo n_i\to 1_\varphi\xo n$ weak* as desired.
Similarly, $-\xo 1_\varphi$ is normal. 
This completes the proof.
\end{proof}

\subsection{The 3-category of 2-categories}

It is well-known that 2-categories form a 3-category $2\Cat$, whose hom 2-categories $2\Cat(\cA\to \cB)$ are given by $\Fun(\cA\to \cB)$.
We now explain 1-composition in this 3-category following \cite[\S5.1]{MR3076451}.
We will then discuss the 3-subcategories $\rm C^*2\Cat$ and $\rm W^*2\Cat$.

\begin{construction}
\label{const:1CompositionIn2Cat}
By \cite[Prop.~5.1]{MR3076451}, given 2-categories $\cA,\cB,\cC$, there is a 2-functor
$$
\xz: 2\Cat(\cB\to \cC) \times 2\Cat(\cA\to \cB) \to 2\Cat(\cA\to \cC).
$$
The 2-functor $\xz$ is the 1-composition in $2\Cat$.
We now describe its definition on 1-morphisms, 2-morphism, and 3-morphisms in $2\Cat$.
\item[\underline{1-composition of 1-morphisms}:]
For $F\in 2\Cat(\cA\to\cB)$ and $G\in 2\Cat(\cB\to\cC)$ the 1-composite 2-functor $G\xz F\in 2\Cat(\cA\to \cC)$ is given by:
\begin{itemize}
\item $(G\xz F)(a)=G(F(a))$ for $a\in\cA$, 
$(G\xz F)(X)=G(F(X))$ for $X\in\cA(a\to b)$, 
and $(G\xz F)(f)=G(F(f))$ for $f\in\cA(X\Rightarrow Y)$.
\item $(G\xz F)^1_a:=G(F^1_a)\xxt G^1_{F(a)}\in \cC(1_{G(F(a))}\Rightarrow G(F(1_a)))$ for $a\in\cA$.
\item $(G\xz F)^2_{X,Y}:=G(F^2_{X,Y})\xxt G^2_{F(X),F(Y)}\in\cC(G(F(X))\xxo G(F(Y))\Rightarrow G(F(X\xxo Y)))$ for $X\in\cA(a\to b)$ and $Y\in\cA(b\to c)$.
\end{itemize}

\item[\underline{1-composition of 2-morphisms}:]
Suppose $F,F'\in 2\Cat(\cA\to\cB)$ and $G,G'\in 2\Cat(\cB\to\cC)$.
In the remainder of this definition, we use the following texture decorations to denote the following composite 2-functors:
$$
\tikzmath{
\filldraw[primedregion=white, rounded corners = 5pt] (0,0) rectangle (.6,.6);
\draw[thin, dotted, rounded corners = 5pt] (0,0) rectangle (.6,.6);
}
=
GF
\qquad\qquad
\tikzmath{
\filldraw[boxregion=white, rounded corners = 5pt] (0,0) rectangle (.6,.6);
\draw[thin, dotted, rounded corners = 5pt] (0,0) rectangle (.6,.6);
}
=
GF'
\qquad\qquad
\tikzmath{
\filldraw[plusregion=white, rounded corners = 5pt] (0,0) rectangle (.6,.6);
\draw[thin, dotted, rounded corners = 5pt] (0,0) rectangle (.6,.6);
}
=
G'F
\qquad\qquad
\tikzmath{
\filldraw[starregion=white, rounded corners = 5pt] (0,0) rectangle (.6,.6);
\draw[thin, dotted, rounded corners = 5pt] (0,0) rectangle (.6,.6);
}
=
G'F'.
$$
Given 2-transformations $\varphi \in 2\Cat(F\Rightarrow F')$ and $\gamma\in2\Cat(G\Rightarrow G')$,
we define $\gamma\xz F\in 2\Cat(G\xz F\Rightarrow G'\xz F)$ component-wise by
\begin{itemize}
\item
For $a\in \cA$, we define $(\gamma \xz F)_a := \gamma_{F(a)}$, and
\item
for $X\in \cA(a\to b)$, we define
$$
(\gamma\xz F)_X:=
\gamma_{F(X)}
=
\tikzmath[scale=.7, transform shape]{
\begin{scope}
\clip[rounded corners = 5] (-.6,0) rectangle (1.8,2.4);
\filldraw[primedregion=\AColor] (0,0) -- (0,.6) .. controls ++(90:.4cm) and ++(-135:.2cm) .. (.6,1.2) .. controls ++(135:.2cm) and ++(270:.4cm) .. (0,1.8) -- (0,3) -- (-.6,3) -- (-.6,0); 
\filldraw[primedregion=\BColor] (1.2,0) -- (1.2,.6) .. controls ++(90:.4cm) and ++(-45:.2cm) .. (.6,1.2) .. controls ++(-135:.2cm) and ++(90:.4cm) .. (0,.6) -- (0,0);
\filldraw[plusregion=\AColor] (0,3) -- (0,1.8) .. controls ++(270:.4cm) and ++(135:.2cm) .. (.6,1.2) .. controls ++(45:.2cm) and ++(270:.4cm) .. (1.2,1.8) -- (1.2,3);
\filldraw[plusregion=\BColor] (1.2,0) -- (1.2,.6) .. controls ++(90:.4cm) and ++(-45:.2cm) .. (.6,1.2) .. controls ++(45:.2cm) and ++(270:.4cm) .. (1.2,1.8) -- (1.2,3) -- (1.8,3) -- (1.8,0);
\end{scope}
\draw[\XColor,thick] (0,0) -- (0,.6) .. controls ++(90:.6cm) and ++(270:.6cm) .. (1.2,1.8) -- (1.2,2.4);
\draw[\gammaColor,thick] (1.2,0) -- (1.2,.6) .. controls ++(90:.6cm) and ++(270:.6cm) .. (0,1.8) -- (0,2.4);
\filldraw[white] (.6,1.2) circle (.1cm);
\draw[thick] (.6,1.2) circle (.1cm); 
\node at (0,-.2) {\scriptsize{$G(F(X))$}};
\node at (1.2,2.6) {\scriptsize{$G'(F(X))$}};
\node at (1.2,-.2) {\scriptsize{$\gamma_{F(b)}$}};
\node at (0,2.6) {\scriptsize{$\gamma_{F(a)}$}};
}
\qquad\qquad
\forall\,
X\in \cA(a\to b).
$$
\end{itemize}
Similarly, we define $G\xz\varphi\in 2\Cat(G\xz F\Rightarrow G\xz F')$ by
\begin{itemize}
\item
For $a\in \cA$, we define $(G\xz\varphi)_a := G(\varphi(a))$, and
\item
for $X\in \cA(a\to b)$, we define
$$
(G\xz\varphi)_X:=
\tikzmath[scale=.7, transform shape]{
\begin{scope}
\clip[rounded corners = 5] (-.6,0) rectangle (1.8,2.4);
\filldraw[primedregion=\AColor] (0,0) -- (0,.6) .. controls ++(90:.4cm) and ++(-135:.2cm) .. (.6,1.2) .. controls ++(135:.2cm) and ++(270:.4cm) .. (0,1.8) -- (0,3) -- (-.6,3) -- (-.6,0); 
\filldraw[primedregion=\BColor] (1.2,0) -- (1.2,.6) .. controls ++(90:.4cm) and ++(-45:.2cm) .. (.6,1.2) .. controls ++(-135:.2cm) and ++(90:.4cm) .. (0,.6) -- (0,0);
\filldraw[boxregion=\AColor] (0,3) -- (0,1.8) .. controls ++(270:.4cm) and ++(135:.2cm) .. (.6,1.2) .. controls ++(45:.2cm) and ++(270:.4cm) .. (1.2,1.8) -- (1.2,3);
\filldraw[boxregion=\BColor] (1.2,0) -- (1.2,.6) .. controls ++(90:.4cm) and ++(-45:.2cm) .. (.6,1.2) .. controls ++(45:.2cm) and ++(270:.4cm) .. (1.2,1.8) -- (1.2,3) -- (1.8,3) -- (1.8,0);
\end{scope}
\draw[\XColor,thick] (0,0) -- (0,.6) .. controls ++(90:.6cm) and ++(270:.6cm) .. (1.2,1.8) -- (1.2,2.4);
\draw[\phiColor,thick] (1.2,0) -- (1.2,.6) .. controls ++(90:.6cm) and ++(270:.6cm) .. (0,1.8) -- (0,2.4);
\filldraw[white] (.6,1.2) circle (.1cm);
\draw[thick] (.6,1.2) circle (.1cm); 
\node at (0,-.2) {\scriptsize{$G(F(X))$}};
\node at (1.2,2.6) {\scriptsize{$G(F'(X))$}};
\node at (1.2,-.2) {\scriptsize{$G(\varphi_b)$}};
\node at (0,2.6) {\scriptsize{$G(\varphi_a)$}};
}
:=
\left(
\resizebox{.35\hsize}{!}{$
\begin{aligned}
G(F(X))\xxo G(\varphi_b)
&\xrightarrow{G^2} G(F(X)\xxo\varphi_b)
\\&\xrightarrow{G(\varphi_X)} G(\varphi_a\xxo F'(X))
\\&\xrightarrow{(G^2)^\dag} G(\varphi_a)\xxo G(F(X))
\end{aligned}
$}
\right)
\qquad\qquad
\forall\,X\in \cA(a\to b).
$$
\end{itemize}

We then use the \emph{cubical convention} to define the 1-composite $\gamma\xz\varphi:=(G\xz\varphi)\xo(\gamma\xz F') \in2\Cat(G\xz F\Rightarrow G'\xz F')$,
whose components are then given by
\[
\tikzmath[scale=.7, transform shape]{
\begin{scope}
\clip[rounded corners = 5] (-.6,0) rectangle (1.8,2.4);
\filldraw[primedregion=\AColor] (0,0) -- (0,.6) .. controls ++(90:.4cm) and ++(-135:.2cm) .. (.6,1.2) .. controls ++(135:.2cm) and ++(270:.4cm) .. (0,1.8) -- (0,3) -- (-.6,3) -- (-.6,0); 
\filldraw[primedregion=\BColor] (1.2,0) -- (1.2,.6) .. controls ++(90:.4cm) and ++(-45:.2cm) .. (.6,1.2) .. controls ++(-135:.2cm) and ++(90:.4cm) .. (0,.6) -- (0,0);
\filldraw[starregion=\AColor] (0,3) -- (0,1.8) .. controls ++(270:.4cm) and ++(135:.2cm) .. (.6,1.2) .. controls ++(45:.2cm) and ++(270:.4cm) .. (1.2,1.8) -- (1.2,3);
\filldraw[starregion=\BColor] (1.2,0) -- (1.2,.6) .. controls ++(90:.4cm) and ++(-45:.2cm) .. (.6,1.2) .. controls ++(45:.2cm) and ++(270:.4cm) .. (1.2,1.8) -- (1.2,3) -- (1.8,3) -- (1.8,0);
\end{scope}
\draw[\XColor,thick] (0,0) -- (0,.6) .. controls ++(90:.6cm) and ++(270:.6cm) .. (1.2,1.8) -- (1.2,2.4);
\draw[\phiColor,thick] (1.15,0) -- (1.15,.6) .. controls ++(90:.6cm) and ++(270:.6cm) .. (-.05,1.8) -- (-.05,2.4);
\draw[\gammaColor,thick] (1.25,0) -- (1.25,.65) .. controls ++(90:.6cm) and ++(270:.6cm) .. (.05,1.85) -- (.05,2.4);
\filldraw[white] (.6,1.2) circle (.1cm);
\draw[thick] (.6,1.2) circle (.1cm); 
\node at (-.1,-.2) {\scriptsize{$G(F(X))$}};
\node at (1.3,2.6) {\scriptsize{$G'(F'(X))$}};
\node at (1.3,-.2) {\scriptsize{$(\gamma\xz\varphi)_b$}};
\node at (-.1,2.6) {\scriptsize{$(\gamma\xz\varphi)_a$}};
} 
\ :=\ 
\tikzmath[scale=.7, transform shape]{
\begin{scope}
\clip[rounded corners = 5] (-1.8,-1.8) rectangle (1.8,1.8);
\filldraw[primedregion=\AColor] (-1.2,-1.8) -- (-1.2,-1.4) .. controls ++(90:.4cm) and ++(-135:.2cm) .. (-.6,-.8) .. controls ++(135:.2cm) and ++(270:.4cm) .. (-1.2,-.2) -- (-1.2,1.8) -- (-.55,1.8) -- (-.55,2.7) -- (-1.8,2.7) -- (-1.8,-1.8);
\filldraw[primedregion=\BColor](0,-1.8) -- (0,-1.4) .. controls ++(90:.4cm) and ++(-45:.2cm) .. (-.6,-.8) .. controls ++(-135:.2cm) and ++(90:.4cm) .. (-1.2,-1.4) -- (-1.2,-1.8);
\filldraw[boxregion=\AColor] (-1.2,1.8) -- (-1.2,-.2) .. controls ++(270:.4cm) and ++(135:.2cm) .. (-.6,-.8) .. controls ++(45:.2cm) and ++(270:.4cm) .. (0,-.2) -- (0,.2) .. controls ++(90:.4cm) and ++(-135:.2cm) .. (.6,.8) .. controls ++(135:.2cm) and ++(270:.4cm) .. (0,1.4) -- (0,1.8);
\filldraw[boxregion=\BColor] (0,-1.8) -- (0,-1.4) .. controls ++(90:.4cm) and ++(-45:.2cm) .. (-.6,-.8) .. controls ++(45:.2cm) and ++(270:.4cm) .. (0,-.2) -- (0,.2) .. controls ++(90:.4cm) and ++(-135:.2cm) .. (.6,.8) .. controls ++(-45:.2cm) and ++(90:.4cm) .. (1.2,.2) -- (1.2,-1.8);
\filldraw[starregion=\AColor] (-.55,2.7) -- (-.55,1.8) -- (0,1.8) -- (0,1.4) .. controls ++(270:.4cm) and ++(135:.2cm) .. (.6,.8) .. controls ++(45:.2cm) and ++(270:.4cm) .. (1.2,1.4) -- (1.2,2.7);
\filldraw[starregion=\BColor] (1.2,-1.8) -- (1.2,.2) .. controls ++(90:.4cm) and ++(-45:.2cm) .. (.6,.8) .. controls ++(45:.2cm) and ++(270:.4cm) .. (1.2,1.4) -- (1.2,2.7) -- (1.8,2.7) -- (1.8,-1.8);
\end{scope}
\draw[\XColor,thick] (-1.2,-1.8) -- (-1.2,-1.4) .. controls ++(90:.6cm) and ++(270:.6cm) .. (0,-.2) -- (0,.2) .. controls ++(90:.6cm) and ++(270:.6cm) .. (1.2,1.4) -- (1.2,1.8);
\draw[\phiColor,thick] (0,-1.8) -- (0,-1.4) .. controls ++(90:.6cm) and ++(270:.6cm) .. (-1.2,-.2) -- (-1.2,1.8);
\draw[\gammaColor,thick] (1.2,-1.8) -- (1.2,.2) .. controls ++(90:.6cm) and ++(270:.6cm) .. (0,1.4) -- (0,1.8);
\filldraw[white] (-.6,-.8) circle (.1cm);
\draw[thick] (-.6,-.8) circle (.1cm); 
\filldraw[white] (.6,.8) circle (.1cm);
\draw[thick] (.6,.8) circle (.1cm); 
\node at (-1.2,-2) {\scriptsize{$G(F(X))$}};
\node at (0,-2) {\scriptsize{$G(\varphi_b)$}};
\node at (1.2,-2) {\scriptsize{$\gamma_{F'(b)}$}};
\node at (-1.2,2) {\scriptsize{$G(\varphi_a)$}};
\node at (0,2) {\scriptsize{$\gamma_{F'(a)}$}};
\node at (1.3,2) {\scriptsize{$G'(F'(X))$}};
}
\qquad\qquad
\begin{aligned}
\tikzmath{
\filldraw[primedregion=white, rounded corners = 5pt] (0,0) rectangle (.6,.6);
\draw[thin, dotted, rounded corners = 5pt] (0,0) rectangle (.6,.6);
}
&=
GF
\\
\tikzmath{
\filldraw[boxregion=white, rounded corners = 5pt] (0,0) rectangle (.6,.6);
\draw[thin, dotted, rounded corners = 5pt] (0,0) rectangle (.6,.6);
}
&=
GF'
\\
\tikzmath{
\filldraw[starregion=white, rounded corners = 5pt] (0,0) rectangle (.6,.6);
\draw[thin, dotted, rounded corners = 5pt] (0,0) rectangle (.6,.6);
}
&=
G'F'.
\end{aligned}
\]

\item[\underline{1-composition of 3-morphisms}:]
Suppose $F,F'\in 2\Cat(\cA\to\cB)$ and $G,G'\in 2\Cat(\cB\to\cC)$ are 2-functors,
$\varphi,\varphi'\in 2\Cat(F\Rightarrow F')$ and $\gamma,\gamma'\in 2\Cat(G\Rightarrow G')$ are 2-transformations, 
and let $n\in 2\Cat(\varphi\Rrightarrow\varphi')$ and $k\in 2\Cat(\gamma\Rrightarrow\gamma')$ be 2-modifications.
We define $k\xz n\in 2\Cat(\gamma\xz\varphi\Rrightarrow\gamma'\xz\varphi')$ component-wise at $a\in \cA$ by
$(k\xz n)_a:= G(n_a)\xxo k_{F(a)}$ as 1-composition of 2-cells in $\cC$.

\item[\underline{Interchanger:}]
For each pair of 1-composable 2-transformations $\varphi,\gamma$,
there is a distinguished invertible modification 
$\chi^{\varphi,\gamma}:(G\xz \varphi)\xo(\gamma\xz F')\Rrightarrow (\gamma\xz F)\xo(G'\xz \varphi)$ 
between the \emph{cubical} and \emph{opcubical} 1-composition conventions for 2-morphisms
called the \emph{interchanger}, 
which is defined component-wise by
\[
\chi^{\varphi,\gamma}_a:=
\tikzmath[scale=.7, transform shape]{
\begin{scope}
\clip[rounded corners = 5] (-.6,0) rectangle (1.8,2.4);
\filldraw[primedregion=\AColor] (0,0) -- (0,.6) .. controls ++(90:.4cm) and ++(-135:.2cm) .. (.6,1.2) .. controls ++(135:.2cm) and ++(270:.4cm) .. (0,1.8) -- (0,3) -- (-.6,3) -- (-.6,0); 
\filldraw[boxregion=\AColor] (1.2,0) -- (1.2,.6) .. controls ++(90:.4cm) and ++(-45:.2cm) .. (.6,1.2) .. controls ++(-135:.2cm) and ++(90:.4cm) .. (0,.6) -- (0,0);
\filldraw[plusregion=\AColor] (0,3) -- (0,1.8) .. controls ++(270:.4cm) and ++(135:.2cm) .. (.6,1.2) .. controls ++(45:.2cm) and ++(270:.4cm) .. (1.2,1.8) -- (1.2,3);
\filldraw[starregion=\AColor] (1.2,0) -- (1.2,.6) .. controls ++(90:.4cm) and ++(-45:.2cm) .. (.6,1.2) .. controls ++(45:.2cm) and ++(270:.4cm) .. (1.2,1.8) -- (1.2,3) -- (1.8,3) -- (1.8,0);
\end{scope}
\draw[\phiColor,thick] (0,0) -- (0,.6) .. controls ++(90:.6cm) and ++(270:.6cm) .. (1.2,1.8) -- (1.2,2.4);
\draw[\gammaColor,thick] (1.2,0) -- (1.2,.6) .. controls ++(90:.6cm) and ++(270:.6cm) .. (0,1.8) -- (0,2.4);
\filldraw[white] (.6,1.2) circle (.1cm);
\draw[thick] (.6,1.2) circle (.1cm); 
\node at (0,-.2) {\scriptsize{$G(\varphi_a)$}};
\node at (1.2,2.6) {\scriptsize{$G'(\varphi_a)$}};
\node at (1.2,-.2) {\scriptsize{$\gamma_{F'(a)}$}};
\node at (0,2.6) {\scriptsize{$\gamma_{F(a)}$}};
} 
=\gamma_{\varphi_a}
\qquad\qquad
\forall\,a\in \cA
\qquad\qquad
\begin{aligned}
\tikzmath{
\filldraw[primedregion=white, rounded corners = 5pt] (0,0) rectangle (.6,.6);
\draw[thin, dotted, rounded corners = 5pt] (0,0) rectangle (.6,.6);
}
&=
GF
&
\tikzmath{
\filldraw[boxregion=white, rounded corners = 5pt] (0,0) rectangle (.6,.6);
\draw[thin, dotted, rounded corners = 5pt] (0,0) rectangle (.6,.6);
}
&=
GF'
\\
\tikzmath{
\filldraw[plusregion=white, rounded corners = 5pt] (0,0) rectangle (.6,.6);
\draw[thin, dotted, rounded corners = 5pt] (0,0) rectangle (.6,.6);
}
&=
G'F
&
\tikzmath{
\filldraw[starregion=white, rounded corners = 5pt] (0,0) rectangle (.6,.6);
\draw[thin, dotted, rounded corners = 5pt] (0,0) rectangle (.6,.6);
}
&=
G'F'.
\end{aligned}
\]
(Recall here that $\varphi_a\in\cB(F(a)\to F'(a))$.)
The interchanger modification is used to prove the \emph{interchange relation} between $\xz, \xo$.
In more detail, given 
$\varphi\in 2\Cat(F\Rightarrow F')$, 
$\varphi'\in 2\Cat(F'\Rightarrow F'')$, 
$\psi\in 2\Cat(G\Rightarrow G')$, 
and $\psi'\in 2\Cat(G'\Rightarrow G'')$,
the interchanger provides an invertible modification 
$$
(\psi\xz\varphi)\xo(\psi'\xz\varphi') 
\Rrightarrow
(\psi\xo\psi')\xz(\varphi\xo\varphi').
$$
We refer the reader to \cite[p.88]{MR3076451} for more details.
\end{construction}

By \cite[p.115]{2002.06055}, 
$\xz$ is strictly associative.
That is, for 
$F\in 2\Cat(\cA\to\cB)$, $G\in2\Cat(\cB\to \cC)$ and $H\in2\Cat(\cC\to\cD)$, then $(H\xz G)\xz F=H\xz (G\xz F):\cA\to\cD$.
By \cite[Props.~5.3 and 5.5]{MR3076451},
we may choose our adjoint equivalences
$a:\xz(\xz\times \mathbf{1})\Rightarrow \xz(\mathbf{1}\times \xz)$,
$\ell:\xz(I_\cA\times \mathbf{1})\Rightarrow \mathbf{1}$,
and
$r:\xz(\mathbf{1}\times I_\cA)\Rightarrow \mathbf{1}$
to be identity transformations, whose inverses are also identity transformations.
Thus by \cite[Thm.~5.7]{MR3076451}, $2\Cat$ is a 3-category.

\begin{defn}
The 3-category $\rm C^*2\Cat$ of $\rm C^*$ 2-categories is the 3-subcategory of $2\Cat$ whose:
\begin{itemize}
\item 
objects are $\rm C^*$ 2-categories,
\item
1-morphisms are $\dag$ 2-functors,
\item
2-morphisms are $\dag$ 2-transformations
\item
3-morphisms are bounded 2-modifications
\end{itemize}
Observe that all higher coherence data in this 3-category is unitary.

The 3-category $\rm W^*2\Cat$ of $\rm W^*$ 2-categories is the full 3-subcategory of $\rm C^*2\Cat$ whose objects are $\rm W^*$ 2-categories and whose 1-morphisms are normal $\dag$ 2-functors.

Observe that $\rm C^*2\Cat$ and $\rm W^*2\Cat$ may be equipped with $\dag$-structures making them into $\dag$ 3-categories.
Indeed, all hom 2-categories are $\rm C^*/W^*$ by Proposition \ref{prop:FunDagC*W*}, 1-composition 2-functors are clearly compatible with the $\dag$-structure, and strictness of associativity of $\xz$ means all coheretors are inherently unitary.
\end{defn}

\subsection{3-endofunctors on \texorpdfstring{$2\Cat$}{2Cat}}

In this section, we give a graphical definition of a (weak) 3-endofunctor $\Phi$ on $2\Cat$.
The definition is considerably easier due to strictness of 1-composition $\xz$.
Our treatment is adapted from \cite[\S4.3]{MR3076451}.

Beyond an assignment of a $k$-morphism in $2\Cat$ for every $k$-morphism in $2\Cat$, $\Phi$ satisfies the following properties:
\begin{itemize}
\item 
$\Phi$ is a 2-functor on all hom 2-categories $2\Cat(\cA\to \cB)=\Fun(\cA\to\cB)$ in $2\Cat$.
That is, for all 
transformations 
$\varphi\in 2\Cat(F\Rightarrow F')$
and
$\psi\in  2\Cat(F'\Rightarrow F'')$
for $F,F',F'': \cA \to \cB$, 
there exist invertible modifications, 
$\Phi^\xo_{\varphi,\psi}:\Phi(\varphi)\xo\Phi(\psi)\Rrightarrow\Phi(\varphi\xo\psi)$ and $\Phi^\xo_{F}:1_{\Phi(F)}\Rrightarrow \Phi(1_{F})$,
which we represent graphically by 
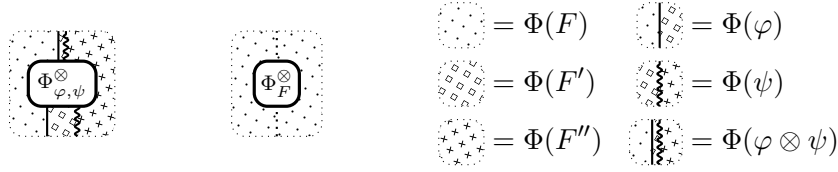
\[
\tikzmath{
\begin{scope}
\clip[rounded corners=5pt] (-.7,.3) rectangle (.7,1.7);
\filldraw[primedregion=white] (-.7,0) rectangle (.7,2);
\filldraw[boxregion=white] (.2,0) -- (.2,1) -- (.05,1) -- (.05,2) -- (-.05,2) -- (-.05,1) -- (-.2,1) -- (-.2,0);
\filldraw[plusregion=white] (.2,0) -- (.2,1) -- (.05,1) -- (.05,2) -- (.7,2) -- (.7,0);
\end{scope}
\draw[\phiColor,thick] (-.2,.3) -- (-.2,1);
\draw[\psiColor,thick] (.2,.3) -- (.2,1);
\draw[\phiColor,thick] (-.05,1) -- (-.05,1.7);
\draw[\psiColor,thick] (.05,1) -- (.05,1.7);
\roundNbox{unshaded}{(0,1)}{.3}{.15}{.15}{\scriptsize{$\Phi^\xo_{\varphi,\psi}$}};
\draw[thin, dotted, rounded corners = 5pt] (-.7,.3) rectangle (.7,1.7);
}
\qquad\qquad
\tikzmath{
\begin{scope}
\clip[rounded corners=5pt] (-.6,.3) rectangle (.6,1.7);
\filldraw[primedregion=white] (-.6,0) rectangle (.6,2);
\end{scope}
\draw[thick,dotted] (0,.3) -- (0,1.7);
\roundNbox{unshaded}{(0,1)}{.3}{0}{0}{\scriptsize{$\Phi^\xo_F$}};
\draw[thin, dotted, rounded corners = 5pt] (-.6,.3) rectangle (.6,1.7);
}
\qquad\qquad
\begin{aligned}
\tikzmath{
\filldraw[primedregion=white, rounded corners = 5pt] (0,0) rectangle (.6,.6);
\draw[thin, dotted, rounded corners = 5pt] (0,0) rectangle (.6,.6);
}
&=
\Phi(F)
&
\tikzmath{
\begin{scope}
\filldraw[primedregion=white, rounded corners = 5pt] (0,0) rectangle (.3,.6);
\filldraw[boxregion=white, rounded corners = 5pt] (.3,0) rectangle (.6,.6);
\end{scope}
\draw[\phiColor,thick] (.3,0) -- (.3,.6);
\draw[thin, dotted, rounded corners = 5pt] (0,0) rectangle (.6,.6);
}
&=
\Phi(\varphi)
\\
\tikzmath{
\filldraw[boxregion=white, rounded corners = 5pt] (0,0) rectangle (.6,.6);
\draw[thin, dotted, rounded corners = 5pt] (0,0) rectangle (.6,.6);
}
&=
\Phi(F')
&
\tikzmath{
\begin{scope}
\filldraw[boxregion=white, rounded corners = 5pt] (0,0) rectangle (.3,.6);
\filldraw[plusregion=white, rounded corners = 5pt] (.3,0) rectangle (.6,.6);
\end{scope}
\draw[\psiColor,thick] (.3,0) -- (.3,.6);
\draw[thin, dotted, rounded corners = 5pt] (0,0) rectangle (.6,.6);
}
&=
\Phi(\psi)
\\
\tikzmath{
\filldraw[plusregion=white, rounded corners = 5pt] (0,0) rectangle (.6,.6);
\draw[thin, dotted, rounded corners = 5pt] (0,0) rectangle (.6,.6);
}
&=
\Phi(F'')
&
\tikzmath{
\begin{scope}
\filldraw[primedregion=white, rounded corners = 5pt] (-.35,0) rectangle (-.05,.6);
\filldraw[boxregion=white, rounded corners = 5pt] (-.05,0) rectangle (.05,.6);
\filldraw[plusregion=white, rounded corners = 5pt] (.05,0) rectangle (.35,.6);
\end{scope}
\draw[\phiColor,thick] (-.05,0) -- (-.05,.6);
\draw[\psiColor,thick] (.05,0) -- (.05,.6);
\draw[thin, dotted, rounded corners = 5pt] (-.35,0) rectangle (.35,.6);
}
&=
\Phi(\varphi\xo\psi)
\end{aligned}
\]
These modifications are subject to the usual associativity and unitality coherence axioms:
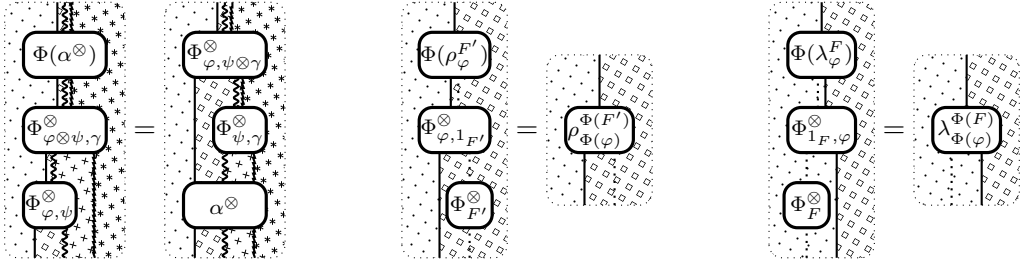
\[
\tikzmath{
\begin{scope}
\clip[rounded corners=5pt] (-.8,.3) rectangle (.8,3.7);
\filldraw[primedregion=white] (-.4,0) -- (-.4,1) -- (-.25,1) -- (-.25,2) -- (-.1,2) -- (-.1,4) -- (-.8,4) -- (-.8,0);
\filldraw[boxregion=white] (0,0) -- (0,1) -- (-.15,1) -- (-.15,2) -- (0,2) -- (0,4) -- (-.1,4) -- (-.1,2) -- (-.25,2) -- (-.25,1) -- (-.4,1) -- (-.4,0);
\filldraw[plusregion=white] (.4,0) -- (.4,2) -- (.1,2) -- (.1,4) -- (0,4) -- (0,2) -- (-.15,2) -- (-.15,1) -- (0,1) -- (0,0);
\filldraw[starregion=white] (.4,0) -- (.4,2) -- (.1,2) -- (.1,4) -- (.8,4) -- (.8,0);
\end{scope}
\draw[\phiColor,thick] (-.4,.3) -- (-.4,1);
\draw[\psiColor,thick] (0,.3) -- (0,1);
\draw[\gammaColor,thick] (.4,.3) -- (.4,2);
\draw[\phiColor,thick] (-.25,1) -- (-.25,2);
\draw[\psiColor,thick] (-.15,1) -- (-.15,2);
\draw[\phiColor,thick] (-.1,2) -- (-.1,3.7);
\draw[\psiColor,thick] (0,2) -- (0,3.7);
\draw[\gammaColor,thick] (.1,2) -- (.1,3.7);
\roundNbox{unshaded}{(-.2,1)}{.3}{.05}{.05}{\scriptsize{$\Phi^\xo_{\varphi,\psi}$}};
\roundNbox{unshaded}{(0,2)}{.3}{.25}{.25}{\scriptsize{$\Phi^\xo_{\varphi\xo \psi,\gamma}$}};
\roundNbox{unshaded}{(0,3)}{.3}{.25}{.25}{\scriptsize{$\Phi(\alpha^{\xo})$}};
\draw[thin, dotted, rounded corners = 5pt] (-.8,.3) rectangle (.8,3.7);
}
=
\tikzmath{
\begin{scope}
\clip[rounded corners=5pt] (-.8,.3) rectangle (.8,3.7);
\filldraw[primedregion=white] (-.4,0) -- (-.4,3) -- (-.1,3) -- (-.1,4) -- (-.8,4) -- (-.8,0);
\filldraw[boxregion=white] (-.4,0) -- (-.4,3) -- (-.1,3) -- (-.1,4) -- (0,4) -- (0,3) -- (.15,3) -- (.15,2) -- (0,2) -- (0,0);
\filldraw[plusregion=white] (.4,0) -- (.4,2) -- (.25,2) -- (.25,3) -- (.1,3) -- (.1,4) -- (0,4) -- (0,3) -- (.15,3) -- (.15,2) -- (0,2) -- (0,0);
\filldraw[starregion=white] (.4,0) -- (.4,2) -- (.25,2) -- (.25,3) -- (.1,3) -- (.1,4) -- (.8,4) -- (.8,0);
\end{scope}
\draw[\phiColor,thick] (-.4,.3) -- (-.4,3);
\draw[\psiColor,thick] (0,.3) -- (0,2);
\draw[\gammaColor,thick] (.4,.3) -- (.4,2);
\draw[\psiColor,thick] (.15,2) -- (.15,3);
\draw[\gammaColor,thick] (.25,2) -- (.25,3);
\draw[\phiColor,thick] (-.1,3) -- (-.1,3.7);
\draw[\psiColor,thick] (0,3) -- (0,3.7);
\draw[\gammaColor,thick] (.1,3) -- (.1,3.7);
\roundNbox{unshaded}{(0,1)}{.3}{.25}{.25}{\scriptsize{$\alpha^\xo$}};
\roundNbox{unshaded}{(.2,2)}{.3}{.05}{.05}{\scriptsize{$\Phi^\xo_{\psi,\gamma}$}};
\roundNbox{unshaded}{(-.1,3)}{.3}{.15}{.35}{\scriptsize{$\Phi^\xo_{\varphi,\psi\xo \gamma}$}};
\draw[thin, dotted, rounded corners = 5pt] (-.8,.3) rectangle (.8,3.7);
}
\qquad\qquad
\tikzmath{
\begin{scope}
\clip[rounded corners=5pt] (-.7,.3) rectangle (.7,3.7);
\filldraw[primedregion=white] (-.2,0) -- (-.2,2) -- (-.05,2) -- (-.05,3) -- (0,3) -- (0,4) -- (-.7,4) -- (-.7,0);
\filldraw[boxregion=white] (-.2,0) -- (-.2,2) -- (-.05,2) -- (-.05,3) -- (0,3) -- (0,4) -- (.7,4) -- (.7,0);
\end{scope}
\draw[\phiColor,thick] (-.2,.3) -- (-.2,2);
\draw[thick,dotted] (.2,.3) -- (.2,2);
\draw[\phiColor,thick] (-.05,2) -- (-.05,3);
\draw[thick,dotted] (.05,2) -- (.05,3);
\draw[\phiColor,thick] (0,3) -- (0,3.7);
\roundNbox{unshaded}{(.2,1)}{.3}{0}{0}{\scriptsize{$\Phi^\xo_{F'}$}};
\roundNbox{unshaded}{(0,2)}{.3}{.17}{.17}{\scriptsize{$\Phi^\xo_{\varphi,1_{F'}}$}};
\roundNbox{unshaded}{(0,3)}{.3}{.17}{.17}{\scriptsize{$\Phi(\rho_\varphi^{F'})$}};
\draw[thin, dotted, rounded corners = 5pt] (-.7,.3) rectangle (.7,3.7);
}
=
\tikzmath{
\begin{scope}
\clip[rounded corners=5pt] (-.7,0) rectangle (.7,2);
\filldraw[primedregion=white] (-.2,0) -- (-.2,1) -- (0,1) -- (0,2) -- (-.7,2) -- (-.7,0);
\filldraw[boxregion=white] (-.2,0) -- (-.2,1) -- (0,1) -- (0,2) -- (.7,2) -- (.7,0);
\end{scope}
\draw[\phiColor,thick] (-.2,0) -- (-.2,1);
\draw[thick,dotted] (.2,0) -- (.2,1);
\draw[\phiColor,thick] (0,1) -- (0,2);
\roundNbox{unshaded}{(0,1)}{.3}{.15}{.15}{\scriptsize{$\rho_{\Phi(\varphi)}^{\Phi(F')}$}};
\draw[thin, dotted, rounded corners = 5pt] (-.7,0) rectangle (.7,2);
}
\qquad\qquad
\tikzmath{
\begin{scope}
\clip[rounded corners=5pt] (-.7,.3) rectangle (.7,3.7);
\filldraw[primedregion=white] (.2,0) -- (.2,2) -- (.05,2) -- (.05,3) -- (0,3) -- (0,4) -- (-.7,4) -- (-.7,0);
\filldraw[boxregion=white] (.2,0) -- (.2,2) -- (.05,2) -- (.05,3) -- (0,3) -- (0,4) -- (.7,4) -- (.7,0);
\end{scope}
\draw[\phiColor,thick] (.2,.3) -- (.2,2);
\draw[thick,dotted] (-.2,.3) -- (-.2,2);
\draw[\phiColor,thick] (.05,2) -- (.05,3);
\draw[thick,dotted] (-.05,2) -- (-.05,3);
\draw[\phiColor,thick] (0,3) -- (0,3.7);
\roundNbox{unshaded}{(-.2,1)}{.3}{0}{0}{\scriptsize{$\Phi^\xo_F$}};
\roundNbox{unshaded}{(0,2)}{.3}{.15}{.15}{\scriptsize{$\Phi^\xo_{1_F,\varphi}$}};
\roundNbox{unshaded}{(0,3)}{.3}{.15}{.15}{\scriptsize{$\Phi(\lambda_\varphi^F)$}};
\draw[thin, dotted, rounded corners = 5pt] (-.7,.3) rectangle (.7,3.7);
}
=
\tikzmath{
\begin{scope}
\clip[rounded corners=5pt] (-.7,0) rectangle (.7,2);
\filldraw[primedregion=white] (.2,0) -- (.2,1) -- (0,1) -- (0,2) -- (-.7,2) -- (-.7,0);
\filldraw[boxregion=white] (.2,0) -- (.2,1) -- (0,1) -- (0,2) -- (.7,2) -- (.7,0);
\end{scope}
\draw[\phiColor,thick] (.2,0) -- (.2,1);
\draw[thick,dotted] (-.2,0) -- (-.2,1);
\draw[\phiColor,thick] (0,1) -- (0,2);
\roundNbox{unshaded}{(0,1)}{.3}{.15}{.15}{\scriptsize{$\lambda_{\Phi(\varphi)}^{\Phi(F)}$}};
\draw[thin, dotted, rounded corners = 5pt] (-.7,0) rectangle (.7,2);
}
\]
\item
We have 1-compositor adjoint equivalence transformations 
$\Phi^\xz_{G,F}:\Phi(G)\xz \Phi(F)\Rightarrow \Phi(G\xz F)$ 
for all $F\in 2\Cat(\cA\to\cB)$ and $G\in 2\Cat(\cB\to\cC)$
and 
$\Phi^\xz_{\cA}:1_{\Phi(\cA)}\Rightarrow \Phi(1_{\cA})$
for all $\cA\in 2\Cat$.
These transformations come equipped with an invertible associator modification
$\omega^{\xz}_{H,G,F}$:
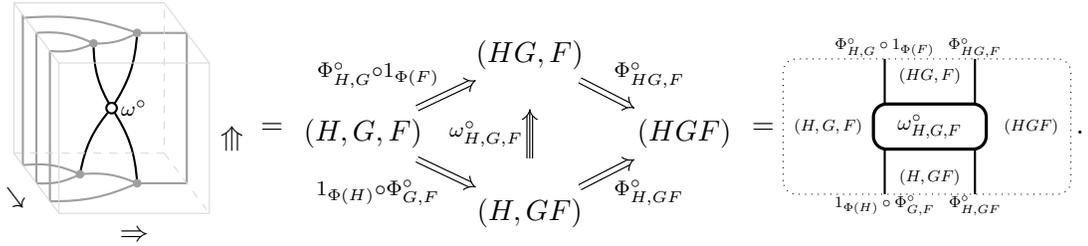
\[
\hspace*{.5cm}
\tikzmath{
\draw[\CColor,thick] (-.12,.16) -- (-.12,2.16);
\draw[\CColor,thick] (-.3,.4) -- (-.3,2.4);
\draw[\CColor,thick] (-.48,.64) -- (-.48,2.64);
\draw[\CColor,thick] (1.7,.4) -- (1.7,2.4);
\draw[\CColor,thick] (-.48,.64) to[bend left=15] (.28,.52);
\draw[\CColor,thick] (-.3,.4) to[bend right=15] (.28,.52);
\draw[\CColor,thick] (.28,.52) to[bend left=14] (1.04,.4);
\draw[\CColor,thick] (-.12,.16) to[bend right=12] (1.04,.4);
\draw[\CColor,thick] (1.04,.4) -- (1.7,.4);
\draw[\CColor,thick] (-.12,2.16) to[bend right=15] (.46,2.26);
\draw[\CColor,thick] (-.3,2.4) to[bend left=15] (.46,2.26);
\draw[\CColor,thick] (.46,2.26) to[bend right=15] (1.04,2.4);
\draw[\CColor,thick] (-.48,2.6) to[bend left=12] (1.04,2.4);
\draw[\CColor,thick] (1.04,2.4) -- (1.7,2.4);
\draw[thick] (.28,.52) to[bend left=15] (.7,1.4);
\draw[thick] (1.04,.4) to[bend right=15] (.7,1.4);
\draw[thick] (.46,2.26) to[bend right=15] (.7,1.4);
\draw[thick] (1.04,2.4) to[bend left=15] (.7,1.4);
\filldraw[white] (.7,1.4) circle (.07cm);
\draw[thick] (.7,1.4) circle (.07cm); %
\node at (1,1.4) {\scriptsize{$\omega^\circ$}};
\filldraw[\CColor] (.28,.52) circle (.05cm);
\filldraw[\CColor] (1.04,.4) circle (.05cm);
\filldraw[\CColor] (.46,2.26) circle (.05cm);
\filldraw[\CColor] (1.04,2.4) circle (.05cm);
\node at (-.6,.2) [rotate=-53] {$\to$};
\node at (1,-.3) {$\Rightarrow$};
\node at (2.3,1) [rotate=90] {$\Rrightarrow$};
\draw[\AColor] (0,0) rectangle (2,2);
\draw[\AColor] (0,0) -- (-.6,.8) -- (-.6,2.8) -- (1.4,2.8) -- (2,2);
\draw[\AColor] (-.6,2.8) -- (0,2);
\draw[\AColor,dashed] (-.6,.8) -- (1.4,.8) -- (1.4,2.8);
\draw[\AColor,dashed] (1.4,.8) -- (2,0);
}
=
\begin{tikzcd}[row sep=1em, column sep=1em]
&
(HG,F)
\arrow[dr,Rightarrow,"\scriptstyle\Phi^\xz_{HG,F}"]
\\
(H,G,F)
\arrow[ur,Rightarrow,"\scriptstyle\Phi^\xz_{H,G}\xz 1_{\Phi(F)}"]
\arrow[dr,Rightarrow,swap,"\scriptstyle1_{\Phi(H)}\xz \Phi^\xz_{G,F}"]
&&
(HGF)
\\
&
(H,GF)
\arrow[ur,Rightarrow,swap,"\scriptstyle\Phi^\xz_{H,GF}"]
\arrow[uu,triplecd,shorten <= 1em, shorten >= 1em,"\scriptstyle\omega^\xz_{H,G,F}"]
\end{tikzcd}
=
\tikzmath[scale=.75, transform shape]{
\draw[thick,\phiColor] (-.8,-1.2) -- (-.8,1.2);
\draw[thick,\phiColor] (.8,-1.2) -- (.8,1.2);
\roundNbox{unshaded}{(0,0)}{.4}{.6}{.6}{\normalsize{$\omega^\xz_{H,G,F}$}};
\node at (-1.8,0) {\scriptsize{$(H,G,F)$}};
\node at (1.8,0) {\scriptsize{$(HGF)$}};
\node at (0,-.9) {\scriptsize{$(H,GF)$}};
\node at (0,.9) {\scriptsize{$(HG,F)$}};
\node at (-.8,-1.4) {\scriptsize{$1_{\Phi(H)}\xz \Phi^\xz_{G,F}$}};
\node at (-.8,1.4) {\scriptsize{$\Phi^\xz_{H,G}\xz 1_{\Phi(F)}$}};
\node at (.8,-1.4) {\scriptsize{$\Phi^\xz_{H,GF}$}};
\node at (.8,1.4) {\scriptsize{$\Phi^\xz_{HG,F}$}};
\draw[thin, dotted, rounded corners = 5pt] (-2.6,-1.2) rectangle (2.5,1.2);
}\,.
\]
Here, we use the abbreviated notation
$(GF):=\Phi(G\xz F)$
and
$(G,F):=\Phi(G)\xz\Phi(F)$,
so that 
$(K,HG,F):=\Phi(K)\xz\Phi(H\xz G)\xz\Phi(F)$ and $\Phi^\xz_{H,GF}:=\Phi^\xz_{H,G\xz F}:(H,GF)\Rightarrow (HGF)$.
The associator $\omega^\xz$ satisfies the coherence axiom
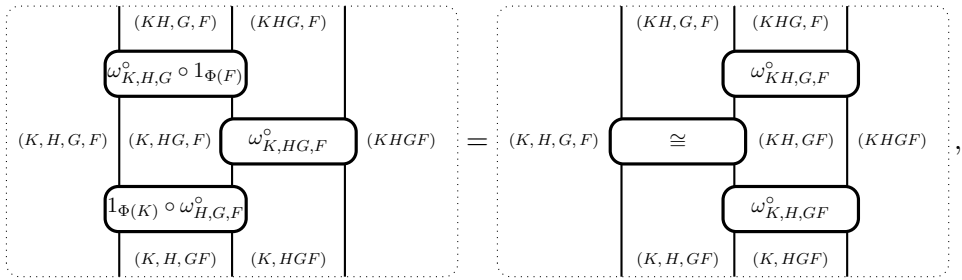
\[
\tikzmath[scale=.75, transform shape]{
\draw[thick,\phiColor] (-2,-2.4) -- (-2,2.4);
\draw[thick,\phiColor] (0,-2.4) -- (0,2.4);
\draw[thick,\phiColor] (2,-2.4) -- (2,2.4);
\roundNbox{unshaded}{(-1,-1.2)}{.4}{.85}{.85}{\normalsize{$1_{\Phi(K)}\xz\omega^\xz_{H,G,F}$}};
\roundNbox{unshaded}{(-1,1.2)}{.4}{.85}{.85}{\normalsize{$\omega^\xz_{K,H,G}\xz 1_{\Phi(F)}$}};
\roundNbox{unshaded}{(1,0)}{.4}{.8}{.8}{\normalsize{$\omega^\xz_{K,HG,F}$}};
\node at (-3,0) {\scriptsize{$(K,H,G,F)$}};
\node at (-1.1,0) {\scriptsize{$(K,HG,F)$}};
\node at (-1,-2.1) {\scriptsize{$(K,H,GF)$}};
\node at (-1,2.1) {\scriptsize{$(KH,G,F)$}};
\node at (1,-2.1) {\scriptsize{$(K,HGF)$}};
\node at (1,2.1) {\scriptsize{$(KHG,F)$}};
\node at (3,0) {\scriptsize{$(KHGF)$}};
\draw[thin, dotted, rounded corners = 5pt] (-4,-2.4) rectangle (4,2.4);
}
=
\tikzmath[scale=.75, transform shape]{
\draw[thick,\phiColor] (-2,-2.4) -- (-2,2.4);
\draw[thick,\phiColor] (0,-2.4) -- (0,2.4);
\draw[thick,\phiColor] (2,-2.4) -- (2,2.4);
\roundNbox{unshaded}{(1,-1.2)}{.4}{.8}{.8}{\normalsize{$\omega^\xz_{K,H,GF}$}};
\roundNbox{unshaded}{(1,1.2)}{.4}{.8}{.8}{\normalsize{$\omega^\xz_{KH,G,F}$}};
\roundNbox{unshaded}{(-1,0)}{.4}{.8}{.8}{\normalsize{$\cong$}};
\node at (-3.2,0) {\scriptsize{$(K,H,G,F)$}};
\node at (1.1,0) {\scriptsize{$(KH,GF)$}};
\node at (-1,-2.1) {\scriptsize{$(K,H,GF)$}};
\node at (-1,2.1) {\scriptsize{$(KH,G,F)$}};
\node at (1,-2.1) {\scriptsize{$(K,HGF)$}};
\node at (1,2.1) {\scriptsize{$(KHG,F)$}};
\node at (2.8,0) {\scriptsize{$(KHGF)$}};
\draw[thin, dotted, rounded corners = 5pt] (-4.2,-2.4) rectangle (3.8,2.4);
}\,,
\]
where the isomorphism on the left of the right hand side is the interchanger from Construction \ref{const:1CompositionIn2Cat}.

Finally, we have invertible unitor modifications $\ell^\xz_F$ and $r^\xz_F$:
\[
\tikzmath{
\draw[\CColor,thick] (-.4,.533) rectangle (1.6,2.533);
\draw[\CColor,thick,dashed] (-.2,2.267) to[bend right=20] (.6,2.533);
\draw[\CColor,thick,dashed] (-.2,2.267) -- (-.2,.267);
\draw[\CColor,thick,dashed] (-.2,.267) -- (.4,.267);
\draw[\CColor,thick,dashed] (.6,2.533) to[bend right=10] (.4,.267);
\filldraw[\CColor] (.6,2.533) circle (.05cm);
\node at (-.6,.2) [rotate=-53] {$\to$};
\node at (1,-.3) {$\Rightarrow$};
\node at (2.3,1) [rotate=90] {$\Rrightarrow$};
\draw[\AColor] (0,0) rectangle (2,2);
\draw[\AColor] (0,0) -- (-.6,.8) -- (-.6,2.8) -- (1.4,2.8) -- (2,2);
\draw[\AColor] (-.6,2.8) -- (0,2);
\draw[\AColor,dashed] (-.6,.8) -- (1.4,.8) -- (1.4,2.8);
\draw[\AColor,dashed] (1.4,.8) -- (2,0);
}
=
\begin{tikzcd}[row sep=1em, column sep=1em]
&
{}
\\
(1_\cB,F)
\arrow[rr,bend left=60,Rightarrow,"\scriptstyle\Phi^\xz_{1_\cB,F}"]
\arrow[rr,bend right=60,Rightarrow,swap,"\scriptstyle\Phi^\xz_\cB\xz 1_{\Phi(F)}"]
&&
(F)
\\
&
{}
\arrow[uu,triplecd,shorten <= 1em, shorten >= 1em,"\scriptstyle\ell^\xz_F"]
\end{tikzcd}
=
\tikzmath{
\draw[thick] (0,-.8) -- (0,.8);
\roundNbox{fill=white}{(0,0)}{.3}{.1}{.1}{\scriptsize{$\ell^\xz_F$}};
\node at (-1,0) {\scriptsize{$(1_\cB,F)$}};
\node at (1,0) {\scriptsize{$(F)$}};
\node at (0,1) {\scriptsize{$\Phi^\xz_{1_\cB,F}$}};
\node at (0,-1) {\scriptsize{$\Phi^\circ_\cB\xz 1_{\Phi(F)}$}};
\draw[thin, dotted, rounded corners = 5pt] (-1.6,-.8) rectangle (1.4,.8);
}
\]
\[
\tikzmath{
\draw[\CColor,thick] (-.2,.267) rectangle (1.8,2.267);
\draw[\CColor,thick,dashed] (-.4,2.533) to[bend left=15] (.8,2.267);
\draw[\CColor,thick,dashed] (-.4,.533) -- (-.4,2.533);
\draw[\CColor,thick,dashed] (-.4,.533) -- (.2,.533);
\draw[\CColor,thick,dashed] (.8,2.267) to[bend right=10] (.2,.533);
\filldraw[\CColor] (.8,2.267) circle (.05cm);
\node at (-.6,.2) [rotate=-53] {$\to$};
\node at (1,-.3) {$\Rightarrow$};
\node at (2.3,1) [rotate=90] {$\Rrightarrow$};
\draw[\AColor] (0,0) rectangle (2,2);
\draw[\AColor] (0,0) -- (-.6,.8) -- (-.6,2.8) -- (1.4,2.8) -- (2,2);
\draw[\AColor] (-.6,2.8) -- (0,2);
\draw[\AColor,dashed] (-.6,.8) -- (1.4,.8) -- (1.4,2.8);
\draw[\AColor,dashed] (1.4,.8) -- (2,0);
}
=
\begin{tikzcd}[row sep=1em, column sep=1em]
&
{}
\\
(G,1_\cB)
\arrow[rr,bend left=60,Rightarrow,"\scriptstyle\Phi^\xz_{G,1_\cB}"]
\arrow[rr,bend right=60,Rightarrow,swap,"\scriptstyle1_{\Phi(G)}\xz\Phi^\xz_\cB"]
&&
(G)
\\
&
{}
\arrow[uu,triplecd,shorten <= 1em, shorten >= 1em,"\scriptstyle r^\xz_G"]
\end{tikzcd}
=
\tikzmath{
\draw[thick] (0,-.8) -- (0,.8);
\roundNbox{fill=white}{(0,0)}{.3}{.1}{.1}{\scriptsize{$r^\xz_G$}};
\node at (-1,0) {\scriptsize{$(G,1_\cB)$}};
\node at (1,0) {\scriptsize{$(G)$}};
\node at (0,1) {\scriptsize{$\Phi^\xz_{G,1_\cB}$}};
\node at (0,-1) {\scriptsize{$1_{\Phi(G)}\xz \Phi^\circ_\cB$}};
\draw[thin, dotted, rounded corners = 5pt] (-1.6,-.8) rectangle (1.4,.8);
}\,.
\]
These unitors satisfy the coherence axiom
\[
\hspace*{1.3cm}
\tikzmath{
\draw[\CColor,thick] (1.04,.4) -- (1.7,.4) -- (1.7,2.4) -- (1.04,2.4);
\draw[\CColor,thick] (-.48,2.6) -- (1.04,2.4);
\draw[\CColor,thick] (-.48,.6) -- (-.48,2.6);
\draw[\CColor,thick] (-.48,.6) -- (1.04,.4);
\draw[\CColor,thick] (-.12,2.2) -- (1.04,2.4);
\draw[\CColor,thick] (-.12,.2) -- (-.12,2.2);
\draw[\CColor,thick] (-.12,.2) -- (1.04,.4);
\draw[\CColor,dashed,thick] (-.3,.4) -- (-.3,2.4);
\draw[\CColor,dashed,thick] (-.3,2.4) to[bend left=12] (.46,2.3);
\draw[\CColor,dashed,thick] (-.3,.4) to[bend right=17] (.28,.5);
\draw[thick] (.28,.5) to[bend left=15] (.7,1.4);
\draw[thick] (1.04,.4) to[bend right=15] (.7,1.4);
\draw[thick] (.46,2.3) to[bend right=15] (.7,1.4);
\draw[thick] (1.04,2.4) to[bend left=15] (.7,1.4);
\filldraw[white] (.7,1.4) circle (.07cm);
\draw[thick] (.7,1.4) circle (.07cm); %
\filldraw[\CColor] (.28,.5) circle (.05cm);
\filldraw[\CColor] (1.04,.4) circle (.05cm);
\filldraw[\CColor] (.46,2.3) circle (.05cm);
\filldraw[\CColor] (1.04,2.4) circle (.05cm);
\draw[\AColor] (0,0) rectangle (2,2);
\draw[\AColor] (0,0) -- (-.6,.8) -- (-.6,2.8) -- (1.4,2.8) -- (2,2);
\draw[\AColor] (-.6,2.8) -- (0,2);
\draw[\AColor,dashed] (-.6,.8) -- (1.4,.8) -- (1.4,2.8);
\draw[\AColor,dashed] (1.4,.8) -- (2,0);
}
=
\tikzmath{
\draw[\CColor,thick] (1.04,.4) rectangle (1.7,2.4);
\draw[\CColor,thick] (-.48,2.6) -- (1.04,2.4);
\draw[\CColor,thick] (-.48,.6) -- (-.48,2.6);
\draw[\CColor,thick] (-.48,.6) -- (1.04,.4);
\draw[\CColor,thick] (-.12,2.2) -- (1.04,2.4);
\draw[\CColor,thick] (-.12,.2) -- (-.12,2.2);
\draw[\CColor,thick] (-.12,.2) -- (1.04,.4);
\draw[\CColor,dashed,thick] (-.3,.4) -- (-.3,2.4);
\draw[\CColor,dashed,thick] (-.3,2.4) to[bend left=12] (.46,2.3);
\draw[\CColor,dashed,thick] (-.3,.4) to[bend right=17] (.28,.5); 
\draw[\CColor,dashed,thick] (.46,2.3) to[bend right=10] (.28,.5);
\draw[\CColor,dashed,thick] (-.3,1.4) -- (.3,1.4);
\filldraw[\CColor] (.28,.5) circle (.05cm);
\filldraw[\CColor] (1.04,.4) circle (.05cm);
\filldraw[\CColor] (.46,2.3) circle (.05cm);
\filldraw[\CColor] (1.04,2.4) circle (.05cm);
\draw[\AColor] (0,0) rectangle (2,2);
\draw[\AColor] (0,0) -- (-.6,.8) -- (-.6,2.8) -- (1.4,2.8) -- (2,2);
\draw[\AColor] (-.6,2.8) -- (0,2);
\draw[\AColor,dashed] (-.6,.8) -- (1.4,.8) -- (1.4,2.8);
\draw[\AColor,dashed] (1.4,.8) -- (2,0);
}
\quad \leftrightarrow\quad
\tikzmath[scale=.75, transform shape]{
\draw[thick,\phiColor] (-.8,-1.2) -- (-.8,1.2);
\draw[thick,\phiColor] (.8,-1.2) -- (.8,1.2);
\roundNbox{unshaded}{(0,0)}{.4}{.6}{.6}{\normalsize{$\omega^\xz_{G,1_\cB,F}$}};
\node at (-1.8,0) {\scriptsize{$(G,1_{\cB},F)$}};
\node at (1.6,0) {\scriptsize{$(GF)$}};
\node at (0,.9) {\scriptsize{$(G1_{\cB},F)$}};
\node at (0,-.9) {\scriptsize{$(G,1_{\cB}F)$}};
\node at (-.8,1.4) {\scriptsize{$\Phi^\xz_{G,1_\cB}\xz 1_{\Phi(F)}$}};
\node at (-.8,-1.4) {\scriptsize{$1_{\Phi(G)}\xz\Phi^\xz_{1_\cB,F}$}};
\node at (.8,-1.4) {\scriptsize{$\Phi^\xz_{G,F}$}};
\node at (.8,1.4) {\scriptsize{$\Phi^\xz_{G,F}$}};
\draw[thin, dotted, rounded corners = 5pt] (-2.7,-1.2) rectangle (2.2,1.2);
}
=
\tikzmath[scale=.75, transform shape]{
\draw[thick,\phiColor] (-.8,-1.9) -- (-.8,1.5);
\draw[thick,\phiColor] (.8,-1.9) -- (.8,1.5);
\roundNbox{unshaded}{(-.8,.7)}{.4}{.6}{.6}{\normalsize{$r^\xz_G \xz 1_{\Phi(F)}$}};
\roundNbox{unshaded}{(-.8,-.7)}{.4}{.6}{.6}{\scriptsize{$1_{\Phi(G)}\xz (\ell^\xz_F)^{-1}$}};
\node at (-1.7,-1.5) {\scriptsize{$(G,1_{\cB},F)$}};
\node at (0,-1.5) {\scriptsize{$(G,F)$}};
\node at (1.5,-1.5) {\scriptsize{$(GF)$}};
\node at (-1.3,0) {\scriptsize{$1_{\Phi(G)}\xz \Phi^{\xz}_{\cB}\xz 1_{\Phi(F)}$}};
\node at (-.8,-2.1) {\scriptsize{$1_{\Phi(G)}\xz\Phi^\xz_{1_\cB,F}$}};
\node at (-.8,1.7) {\scriptsize{$\Phi^\xz_{G,1_\cB}\xz 1_{\Phi(F)}$}};
\node at (.8,-2.1) {\scriptsize{$\Phi^\xz_{G,F}$}};
\draw[thin, dotted, rounded corners = 5pt] (-2.7,-1.9) rectangle (2.2,1.5);
}\,.
\]
Here, we note that $F\xz 1_{\cA}=F=1_{\cB}\xz F$, so $(G1_{\cB}F)=(GF)$, $(G1_{\cB},F)=(G,F) =(G,1_{\cB}F)$ and $(G,F)=\Phi(G)\xz 1_{\Phi(\cB)}\xz \Phi(F)$.
\end{itemize}

Given a weak 3-functor $\Phi$ on $2\Cat$ which preserves the 3-subcategories $\rm C^*2\Cat$ and $\rm W^*2\Cat$, we can ask whether $\Phi$ restricts to a $\dag$ 3-functor.
This consists of the following conditions:
\begin{itemize}
\item 
$\Phi(n^\dag)=\Phi(n)^\dag$ for all bounded 2-modifications $n$,
\item
the coheretors 
$\Phi^\xo_{\varphi,\psi}$
and
$\Phi^\xo_{F}$
are unitary,
\item 
$\Phi^\xz_{G,F}$ and $\Phi^\xz_{\cA}$ are unitary adjoint equivalences, and
\item
the associators $\omega^\xz_{H,G,F}$
and unitors $\ell^\xz_{F},r^\xz_F$ are unitary.
\end{itemize}

\section{Q-system completion is a 3-functor}

In this section, we rapidly recall the definition of Q-system completion for a $\rm C^*/W^*$ 2-category from \cite[\S3]{2105.12010}, and we prove Theorem \ref{thm:QSys3Functor} that Q-system completion is a 3-functor.

\subsection{Graphical calculus for Q-systems and their bimodules}

Q-systems were first defined in \cite{MR1257245}, and were subsequently studied in \cite{MR1444286,MR2298822,MR3308880}.
For this section, we fix a $\rm C^*/W^*$ 2-category $\cC$ which we assume is locally unitarily Cauchy complete, i.e., every hom 1-category has orthogonal direct sums and all orthogonal projections split orthogonally.

\begin{defn}
A \emph{Q-system} in $\cC$ consists of a triple $(Q,m,i)$
where $Q\in \cC(b\to b)$, $m\in \cC(Q\xxo Q \Rightarrow Q)$, and $i\in \cC(1_b \Rightarrow Q)$, which satisfy certain axioms.
We represent $b,Q,m,i$ and the adjoints $m^\dag,i^\dag$ graphically as follows:
$$
\tikzmath{\filldraw[\BColor, rounded corners=5, very thin, baseline=1cm] (0,0) rectangle (.6,.6);}=b
\qquad
\tikzmath{
\fill[\BColor, rounded corners=5pt ] (0,0) rectangle (.6,.6);
\draw[\QsColor,thick] (.3,0) -- (.3,.6);
}={}_bQ_b.
\qquad
\tikzmath{
\fill[\BColor, rounded corners=5pt] (-.3,0) rectangle (.9,.6);
\draw[\QsColor,thick] (0,0) arc (180:0:.3cm);
\draw[\QsColor,thick] (.3,.3) -- (.3,.6);
\filldraw[\QsColor] (.3,.3) circle (.05cm);
}=m
\qquad
\tikzmath{
\fill[\BColor, rounded corners=5pt] (-.3,0) rectangle (.9,-.6);
\draw[\QsColor,thick] (0,0) arc (-180:0:.3cm);
\draw[\QsColor,thick] (.3,-.3) -- (.3,-.6);
\filldraw[\QsColor] (.3,-.3) circle (.05cm);
}=m^\dag
\qquad
\tikzmath{
\fill[\BColor, rounded corners=5pt] (0,0) rectangle (.6,.6);
\draw[\QsColor,thick] (.3,.3) -- (.3,.6);
\filldraw[\QsColor] (.3,.3) circle (.05cm);
}=i
\qquad
\tikzmath{
\fill[\BColor, rounded corners=5pt] (0,0) rectangle (.6,-.6);
\draw[\QsColor,thick] (.3,-.3) -- (.3,-.6);
\filldraw[\QsColor] (.3,-.3) circle (.05cm);
}=i^\dag.
$$
The Q-system axioms are as follows:
\begin{enumerate}[label=(Q\arabic*)]
\item 
\label{Q:associativity}
(associativity)
$\tikzmath{
\fill[\BColor, rounded corners=5pt] (-.3,-.3) rectangle (1.2,.6);
\draw[\QsColor,thick] (0,-.3) -- (0,0) arc (180:0:.3cm);
\draw[\QsColor,thick] (.3,-.3) arc (180:0:.3cm);
\draw[\QsColor,thick] (.3,.3) -- (.3,.6);
\filldraw[\QsColor] (.3,.3) circle (.05cm);
\filldraw[\QsColor] (.6,0) circle (.05cm);
}
=
\tikzmath{
\fill[\BColor, rounded corners=5pt] (-.6,-.3) rectangle (.9,.6);
\draw[\QsColor,thick] (0,0) arc (180:0:.3cm) -- (.6,-.3);
\draw[\QsColor,thick] (-.3,-.3) arc (180:0:.3cm);
\draw[\QsColor,thick] (.3,.3) -- (.3,.6);
\filldraw[\QsColor] (.3,.3) circle (.05cm);
\filldraw[\QsColor] (0,0) circle (.05cm);
}$
\item
\label{Q:unitality}
(unitality)
$\tikzmath{
\fill[\BColor, rounded corners=5pt] (-.3,-.3) rectangle (.9,.6);
\draw[\QsColor,thick] (0,-.1) -- (0,0) arc (180:0:.3cm) -- (.6,-.3);
\draw[\QsColor,thick] (.3,.3) -- (.3,.6);
\filldraw[\QsColor] (.3,.3) circle (.05cm);
\filldraw[\QsColor] (0,-.1) circle (.05cm);
}
=
\tikzmath{
\fill[\BColor, rounded corners=5pt ] (0,-.3) rectangle (.6,.6);
\draw[\QsColor,thick] (.3,-.3) -- (.3,.6);
}
=
\tikzmath{
\fill[\BColor, rounded corners=5pt] (-.3,-.3) rectangle (.9,.6);
\draw[\QsColor,thick] (0,-.3) -- (0,0) arc (180:0:.3cm) -- (.6,-.1);
\draw[\QsColor,thick] (.3,.3) -- (.3,.6);
\filldraw[\QsColor] (.3,.3) circle (.05cm);
\filldraw[\QsColor] (.6,-.1) circle (.05cm);
}$
\item
\label{Q:Frobenius}
(Frobenius)
$
\tikzmath{
\fill[\BColor, rounded corners=5pt] (-.3,-.6) rectangle (1.5,.6);
\draw[\QsColor,thick] (0,-.6) -- (0,0) arc (180:0:.3cm) arc (-180:0:.3cm) -- (1.2,.6);
\draw[\QsColor,thick] (.3,.3) -- (.3,.6);
\draw[\QsColor,thick] (.9,-.3) -- (.9,-.6);
\filldraw[\QsColor] (.3,.3) circle (.05cm);
\filldraw[\QsColor] (.9,-.3) circle (.05cm);
}
=
\tikzmath{
\fill[\BColor, rounded corners=5pt] (-.3,0) rectangle (.9,1.2);
\draw[\QsColor,thick] (0,0) arc (180:0:.3cm);
\draw[\QsColor,thick] (0,1.2) arc (-180:0:.3cm);
\draw[\QsColor,thick] (.3,.3) -- (.3,.9);
\filldraw[\QsColor] (.3,.3) circle (.05cm);
\filldraw[\QsColor] (.3,.9) circle (.05cm);
}
=
\tikzmath{
\fill[\BColor, rounded corners=5pt] (-.3,.6) rectangle (1.5,-.6);
\draw[\QsColor,thick] (0,.6) -- (0,0) arc (-180:0:.3cm) arc (180:0:.3cm) -- (1.2,-.6);
\draw[\QsColor,thick] (.3,-.3) -- (.3,-.6);
\draw[\QsColor,thick] (.9,.3) -- (.9,.6);
\filldraw[\QsColor] (.3,-.3) circle (.05cm);
\filldraw[\QsColor] (.9,.3) circle (.05cm);
}
$
\item
\label{Q:separable}
(separable)
$
\tikzmath{
\fill[\BColor, rounded corners=5pt] (-.3,0) rectangle (.9,1.2);
\draw[\QsColor,thick] (0,.6) arc (180:-180:.3cm);
\draw[\QsColor,thick] (.3,1.2) -- (.3,.9);
\draw[\QsColor,thick] (.3,0) -- (.3,.3);
\filldraw[\QsColor] (.3,.3) circle (.05cm);
\filldraw[\QsColor] (.3,.9) circle (.05cm);
}
=
\tikzmath{
\fill[\BColor, rounded corners=5pt ] (0,0) rectangle (.6,1.2);
\draw[\QsColor,thick] (.3,0) -- (.3,1.2);
}
$
\end{enumerate}
We refer the reader to 
\cite[Prop.~5.17]{MR2298822}
or
\cite[Facts~3.4]{2105.12010}
for various dependencies amongst these axioms.
\end{defn}

\begin{defn}
Suppose $P\in \cC(a\to a)$ and $Q\in \cC(b\to b)$ are Q-systems.
A $P-Q$ bimodule is a triple $(X,\lambda_X, \rho_X)$
consisting of
$X\in \cC(a\to b)$,
$\lambda_X \in \cC( P\xxo X \Rightarrow X)$,
and $\rho_X\in \cC(X\xxo Q\Rightarrow X)$,
again satisfying certain properties.
We represent $a,b,X,P,Q$ graphically by 
$$
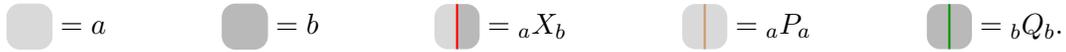
\tikzmath{\filldraw[\AColor, rounded corners=5, very thin, baseline=1cm] (0,0) rectangle (.6,.6);}=a
\qquad\qquad
\tikzmath{\filldraw[\BColor, rounded corners=5, very thin, baseline=1cm] (0,0) rectangle (.6,.6);}=b
\qquad\qquad
\tikzmath{
\begin{scope}
\clip[rounded corners=5pt] (-.3,0) rectangle (.3,.6);
\fill[\AColor] (0,0) rectangle (-.3,.6);
\fill[\BColor] (0,0) rectangle (.3,.6);
\end{scope}
\draw[thick, \XColor] (0,0) -- (0,.6);
}={}_aX_b
\qquad\qquad
\tikzmath{
\fill[\AColor, rounded corners=5pt ] (0,0) rectangle (.6,.6);
\draw[\PsColor,thick] (.3,0) -- (.3,.6);
}={}_aP_a
\qquad\qquad
\tikzmath{
\fill[\BColor, rounded corners=5pt ] (0,0) rectangle (.6,.6);
\draw[\QsColor,thick] (.3,0) -- (.3,.6);
}={}_bQ_b.
$$
We denote $\lambda_X,\rho_X$ and $\lambda_X^\dag,\rho_X^\dag$ by trivalent vertices:
$$
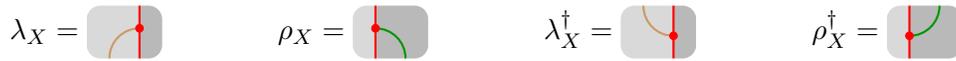
\lambda_X =
\tikzmath{
\begin{scope}
\clip[rounded corners = 5pt] (-.7,-.2) rectangle (.3,.5);
\filldraw[\AColor] (-.7,-.2) rectangle (0,.5);
\filldraw[\BColor] (0,-.2) rectangle (.3,.5);
\end{scope}
\draw[\XColor,thick] (0,-.2) -- (0,.5);
\draw[\PsColor,thick] (-.4,-.2) arc (180:90:.4cm);
\filldraw[\XColor] (0,.2) circle (.05cm);
}
\qquad\qquad
\rho_X =
\tikzmath{
\begin{scope}
\clip[rounded corners = 5pt] (-.3,-.2) rectangle (.7,.5);
\filldraw[\AColor] (-.3,-.2) rectangle (0,.5);
\filldraw[\BColor] (0,-.2) rectangle (.7,.5);
\end{scope}
\draw[\XColor,thick] (0,-.2) -- (0,.5);
\draw[\QsColor,thick] (.4,-.2) arc (0:90:.4cm);
\filldraw[\XColor] (0,.2) circle (.05cm);
}
\qquad\qquad
\lambda_X^\dag =
\tikzmath{
\begin{scope}
\clip[rounded corners = 5pt] (-.7,-.5) rectangle (.3,.2);
\filldraw[\AColor] (-.7,-.5) rectangle (0,.2);
\filldraw[\BColor] (0,-.5) rectangle (.3,.2);
\end{scope}
\draw[\XColor,thick] (0,-.5) -- (0,.2);
\draw[\PsColor,thick] (-.4,.2) arc (180:270:.4cm);
\filldraw[\XColor] (0,-.2) circle (.05cm);
}
\qquad\qquad
\rho_X^\dag =
\tikzmath{
\begin{scope}
\clip[rounded corners = 5pt] (-.3,-.5) rectangle (.7,.2);
\filldraw[\AColor] (-.3,-.5) rectangle (0,.2);
\filldraw[\BColor] (0,-.5) rectangle (.7,.2);
\end{scope}
\draw[\XColor,thick] (0,-.5) -- (0,.2);
\draw[\QsColor,thick] (.4,.2) arc (0:-90:.4cm);
\filldraw[\XColor] (0,-.2) circle (.05cm);
}
$$
The bimodule axioms are as follows:
\begin{enumerate}[label=(B\arabic*)]
\item 
\label{M:associativity}
(associativity)
$
\tikzmath{
\begin{scope}
\clip[rounded corners = 5pt] (-.9,-.6) rectangle (.3,.5);
\filldraw[\AColor] (-.9,-.6) rectangle (0,.5);
\filldraw[\BColor] (0,-.6) rectangle (.3,.5);
\end{scope}
\draw[\XColor,thick] (0,-.6) -- (0,.5);
\draw[\PsColor,thick] (-.6,-.6) -- (-.6,-.4) arc (180:90:.6cm);
\draw[\PsColor,thick] (-.3,-.6) -- (-.3,-.4) arc (180:90:.3cm);
\filldraw[\XColor] (0,.2) circle (.05cm);
\filldraw[\XColor] (0,-.1) circle (.05cm);
}
=
\tikzmath{
\begin{scope}
\clip[rounded corners = 5pt] (-.9,-.6) rectangle (.3,.5);
\filldraw[\AColor] (-.9,-.6) rectangle (0,.5);
\filldraw[\BColor] (0,-.6) rectangle (.3,.5);
\end{scope}
\draw[\XColor,thick] (0,-.6) -- (0,.5);
\draw[\PsColor,thick] (-.4,-.2) arc (180:90:.4cm);
\draw[\PsColor,thick] (-.6,-.6) -- (-.6,-.4)  arc (180:0:.2cm) -- (-.2,-.6);
\filldraw[\XColor] (0,.2) circle (.05cm);
\filldraw[\PsColor] (-.4,-.2) circle (.05cm);
}
$,
$
\tikzmath{
\begin{scope}
\clip[rounded corners = 5pt] (-.3,-.6) rectangle (.9,.5);
\filldraw[\AColor] (-.3,-.6) rectangle (0,.5);
\filldraw[\BColor] (0,-.6) rectangle (.9,.5);
\end{scope}
\draw[\XColor,thick] (0,-.6) -- (0,.5);
\draw[\QsColor,thick] (.6,-.6) -- (.6,-.4) arc (0:90:.6cm);
\draw[\QsColor,thick] (.3,-.6) -- (.3,-.4) arc (0:90:.3cm);
\filldraw[\XColor] (0,.2) circle (.05cm);
\filldraw[\XColor] (0,-.1) circle (.05cm);
}
=
\tikzmath{
\begin{scope}
\clip[rounded corners = 5pt] (-.3,-.6) rectangle (.9,.5);
\filldraw[\AColor] (-.3,-.6) rectangle (0,.5);
\filldraw[\BColor] (0,-.6) rectangle (.9,.5);
\end{scope}
\draw[\XColor,thick] (0,-.6) -- (0,.5);
\draw[\QsColor,thick] (.4,-.2) arc (0:90:.4cm);
\draw[\QsColor,thick] (.6,-.6) -- (.6,-.4)  arc (0:180:.2cm) -- (.2,-.6);
\filldraw[\XColor] (0,.2) circle (.05cm);
\filldraw[\QsColor] (.4,-.2) circle (.05cm);
}
$, and
$
\tikzmath{
\begin{scope}
\clip[rounded corners = 5pt] (-.7,-.5) rectangle (.7,.5);
\filldraw[\AColor] (-.7,-.5) rectangle (0,.5);
\filldraw[\BColor] (0,-.5) rectangle (.7,.5);
\end{scope}
\draw[\XColor,thick] (0,-.5) -- (0,.5);
\draw[\PsColor,thick] (-.4,-.5) -- (-.4,-.2) arc (180:90:.4cm);
\draw[\QsColor,thick] (.4,-.5) arc (0:90:.4cm);
\filldraw[\XColor] (0,.2) circle (.05cm);
\filldraw[\XColor] (0,-.1) circle (.05cm);
}
=
\tikzmath{
\begin{scope}
\clip[rounded corners = 5pt] (-.7,-.5) rectangle (.7,.5);
\filldraw[\AColor] (-.7,-.5) rectangle (0,.5);
\filldraw[\BColor] (0,-.5) rectangle (.7,.5);
\end{scope}
\draw[\XColor,thick] (0,-.5) -- (0,.5);
\draw[\PsColor,thick] (-.4,-.5) arc (180:90:.4cm);
\draw[\QsColor,thick] (.4,-.5) -- (.4,-.2) arc (0:90:.4cm);
\filldraw[\XColor] (0,.2) circle (.05cm);
\filldraw[\XColor] (0,-.1) circle (.05cm);
}
$
\item
\label{M:unitality}
(unitality)
$
\tikzmath{
\begin{scope}
\clip[rounded corners = 5pt] (-.7,-.5) rectangle (.3,.5);
\filldraw[\AColor] (-.7,-.5) rectangle (0,.5);
\filldraw[\BColor] (0,-.5) rectangle (.3,.5);
\end{scope}
\draw[\XColor,thick] (0,-.5) -- (0,.5);
\draw[\PsColor,thick] (-.4,-.2) arc (180:90:.4cm);
\filldraw[\XColor] (0,.2) circle (.05cm);
\filldraw[\PsColor] (-.4,-.2) circle (.05cm);
}
=
\tikzmath{
\begin{scope}
\clip[rounded corners = 5pt] (-.3,-.5) rectangle (.3,.5);
\filldraw[\AColor] (-.3,-.5) rectangle (0,.5);
\filldraw[\BColor] (0,-.5) rectangle (.3,.5);
\end{scope}
\draw[\XColor,thick] (0,-.5) -- (0,.5);
}
$ 
=
$
\tikzmath{
\begin{scope}
\clip[rounded corners = 5pt] (-.3,-.5) rectangle (.7,.5);
\filldraw[\AColor] (-.3,-.5) rectangle (0,.5);
\filldraw[\BColor] (0,-.5) rectangle (.7,.5);
\end{scope}
\draw[\XColor,thick] (0,-.5) -- (0,.5);
\draw[\QsColor,thick] (.4,-.2) arc (0:90:.4cm);
\filldraw[\XColor] (0,.2) circle (.05cm);
\filldraw[\QsColor] (.4,-.2) circle (.05cm);
}
$
\item
\label{M:Frobenius}
(Frobenius)
$
\tikzmath{
\begin{scope}
\clip[rounded corners = 5pt] (-1.3,-.5) rectangle (.3,.8);
\filldraw[\AColor] (-1.3,-.5) rectangle (0,.8);
\filldraw[\BColor] (0,-.5) rectangle (.3,.8);
\end{scope}
\draw[\XColor,thick] (0,-.5) -- (0,.8);
\draw[\PsColor,thick] (-1,-.5) -- (-1,.2) arc (180:0:.3cm) arc (180:270:.4cm);
\draw[\PsColor,thick] (-.7,.5) -- (-.7,.8);
\filldraw[\XColor] (0,-.2) circle (.05cm);
\filldraw[\PsColor] (-.7,.5) circle (.05cm);
}
=
\tikzmath{
\begin{scope}
\clip[rounded corners = 5pt] (-.7,-.2) rectangle (.3,1.1);
\filldraw[\AColor] (-.7,-.2) rectangle (0,1.1);
\filldraw[\BColor] (0,-.2) rectangle (.3,1.1);
\end{scope}
\draw[\XColor,thick] (0,-.2) -- (0,1.1);
\draw[\PsColor,thick] (-.4,-.2) arc (180:90:.4cm);
\draw[\PsColor,thick] (-.4,1.1) arc (180:270:.4cm);
\filldraw[\XColor] (0,.2) circle (.05cm);
\filldraw[\XColor] (0,.7) circle (.05cm);
}
=
\tikzmath{
\begin{scope}
\clip[rounded corners = 5pt] (-1.3,-.8) rectangle (.3,.5);
\filldraw[\AColor] (-1.3,-.8) rectangle (0,.5);
\filldraw[\BColor] (0,-.8) rectangle (.3,.5);
\end{scope}
\draw[\XColor,thick] (0,-.8) -- (0,.5);
\draw[\PsColor,thick] (-1,.5) -- (-1,-.2) arc (-180:0:.3cm) arc (180:90:.4cm);
\draw[\PsColor,thick] (-.7,-.5) -- (-.7,-.8);
\filldraw[\XColor] (0,.2) circle (.05cm);
\filldraw[\PsColor] (-.7,-.5) circle (.05cm);
}
$ and $
\tikzmath{
\begin{scope}
\clip[rounded corners = 5pt] (-.3,-.5) rectangle (1.3,.8);
\filldraw[\AColor] (-.3,-.5) rectangle (0,.8);
\filldraw[\BColor] (0,-.5) rectangle (1.3,.8);
\end{scope}
\draw[\XColor,thick] (0,-.5) -- (0,.8);
\draw[\QsColor,thick] (1,-.5) -- (1,.2) arc (0:180:.3cm) arc (0:-90:.4cm);
\draw[\QsColor,thick] (.7,.5) -- (.7,.8);
\filldraw[\XColor] (0,-.2) circle (.05cm);
\filldraw[\QsColor] (.7,.5) circle (.05cm);
}
=
\tikzmath{
\begin{scope}
\clip[rounded corners = 5pt] (-.3,-.2) rectangle (.7,1.1);
\filldraw[\AColor] (-.3,-.2) rectangle (0,1.1);
\filldraw[\BColor] (0,-.2) rectangle (.7,1.1);
\end{scope}
\draw[\XColor,thick] (0,-.2) -- (0,1.1);
\draw[\QsColor,thick] (.4,-.2) arc (0:90:.4cm);
\draw[\QsColor,thick] (.4,1.1) arc (0:-90:.4cm);
\filldraw[\XColor] (0,.2) circle (.05cm);
\filldraw[\XColor] (0,.7) circle (.05cm);
}
=
\tikzmath{
\begin{scope}
\clip[rounded corners = 5pt] (-.3,-.8) rectangle (1.3,.5);
\filldraw[\AColor] (-.3,-.8) rectangle (0,.5);
\filldraw[\BColor] (0,-.8) rectangle (1.3,.5);
\end{scope}
\draw[\XColor,thick] (0,-.8) -- (0,.5);
\draw[\QsColor,thick] (1,.5) -- (1,-.2) arc (0:-180:.3cm) arc (0:90:.4cm);
\draw[\QsColor,thick] (.7,-.5) -- (.7,-.8);
\filldraw[\XColor] (0,.2) circle (.05cm);
\filldraw[\QsColor] (.7,-.5) circle (.05cm);
}
$
\item
\label{M:separable}
(separable)
$
\tikzmath{
\begin{scope}
\clip[rounded corners = 5pt] (-.5,-.5) rectangle (.3,.5);
\filldraw[\AColor] (-.6,-.5) rectangle (0,.5);
\filldraw[\BColor] (0,-.5) rectangle (.3,.5);
\end{scope}
\draw[\XColor,thick] (0,-.5) -- (0,.5);
\draw[\PsColor,thick] (0,-.3) arc (270:90:.3cm);
\filldraw[\XColor] (0,.3) circle (.05cm);
\filldraw[\XColor] (0,-.3) circle (.05cm);
}
=
\tikzmath{
\begin{scope}
\clip[rounded corners = 5pt] (-.3,-.5) rectangle (.3,.5);
\filldraw[\AColor] (-.3,-.5) rectangle (0,.5);
\filldraw[\BColor] (0,-.5) rectangle (.3,.5);
\end{scope}
\draw[\XColor,thick] (0,-.5) -- (0,.5);
}
=
\tikzmath{
\begin{scope}
\clip[rounded corners = 5pt] (-.3,-.5) rectangle (.5,.5);
\filldraw[\AColor] (-.3,-.5) rectangle (0,.5);
\filldraw[\BColor] (0,-.5) rectangle (.6,.5);
\end{scope}
\draw[\XColor,thick] (0,-.5) -- (0,.5);
\draw[\QsColor,thick] (0,-.3) arc (-90:90:.3cm);
\filldraw[\XColor] (0,.3) circle (.05cm);
\filldraw[\XColor] (0,-.3) circle (.05cm);
}
$
\end{enumerate}
We refer the reader to 
\cite[Facts~3.16]{2105.12010}
for various dependencies amongst these axioms.
\end{defn}

\begin{defn}
For $\cC$ a $\rm C^*/\rm W^*$ 2-category,
its \emph{Q-system completion} is the $\rm C^*/\rm W^*$ 2-category $\QSys(\cC)$ whose:
\begin{itemize}
\item 
0-cells are Q-systems $(Q,m,i)\in \cC(b\to b)$,

\item
1-cells between Q-systems $P\in \cC(a\to a)$ and $Q\in \cC(b\to b)$ are (unital Frobenius) bimodules $({}_aX_b,\lambda_X,\rho_X)\in \cC(a\to b)$, and

\item
2-cells are bimodule intertwiners, i.e., given Q-systems
${}_aP_a,{}_bQ_b$ and $P-Q$ bimodules ${}_aX_b, {}_aY_b$,
$\QSys(\cC)({}_PX_Q\Rightarrow {}_PY_Q)$ is the set of $f\in \cC({}_aX_b\Rightarrow {}_aY_b)$ such that
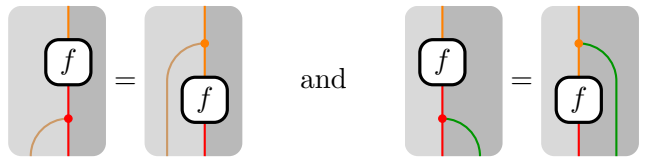
\begin{equation*}
\tikzmath{
\begin{scope}
\clip[rounded corners=5pt] (-.8,.5) rectangle (.5,-1.5);
\fill[\AColor] (-.8,.5) rectangle (0,-1.5);
\fill[\BColor] (0,.5) rectangle (.5,-1.5);
\end{scope}
\draw[\YColor,thick] (0,.5) -- (0,-.25);
\draw[\XColor,thick] (0,-.25) -- (0,-1.5);
\draw[\PsColor,thick] (-.5,-1.5) arc (180:90:.5cm);
\filldraw[\XColor] (0,-1) circle (.05cm);
\roundNbox{unshaded}{(0,-.25)}{.3}{0}{0}{$f$};
}
=
\tikzmath{
\begin{scope}
\clip[rounded corners=5pt] (-.8,-.5) rectangle (.5,-2.5);
\fill[\AColor] (-.8,-.5) rectangle (0,-2.5);
\fill[\BColor] (0,-.5) rectangle (.5,-2.5);
\end{scope}
\draw[\XColor,thick] (0,-1.75) -- (0,-2.5);
\draw[\YColor,thick] (0,-0.5) -- (0,-1.75);
\draw[\PsColor,thick] (-.5,-2.5) -- (-.5,-1.5) arc (180:90:.5cm);
\filldraw[\YColor] (0,-1) circle (.05cm);
\roundNbox{unshaded}{(0,-1.75)}{.3}{0}{0}{$f$};
}
\qquad\text{and}\qquad
\tikzmath{
\begin{scope}
\clip[rounded corners=5pt] (.8,.5) rectangle (-.5,-1.5);
\fill[\BColor] (.8,.5) rectangle (0,-1.5);
\fill[\AColor] (0,.5) rectangle (-.5,-1.5);
\end{scope}
\draw[\YColor,thick] (0,.5) -- (0,-.25);
\draw[\XColor,thick] (0,-.25) -- (0,-1.5);
\draw[\QsColor,thick] (.5,-1.5) arc (0:90:.5cm);
\filldraw[\XColor] (0,-1) circle (.05cm);
\roundNbox{unshaded}{(0,-.25)}{.3}{0}{0}{$f$};
}
=
\tikzmath{
\begin{scope}
\clip[rounded corners=5pt] (.8,-.5) rectangle (-.5,-2.5);
\fill[\BColor] (.8,-.5) rectangle (0,-2.5);
\fill[\AColor] (0,-.5) rectangle (-.5,-2.5);
\end{scope}
\draw[\XColor,thick] (0,-1.75) -- (0,-2.5);
\draw[\YColor,thick] (0,-0.5) -- (0,-1.75);
\draw[\QsColor,thick] (.5,-2.5) -- (.5,-1.5) arc (0:90:.5cm);
\filldraw[\YColor] (0,-1) circle (.05cm);
\roundNbox{unshaded}{(0,-1.75)}{.3}{0}{0}{$f$};
}\,.
\end{equation*}
\item
1-composition in $\QSys(\cC)$ is performed by orthogonally splitting the \emph{separability projector}
\begin{equation}
\label{eq:SeparabilityProjector}
p^Q_{X,Y}
:=
\tikzmath{
\begin{scope}
\clip[rounded corners = 5pt] (-.5,-.5) rectangle (.5,.5);
\filldraw[\AColor] (-.5,-.5) rectangle (-.2,.5);
\filldraw[\BColor] (-.2,-.5) rectangle (.2,.5);
\filldraw[\CColor] (.2,-.5) rectangle (.5,.5);
\end{scope}
\draw[thick, \XColor] (-.2,-.5) -- (-.2,.5);
\draw[thick,\QsColor] (-.2,0) -- (.2,0);
\draw[thick, \YColor] (.2,-.5) -- (.2,.5);
}
:=
\tikzmath{
\begin{scope}
\clip[rounded corners = 5pt] (-.6,-.6) rectangle (.6,.6);
\filldraw[\AColor] (-.6,-.6) rectangle (-.3,.6);
\filldraw[\BColor] (-.3,-.6) rectangle (.3,.6);
\filldraw[\CColor] (.3,-.6) rectangle (.6,.6);
\end{scope}
\draw[thick,\QsColor] (-.3,-.3) arc (-90:0:.3cm) arc (180:90:.3cm);
\draw[thick, \XColor] (-.3,-.6) -- (-.3,.6);
\draw[thick, \YColor] (.3,-.6) -- (.3,.6);
\filldraw[\XColor] (-.3,-.3) circle (.05cm);
\filldraw[\YColor] (.3,.3) circle (.05cm);
}
=
\tikzmath{
\begin{scope}
\clip[rounded corners = 5pt] (-.7,-.7) rectangle (.7,.5);
\filldraw[\AColor] (-.7,-.7) rectangle (-.4,.5);
\filldraw[\BColor] (-.4,-.7) rectangle (.4,.5);
\filldraw[\CColor] (.4,-.7) rectangle (.7,.5);
\end{scope}
\draw[thick,\QsColor] (-.4,.2) arc (90:0:.2cm) arc (-180:0:.2cm) arc (180:90:.2cm);
\draw[thick,\QsColor] (0,-.4) -- (0,-.2);
\draw[thick, \XColor] (-.4,-.7) -- (-.4,.5);
\draw[thick, \YColor] (.4,-.7) -- (.4,.5);
\filldraw[\XColor] (-.4,.2) circle (.05cm);
\filldraw[\YColor] (.4,.2) circle (.05cm);
\filldraw[\QsColor] (0,-.2) circle (.05cm);
\filldraw[\QsColor] (0,-.4) circle (.05cm);
}
\end{equation}
The object ${}_aX\xxo_Q Y_b \in \QSys(\cC)(P \to R)$ and a $P-R$ bimodular coisometry $u_{X,Y}^Q: X\xxo_b Y \to X\xxo_Q Y$, unique up to canonical unitary, such that $p_{X,Y}^Q=(u_{X,Y}^Q)^\dag\xxt u_{X,Y}^Q$. 
\end{itemize}
We refer the reader to \cite[\S3.2]{2105.12010} for the full details that $\QSys(\cC)$ is a $\dag$ 2-category, which is $\rm C^*/W^*$ whenever $\cC$ is respectively.
\end{defn}

\begin{nota}
\label{nota:QSys(C)andOther}
We use the graphical notation for $\QSys(\cC)$ from \cite[\S3.3]{2105.12010}, where shaded regions for Q-systems are denoted by colored regions, but trivial Q-systems are still represented in gray-scale:
$$
\tikzmath{\filldraw[\PrColor, rounded corners=5, very thin, baseline=1cm] (0,0) rectangle (.6,.6);}=P
\qquad\qquad
\tikzmath{\filldraw[\QrColor, rounded corners=5, very thin, baseline=1cm] (0,0) rectangle (.6,.6);}=Q
\qquad\qquad
\tikzmath{\filldraw[\AColor, rounded corners=5, very thin, baseline=1cm] (0,0) rectangle (.6,.6);}=1_a
\qquad\qquad
\tikzmath{\filldraw[\BColor, rounded corners=5, very thin, baseline=1cm] (0,0) rectangle (.6,.6);}=1_b.
$$
If ${}_aP_a, {}_bQ_b\in \QSys(\cC)$ are Q-systems and $X\in \QSys(\cC)(P\to Q)$,
then $X$ may be also viewed as a $1_a-Q$, $P-1_b$, and a $1_a-1_b$ bimodule;
we represent these four possibilities by varying the shadings:
$$
\tikzmath{
\begin{scope}
\clip[rounded corners=5pt] (-.3,0) rectangle (.3,.6);
\fill[\PrColor] (0,0) rectangle (-.3,.6);
\fill[\QrColor] (0,0) rectangle (.3,.6);
\end{scope}
\draw[thick, \XColor] (0,0) -- (0,.6);
}={}_PX_Q
\qquad\qquad
\tikzmath{
\begin{scope}
\clip[rounded corners=5pt] (-.3,0) rectangle (.3,.6);
\fill[\AColor] (0,0) rectangle (-.3,.6);
\fill[\QrColor] (0,0) rectangle (.3,.6);
\end{scope}
\draw[thick, \XColor] (0,0) -- (0,.6);
}={}_{1_a}X_Q
\qquad\qquad
\tikzmath{
\begin{scope}
\clip[rounded corners=5pt] (-.3,0) rectangle (.3,.6);
\fill[\PrColor] (0,0) rectangle (-.3,.6);
\fill[\BColor] (0,0) rectangle (.3,.6);
\end{scope}
\draw[thick, \XColor] (0,0) -- (0,.6);
}={}_PX_{1_b}
\qquad\qquad
\tikzmath{
\begin{scope}
\clip[rounded corners=5pt] (-.3,0) rectangle (.3,.6);
\fill[\AColor] (0,0) rectangle (-.3,.6);
\fill[\BColor] (0,0) rectangle (.3,.6);
\end{scope}
\draw[thick, \XColor] (0,0) -- (0,.6);
}={}_{1_a}X_{1_b}.
$$
We use a similar convention for intertwiners of bimodules.
We often suppress the external shading when drawing 2-cells in $\QSys(\cC)$; when we do so, it should be inferred that the diagram/relation depicted holds for any consistent external shading applied to the diagram(s).

Given $X\in \QSys(\cC)(P\to Q)$ and $Y\in \QSys(\cC)(Q\to R)$,
we denote the coisometry $u_{X,Y}^Q$ and its adjoint in the graphical calculus of $\QSys(\cC)$ by
\[
u^Q_{X,Y}
:=
\tikzmath{
\begin{scope}
\clip[rounded corners = 5pt] (-.5,-.5) rectangle (.5,.5);
\filldraw[\QrColor] (-.2,0) rectangle (.2,.5);
\filldraw[\BColor] (-.2,-.5) rectangle (.2,0);
\end{scope}
\draw[thick, \XColor] (-.2,-.5) -- (-.2,.5);
\draw[thick, \QsColor] (-.2,0) -- (.2,0);
\draw[thick, \YColor] (.2,-.5) -- (.2,.5);
}
: X\xxo_b Y \to X\xxo_Q Y
\qquad\text{and}\qquad
(u^Q_{X,Y})^\dag
=
\tikzmath{
\begin{scope}
\clip[rounded corners = 5pt] (-.5,-.5) rectangle (.5,.5);
\filldraw[\QrColor] (-.2,0) rectangle (.2,-.5);
\filldraw[\BColor] (-.2,.5) rectangle (.2,0);
\end{scope}
\draw[thick, \XColor] (-.2,-.5) -- (-.2,.5);
\draw[thick, \QsColor] (-.2,0) -- (.2,0);
\draw[thick, \YColor] (.2,-.5) -- (.2,.5);
}\,.
\]
We thus get the following relations:
\[
u^Q_{X,Y}\xxt (u^Q_{X,Y})^\dag 
=
\tikzmath{
\begin{scope}
\clip[rounded corners = 5pt] (-.5,-.5) rectangle (.5,1);
\filldraw[\QrColor] (-.2,0) rectangle (.2,-.5);
\filldraw[\BColor] (-.2,.5) rectangle (.2,0);
\filldraw[\QrColor] (-.2,1) rectangle (.2,.5);
\end{scope}
\draw[thick, \XColor] (-.2,-.5) -- (-.2,1);
\draw[thick, \QsColor] (-.2,0) -- (.2,0);
\draw[thick, \QsColor] (-.2,.5) -- (.2,.5);
\draw[thick, \YColor] (.2,-.5) -- (.2,1);
}
=
\tikzmath{
\begin{scope}
\clip[rounded corners = 5pt] (-.5,-.5) rectangle (.5,1);
\filldraw[\QrColor] (-.2,1) rectangle (.2,-.5);
\end{scope}
\draw[thick, \XColor] (-.2,-.5) -- (-.2,1);
\draw[thick, \YColor] (.2,-.5) -- (.2,1);
}
=
\id_{X\xxo_Q Y}
\qquad
(u^Q_{X,Y})^\dag \xxt u^Q_{X,Y} 
=
\tikzmath{
\begin{scope}
\clip[rounded corners = 5pt] (-.5,-.5) rectangle (.5,1);
\filldraw[\BColor] (-.2,0) rectangle (.2,-.5);
\filldraw[\QrColor] (-.2,.5) rectangle (.2,0);
\filldraw[\BColor] (-.2,1) rectangle (.2,.5);
\end{scope}
\draw[thick, \XColor] (-.2,-.5) -- (-.2,1);
\draw[thick, \QsColor] (-.2,0) -- (.2,0);
\draw[thick, \QsColor] (-.2,.5) -- (.2,.5);
\draw[thick, \YColor] (.2,-.5) -- (.2,1);
}
=
\tikzmath{
\begin{scope}
\clip[rounded corners = 5pt] (-.5,-.5) rectangle (.5,1);
\filldraw[\BColor] (-.2,1) rectangle (.2,-.5);
\end{scope}
\draw[thick, \XColor] (-.2,-.5) -- (-.2,1);
\draw[thick, \QsColor] (-.2,.25) -- (.2,.25);
\draw[thick, \YColor] (.2,-.5) -- (.2,1);
}
=
p^Q_{X, Y}.
\]
We define canonical unitor trivalent vertices by
\[
\lambda^P_X=
\tikzmath{
\begin{scope}
\clip[rounded corners = 5pt] (-.7,-.2) rectangle (.3,.5);
\filldraw[\AColor] (-.7,-.2) rectangle (0,.5);
\filldraw[\PrColor,thick] (-.4,-.2) arc (180:90:.4cm) -- (0,-.2);
\filldraw[\BColor] (0,-.2) rectangle (.3,.5);
\end{scope}
\draw[\XColor,thick] (0,-.2) -- (0,.5);
\draw[\PsColor,thick] (-.4,-.2) arc (180:90:.4cm);
\filldraw[\XColor] (0,.2) circle (.05cm);
}
:=
\tikzmath{
\begin{scope}
\clip[rounded corners = 5pt] (-.7,-.6) rectangle (.3,.5);
\filldraw[\AColor] (-.7,-.6) rectangle (0,.5);
\filldraw[\PrColor,thick] (-.4,-.6) rectangle (0,-.2);
\filldraw[\BColor] (0,-.6) rectangle (.3,.5);
\end{scope}
\draw[\XColor,thick] (0,-.6) -- (0,.5);
\draw[\PsColor,thick] (-.4,-.6) -- (-.4,-.2) arc (180:90:.4cm);
\draw[\PsColor,thick] (-.4,-.2) -- (0,-.2);
\filldraw[\XColor] (0,.2) circle (.05cm);
}
=\lambda_X\xxt (u^P_{P,X})^\dag
\qquad\text{and}\qquad
\rho_X^Q=
\tikzmath{
\begin{scope}
\clip[rounded corners = 5pt] (-.3,-.2) rectangle (.7,.5);
\filldraw[\AColor] (-.3,-.2) rectangle (0,.5);
\filldraw[\BColor] (0,-.2) rectangle (.7,.5);
\filldraw[\QrColor] (.4,-.2) arc (0:90:.4cm) -- (0,-.2);
\end{scope}
\draw[\XColor,thick] (0,-.2) -- (0,.5);
\draw[\QsColor,thick] (.4,-.2) arc (0:90:.4cm);
\filldraw[\XColor] (0,.2) circle (.05cm);
}
:=
\tikzmath{
\begin{scope}
\clip[rounded corners = 5pt] (-.3,-.6) rectangle (.7,.5);
\filldraw[\AColor] (-.3,-.6) rectangle (0,.5);
\filldraw[\BColor] (0,-.6) rectangle (.7,.5);
\filldraw[\QrColor] (.4,-.6) rectangle (0,-.2);
\end{scope}
\draw[\XColor,thick] (0,-.6) -- (0,.5);
\draw[\QsColor,thick] (.4,-.6) -- (.4,-.2) arc (0:90:.4cm);
\draw[\QsColor,thick] (0,-.2) -- (.4,-.2);
\filldraw[\XColor] (0,.2) circle (.05cm);
}
= \rho_X\xxt (u_{X,Q}^Q)^\dag.
\]
It is straightforward to verify that $\lambda^P_X$ and $\rho^Q_X$ are unitaries (see \cite[\S3.3]{2105.12010}).
In this graphical notation, the associator of $\QSys(\cC)$ is uniquely determined by the formula on the left hand side:
\[
\tikzmath{
\filldraw[\BColor] (-.4,.5) rectangle (0,0);
\filldraw[\QrColor] (-.4,.5) rectangle (0,2);
\filldraw[\CColor] (.4,.2) rectangle (0,0);
\filldraw[\RrColor] (.4,.2) rectangle (0,2);
\draw[\QsColor,thick] (-.4,.5) -- (0,.5);
\draw[\RsColor,thick] (.4,.2) -- (0,.2);
\draw[\XColor,thick] (-.4,0) -- (-.4,2);
\draw[\YColor,thick] (0,0) -- (0,2);
\draw[\ZColor,thick] (.4,0) -- (.4,2);
\roundNbox{unshaded}{(0,1)}{.3}{.3}{.3}{\scriptsize{$\alpha^{\QSys(\cC)}_{X,Y,Z}$}};
}
=
\tikzmath{
\filldraw[\BColor] (-.4,0) rectangle (0,1.8);
\filldraw[\QrColor] (-.4,2) rectangle (0,1.8);
\filldraw[\CColor] (.4,0) rectangle (0,1.5);
\filldraw[\RrColor] (.4,2) rectangle (0,1.5);
\draw[\QsColor,thick] (-.4,1.8) -- (0,1.8);
\draw[\RsColor,thick] (.4,1.5) -- (0,1.5);
\draw[\XColor,thick] (-.4,0) -- (-.4,2);
\draw[\YColor,thick] (0,0) -- (0,2);
\draw[\ZColor,thick] (.4,0) -- (.4,2);
\roundNbox{unshaded}{(0,1)}{.3}{.3}{.3}{\scriptsize{$\alpha^{\cC}_{X,Y,Z}$}};
}
\qquad\qquad
\Longrightarrow
\qquad\qquad
\underset{
(X\xxo_Q Y)\xxo_R Z \Rightarrow X\xxo_Q (Y\xxo_R Z)
}{
\tikzmath{
\filldraw[\QrColor] (-.4,0) rectangle (0,2);
\filldraw[\RrColor] (.4,0) rectangle (0,2);
\draw[\XColor,thick] (-.4,0) -- (-.4,2);
\draw[\YColor,thick] (0,0) -- (0,2);
\draw[\ZColor,thick] (.4,0) -- (.4,2);
\roundNbox{unshaded}{(0,1)}{.3}{.3}{.3}{\scriptsize{$\alpha^{\QSys(\cC)}_{X,Y,Z}$}};
}
=
\tikzmath{
\filldraw[\QrColor] (-.4,0) rectangle (0,.5);
\filldraw[\BColor] (-.4,1.8) rectangle (0,.5);
\filldraw[\QrColor] (-.4,2) rectangle (0,1.8);
\filldraw[\RrColor] (.4,0) rectangle (0,.2);
\filldraw[\CColor] (.4,1.5) rectangle (0,.2);
\filldraw[\RrColor] (.4,2) rectangle (0,1.5);
\draw[\QsColor,thick] (-.4,.5) -- (0,.5);
\draw[\RsColor,thick] (.4,.2) -- (0,.2);
\draw[\QsColor,thick] (-.4,1.8) -- (0,1.8);
\draw[\RsColor,thick] (.4,1.5) -- (0,1.5);
\draw[\XColor,thick] (-.4,0) -- (-.4,2);
\draw[\YColor,thick] (0,0) -- (0,2);
\draw[\ZColor,thick] (.4,0) -- (.4,2);
\roundNbox{unshaded}{(0,1)}{.3}{.3}{.3}{\scriptsize{$\alpha^{\cC}_{X,Y,Z}$}};
}
}
.
\]
\end{nota}

\subsection{Constructions on 1-morphisms, 2-morphisms, and 3-morphisms in \texorpdfstring{$2\Cat$}{2Cat}}

For this section, we fix two $\rm C^*/W^*$ 2-categories $\cC,\cD$.

\begin{construction}[{\cite[Const.~3.29]{2105.12010}}]
\label{construction:Qsys(F)}
A $\dag$ 2-functor $F: \cC \to \cD$ between $\rm C^*/W^*$ 2-categories induces a $\dag$-2-functor $\QSys(F): \QSys(\cC) \to \QSys(\cD)$.
\begin{itemize}
\item 
For $({}_bQ_b,m,i)\in \QSys(\cC)$,
we define 
$$
\QSys(F)({}_bQ_b):=({}_{F(b)}F(Q)_{F(b)},F(m)\xxt F^2_{Q,Q},F(i)\xxt F^1_b)\in \QSys(\cD).
$$

\item
For $({}_PX_Q,\lambda, \rho)\in \QSys(\cC)(P\to Q)$,
we define
$$
\QSys(F)({}_PX_Q)
:=
(F(X),F(\lambda)\xxt F^2_{P,X},F(\rho)\xxt F^2_{X,Q})\in \QSys(\cD)(F(P)\to F(Q))
$$

\item
For $f \in \QSys(\cC)({}_PX_Q \Rightarrow {}_PY_Q)$
we define
$$
\QSys(F)(f):=
F(f)\in \QSys(\cD)({}_{F(P)}F(X)_{F(Q)} \Rightarrow {}_{F(P)}F(Y)_{F(Q)}).
$$
Since $F$ is a $\dag$ 2-functor, $\QSys(F)$ will be as well.
Moreover, when $\cA,\cB$ are $\rm W^*$ and $F: \cA\to \cB$ is normal, so is $\QSys(F)$.

\item
For ${}_P X_Q \in \QSys(\cC)(P \to Q)$ and ${}_Q Y_R \in \QSys(\cC)(Q\to R)$, 
we define
\begin{equation}
\label{eq:DefOfQSysF2}
\QSys(F)_{X,Y}^{2} :=F(u_{X,Y}) \xxt F^2_{X,Y} \xxt u^\dag_{F(X), F(Y)}
\in \QSys(\cD)(F(X)\xxo_{F(Q)}F(Y)\Rightarrow F(X\xxo_Q Y)).
\end{equation}
Finally, for a  Q-system $Q\in \cC(b\to b)$, 
we define 
$$
\QSys(F)^1_{F(Q)}:=\id\in \QSys(\cD)(1_{F(Q)}\Rightarrow F(1_Q)).
$$
\end{itemize}
For convenience of the reader, we provide a diagrammatic proof below that $\QSys(F)$ is a $\dag$ 2-functor.
We graphically represent
\begin{align*}
\tikzmath{
\filldraw[primedregion=white, rounded corners = 5pt] (0,0) rectangle (.6,.6);
\draw[thin, dotted, rounded corners = 5pt] (0,0) rectangle (.6,.6);
}
&=
F
&
\tikzmath{
\begin{scope}
\clip[rounded corners=5pt] (-.5,-.3) rectangle (.5,.3);
\filldraw[primedregion=white] (-.5,-.3) rectangle (-.2,.3);
\filldraw[primedregion=\QrColor] (-.2,-.3) rectangle (.2,.3);
\filldraw[primedregion=white] (.5,-.3) rectangle (.2,.3);
\end{scope}
\draw[\XColor,thick] (-.2,-.3) -- (-.2,.3);
\draw[\YColor,thick] (.2,-.3) -- (.2,.3);
\draw[thin, dotted, rounded corners = 5pt] (-.5,-.3) rectangle (.5,.3);
}
&=
F(X)\xxo_{F(Q)} F(Y)
&
\tikzmath{
\begin{scope}
\clip[rounded corners=5pt] (-.35,-.3) rectangle (.35,.3);
\filldraw[primedregion=white] (-.05,-.3) rectangle (-.35,.3);
\filldraw[primedregion=\QrColor] (-.05,-.3) rectangle (.05,.3);
\filldraw[primedregion=white] (.05,-.3) rectangle (.35,.3);
\end{scope}
\draw[\XColor,thick] (-.05,-.3) -- (-.05,.3);
\draw[\YColor,thick] (.05,-.3) -- (.05,.3);
\draw[thin, dotted, rounded corners = 5pt] (-.35,-.3) rectangle (.35,.3);
}
&=
F(X\xxo_Q Y)
\\
\tikzmath{
\begin{scope}
\clip[rounded corners=5pt] (-.5,-.3) rectangle (.5,.3);
\filldraw[primedregion=white] (-.5,-.3) rectangle (-.2,.3);
\filldraw[primedregion=\QrColor] (-.2,0) rectangle (.2,.3);
\filldraw[primedregion=\BColor] (-.2,0) rectangle (.2,-.3);
\filldraw[primedregion=white] (.5,-.3) rectangle (.2,.3);
\end{scope}
\draw[\XColor,thick] (-.2,-.3) -- (-.2,.3);
\draw[\YColor,thick] (.2,-.3) -- (.2,.3);
\draw[\QsColor,thick] (-.2,0) -- (.2,0);
\draw[thin, dotted, rounded corners = 5pt] (-.5,-.3) rectangle (.5,.3);
}
&=
u^{F(Q)}_{F(X),F(Y)}
&
\tikzmath{
\begin{scope}
\clip[rounded corners=5pt] (-.35,-.3) rectangle (.35,.3);
\filldraw[primedregion=white] (-.05,-.3) rectangle (-.35,.3);
\filldraw[primedregion=\BColor] (-.05,.3) rectangle (.05,0);
\filldraw[primedregion=\QrColor] (-.05,-.3) rectangle (.05,0);
\filldraw[primedregion=white] (.05,-.3) rectangle (.35,.3);
\end{scope}
\draw[\XColor,thick] (-.05,-.3) -- (-.05,.3);
\draw[\YColor,thick] (.05,-.3) -- (.05,.3);
\draw[\QsColor,thick] (-.05,0) -- (.05,0);
\draw[thin, dotted, rounded corners = 5pt] (-.35,-.3) rectangle (.35,.3);
}
&=
F(u^Q_{X,Y})^\dag
&
\tikzmath{
\begin{scope}
\clip[rounded corners=5pt] (-.35,-.3) rectangle (.35,.3);
\filldraw[primedregion=white] (-.05,-.3) rectangle (-.35,.3);
\filldraw[primedregion=\BColor] (-.05,-.3) rectangle (.05,.3);
\filldraw[primedregion=white] (.05,-.3) rectangle (.35,.3);
\end{scope}
\draw[\XColor,thick] (-.05,-.3) -- (-.05,.3);
\draw[\YColor,thick] (.05,-.3) -- (.05,.3);
\draw[\QsColor,thick] (-.05,0) -- (.05,0);
\draw[thin, dotted, rounded corners = 5pt] (-.35,-.3) rectangle (.35,.3);
}
&=
F(p^Q_{X,Y}).
\end{align*}
We then define
\[
\tikzmath[scale=.5, transform shape]{
\begin{scope}
\clip[rounded corners=5pt] (-1.4,-2) rectangle (1.4,2);
\filldraw[primedregion=white] (-1.4,-2) rectangle (1.4,2);
\filldraw[primedregion=\QrColor] (-.4,-2) rectangle (.4,0);
\filldraw[primedregion=\QrColor] (-.1,2) rectangle (.1,0);
\end{scope}
\draw[\XColor,thick] (-.4,-2) -- (-.4,0);
\draw[\YColor,thick] (.4,-2) -- (.4,0);
\draw[\XColor,thick] (-.1,2) -- (-.1,0);
\draw[\YColor,thick] (.1,2) -- (.1,0);
\roundNbox{unshaded}{(0,0)}{.6}{.45}{.45}{\small{$\QSys(F)^2_{X,Y}$}};
\draw[thin, dotted, rounded corners = 5pt] (-1.4,-2) rectangle (1.4,2);
}
:=
\tikzmath{
\begin{scope}
\clip[rounded corners=5pt] (-.7,-1) rectangle (.7,1);
\filldraw[primedregion=white] (-.7,-1) rectangle (.7,1);
\filldraw[primedregion=\QrColor] (-.2,-1) rectangle (.2,-.7);
\filldraw[primedregion=\BColor] (-.2,-.3) rectangle (.2,-.7);
\filldraw[primedregion=\QrColor] (-.05,1) rectangle (.05,.7);
\filldraw[primedregion=\BColor] (-.05,.3) rectangle (.05,.7);
\end{scope}
\draw[\QsColor,thick] (-.2,-.7) -- (.2,-.7);
\draw[\QsColor,thick] (-.05,.7) -- (.05,.7);
\draw[\XColor,thick] (-.2,-1) -- (-.2,0);
\draw[\YColor,thick] (.2,-1) -- (.2,0);
\draw[\XColor,thick] (-.05,1) -- (-.05,0);
\draw[\YColor,thick] (.05,1) -- (.05,0);
\roundNbox{unshaded}{(0,0)}{.3}{.2}{.2}{\scriptsize{$F^2_{X,Y}$}};
\draw[thin, dotted, rounded corners = 5pt] (-.7,-1) rectangle (.7,1);
}\,.
\]
By definition of the separability projector \eqref{eq:SeparabilityProjector} for $F(X)\xxo_{F(Q)}F(Y)$, we have
\[
p_{F(X),F(Y)} =
\tikzmath{
\begin{scope}
\clip[rounded corners=5pt] (-.5,-1) rectangle (.5,1);
\filldraw[primedregion=white] (-.5,-1) rectangle (.5,1);
\filldraw[primedregion=\BColor] (-.2,-1) rectangle (.2,1);
\end{scope}
\draw[\QsColor,thick] (-.2,0) -- (.2,0);
\draw[\XColor,thick] (-.2,-1) -- (-.2,1);
\draw[\YColor,thick] (.2,-1) -- (.2,1);
\draw[thin, dotted, rounded corners = 5pt] (-.5,-1) rectangle (.5,1);
}
:=
\tikzmath{
\begin{scope}
\clip[rounded corners=5pt] (-.7,-1.5) rectangle (.7,1.5);
\filldraw[primedregion=white] (-.7,-1.5) rectangle (.7,1.5);
\filldraw[primedregion=\BColor] (-.2,-1) rectangle (.2,-1.5);
\filldraw[primedregion=\BColor] (-.05,-1) rectangle (.05,1);
\filldraw[primedregion=\BColor] (-.2,1) rectangle (.2,1.5);
\end{scope}
\draw[\QsColor,thick] (-.05,0) -- (.05,0);
\draw[\XColor,thick] (-.2,-1.5) -- (-.2,-1);
\draw[\YColor,thick] (.2,-1.5) -- (.2,-1);
\draw[\XColor,thick] (-.05,-1) -- (-.05,1);
\draw[\YColor,thick] (.05,-1) -- (.05,1);
\draw[\XColor,thick] (-.2,1.5) -- (-.2,1);
\draw[\YColor,thick] (.2,1.5) -- (.2,1);
\roundNbox{unshaded}{(0,-.7)}{.3}{.2}{.2}{\scriptsize{$F^2_{X,Y}$}};
\roundNbox{unshaded}{(0,.7)}{.3}{.25}{.25}{\tiny{$(F^2_{X,Y})^\dag$}};
\draw[thin, dotted, rounded corners = 5pt] (-.7,-1.5) rectangle (.7,1.5);
}
\qquad\qquad
\Longrightarrow
\qquad\qquad
\tikzmath{
\begin{scope}
\clip[rounded corners=5pt] (-.7,-.8) rectangle (.7,1.1);
\filldraw[primedregion=white] (-.7,-.8) rectangle (.7,1.1);
\filldraw[primedregion=\BColor] (-.2,-.8) rectangle (.2,-.3);
\filldraw[primedregion=\BColor] (-.05,.3) rectangle (.05,1.1);
\end{scope}
\draw[\QsColor,thick] (-.05,.7) -- (.05,.7);
\draw[\XColor,thick] (-.2,-.8) -- (-.2,-.3);
\draw[\YColor,thick] (.2,-.8) -- (.2,-.3);
\draw[\XColor,thick] (-.05,.3) -- (-.05,1.1);
\draw[\YColor,thick] (.05,.3) -- (.05,1.1);
\roundNbox{unshaded}{(0,0)}{.3}{.2}{.2}{\scriptsize{$F^2_{X,Y}$}};
\draw[thin, dotted, rounded corners = 5pt] (-.7,-.8) rectangle (.7,1.1);
}
=
\tikzmath{
\begin{scope}
\clip[rounded corners=5pt] (-.7,-1.1) rectangle (.7,.8);
\filldraw[primedregion=white] (-.7,-1.1) rectangle (.7,.8);
\filldraw[primedregion=\BColor] (-.2,-1.1) rectangle (.2,-.3);
\filldraw[primedregion=\BColor] (-.05,.3) rectangle (.05,.8);
\end{scope}
\draw[\QsColor,thick] (-.2,-.7) -- (.2,-.7);
\draw[\XColor,thick] (-.2,-1.1) -- (-.2,-.3);
\draw[\YColor,thick] (.2,-1.1) -- (.2,-.3);
\draw[\XColor,thick] (-.05,.3) -- (-.05,.8);
\draw[\YColor,thick] (.05,.3) -- (.05,.8);
\roundNbox{unshaded}{(0,0)}{.3}{.2}{.2}{\scriptsize{$F^2_{X,Y}$}};
\draw[thin, dotted, rounded corners = 5pt] (-.7,-1.1) rectangle (.7,.8);
}\,.
\]
This formula for $p_{F(X), F(Y)}$ immediately implies 
$\QSys(F)^2_{X,Y}$ is unitary:
\[
\tikzmath[scale=.5, transform shape]{
\begin{scope}
\clip[rounded corners=5pt] (-1.6,-3) rectangle (1.6,3);
\filldraw[primedregion=white] (-1.6,-3) rectangle (1.6,3);
\filldraw[primedregion=\QrColor] (-.4,-2) rectangle (.4,-3);
\filldraw[primedregion=\QrColor] (-.4,2) rectangle (.4,3);
\filldraw[primedregion=\QrColor] (-.1,-2) rectangle (.1,2);
\end{scope}
\draw[\XColor,thick] (-.4,-3) -- (-.4,-2);
\draw[\YColor,thick] (.4,-3) -- (.4,-2);
\draw[\XColor,thick] (-.1,-2) -- (-.1,2);
\draw[\YColor,thick] (.1,-2) -- (.1,2);
\draw[\XColor,thick] (-.4,3) -- (-.4,2);
\draw[\YColor,thick] (.4,3) -- (.4,2);
\roundNbox{unshaded}{(0,-1.4)}{.6}{.7}{.7}{\small{$\QSys(F)^2_{X,Y}$}};
\roundNbox{unshaded}{(0,1.4)}{.6}{.7}{.7}{\small{$(\QSys(F)^2_{X,Y})^\dag$}};
\draw[thin, dotted, rounded corners = 5pt] (-1.6,-3) rectangle (1.6,3);
}
=
\tikzmath{
\begin{scope}
\clip[rounded corners=5pt] (-.7,-1.5) rectangle (.7,1.5);
\filldraw[primedregion=white] (-.7,-1.5) rectangle (.7,1.5);
\filldraw[primedregion=\QrColor] (-.2,1.5) rectangle (.2,1.2);
\filldraw[primedregion=\BColor] (-.2,.7) rectangle (.2,1.2);
\filldraw[primedregion=\BColor] (-.05,.2) rectangle (.05,.7);
\filldraw[primedregion=\QrColor] (-.05,-.2) rectangle (.05,.2);
\filldraw[primedregion=\BColor] (-.05,-.2) rectangle (.05,-.7);
\filldraw[primedregion=\BColor] (-.2,-.7) rectangle (.2,-1.2);
\filldraw[primedregion=\QrColor] (-.2,-1.5) rectangle (.2,-1.2);
\end{scope}
\draw[\QsColor,thick] (-.2,1.2) -- (.2,1.2);
\draw[\QsColor,thick] (-.05,.2) -- (.05,.2);
\draw[\QsColor,thick] (-.05,-.2) -- (.05,-.2);
\draw[\QsColor,thick] (-.2,-1.2) -- (.2,-1.2);
\draw[\XColor,thick] (-.2,-1.5) -- (-.2,-1);
\draw[\YColor,thick] (.2,-1.5) -- (.2,-1);
\draw[\XColor,thick] (-.05,-1) -- (-.05,1);
\draw[\YColor,thick] (.05,-1) -- (.05,1);
\draw[\XColor,thick] (-.2,1.5) -- (-.2,1);
\draw[\YColor,thick] (.2,1.5) -- (.2,1);
\roundNbox{unshaded}{(0,-.7)}{.3}{.2}{.2}{\scriptsize{$F^2_{X,Y}$}};
\roundNbox{unshaded}{(0,.7)}{.3}{.25}{.25}{\tiny{$(F^2_{X,Y})^\dag$}};
\draw[thin, dotted, rounded corners = 5pt] (-.7,-1.5) rectangle (.7,1.5);
}
=
\tikzmath{
\begin{scope}
\clip[rounded corners=5pt] (-.7,-1.5) rectangle (.7,1.5);
\filldraw[primedregion=white] (-.7,-1.5) rectangle (.7,1.5);
\filldraw[primedregion=\QrColor] (-.2,1.5) rectangle (.2,1.2);
\filldraw[primedregion=\BColor] (-.2,.7) rectangle (.2,1.2);
\filldraw[primedregion=\BColor] (-.05,.7) rectangle (.05,-.7);
\filldraw[primedregion=\BColor] (-.2,-.7) rectangle (.2,-1.2);
\filldraw[primedregion=\QrColor] (-.2,-1.5) rectangle (.2,-1.2);
\end{scope}
\draw[\QsColor,thick] (-.2,1.2) -- (.2,1.2);
\draw[\QsColor,thick] (-.05,0) -- (.05,0);
\draw[\QsColor,thick] (-.2,-1.2) -- (.2,-1.2);
\draw[\XColor,thick] (-.2,-1.5) -- (-.2,-1);
\draw[\YColor,thick] (.2,-1.5) -- (.2,-1);
\draw[\XColor,thick] (-.05,-1) -- (-.05,1);
\draw[\YColor,thick] (.05,-1) -- (.05,1);
\draw[\XColor,thick] (-.2,1.5) -- (-.2,1);
\draw[\YColor,thick] (.2,1.5) -- (.2,1);
\roundNbox{unshaded}{(0,-.7)}{.3}{.2}{.2}{\scriptsize{$F^2_{X,Y}$}};
\roundNbox{unshaded}{(0,.7)}{.3}{.25}{.25}{\tiny{$(F^2_{X,Y})^\dag$}};
\draw[thin, dotted, rounded corners = 5pt] (-.7,-1.5) rectangle (.7,1.5);
}
=
\tikzmath{
\begin{scope}
\clip[rounded corners=5pt] (-.7,-1.5) rectangle (.7,1.5);
\filldraw[primedregion=white] (-.7,-1.5) rectangle (.7,1.5);
\filldraw[primedregion=\QrColor] (-.2,1.5) rectangle (.2,1.2);
\filldraw[primedregion=\BColor] (-.2,.7) rectangle (.2,1.2);
\filldraw[primedregion=\BColor] (-.05,.7) rectangle (.05,-.7);
\filldraw[primedregion=\BColor] (-.2,-.7) rectangle (.2,-1.2);
\filldraw[primedregion=\QrColor] (-.2,-1.5) rectangle (.2,-1.2);
\end{scope}
\draw[\QsColor,thick] (-.2,1.2) -- (.2,1.2);
\draw[\QsColor,thick] (-.2,-1.2) -- (.2,-1.2);
\draw[\XColor,thick] (-.2,-1.5) -- (-.2,-1);
\draw[\YColor,thick] (.2,-1.5) -- (.2,-1);
\draw[\XColor,thick] (-.05,-1) -- (-.05,1);
\draw[\YColor,thick] (.05,-1) -- (.05,1);
\draw[\XColor,thick] (-.2,1.5) -- (-.2,1);
\draw[\YColor,thick] (.2,1.5) -- (.2,1);
\roundNbox{unshaded}{(0,-.7)}{.3}{.2}{.2}{\scriptsize{$F^2_{X,Y}$}};
\roundNbox{unshaded}{(0,.7)}{.3}{.25}{.25}{\tiny{$(F^2_{X,Y})^\dag$}};
\draw[thin, dotted, rounded corners = 5pt] (-.7,-1.5) rectangle (.7,1.5);
}
=
\tikzmath{
\begin{scope}
\clip[rounded corners=5pt] (-.5,-1) rectangle (.5,1);
\filldraw[primedregion=white] (-.5,-1) rectangle (.5,1);
\filldraw[primedregion=\QrColor] (-.2,1) rectangle (.2,.5);
\filldraw[primedregion=\BColor] (-.2,-.5) rectangle (.2,.5);
\filldraw[primedregion=\QrColor] (-.2,-1) rectangle (.2,-.5);
\end{scope}
\draw[\QsColor,thick] (-.2,-.5) -- (.2,-.5);
\draw[\QsColor,thick] (-.2,.5) -- (.2,.5);
\draw[\XColor,thick] (-.2,-1) -- (-.2,1);
\draw[\YColor,thick] (.2,-1) -- (.2,1);
\draw[thin, dotted, rounded corners = 5pt] (-.5,-1) rectangle (.5,1);
}
=
\tikzmath{
\begin{scope}
\clip[rounded corners=5pt] (-.5,-1) rectangle (.5,1);
\filldraw[primedregion=white] (-.5,-1) rectangle (.5,1);
\filldraw[primedregion=\QrColor] (-.2,1) rectangle (.2,-1);
\end{scope}
\draw[\XColor,thick] (-.2,-1) -- (-.2,1);
\draw[\YColor,thick] (.2,-1) -- (.2,1);
\draw[thin, dotted, rounded corners = 5pt] (-.5,-1) rectangle (.5,1);
}\,;
\quad\text{similarly,}\quad
\tikzmath[scale=.5, transform shape]{
\begin{scope}
\clip[rounded corners=5pt] (-1.6,-3) rectangle (1.6,3);
\filldraw[primedregion=white] (-1.6,-3) rectangle (1.6,3);
\filldraw[primedregion=\QrColor] (-.1,-2) rectangle (.1,-3);
\filldraw[primedregion=\QrColor] (-.1,2) rectangle (.1,3);
\filldraw[primedregion=\QrColor] (-.4,-2) rectangle (.4,2);
\end{scope}
\draw[\XColor,thick] (-.1,-3) -- (-.1,-2);
\draw[\YColor,thick] (.1,-3) -- (.1,-2);
\draw[\XColor,thick] (-.4,-2) -- (-.4,2);
\draw[\YColor,thick] (.4,-2) -- (.4,2);
\draw[\XColor,thick] (-.1,3) -- (-.1,2);
\draw[\YColor,thick] (.1,3) -- (.1,2);
\roundNbox{unshaded}{(0,1.4)}{.6}{.7}{.7}{\small{$\QSys(F)^2_{X,Y}$}};
\roundNbox{unshaded}{(0,-1.4)}{.6}{.7}{.7}{\small{$(\QSys(F)^2_{X,Y})^\dag$}};
\draw[thin, dotted, rounded corners = 5pt] (-1.6,-3) rectangle (1.6,3);
}
=
\tikzmath{
\begin{scope}
\clip[rounded corners=5pt] (-.35,-1) rectangle (.35,1);
\filldraw[primedregion=white] (-.35,-1) rectangle (.35,1);
\filldraw[primedregion=\QrColor] (-.05,1) rectangle (.05,-1);
\end{scope}
\draw[\XColor,thick] (-.05,-1) -- (-.05,1);
\draw[\YColor,thick] (.05,-1) -- (.05,1);
\draw[thin, dotted, rounded corners = 5pt] (-.35,-1) rectangle (.35,1);
}\,.
\]
Using \eqref{eq:DefOfQSysF2}, unitarity of $\QSys(F)^2$, and that $u$ is a coisometry, we have
\[
\tikzmath[scale=.5, transform shape]{
\begin{scope}
\clip[rounded corners=5pt] (-1.4,-2) rectangle (1.4,2);
\filldraw[primedregion=white] (-1.4,-2) rectangle (1.4,2);
\filldraw[primedregion=\BColor] (-.4,-1.4) rectangle (.4,-2);
\filldraw[primedregion=\QrColor] (-.4,-1.4) rectangle (.4,0);
\filldraw[primedregion=\QrColor] (-.1,2) rectangle (.1,0);
\end{scope}
\draw[\QsColor,thick] (-.4,-1.4) -- (.4,-1.4);
\draw[\XColor,thick] (-.4,-2) -- (-.4,0);
\draw[\YColor,thick] (.4,-2) -- (.4,0);
\draw[\XColor,thick] (-.1,2) -- (-.1,0);
\draw[\YColor,thick] (.1,2) -- (.1,0);
\roundNbox{unshaded}{(0,0)}{.6}{.45}{.45}{\small{$\QSys(F)^2_{X,Y}$}};
\draw[thin, dotted, rounded corners = 5pt] (-1.4,-2) rectangle (1.4,2);
}
=
\tikzmath{
\begin{scope}
\clip[rounded corners=5pt] (-.7,-1) rectangle (.7,1);
\filldraw[primedregion=white] (-.7,-1) rectangle (.7,1);
\filldraw[primedregion=\BColor] (-.2,-1) rectangle (.2,0);
\filldraw[primedregion=\BColor] (-.05,0) rectangle (.05,.7);
\filldraw[primedregion=\QrColor] (-.05,1) rectangle (.05,.7);
\end{scope}
\draw[\QsColor,thick] (-.05,.7) -- (.05,.7);
\draw[\XColor,thick] (-.2,-1) -- (-.2,0);
\draw[\YColor,thick] (.2,-1) -- (.2,0);
\draw[\XColor,thick] (-.05,1) -- (-.05,0);
\draw[\YColor,thick] (.05,1) -- (.05,0);
\roundNbox{unshaded}{(0,0)}{.3}{.2}{.2}{\scriptsize{$F^2_{X,Y}$}};
\draw[thin, dotted, rounded corners = 5pt] (-.7,-1) rectangle (.7,1);
}
\qquad\text{and}\qquad
\tikzmath[scale=.5, transform shape]{
\begin{scope}
\clip[rounded corners=5pt] (-1.4,-2) rectangle (1.4,2);
\filldraw[primedregion=white] (-1.4,-2) rectangle (1.4,2);
\filldraw[primedregion=\QrColor] (-.4,-2) rectangle (.4,0);
\filldraw[primedregion=\QrColor] (-.1,1.4) rectangle (.1,0);
\filldraw[primedregion=\BColor] (-.1,1.4) rectangle (.1,2);
\end{scope}
\draw[\QsColor,thick] (-.1,1.4) -- (.1,1.4);
\draw[\XColor,thick] (-.4,-2) -- (-.4,0);
\draw[\YColor,thick] (.4,-2) -- (.4,0);
\draw[\XColor,thick] (-.1,2) -- (-.1,0);
\draw[\YColor,thick] (.1,2) -- (.1,0);
\roundNbox{unshaded}{(0,0)}{.6}{.45}{.45}{\small{$\QSys(F)^2_{X,Y}$}};
\draw[thin, dotted, rounded corners = 5pt] (-1.4,-2) rectangle (1.4,2);
}
=
\tikzmath{
\begin{scope}
\clip[rounded corners=5pt] (-.7,-1) rectangle (.7,1);
\filldraw[primedregion=white] (-.7,-1) rectangle (.7,1);
\filldraw[primedregion=\QrColor] (-.2,-1) rectangle (.2,-.7);
\filldraw[primedregion=\BColor] (-.2,0) rectangle (.2,-.7);
\filldraw[primedregion=\BColor] (-.05,1) rectangle (.05,0);
\end{scope}
\draw[\QsColor,thick] (-.2,-.7) -- (.2,-.7);
\draw[\XColor,thick] (-.2,-1) -- (-.2,0);
\draw[\YColor,thick] (.2,-1) -- (.2,0);
\draw[\XColor,thick] (-.05,1) -- (-.05,0);
\draw[\YColor,thick] (.05,1) -- (.05,0);
\roundNbox{unshaded}{(0,0)}{.3}{.2}{.2}{\scriptsize{$F^2_{X,Y}$}};
\draw[thin, dotted, rounded corners = 5pt] (-.7,-1) rectangle (.7,1);
}
\]
By naturality, we have
\[
\tikzmath{
\begin{scope}
\clip[rounded corners=5pt] (-.85,-1) rectangle (.85,1);
\filldraw[primedregion=white] (-.85,-1) rectangle (.85,1);
\filldraw[primedregion=\QrColor] (-.3,-1) rectangle (-.2,0);
\filldraw[primedregion=\QrColor] (-.1,0) rectangle (0,.7);
\filldraw[primedregion=\BColor] (-.1,.7) rectangle (0,1);
\filldraw[primedregion=\CColor] (-.2,-1) rectangle (.2,0);
\filldraw[primedregion=\CColor] (0,0) rectangle (.1,1);
\end{scope}
\draw[\QsColor,thick] (-.1,.7) -- (0,.7);
\draw[\XColor,thick] (-.3,-1) -- (-.3,0);
\draw[\YColor,thick] (-.2,-1) -- (-.2,0);
\draw[\ZColor,thick] (.2,-1) -- (.2,0);
\draw[\XColor,thick] (-.1,1) -- (-.1,0);
\draw[\YColor,thick] (0,1) -- (0,0);
\draw[\ZColor,thick] (.1,1) -- (.1,0);
\roundNbox{unshaded}{(0,0)}{.3}{.35}{.35}{\scriptsize{$F^2_{X\xxo_Q Y,Z}$}};
\draw[thin, dotted, rounded corners = 5pt] (-.85,-1) rectangle (.85,1);
}
=
\tikzmath{
\begin{scope}
\clip[rounded corners=5pt] (-.85,-1) rectangle (.85,1);
\filldraw[primedregion=white] (-.85,-1) rectangle (.85,1);
\filldraw[primedregion=\BColor] (-.1,0) rectangle (0,1);
\filldraw[primedregion=\BColor] (-.3,-.7) rectangle (-.2,0);
\filldraw[primedregion=\QrColor] (-.3,-1) rectangle (-.2,-.7);
\filldraw[primedregion=\CColor] (-.2,-1) rectangle (.2,0);
\filldraw[primedregion=\CColor] (0,0) rectangle (.1,1);
\end{scope}
\draw[\QsColor,thick] (-.3,-.7) -- (-.2,-.7);
\draw[\XColor,thick] (-.3,-1) -- (-.3,0);
\draw[\YColor,thick] (-.2,-1) -- (-.2,0);
\draw[\ZColor,thick] (.2,-1) -- (.2,0);
\draw[\XColor,thick] (-.1,1) -- (-.1,0);
\draw[\YColor,thick] (0,1) -- (0,0);
\draw[\ZColor,thick] (.1,1) -- (.1,0);
\roundNbox{unshaded}{(0,0)}{.3}{.35}{.35}{\scriptsize{$F^2_{X\xxo Y,Z}$}};
\draw[thin, dotted, rounded corners = 5pt] (-.85,-1) rectangle (.85,1);
}
: F(X\xxo_Q Y)\xxo F(Z)\to F(X\xxo Y)\xxo F(Z).
\]
These identities are used to prove the hexagon associativity coherence for $\QSys(F)^2$ and the triangle unit coherences for $\QSys(F)^1$:
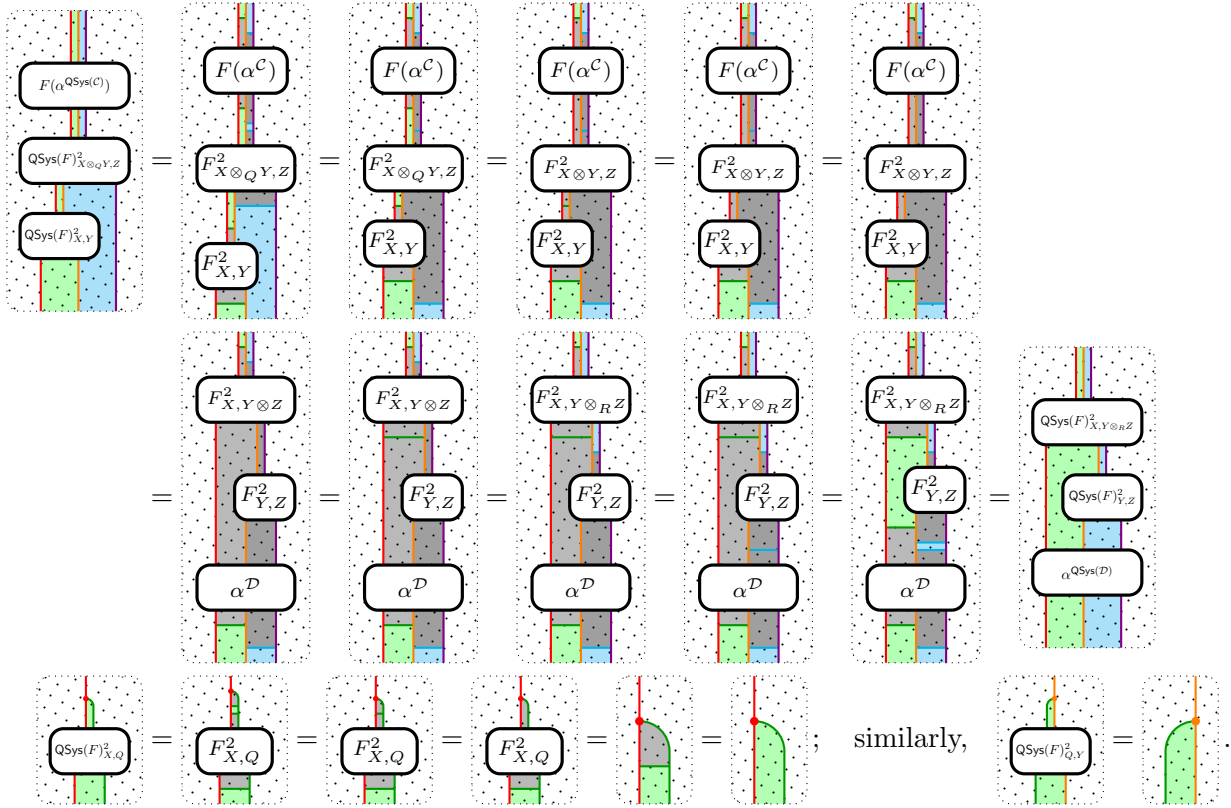
\begin{align*}
\tikzmath[scale=.5, transform shape]{
\begin{scope}
\clip[rounded corners=5pt] (-1.9,0) rectangle (1.7,8);
\filldraw[primedregion=white] (-1.9,0) rectangle (1.7,8);
\filldraw[primedregion=\RrColor] (1,0) rectangle (-.4,4);
\filldraw[primedregion=\QrColor] (-1,0) rectangle (0,2);
\filldraw[primedregion=\QrColor] (-.6,2) rectangle (-.4,4);
\filldraw[primedregion=\QrColor] (-.2,4) rectangle (0,8);
\filldraw[primedregion=\RrColor] (0,4) rectangle (.2,8);
\end{scope}
\draw[\XColor,thick] (-1,0) -- (-1,2);
\draw[\YColor,thick] (0,0) -- (0,2);
\draw[\XColor,thick] (-.6,2) -- (-.6,4);
\draw[\YColor,thick] (-.4,2) -- (-.4,4);
\draw[\ZColor,thick] (1,0) -- (1,4);
\draw[\XColor,thick] (-.2,4) -- (-.2,8);
\draw[\YColor,thick] (0,4) -- (0,8);
\draw[\ZColor,thick] (.2,4) -- (.2,8);
\roundNbox{unshaded}{(-.5,2)}{.6}{.45}{.45}{\small{$\QSys(F)^2_{X,Y}$}};
\roundNbox{unshaded}{(0,4)}{.6}{.95}{.75}{\small{$\QSys(F)^2_{X\xxo_Q Y,Z}$}};
\roundNbox{unshaded}{(0,6)}{.6}{.95}{.75}{\normalsize{$F(\alpha^{\QSys(\cC)})$}};
\draw[thin, dotted, rounded corners = 5pt] (-1.9,0) rectangle (1.7,8);
}
&=
\tikzmath{
\begin{scope}
\clip[rounded corners=5pt] (-.85,0) rectangle (.85,4.2);
\filldraw[primedregion=white] (-.85,0) rectangle (.85,4.2);
\filldraw[primedregion=\QrColor] (-.4,0) rectangle (0,.2);
\filldraw[primedregion=\QrColor] (-.25,1.2) rectangle (-.15,2);
\filldraw[primedregion=\QrColor] (-.1,2) rectangle (0,2.8);
\filldraw[primedregion=\QrColor] (-.1,4) rectangle (0,4.2);
\filldraw[primedregion=\BColor] (-.4,.2) rectangle (0,.7);
\filldraw[primedregion=\BColor] (-.25,.7) rectangle (-.15,1.2);
\filldraw[primedregion=\BColor] (-.1,2.8) rectangle (0,4);
\filldraw[primedregion=\RrColor] (0,0) rectangle (.4,1.5);
\filldraw[primedregion=\RrColor] (-.15,1) rectangle (.4,1.5);
\filldraw[primedregion=\RrColor] (0,2.5) rectangle (.1,2.6);
\filldraw[primedregion=\RrColor] (0,3.8) rectangle (.1,4.2);
\filldraw[primedregion=\CColor] (-.15,1.5) rectangle (.4,2);
\filldraw[primedregion=\CColor] (0,2) rectangle (.1,2.5);
\filldraw[primedregion=\CColor] (0,2.6) rectangle (.1,3.8);
\end{scope}
\draw[\QsColor,thick] (-.4,.2) -- (0,.2);
\draw[\QsColor,thick] (-.25,1.2) -- (-.15,1.2);
\draw[\QsColor,thick] (-.1,2.8) -- (0,2.8);
\draw[\QsColor,thick] (-.1,4) -- (0,4);
\draw[\RsColor,thick] (-.15,1.5) -- (.4,1.5);
\draw[\RsColor,thick] (0,2.6) -- (.1,2.6);
\draw[\RsColor,thick] (0,2.5) -- (.1,2.5);
\draw[\RsColor,thick] (0,3.8) -- (.1,3.8);
\draw[\XColor,thick] (-.4,0) -- (-.4,1);
\draw[\YColor,thick] (0,0) -- (0,1);
\draw[\ZColor,thick] (.4,0) -- (.4,2);
\draw[\XColor,thick] (-.25,1) -- (-.25,2);
\draw[\YColor,thick] (-.15,1) -- (-.15,2);
\draw[\XColor,thick] (-.1,2) -- (-.1,4.2);
\draw[\YColor,thick] (0,2) -- (0,4.2);
\draw[\ZColor,thick] (.1,2) -- (.1,4.2);
\roundNbox{unshaded}{(-.2,.7)}{.3}{.15}{.05}{\scriptsize{$F^2_{X,Y}$}};
\roundNbox{unshaded}{(0,2)}{.3}{.35}{.35}{\tiny{$F^2_{X\xxo_Q Y,Z}$}};
\roundNbox{unshaded}{(0,3.3)}{.3}{.25}{.25}{\scriptsize{$F(\alpha^{\cC})$}};
\draw[thin, dotted, rounded corners = 5pt] (-.85,0) rectangle (.85,4.2);
}
=
\tikzmath{
\begin{scope}
\clip[rounded corners=5pt] (-.85,0) rectangle (.85,4.2);
\filldraw[primedregion=white] (-.85,0) rectangle (.85,4.2);
\filldraw[primedregion=\QrColor] (-.4,0) rectangle (0,.5);
\filldraw[primedregion=\QrColor] (-.25,1.5) rectangle (-.15,2);
\filldraw[primedregion=\QrColor] (-.1,2) rectangle (0,2.8);
\filldraw[primedregion=\QrColor] (-.1,4) rectangle (0,4.2);
\filldraw[primedregion=\BColor] (-.4,.5) rectangle (0,1);
\filldraw[primedregion=\BColor] (-.25,1) rectangle (-.15,1.5);
\filldraw[primedregion=\BColor] (-.1,2.8) rectangle (0,4);
\filldraw[primedregion=\RrColor] (0,0) rectangle (.4,.2);
\filldraw[primedregion=\RrColor] (0,3.8) rectangle (.1,4.2);
\filldraw[primedregion=\CColor] (0,2) rectangle (.1,3.8);
\filldraw[primedregion=\CColor] (-.15,2) -- (-.15,1) -- (0,1) -- (0,.2) -- (.4,.2) -- (.4,2);
\end{scope}
\draw[\QsColor,thick] (-.4,.5) -- (0,.5);
\draw[\QsColor,thick] (-.25,1.5) -- (-.15,1.5);
\draw[\QsColor,thick] (-.1,2.8) -- (0,2.8);
\draw[\QsColor,thick] (-.1,4) -- (0,4);
\draw[\RsColor,thick] (0,.2) -- (.4,.2);
\draw[\RsColor,thick] (0,2.5) -- (.1,2.5);
\draw[\RsColor,thick] (0,3.8) -- (.1,3.8);
\draw[\XColor,thick] (-.4,0) -- (-.4,1);
\draw[\YColor,thick] (0,0) -- (0,1);
\draw[\ZColor,thick] (.4,0) -- (.4,2);
\draw[\XColor,thick] (-.25,1) -- (-.25,2);
\draw[\YColor,thick] (-.15,1) -- (-.15,2);
\draw[\XColor,thick] (-.1,2) -- (-.1,4.2);
\draw[\YColor,thick] (0,2) -- (0,4.2);
\draw[\ZColor,thick] (.1,2) -- (.1,4.2);
\roundNbox{unshaded}{(-.2,1)}{.3}{.15}{.05}{\scriptsize{$F^2_{X,Y}$}};
\roundNbox{unshaded}{(0,2)}{.3}{.35}{.35}{\tiny{$F^2_{X\xxo_Q Y,Z}$}};
\roundNbox{unshaded}{(0,3.3)}{.3}{.25}{.25}{\scriptsize{$F(\alpha^{\cC})$}};
\draw[thin, dotted, rounded corners = 5pt] (-.85,0) rectangle (.85,4.2);
}
=
\tikzmath{
\begin{scope}
\clip[rounded corners=5pt] (-.85,0) rectangle (.85,4.2);
\filldraw[primedregion=white] (-.85,0) rectangle (.85,4.2);
\filldraw[primedregion=\QrColor] (-.4,0) rectangle (0,.5);
\filldraw[primedregion=\QrColor] (-.1,4) rectangle (0,4.2);
\filldraw[primedregion=\BColor] (-.4,.5) rectangle (0,1);
\filldraw[primedregion=\BColor] (-.25,1) rectangle (-.15,2);
\filldraw[primedregion=\BColor] (-.1,2) rectangle (0,4);
\filldraw[primedregion=\RrColor] (0,0) rectangle (.4,.2);
\filldraw[primedregion=\RrColor] (0,3.8) rectangle (.1,4.2);
\filldraw[primedregion=\CColor] (0,2) rectangle (.1,3.8);
\filldraw[primedregion=\CColor] (-.15,2) -- (-.15,1) -- (0,1) -- (0,.2) -- (.4,.2) -- (.4,2);
\end{scope}
\draw[\QsColor,thick] (-.4,.5) -- (0,.5);
\draw[\QsColor,thick] (-.25,1.5) -- (-.15,1.5);
\draw[\QsColor,thick] (-.1,4) -- (0,4);
\draw[\RsColor,thick] (0,.2) -- (.4,.2);
\draw[\RsColor,thick] (0,2.5) -- (.1,2.5);
\draw[\RsColor,thick] (0,3.8) -- (.1,3.8);
\draw[\XColor,thick] (-.4,0) -- (-.4,1);
\draw[\YColor,thick] (0,0) -- (0,1);
\draw[\ZColor,thick] (.4,0) -- (.4,2);
\draw[\XColor,thick] (-.25,1) -- (-.25,2);
\draw[\YColor,thick] (-.15,1) -- (-.15,2);
\draw[\XColor,thick] (-.1,2) -- (-.1,4.2);
\draw[\YColor,thick] (0,2) -- (0,4.2);
\draw[\ZColor,thick] (.1,2) -- (.1,4.2);
\roundNbox{unshaded}{(-.2,1)}{.3}{.15}{.05}{\scriptsize{$F^2_{X,Y}$}};
\roundNbox{unshaded}{(0,2)}{.3}{.35}{.35}{\tiny{$F^2_{X\xxo Y,Z}$}};
\roundNbox{unshaded}{(0,3.3)}{.3}{.25}{.25}{\scriptsize{$F(\alpha^{\cC})$}};
\draw[thin, dotted, rounded corners = 5pt] (-.85,0) rectangle (.85,4.2);
}
=
\tikzmath{
\begin{scope}
\clip[rounded corners=5pt] (-.85,0) rectangle (.85,4.2);
\filldraw[primedregion=white] (-.85,0) rectangle (.85,4.2);
\filldraw[primedregion=\QrColor] (-.4,0) rectangle (0,.5);
\filldraw[primedregion=\QrColor] (-.1,4) rectangle (0,4.2);
\filldraw[primedregion=\BColor] (-.4,.5) rectangle (0,1);
\filldraw[primedregion=\BColor] (-.25,1) rectangle (-.15,2);
\filldraw[primedregion=\BColor] (-.1,2) rectangle (0,4);
\filldraw[primedregion=\RrColor] (0,0) rectangle (.4,.2);
\filldraw[primedregion=\RrColor] (0,3.8) rectangle (.1,4.2);
\filldraw[primedregion=\CColor] (0,2) rectangle (.1,3.8);
\filldraw[primedregion=\CColor] (-.15,2) -- (-.15,1) -- (0,1) -- (0,.2) -- (.4,.2) -- (.4,2);
\end{scope}
\draw[\QsColor,thick] (-.4,.5) -- (0,.5);
\draw[\QsColor,thick] (-.1,4) -- (0,4);
\draw[\RsColor,thick] (0,.2) -- (.4,.2);
\draw[\RsColor,thick] (0,2.5) -- (.1,2.5);
\draw[\RsColor,thick] (0,3.8) -- (.1,3.8);
\draw[\XColor,thick] (-.4,0) -- (-.4,1);
\draw[\YColor,thick] (0,0) -- (0,1);
\draw[\ZColor,thick] (.4,0) -- (.4,2);
\draw[\XColor,thick] (-.25,1) -- (-.25,2);
\draw[\YColor,thick] (-.15,1) -- (-.15,2);
\draw[\XColor,thick] (-.1,2) -- (-.1,4.2);
\draw[\YColor,thick] (0,2) -- (0,4.2);
\draw[\ZColor,thick] (.1,2) -- (.1,4.2);
\roundNbox{unshaded}{(-.2,1)}{.3}{.15}{.05}{\scriptsize{$F^2_{X,Y}$}};
\roundNbox{unshaded}{(0,2)}{.3}{.35}{.35}{\tiny{$F^2_{X\xxo Y,Z}$}};
\roundNbox{unshaded}{(0,3.3)}{.3}{.25}{.25}{\scriptsize{$F(\alpha^{\cC})$}};
\draw[thin, dotted, rounded corners = 5pt] (-.85,0) rectangle (.85,4.2);
}
=
\tikzmath{
\begin{scope}
\clip[rounded corners=5pt] (-.85,0) rectangle (.85,4.2);
\filldraw[primedregion=white] (-.85,0) rectangle (.85,4.2);
\filldraw[primedregion=\QrColor] (-.4,0) rectangle (0,.5);
\filldraw[primedregion=\QrColor] (-.1,4) rectangle (0,4.2);
\filldraw[primedregion=\BColor] (-.4,.5) rectangle (0,1);
\filldraw[primedregion=\BColor] (-.25,1) rectangle (-.15,2);
\filldraw[primedregion=\BColor] (-.1,2) rectangle (0,4);
\filldraw[primedregion=\RrColor] (0,0) rectangle (.4,.2);
\filldraw[primedregion=\RrColor] (0,3.8) rectangle (.1,4.2);
\filldraw[primedregion=\CColor] (0,2) rectangle (.1,3.8);
\filldraw[primedregion=\CColor] (-.15,2) -- (-.15,1) -- (0,1) -- (0,.2) -- (.4,.2) -- (.4,2);
\end{scope}
\draw[\QsColor,thick] (-.4,.5) -- (0,.5);
\draw[\QsColor,thick] (-.1,4) -- (0,4);
\draw[\RsColor,thick] (0,.2) -- (.4,.2);
\draw[\RsColor,thick] (0,3.8) -- (.1,3.8);
\draw[\XColor,thick] (-.4,0) -- (-.4,1);
\draw[\YColor,thick] (0,0) -- (0,1);
\draw[\ZColor,thick] (.4,0) -- (.4,2);
\draw[\XColor,thick] (-.25,1) -- (-.25,2);
\draw[\YColor,thick] (-.15,1) -- (-.15,2);
\draw[\XColor,thick] (-.1,2) -- (-.1,4.2);
\draw[\YColor,thick] (0,2) -- (0,4.2);
\draw[\ZColor,thick] (.1,2) -- (.1,4.2);
\roundNbox{unshaded}{(-.2,1)}{.3}{.15}{.05}{\scriptsize{$F^2_{X,Y}$}};
\roundNbox{unshaded}{(0,2)}{.3}{.35}{.35}{\tiny{$F^2_{X\xxo Y,Z}$}};
\roundNbox{unshaded}{(0,3.3)}{.3}{.25}{.25}{\scriptsize{$F(\alpha^{\cC})$}};
\draw[thin, dotted, rounded corners = 5pt] (-.85,0) rectangle (.85,4.2);
}
\\
&=
\tikzmath{
\begin{scope}
\clip[rounded corners = 5] (-.85,0) rectangle (.85,4.4);
\filldraw[primedregion=white] (-.85,0) rectangle (.85,4.4);
\filldraw[primedregion=\QrColor] (-.4,0) rectangle (0,.5);
\filldraw[primedregion=\QrColor] (-.1,4.2) rectangle (0,4.4);
\filldraw[primedregion=\BColor] (-.4,.5) -- (-.4,3.5) -- (.15,3.5) -- (.15,2.3) -- (0,2.3) -- (0,.5);
\filldraw[primedregion=\BColor] (-.1,3.5) rectangle (0,4.2);
\filldraw[primedregion=\RrColor] (0,0) rectangle (.4,.2);
\filldraw[primedregion=\RrColor] (0,4) rectangle (.1,4.4);
\filldraw[primedregion=\CColor] (0,.2) rectangle (.4,2.3);
\filldraw[primedregion=\CColor] (.15,2.3) rectangle (.25,3.5);
\filldraw[primedregion=\CColor] (0,3.5) rectangle (.1,4);
\end{scope}
\draw[\QsColor,thick] (-.4,.5) -- (0,.5);
\draw[\QsColor,thick] (-.1,4.2) -- (0,4.2);
\draw[\RsColor,thick] (0,.2) -- (.4,.2);
\draw[\RsColor,thick] (0,4) -- (.1,4);
\draw[\XColor,thick] (-.4,0) -- (-.4,3.5);
\draw[\YColor,thick] (0,0) -- (0,2);
\draw[\ZColor,thick] (.4,0) -- (.4,2);
\draw[\YColor,thick] (.15,2) -- (.15,3.5);
\draw[\ZColor,thick] (.25,2) -- (.25,3.5);
\draw[\XColor,thick] (-.1,3.5) -- (-.1,4.4);
\draw[\YColor,thick] (0,3.5) -- (0,4.4);
\draw[\ZColor,thick] (.1,3.5) -- (.1,4.4);
\roundNbox{unshaded}{(0,1)}{.3}{.35}{.35}{\scriptsize{$\alpha^\cD$}};
\roundNbox{unshaded}{(.2,2.2)}{.3}{.05}{.15}{\scriptsize{$F^2_{Y,Z}$}};
\roundNbox{unshaded}{(0,3.5)}{.3}{.35}{.35}{\tiny{$F^2_{X,Y\xxo Z}$}};
\draw[thin, dotted, rounded corners = 5pt] (-.85,0) rectangle (.85,4.4);
}
=
\tikzmath{
\begin{scope}
\clip[rounded corners = 5] (-.85,0) rectangle (.85,4.4);
\filldraw[primedregion=white] (-.85,0) rectangle (.85,4.4);
\filldraw[primedregion=\QrColor] (-.4,0) rectangle (0,.5);
\filldraw[primedregion=\QrColor] (-.1,4.2) rectangle (0,4.4);
\filldraw[primedregion=\BColor] (-.4,.5) -- (-.4,3.5) -- (.15,3.5) -- (.15,2.3) -- (0,2.3) -- (0,.5);
\filldraw[primedregion=\BColor] (-.1,3.5) rectangle (0,4.2);
\filldraw[primedregion=\RrColor] (0,0) rectangle (.4,.2);
\filldraw[primedregion=\RrColor] (0,4) rectangle (.1,4.4);
\filldraw[primedregion=\CColor] (0,.2) rectangle (.4,2.3);
\filldraw[primedregion=\CColor] (.15,2.3) rectangle (.25,3.5);
\filldraw[primedregion=\CColor] (0,3.5) rectangle (.1,4);
\end{scope}
\draw[\QsColor,thick] (-.4,.5) -- (0,.5);
\draw[\QsColor,thick] (-.4,3) -- (.15,3);
\draw[\QsColor,thick] (-.1,4.2) -- (0,4.2);
\draw[\RsColor,thick] (0,.2) -- (.4,.2);
\draw[\RsColor,thick] (0,4) -- (.1,4);
\draw[\XColor,thick] (-.4,0) -- (-.4,3.5);
\draw[\YColor,thick] (0,0) -- (0,2);
\draw[\ZColor,thick] (.4,0) -- (.4,2);
\draw[\YColor,thick] (.15,2) -- (.15,3.5);
\draw[\ZColor,thick] (.25,2) -- (.25,3.5);
\draw[\XColor,thick] (-.1,3.5) -- (-.1,4.4);
\draw[\YColor,thick] (0,3.5) -- (0,4.4);
\draw[\ZColor,thick] (.1,3.5) -- (.1,4.4);
\roundNbox{unshaded}{(0,1)}{.3}{.35}{.35}{\scriptsize{$\alpha^\cD$}};
\roundNbox{unshaded}{(.2,2.2)}{.3}{.05}{.15}{\scriptsize{$F^2_{Y,Z}$}};
\roundNbox{unshaded}{(0,3.5)}{.3}{.35}{.35}{\tiny{$F^2_{X,Y\xxo Z}$}};
\draw[thin, dotted, rounded corners = 5pt] (-.85,0) rectangle (.85,4.4);
}
=
\tikzmath{
\begin{scope}
\clip[rounded corners = 5] (-.85,0) rectangle (.85,4.4);
\filldraw[primedregion=white] (-.85,0) rectangle (.85,4.4);
\filldraw[primedregion=\QrColor] (-.4,0) rectangle (0,.5);
\filldraw[primedregion=\QrColor] (-.1,4.2) rectangle (0,4.4);
\filldraw[primedregion=\BColor] (-.4,.5) -- (-.4,3.5) -- (.15,3.5) -- (.15,2.3) -- (0,2.3) -- (0,.5);
\filldraw[primedregion=\BColor] (-.1,3.5) rectangle (0,4.2);
\filldraw[primedregion=\RrColor] (0,0) rectangle (.4,.2);
\filldraw[primedregion=\RrColor] (.15,2.8) rectangle (.25,3.5);
\filldraw[primedregion=\RrColor] (0,3.5) rectangle (.1,4.4);
\filldraw[primedregion=\CColor] (0,.2) rectangle (.4,2.3);
\filldraw[primedregion=\CColor] (.15,2.3) rectangle (.25,2.8);
\end{scope}
\draw[\QsColor,thick] (-.4,.5) -- (0,.5);
\draw[\QsColor,thick] (-.4,3) -- (.15,3);
\draw[\QsColor,thick] (-.1,4.2) -- (0,4.2);
\draw[\RsColor,thick] (0,.2) -- (.4,.2);
\draw[\RsColor,thick] (.15,2.8) -- (.25,2.8);
\draw[\XColor,thick] (-.4,0) -- (-.4,3.5);
\draw[\YColor,thick] (0,0) -- (0,2);
\draw[\ZColor,thick] (.4,0) -- (.4,2);
\draw[\YColor,thick] (.15,2) -- (.15,3.5);
\draw[\ZColor,thick] (.25,2) -- (.25,3.5);
\draw[\XColor,thick] (-.1,3.5) -- (-.1,4.4);
\draw[\YColor,thick] (0,3.5) -- (0,4.4);
\draw[\ZColor,thick] (.1,3.5) -- (.1,4.4);
\roundNbox{unshaded}{(0,1)}{.3}{.35}{.35}{\scriptsize{$\alpha^\cD$}};
\roundNbox{unshaded}{(.2,2.2)}{.3}{.05}{.15}{\scriptsize{$F^2_{Y,Z}$}};
\roundNbox{unshaded}{(0,3.5)}{.3}{.35}{.35}{\tiny{$F^2_{X,Y\xxo_R Z}$}};
\draw[thin, dotted, rounded corners = 5pt] (-.85,0) rectangle (.85,4.4);
}
=
\tikzmath{
\begin{scope}
\clip[rounded corners = 5] (-.85,0) rectangle (.85,4.4);
\filldraw[primedregion=white] (-.85,0) rectangle (.85,4.4);
\filldraw[primedregion=\QrColor] (-.4,0) rectangle (0,.5);
\filldraw[primedregion=\QrColor] (-.1,4.2) rectangle (0,4.4);
\filldraw[primedregion=\BColor] (-.4,.5) -- (-.4,3.5) -- (.15,3.5) -- (.15,2.3) -- (0,2.3) -- (0,.5);
\filldraw[primedregion=\BColor] (-.1,3.5) rectangle (0,4.2);
\filldraw[primedregion=\RrColor] (0,0) rectangle (.4,.2);
\filldraw[primedregion=\RrColor] (.15,2.8) rectangle (.25,3.5);
\filldraw[primedregion=\RrColor] (0,3.5) rectangle (.1,4.4);
\filldraw[primedregion=\CColor] (0,.2) rectangle (.4,2.3);
\filldraw[primedregion=\CColor] (.15,2.3) rectangle (.25,2.8);
\end{scope}
\draw[\QsColor,thick] (-.4,.5) -- (0,.5);
\draw[\QsColor,thick] (-.4,3) -- (.15,3);
\draw[\QsColor,thick] (-.1,4.2) -- (0,4.2);
\draw[\RsColor,thick] (0,.2) -- (.4,.2);
\draw[\RsColor,thick] (0,1.5) -- (.4,1.5);
\draw[\RsColor,thick] (.15,2.8) -- (.25,2.8);
\draw[\XColor,thick] (-.4,0) -- (-.4,3.5);
\draw[\YColor,thick] (0,0) -- (0,2);
\draw[\ZColor,thick] (.4,0) -- (.4,2);
\draw[\YColor,thick] (.15,2) -- (.15,3.5);
\draw[\ZColor,thick] (.25,2) -- (.25,3.5);
\draw[\XColor,thick] (-.1,3.5) -- (-.1,4.4);
\draw[\YColor,thick] (0,3.5) -- (0,4.4);
\draw[\ZColor,thick] (.1,3.5) -- (.1,4.4);
\roundNbox{unshaded}{(0,1)}{.3}{.35}{.35}{\scriptsize{$\alpha^\cD$}};
\roundNbox{unshaded}{(.2,2.2)}{.3}{.05}{.15}{\scriptsize{$F^2_{Y,Z}$}};
\roundNbox{unshaded}{(0,3.5)}{.3}{.35}{.35}{\tiny{$F^2_{X,Y\xxo_R Z}$}};
\draw[thin, dotted, rounded corners = 5pt] (-.85,0) rectangle (.85,4.4);
}
=
\tikzmath{
\begin{scope}
\clip[rounded corners = 5] (-.85,0) rectangle (.85,4.4);
\filldraw[primedregion=white] (-.85,0) rectangle (.85,4.4);
\filldraw[primedregion=\QrColor] (-.4,0) rectangle (0,.5);
\filldraw[primedregion=\QrColor] (-.4,1.8) rectangle (0,3);
\filldraw[primedregion=\QrColor] (-.4,2) rectangle (.15,3);
\filldraw[primedregion=\QrColor] (-.1,4.2) rectangle (0,4.4);
\filldraw[primedregion=\BColor] (-.4,.5) rectangle (0,1.8);
\filldraw[primedregion=\BColor] (-.4,3) rectangle (.15,3.5);
\filldraw[primedregion=\BColor] (-.1,3.5) rectangle (0,4.2);
\filldraw[primedregion=\RrColor] (0,0) rectangle (.4,.2);
\filldraw[primedregion=\RrColor] (0,1.5) rectangle (.4,1.6);
\filldraw[primedregion=\RrColor] (.15,2.8) rectangle (.25,3.5);
\filldraw[primedregion=\RrColor] (0,3.5) rectangle (.1,4.4);
\filldraw[primedregion=\CColor] (0,.2) rectangle (.4,1.5);
\filldraw[primedregion=\CColor] (0,1.6) rectangle (.4,2.3);
\filldraw[primedregion=\CColor] (.15,2.3) rectangle (.25,2.8);
\end{scope}
\draw[\QsColor,thick] (-.4,.5) -- (0,.5);
\draw[\QsColor,thick] (-.4,1.8) -- (0,1.8);
\draw[\QsColor,thick] (-.4,3) -- (.15,3);
\draw[\QsColor,thick] (-.1,4.2) -- (0,4.2);
\draw[\RsColor,thick] (0,.2) -- (.4,.2);
\draw[\RsColor,thick] (0,1.5) -- (.4,1.5);
\draw[\RsColor,thick] (0,1.6) -- (.4,1.6);
\draw[\RsColor,thick] (.15,2.8) -- (.25,2.8);
\draw[\XColor,thick] (-.4,0) -- (-.4,3.5);
\draw[\YColor,thick] (0,0) -- (0,2);
\draw[\ZColor,thick] (.4,0) -- (.4,2);
\draw[\YColor,thick] (.15,2) -- (.15,3.5);
\draw[\ZColor,thick] (.25,2) -- (.25,3.5);
\draw[\XColor,thick] (-.1,3.5) -- (-.1,4.4);
\draw[\YColor,thick] (0,3.5) -- (0,4.4);
\draw[\ZColor,thick] (.1,3.5) -- (.1,4.4);
\roundNbox{unshaded}{(0,1)}{.3}{.35}{.35}{\scriptsize{$\alpha^\cD$}};
\roundNbox{unshaded}{(.2,2.3)}{.3}{.05}{.15}{\scriptsize{$F^2_{Y,Z}$}};
\roundNbox{unshaded}{(0,3.5)}{.3}{.35}{.35}{\tiny{$F^2_{X,Y\xxo_R Z}$}};
\draw[thin, dotted, rounded corners = 5pt] (-.85,0) rectangle (.85,4.4);
}
=
\tikzmath[scale=.5, transform shape]{
\begin{scope}
\clip[rounded corners = 5] (-1.7,0) rectangle (1.9,8);
\filldraw[primedregion=white] (-1.7,0) rectangle (1.9,8);
\filldraw[primedregion=\QrColor] (-1,0) rectangle (.4,6);
\filldraw[primedregion=\RrColor] (0,0) rectangle (1,4);
\filldraw[primedregion=\RrColor] (.6,4) rectangle (.4,6);
\filldraw[primedregion=\QrColor] (-.2,6) rectangle (0,8);
\filldraw[primedregion=\RrColor] (0,6) rectangle (.2,8);
\end{scope}
\draw[\XColor,thick] (-1,0) -- (-1,6);
\draw[\YColor,thick] (0,0) -- (0,4);
\draw[\ZColor,thick] (1,0) -- (1,4);
\draw[\YColor,thick] (.4,4) -- (.4,6);
\draw[\ZColor,thick] (.6,4) -- (.6,6);
\draw[\XColor,thick] (-.2,6) -- (-.2,8);
\draw[\YColor,thick] (0,6) -- (0,8);
\draw[\ZColor,thick] (.2,6) -- (.2,8);
\roundNbox{unshaded}{(0,2)}{.6}{.75}{.95}{\normalsize{$\alpha^{\QSys(\cD)}$}};
\roundNbox{unshaded}{(.5,4)}{.6}{.45}{.45}{\small{$\QSys(F)^2_{Y,Z}$}};
\roundNbox{unshaded}{(0,6)}{.6}{.75}{.95}{\small{$\QSys(F)^2_{X,Y\xxo_R Z}$}};
\draw[thin, dotted, rounded corners = 5pt] (-1.7,0) rectangle (1.9,8);
}
\\
\tikzmath[scale=.5, transform shape]{
\begin{scope}
\clip[rounded corners=5pt] (-1.4,-1.4) rectangle (1.4,2);
\filldraw[primedregion=white] (-1.4,-1.4) rectangle (1.4,2);
\filldraw[primedregion=\QrColor] (-.4,-1.4) rectangle (.4,0);
\filldraw[primedregion=\QrColor] (.1,0) -- (.1,1.2) arc (0:90:.2cm) -- (-.1,0);
\end{scope}
\draw[\XColor,thick] (-.4,-1.4) -- (-.4,0);
\draw[\QsColor,thick] (.4,-1.4) -- (.4,0);
\draw[\XColor,thick] (-.1,2) -- (-.1,0);
\draw[\QsColor,thick] (.1,0) -- (.1,1.2) arc (0:90:.2cm);
\filldraw[\XColor] (-.1,1.4) circle (.06cm);
\roundNbox{unshaded}{(0,0)}{.6}{.45}{.45}{\small{$\QSys(F)^2_{X,Q}$}};
\draw[thin, dotted, rounded corners = 5pt] (-1.4,-1.4) rectangle (1.4,2);
}
&=
\tikzmath{
\begin{scope}
\clip[rounded corners=5pt] (-.7,-.7) rectangle (.7,1);
\filldraw[primedregion=white] (-.7,-.7) rectangle (.7,1);
\filldraw[primedregion=\BColor] (-.05,0) -- (.05,0) -- (.05,.7) arc (0:90:.1cm);
\filldraw[primedregion=\QrColor] (-.05,.5) rectangle (.05,.6);
\filldraw[primedregion=\BColor] (-.2,0) rectangle (.2,-.5);
\filldraw[primedregion=\QrColor] (-.2,-.7) rectangle (.2,-.5);
\end{scope}
\draw[\QsColor,thick] (-.2,-.5) -- (.2,-.5);
\draw[\QsColor,thick] (-.05,.5) -- (.05,.5);
\draw[\QsColor,thick] (-.05,.6) -- (.05,.6);
\draw[\XColor,thick] (-.2,-.7) -- (-.2,0);
\draw[\QsColor,thick] (.2,-.7) -- (.2,0);
\draw[\XColor,thick] (-.05,1) -- (-.05,0);
\draw[\QsColor,thick] (.05,0) -- (.05,.7) arc (0:90:.1cm);
\filldraw[\XColor] (-.05,.8) circle (.03cm);
\roundNbox{unshaded}{(0,0)}{.3}{.2}{.2}{\scriptsize{$F^2_{X,Q}$}};
\draw[thin, dotted, rounded corners = 5pt] (-.7,-.7) rectangle (.7,1);
}
=
\tikzmath{
\begin{scope}
\clip[rounded corners=5pt] (-.7,-.7) rectangle (.7,1);
\filldraw[primedregion=white] (-.7,-.7) rectangle (.7,1);
\filldraw[primedregion=\BColor] (-.05,0) -- (.05,0) -- (.05,.6) arc (0:90:.1cm);
\filldraw[primedregion=\QrColor] (-.2,-.7) rectangle (.2,-.5);
\filldraw[primedregion=\BColor] (-.2,0) rectangle (.2,-.5);
\end{scope}
\draw[\QsColor,thick] (-.2,-.5) -- (.2,-.5);
\draw[\QsColor,thick] (-.05,.5) -- (.05,.5);
\draw[\XColor,thick] (-.2,-.7) -- (-.2,0);
\draw[\QsColor,thick] (.2,-.7) -- (.2,0);
\draw[\XColor,thick] (-.05,1) -- (-.05,0);
\draw[\QsColor,thick] (.05,0) -- (.05,.6) arc (0:90:.1cm);
\filldraw[\XColor] (-.05,.7) circle (.03cm);
\roundNbox{unshaded}{(0,0)}{.3}{.2}{.2}{\scriptsize{$F^2_{X,Q}$}};
\draw[thin, dotted, rounded corners = 5pt] (-.7,-.7) rectangle (.7,1);
}
=
\tikzmath{
\begin{scope}
\clip[rounded corners=5pt] (-.7,-.7) rectangle (.7,1);
\filldraw[primedregion=white] (-.7,-.7) rectangle (.7,1);
\filldraw[primedregion=\BColor] (-.05,0) -- (.05,0) -- (.05,.6) arc (0:90:.1cm);
\filldraw[primedregion=\QrColor] (-.2,-.7) rectangle (.2,-.5);
\filldraw[primedregion=\BColor] (-.2,0) rectangle (.2,-.5);
\end{scope}
\draw[\QsColor,thick] (-.2,-.5) -- (.2,-.5);
%
\draw[\XColor,thick] (-.2,-.7) -- (-.2,0);
\draw[\QsColor,thick] (.2,-.7) -- (.2,0);
\draw[\XColor,thick] (-.05,1) -- (-.05,0);
\draw[\QsColor,thick] (.05,0) -- (.05,.6) arc (0:90:.1cm);
\filldraw[\XColor] (-.05,.7) circle (.03cm);
\roundNbox{unshaded}{(0,0)}{.3}{.2}{.2}{\scriptsize{$F^2_{X,Q}$}};
\draw[thin, dotted, rounded corners = 5pt] (-.7,-.7) rectangle (.7,1);
}
=
\tikzmath{
\begin{scope}
\clip[rounded corners=5pt] (-.3,0) rectangle (.7,1.7);
\filldraw[primedregion=white] (-.3,0) rectangle (.7,1.7);
\filldraw[primedregion=\BColor] (0,.5) -- (.4,.5) -- (.4,.7) arc (0:90:.4cm);
\filldraw[primedregion=\QrColor] (0,0) rectangle (.4,.5);
\end{scope}
\draw[\QsColor,thick] (0,.5) -- (.4,.5);
\draw[\XColor,thick] (0,0) -- (0,1.7);
\draw[\QsColor,thick] (.4,0) -- (.4,.7) arc (0:90:.4cm);
\filldraw[\XColor] (0,1.1) circle (.05cm);
\draw[thin, dotted, rounded corners = 5pt] (-.3,0) rectangle (.7,1.7);
}
=
\tikzmath{
\begin{scope}
\clip[rounded corners=5pt] (-.3,0) rectangle (.7,1.7);
\filldraw[primedregion=white] (-.3,0) rectangle (.7,1.7);
\filldraw[primedregion=\QrColor] (.4,0) -- (.4,.7) arc (0:90:.4cm) -- (0,0);
\end{scope}
\draw[\XColor,thick] (0,0) -- (0,1.7);
\draw[\QsColor,thick] (.4,0) -- (.4,.7) arc (0:90:.4cm);
\filldraw[\XColor] (0,1.1) circle (.05cm);
\draw[thin, dotted, rounded corners = 5pt] (-.3,0) rectangle (.7,1.7);
}\,;
\quad\text{similarly,}\quad
\tikzmath[scale=.5, transform shape]{
\begin{scope}
\clip[rounded corners=5pt] (-1.4,-1.4) rectangle (1.4,2);
\filldraw[primedregion=white] (-1.4,-1.4) rectangle (1.4,2);
\filldraw[primedregion=\QrColor] (-.4,-1.4) rectangle (.4,0);
\filldraw[primedregion=\QrColor] (-.1,0) -- (-.1,1.2) arc (180:90:.2cm) -- (.1,0);
\end{scope}
\draw[\QsColor,thick] (-.4,-1.4) -- (-.4,0);
\draw[\YColor,thick] (.4,-1.4) -- (.4,0);
\draw[\QsColor,thick] (-.1,0) -- (-.1,1.2) arc (180:90:.2cm);
\draw[\YColor,thick] (.1,0) -- (.1,2);
\filldraw[\YColor] (.1,1.4) circle (.06cm);
\roundNbox{unshaded}{(0,0)}{.6}{.45}{.45}{\small{$\QSys(F)^2_{Q,Y}$}};
\draw[thin, dotted, rounded corners = 5pt] (-1.4,-1.4) rectangle (1.4,2);
}
=
\tikzmath{
\begin{scope}
\clip[rounded corners=5pt] (-.3,0) rectangle (.7,1.7);
\filldraw[primedregion=white] (-.3,0) rectangle (.7,1.7);
\filldraw[primedregion=\QrColor] (0,0) -- (0,.7) arc (180:90:.4cm) -- (.4,0);
\end{scope}
\draw[\QsColor,thick] (0,0) -- (0,.7) arc (180:90:.4cm);
\draw[\YColor,thick] (.4,0) -- (.4,1.7);
\filldraw[\YColor] (.4,1.1) circle (.05cm);
\draw[thin, dotted, rounded corners = 5pt] (-.3,0) rectangle (.7,1.7);
}\,.
\end{align*}
\end{construction}

For the rest of this section, we fix two $\dag$ 2-functors $F,G: \cC \to \cD$.

\begin{rem}
Suppose $F,G: \cC\to \cD$ are two $\dag$ 2-functors and
$\varphi: F\Rightarrow G$ is a $\dag$ 2-transformation.
For any dualizable 1-cell $X\in\cC(a\to b)$, we have
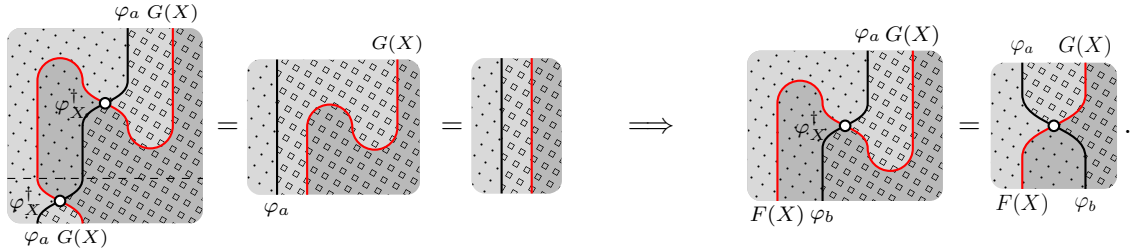
\[
\tikzmath{
\begin{scope}
\clip[rounded corners = 5] (-1.3,-1.6) rectangle (1.3,1);
\filldraw[primedregion=\BColor] (-.9,-1.6) .. controls ++(90:.3cm) and ++(270:.3cm) .. (-.3,-1) -- (-.3,-.3) .. controls ++(90:.3cm) and ++(270:.3cm) .. (.3,.3) -- (.3,1) -- (-1.3,1) -- (-1.3,-1.6);
\filldraw[primedregion=\AColor] (-.9,-1.6) .. controls ++(90:.2cm) and ++(-135:.1cm) .. (-.6,-1.3) .. controls ++(135:.1cm) and ++(270:.2cm) .. (-.9,-1) -- (-.9,.3) arc (180:0:.3cm) .. controls ++(270:.2cm) and ++(135:.1cm) .. (0,0) .. controls ++(45:.1cm) and ++(270:.2cm) .. (.3,.3) -- (.3,1) -- (-1.3,1) -- (-1.3,-1.6);
\filldraw[boxregion=\BColor] (-.9,-1.6) .. controls ++(90:.3cm) and ++(270:.3cm) .. (-.3,-1) -- (-.3,-.3) .. controls ++(90:.3cm) and ++(270:.3cm) .. (.3,.3) -- (.3,1) -- (1.3,1) -- (1.3,-1.6);
\filldraw[boxregion=\AColor] (-.9,-1.6) .. controls ++(90:.2cm) and ++(-135:.1cm) .. (-.6,-1.3) .. controls ++(-45:.1cm) and ++(90:.2cm) .. (-.3,-1.6);
\filldraw[boxregion=\AColor] (.9,1) -- (.9,-.3) arc (0:-180:.3cm) .. controls ++(90:.2cm) and ++(-45:.1cm) .. (0,0) .. controls ++(45:.1cm) and ++(270:.2cm) .. (.3,.3) -- (.3,1);
\end{scope}
\draw[dashed] (-1.3,-1) -- (1.3,-1);
\draw[\phiColor,thick] (-.9,-1.6) .. controls ++(90:.3cm) and ++(270:.3cm) .. (-.3,-1) -- (-.3,-.3) .. controls ++(90:.3cm) and ++(270:.3cm) .. (.3,.3) -- (.3,1);
\draw[\XColor,thick] (.9,1) -- (.9,-.3) arc (0:-180:.3cm) .. controls ++(90:.3cm) and ++(270:.3cm) .. (-.3,.3) arc (0:180:.3cm) -- (-.9,-1);
\draw[\XColor,thick] (-.3,-1.6) .. controls ++(90:.3cm) and ++(270:.3cm) .. (-.9,-1);
\filldraw[white] (0,0) circle (.07cm);
\draw[thick] (0,0) circle (.07cm); 
\filldraw[white] (-.6,-1.3) circle (.07cm);
\draw[thick] (-.6,-1.3) circle (.07cm); 
\node at (.9,1.2) {\tiny{$G(X)$}};
\node at (.3,1.2) {\scriptsize{$\varphi_a$}};
\node at (-.3,-1.8) {\tiny{$G(X)$}};
\node at (-.9,-1.8) {\scriptsize{$\varphi_a$}};
\node at (-.45,0) {\scriptsize{$\varphi_X^\dag$}};
\node at (-1.05,-1.3) {\scriptsize{$\varphi_X^\dag$}};
}
=
\tikzmath{
\begin{scope}
\clip[rounded corners = 5] (-.8,-.9) rectangle (1.5,.9);
\filldraw[primedregion=\AColor] (-.8,-.9) rectangle (-.4,.9);
\filldraw[boxregion=\AColor] (0,-.9) -- (0,0) arc (180:0:.3cm) arc (-180:0:.3cm) -- (1.2,.9) -- (-.4,.9) -- (-.4,-.9);
\filldraw[boxregion=\BColor] (0,-.9) -- (0,0) arc (180:0:.3cm) arc (-180:0:.3cm) -- (1.2,.9) -- (1.5,.9) -- (1.5,-.9);
\end{scope}
\draw[\phiColor,thick] (-.4,-.9) -- (-.4,.9);
\draw[\XColor,thick] (0,-.9) -- (0,0) arc (180:0:.3cm) arc (-180:0:.3cm) -- (1.2,.9);
\node at (-.4,-1.1) {\scriptsize{$\varphi_a$}};
\node at (1.2,1.1) {\tiny{$G(X)$}};
}
=
\tikzmath{
\begin{scope}
\clip[rounded corners = 5] (-.6,-.9) rectangle (.6,.9);
\filldraw[primedregion=\AColor] (-.6,-.9) rectangle (-.2,.9);
\filldraw[boxregion=\AColor] (-.2,-.9) rectangle (.2,.9);
\filldraw[boxregion=\BColor] (.2,-.9) rectangle (.6,.9);
\end{scope}
\draw[\phiColor,thick] (-.2,-.9) -- (-.2,.9);
\draw[\XColor,thick] (.2,-.9) -- (.2,.9);
}
\qquad
\Longrightarrow
\qquad
\tikzmath{
\begin{scope}
\clip[rounded corners = 5] (-1.3,-1) rectangle (1.3,1);
\filldraw[primedregion=\BColor] (-.3,-1) -- (-.3,-.3) .. controls ++(90:.3cm) and ++(270:.3cm) .. (.3,.3) -- (.3,1) -- (-1.3,1) -- (-1.3,-1);
\filldraw[primedregion=\AColor] (-.9,-1) -- (-.9,.3) arc (180:0:.3cm) .. controls ++(270:.2cm) and ++(135:.1cm) .. (0,0) .. controls ++(45:.1cm) and ++(270:.2cm) .. (.3,.3) -- (.3,1) -- (-1.3,1) -- (-1.3,-1);
\filldraw[boxregion=\BColor] (-.3,-1) -- (-.3,-.3) .. controls ++(90:.3cm) and ++(270:.3cm) .. (.3,.3) -- (.3,1) -- (1.3,1) -- (1.3,-1);
\filldraw[boxregion=\AColor] (.9,1) -- (.9,-.3) arc (0:-180:.3cm) .. controls ++(90:.2cm) and ++(-45:.1cm) .. (0,0) .. controls ++(45:.1cm) and ++(270:.2cm) .. (.3,.3) -- (.3,1);
\end{scope}
\draw[\phiColor,thick] (-.3,-1) -- (-.3,-.3) .. controls ++(90:.3cm) and ++(270:.3cm) .. (.3,.3) -- (.3,1);
\draw[\XColor,thick] (.9,1) -- (.9,-.3) arc (0:-180:.3cm) .. controls ++(90:.3cm) and ++(270:.3cm) .. (-.3,.3) arc (0:180:.3cm) -- (-.9,-1);
\filldraw[white] (0,0) circle (.07cm);
\draw[thick] (0,0) circle (.07cm); 
\node at (.9,1.2) {\scriptsize{$G(X)$}};
\node at (.3,1.2) {\scriptsize{$\varphi_a$}};
\node at (-.3,-1.2) {\scriptsize{$\varphi_b$}};
\node at (-.9,-1.2) {\scriptsize{$F(X)$}};
\node at (-.45,0) {\scriptsize{$\varphi_X^\dag$}};
}
=
\tikzmath[scale=.7, transform shape]{
\begin{scope}
\clip[rounded corners = 5] (-.6,0) rectangle (1.8,2.4);
\filldraw[primedregion=\AColor] (0,0) -- (0,.6) .. controls ++(90:.4cm) and ++(-135:.2cm) .. (.6,1.2) .. controls ++(135:.2cm) and ++(270:.4cm) .. (0,1.8) -- (0,3) -- (-.6,3) -- (-.6,0); 
\filldraw[primedregion=\BColor] (1.2,0) -- (1.2,.6) .. controls ++(90:.4cm) and ++(-45:.2cm) .. (.6,1.2) .. controls ++(-135:.2cm) and ++(90:.4cm) .. (0,.6) -- (0,0);
\filldraw[boxregion=\AColor] (0,3) -- (0,1.8) .. controls ++(270:.4cm) and ++(135:.2cm) .. (.6,1.2) .. controls ++(45:.2cm) and ++(270:.4cm) .. (1.2,1.8) -- (1.2,3);
\filldraw[boxregion=\BColor] (1.2,0) -- (1.2,.6) .. controls ++(90:.4cm) and ++(-45:.2cm) .. (.6,1.2) .. controls ++(45:.2cm) and ++(270:.4cm) .. (1.2,1.8) -- (1.2,3) -- (1.8,3) -- (1.8,0);
\end{scope}
\draw[\XColor,thick] (0,0) -- (0,.6) .. controls ++(90:.6cm) and ++(270:.6cm) .. (1.2,1.8) -- (1.2,2.4);
\draw[\phiColor,thick] (1.2,0) -- (1.2,.6) .. controls ++(90:.6cm) and ++(270:.6cm) .. (0,1.8) -- (0,2.4);
\filldraw[white] (.6,1.2) circle (.1cm);
\draw[thick] (.6,1.2) circle (.1cm); 
\node at (0,-.3) {\large{$F(X)$}};
\node at (1.2,2.7) {\large{$G(X)$}};
\node at (1.2,-.3) {\large{$\varphi_b$}};
\node at (0,2.7) {\large{$\varphi_a$}};
}\,.
\]
\end{rem}

\begin{construction}
\label{construction:QSys(phi)}
Given a $\dag$-transformation $\varphi : F\Rightarrow G$, 
we define a $\dag$-transformation $\QSys(\varphi):\QSys(F)\Rightarrow \QSys(G)$. 
In the diagrams below, we suppress all coherence isomorphisms for $F$ and $G$.

For a Q-system $({}_bQ_b,m,i)\in\QSys(\cC)$, 
we define $\QSys(\varphi)_Q$ by orthogonally splitting the orthgonal projection which appears as the second diagram below.
The other diagrams in the next two rows prove that this second diagram is an orthogonal projection.
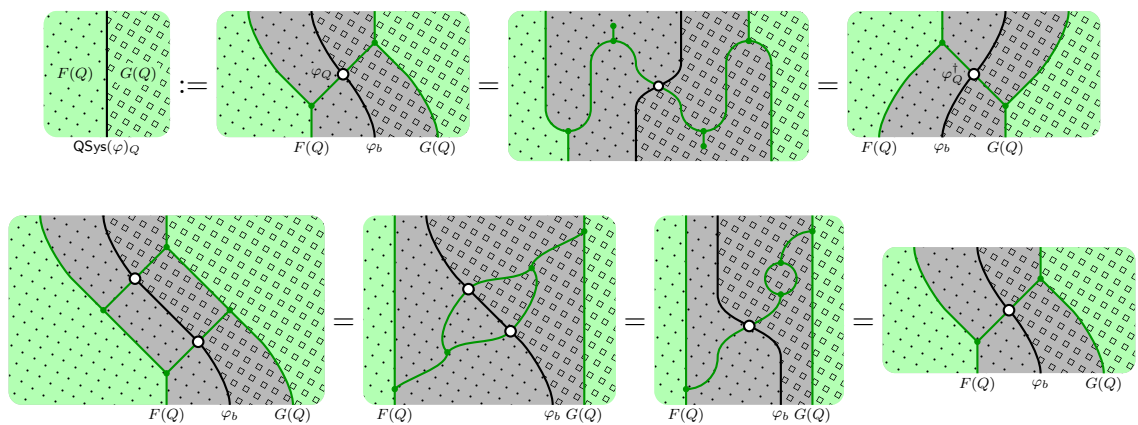
\[
\tikzmath[scale=.7, transform shape]{
\begin{scope}
\clip[rounded corners = 5] (-1.2,-1.2) rectangle (1.2,1.2);
\filldraw[primedregion=\QrColor] (-1.2,-1.2) rectangle (0,1.2);
\filldraw[boxregion=\QrColor] (1.2,1.2) rectangle (0,-1.2);
\filldraw[\QrColor] (-1,-.2) rectangle (-.2,.2);
\filldraw[\QrColor] (.2,-.2) rectangle (1,.2);
\end{scope}
\draw[\phiColor,thick] (0,-1.2) -- (0,1.2);
\node at (.6,0) {\scriptsize{$G(Q)$}};
\node at (-.6,0) {\scriptsize{$F(Q)$}};
\node at (0,-1.4) {\scriptsize{$\QSys(\varphi)_Q$}};
}
:=
\tikzmath[scale=.7, transform shape]{
\begin{scope}
\clip[rounded corners = 5] (-2.4,-1.2) rectangle (2.4,1.2);
\filldraw[primedregion=\BColor] (.6,-1.2) .. controls ++(90:.4cm) and ++(-45:.2cm) .. (0,0) .. controls ++(135:.2cm) and ++(270:.4cm) .. (-.6,1.2) -- (-2.4,1.2) -- (-2.4,-1.2);
\filldraw[boxregion=\BColor] (.6,-1.2) .. controls ++(90:.4cm) and ++(-45:.2cm) .. (0,0) .. controls ++(135:.2cm) and ++(270:.4cm) .. (-.6,1.2) -- (2.4,1.2) -- (2.4,-1.2);
\filldraw[primedregion=\QrColor] (-.6,-1.2) -- (-.6,-.6) -- (-1.2,0) .. controls ++(135:.2cm) and ++(270:.4cm) .. (-1.8,1.2) -- (-2.4,1.2) -- (-2.4,-1.2);
\filldraw[boxregion=\QrColor] (1.8,-1.2) .. controls ++(90:.4cm) and ++(-45:.2cm) .. (1.2,0) -- (.6,.6) -- (.6,1.2) -- (2.4,1.2) -- (2.4,-1.2);
\filldraw[\BColor] (-.6,-.2) rectangle (-.2,.2);
\end{scope}
\draw[\phiColor,thick] (.6,-1.2) .. controls ++(90:.4cm) and ++(-45:.2cm) .. (0,0) .. controls ++(135:.2cm) and ++(270:.4cm) .. (-.6,1.2);
\draw[\QsColor,thick] (-.6,-1.2) -- (-.6,-.6) -- (-1.2,0) .. controls ++(135:.2cm) and ++(270:.4cm) .. (-1.8,1.2);
\draw[\QsColor,thick] (1.8,-1.2) .. controls ++(90:.4cm) and ++(-45:.2cm) .. (1.2,0) -- (.6,.6) -- (.6,1.2);
\draw[\QsColor,thick] (-.6,-.6) -- (.6,.6);
\filldraw[\QsColor] (.6,.6) circle (.05cm);
\filldraw[\QsColor] (-.6,-.6) circle (.05cm);
\filldraw[white] (0,0) circle (.1cm);
\draw[thick] (0,0) circle (.1cm); 
\node at (-.6,-1.4) {\scriptsize{$F(Q)$}};
\node at (.6,-1.4) {\scriptsize{$\varphi_b$}};
\node at (1.8,-1.4) {\scriptsize{$G(Q)$}};
\node at (-.4,0) {\scriptsize{$\varphi_Q$}};
}
=
\tikzmath[scale=.5, transform shape]{
\begin{scope}
\clip[rounded corners = 5] (-4,-2) rectangle (4,2);
\filldraw[primedregion=\BColor] (-.6,-2) -- (-.6,-.6) .. controls ++(90:.6cm) and ++(270:.6cm) .. (.6,.6) -- (.6,2) -- (-4,2) -- (-4,-2);
\filldraw[boxregion=\BColor] (-.6,-2) -- (-.6,-.6) .. controls ++(90:.6cm) and ++(270:.6cm) .. (.6,.6) -- (.6,2) -- (4,2) -- (4,-2);
\filldraw[primedregion=\QrColor] (-2.4,-2) -- (-2.4,-1.2) arc (-90:-180:.6cm) -- (-3,2) --  (-4,2) -- (-4,-2);
\filldraw[boxregion=\QrColor] (2.4,2) -- (2.4,1.2) arc (90:0:.6cm) -- (3,-2) --  (4,-2) -- (4,2); 
\end{scope}
\draw[\phiColor,thick] (-.6,-2) -- (-.6,-.6) .. controls ++(90:.6cm) and ++(270:.6cm) .. (.6,.6) -- (.6,2); 
\draw[\QsColor,thick] (3,-2) -- (3,.6) arc (0:180:.6cm) -- (1.8,-.6) arc (0:-180:.6cm) .. controls ++(90:.6cm) and ++(270:.6cm) .. (-.6,.6) arc (0:180:.6cm) -- (-1.8,-.6) arc (0:-180:.6cm) -- (-3,2);
\draw[\QsColor,thick] (2.4,1.2) -- (2.4,2);
\draw[\QsColor,thick] (-2.4,-1.2) -- (-2.4,-2);
\draw[\QsColor,thick] (1.2,-1.2) -- (1.2,-1.6);
\draw[\QsColor,thick] (-1.2,1.2) -- (-1.2,1.6);
\filldraw[\QsColor] (2.4,1.2) circle (.07cm);
\filldraw[\QsColor] (-2.4,-1.2) circle (.07cm);
\filldraw[\QsColor] (1.2,-1.2) circle (.07cm);
\filldraw[\QsColor] (1.2,-1.6) circle (.07cm);
\filldraw[\QsColor] (-1.2,1.2) circle (.07cm);
\filldraw[\QsColor] (-1.2,1.6) circle (.07cm);
\filldraw[white] (0,0) circle (.12cm);
\draw[thick] (0,0) circle (.12cm); 
}
=
\tikzmath[scale=.7, transform shape]{
\begin{scope}
\clip[rounded corners = 5] (-2.4,-1.2) rectangle (2.4,1.2);
\filldraw[primedregion=\BColor] (-.6,-1.2) .. controls ++(90:.4cm) and ++(-135:.2cm) .. (0,0) .. controls ++(45:.2cm) and ++(270:.4cm) .. (.6,1.2) -- (-2.4,1.2) -- (-2.4,-1.2);
\filldraw[boxregion=\BColor] (-.6,-1.2) .. controls ++(90:.4cm) and ++(-135:.2cm) .. (0,0) .. controls ++(45:.2cm) and ++(270:.4cm) .. (.6,1.2) -- (2.4,1.2) -- (2.4,-1.2);
\filldraw[primedregion=\QrColor] (-1.8,-1.2) .. controls ++(90:.4cm) and ++(-135:.2cm) .. (-1.2,0) -- (-.6,.6) -- (-.6,1.2) -- (-2.4,1.2) -- (-2.4,-1.2);
\filldraw[boxregion=\QrColor] (.6,-1.2) -- (.6,-.6) -- (1.2,0) .. controls ++(45:.2cm) and ++(270:.4cm) .. (1.8,1.2) -- (2.4,1.2) -- (2.4,-1.2);
\filldraw[\BColor] (-.6,-.2) rectangle (-.2,.2);
\end{scope}
\draw[\phiColor,thick] (-.6,-1.2) .. controls ++(90:.4cm) and ++(-135:.2cm) .. (0,0) .. controls ++(45:.2cm) and ++(270:.4cm) .. (.6,1.2);
\draw[\QsColor,thick] (.6,-1.2) -- (.6,-.6) -- (1.2,0) .. controls ++(45:.2cm) and ++(270:.4cm) .. (1.8,1.2);
\draw[\QsColor,thick] (-1.8,-1.2) .. controls ++(90:.4cm) and ++(-135:.2cm) .. (-1.2,0) -- (-.6,.6) -- (-.6,1.2);
\draw[\QsColor,thick] (.6,-.6) -- (-.6,.6);
\filldraw[\QsColor] (.6,-.6) circle (.05cm);
\filldraw[\QsColor] (-.6,.6) circle (.05cm);
\filldraw[white] (0,0) circle (.1cm);
\draw[thick] (0,0) circle (.1cm); 
\node at (-1.8,-1.4) {\scriptsize{$F(Q)$}};
\node at (-.6,-1.4) {\scriptsize{$\varphi_b$}};
\node at (.6,-1.4) {\scriptsize{$G(Q)$}};
\node at (-.4,0) {\scriptsize{$\varphi_Q^\dag$}};
}
\]
\[
\tikzmath[scale=.7, transform shape]{
\begin{scope}
\clip[rounded corners = 5] (-3,-1.8) rectangle (3,1.8);
\filldraw[primedregion=\BColor] (1.2,-1.8) .. controls ++(90:.4cm) and ++(-45:.2cm) .. (.6,-.6) -- (-.6,.6) .. controls ++(135:.2cm) and ++(270:.4cm) .. (-1.2,1.8) -- (-3,1.8) -- (-3,-1.8);
\filldraw[boxregion=\BColor] (1.2,-1.8) .. controls ++(90:.4cm) and ++(-45:.2cm) .. (.6,-.6) -- (-.6,.6) .. controls ++(135:.2cm) and ++(270:.4cm) .. (-1.2,1.8) -- (3,1.8) -- (3,-1.8);
\filldraw[primedregion=\QrColor] (0,-1.8) -- (0,-1.2) -- (-1.8,.6) .. controls ++(135:.2cm) and ++(270:.4cm) .. (-2.4,1.8) -- (-3,1.8) -- (-3,-1.8);
\filldraw[boxregion=\QrColor] (2.4,-1.8) .. controls ++(90:.4cm) and ++(-45:.2cm) .. (1.8,-.6) -- (0,1.2) -- (0,1.8) -- (3,1.8) -- (3,-1.8);
\end{scope} 
\draw[\phiColor,thick] (1.2,-1.8) .. controls ++(90:.4cm) and ++(-45:.2cm) .. (.6,-.6) -- (-.6,.6) .. controls ++(135:.2cm) and ++(270:.4cm) .. (-1.2,1.8);
\draw[\QsColor,thick] (0,-1.8) -- (0,-1.2) -- (-1.8,.6) .. controls ++(135:.2cm) and ++(270:.4cm) .. (-2.4,1.8);
\draw[\QsColor,thick] (2.4,-1.8) .. controls ++(90:.4cm) and ++(-45:.2cm) .. (1.8,-.6) -- (0,1.2) -- (0,1.8);
\draw[\QsColor,thick] (0,-1.2) -- (1.2,0);
\draw[\QsColor,thick] (0,1.2) -- (-1.2,0);
\filldraw[\QsColor] (0,-1.2) circle (.05cm);
\filldraw[\QsColor] (-1.2,0) circle (.05cm);
\filldraw[\QsColor] (0,1.2) circle (.05cm);
\filldraw[\QsColor] (1.2,0) circle (.05cm);
\filldraw[white] (.6,-.6) circle (.1cm);
\draw[thick] (.6,-.6) circle (.1cm); 
\filldraw[white] (-.6,.6) circle (.1cm);
\draw[thick] (-.6,.6) circle (.1cm); 
\node at (0,-2) {\scriptsize{$F(Q)$}};
\node at (1.2,-2) {\scriptsize{$\varphi_b$}};
\node at (2.4,-2) {\scriptsize{$G(Q)$}};
}
=
\tikzmath[scale=.7, transform shape]{
\begin{scope}
\clip[rounded corners = 5] (-2.4,-1.8) rectangle (2.4,1.8);
\filldraw[primedregion=\BColor] (1.2,-1.8) .. controls ++(90:.4cm) and ++(-45:.2cm) .. (.6,-.6) -- (-.6,.6) .. controls ++(135:.2cm) and ++(270:.4cm) .. (-1.2,1.8) -- (-1.8,1.8) -- (-1.8,-1.8);
\filldraw[boxregion=\BColor] (1.2,-1.8) .. controls ++(90:.4cm) and ++(-45:.2cm) .. (.6,-.6) -- (-.6,.6) .. controls ++(135:.2cm) and ++(270:.4cm) .. (-1.2,1.8) -- (1.8,1.8) -- (1.8,-1.8);
\filldraw[primedregion=\QrColor] (-1.8,1.8) -- (-1.8,-1.8) -- (-2.4,-1.8) -- (-2.4,1.8);
\filldraw[boxregion=\QrColor] (1.8,-1.8) -- (1.8,1.8) -- (2.4,1.8) -- (2.4,-1.8);
\end{scope}
\coordinate (a) at (.8,.8);
\coordinate (b) at (-.8,-.8);
\draw[\phiColor,thick] (1.2,-1.8) .. controls ++(90:.4cm) and ++(-45:.2cm) .. (.6,-.6) -- (-.6,.6) .. controls ++(135:.2cm) and ++(270:.4cm) .. (-1.2,1.8);
\draw[\QsColor,thick] (1.8,-1.8) -- (1.8,1.8);
\draw[\QsColor,thick] (-1.8,-1.8) -- (-1.8,1.8);
\draw[\QsColor,thick] (-.4,.4) .. controls ++(45:.4cm) and ++(-135:.4cm) .. (a);
\draw[\QsColor,thick] (.4,-.4) .. controls ++(45:.4cm) and ++(-45:.4cm) .. (a);
\draw[\QsColor,thick] (a) .. controls ++(90:.3cm) and ++(-135:.3cm) .. (1.8,1.5);
\draw[\QsColor,thick] (-.4,.4) .. controls ++(-135:.4cm) and ++(135:.4cm) .. (b);
\draw[\QsColor,thick] (.4,-.4) .. controls ++(-135:.4cm) and ++(45:.4cm) .. (b);
\draw[\QsColor,thick] (b) .. controls ++(270:.3cm) and ++(45:.3cm) .. (-1.8,-1.5);
\filldraw[\QsColor] (a) circle (.05cm);
\filldraw[\QsColor] (b) circle (.05cm);
\filldraw[\QsColor] (-1.8,-1.5) circle (.05cm);
\filldraw[\QsColor] (1.8,1.5) circle (.05cm);
\filldraw[fill=white, thick] (.4,-.4) circle (.1cm);
\filldraw[fill=white, thick] (-.4,.4) circle (.1cm);
\node at (-1.8,-2) {\scriptsize{$F(Q)$}};
\node at (1.2,-2) {\scriptsize{$\varphi_b$}};
\node at (1.8,-2) {\scriptsize{$G(Q)$}};
}
=
\tikzmath[scale=.7, transform shape]{
\begin{scope}
\clip[rounded corners = 5] (-1.8,-1.5) rectangle (1.8,2.1);
\filldraw[primedregion=\BColor] (.6,-1.5) -- (.6,-.6) .. controls ++(90:.6cm) and ++(270:.6cm) .. (-.6,.6) -- (-.6,2.1) -- (-1.2,2.1) -- (-1.2,-1.5);
\filldraw[boxregion=\BColor] (.6,-1.5) -- (.6,-.6) .. controls ++(90:.6cm) and ++(270:.6cm) .. (-.6,.6) -- (-.6,2.1) -- (1.2,2.1) -- (1.2,-1.5);
\filldraw[primedregion=\QrColor] (-1.8,2.1) -- (-1.8,-1.5) -- (-1.2,-1.5) -- (-1.2,2.1);
\filldraw[boxregion=\QrColor] (1.8,2.1) -- (1.8,-1.5) -- (1.2,-1.5) -- (1.2,2.1);
\end{scope}
\draw[\phiColor,thick] (.6,-1.5) -- (.6,-.6) .. controls ++(90:.6cm) and ++(270:.6cm) .. (-.6,.6) -- (-.6,2.1);
\draw[\QsColor,thick] (-1.2,-1.2) arc (-90:0:.6cm) .. controls ++(90:.6cm) and ++(270:.6cm) .. (.6,.6);
\draw[\QsColor,thick] (.6,1.2) arc (90:450:.3cm) arc (180:90:.6cm);
\draw[\QsColor,thick] (-1.2,-1.5) -- (-1.2,2.1);
\draw[\QsColor,thick] (1.2,-1.5) -- (1.2,2.1);
\filldraw[\QsColor] (.6,.6) circle (.05cm);
\filldraw[\QsColor] (.6,1.2) circle (.05cm);
\filldraw[\QsColor] (1.2,1.8) circle (.05cm);
\filldraw[\QsColor] (-1.2,-1.2) circle (.05cm);
\filldraw[white] (0,0) circle (.1cm);
\draw[thick] (0,0) circle (.1cm); 
\node at (-1.2,-1.7) {\scriptsize{$F(Q)$}};
\node at (.6,-1.7) {\scriptsize{$\varphi_b$}};
\node at (1.2,-1.7) {\scriptsize{$G(Q)$}};
}
=
\tikzmath[scale=.7, transform shape]{
\begin{scope}
\clip[rounded corners = 5] (-2.4,-1.2) rectangle (2.4,1.2);
\filldraw[primedregion=\BColor] (.6,-1.2) .. controls ++(90:.4cm) and ++(-45:.2cm) .. (0,0) .. controls ++(135:.2cm) and ++(270:.4cm) .. (-.6,1.2) -- (-2.4,1.2) -- (-2.4,-1.2);
\filldraw[boxregion=\BColor] (.6,-1.2) .. controls ++(90:.4cm) and ++(-45:.2cm) .. (0,0) .. controls ++(135:.2cm) and ++(270:.4cm) .. (-.6,1.2) -- (2.4,1.2) -- (2.4,-1.2);
\filldraw[primedregion=\QrColor] (-.6,-1.2) -- (-.6,-.6) -- (-1.2,0) .. controls ++(135:.2cm) and ++(270:.4cm) .. (-1.8,1.2) -- (-2.4,1.2) -- (-2.4,-1.2);
\filldraw[boxregion=\QrColor] (1.8,-1.2) .. controls ++(90:.4cm) and ++(-45:.2cm) .. (1.2,0) -- (.6,.6) -- (.6,1.2) -- (2.4,1.2) -- (2.4,-1.2);
\end{scope}
\draw[\phiColor,thick] (.6,-1.2) .. controls ++(90:.4cm) and ++(-45:.2cm) .. (0,0) .. controls ++(135:.2cm) and ++(270:.4cm) .. (-.6,1.2);
\draw[\QsColor,thick] (-.6,-1.2) -- (-.6,-.6) -- (-1.2,0) .. controls ++(135:.2cm) and ++(270:.4cm) .. (-1.8,1.2);
\draw[\QsColor,thick] (1.8,-1.2) .. controls ++(90:.4cm) and ++(-45:.2cm) .. (1.2,0) -- (.6,.6) -- (.6,1.2);
\draw[\QsColor,thick] (-.6,-.6) -- (.6,.6);
\filldraw[\QsColor] (.6,.6) circle (.05cm);
\filldraw[\QsColor] (-.6,-.6) circle (.05cm);
\filldraw[white] (0,0) circle (.1cm);
\draw[thick] (0,0) circle (.1cm); 
\node at (-.6,-1.4) {\scriptsize{$F(Q)$}};
\node at (.6,-1.4) {\scriptsize{$\varphi_b$}};
\node at (1.8,-1.4) {\scriptsize{$G(Q)$}};
}
\]

For a 1-cell $({}_P X_Q,\lambda,\rho)$, we define $\QSys(\varphi)_X: F(X)\xxo_{F(Q)} \QSys(\varphi)_Q \Rightarrow \QSys(\varphi)_P\xxo_{G(P)}G(X)$ by
\[
\QSys(\varphi)_X=
\tikzmath[scale=.75, transform shape]{
\begin{scope}
\clip[rounded corners = 5] (-.6,0) rectangle (1.8,2.4);
\filldraw[primedregion=\PrColor] (0,0) -- (0,.6) .. controls ++(90:.4cm) and ++(-135:.2cm) .. (.6,1.2) .. controls ++(135:.2cm) and ++(270:.4cm) .. (0,1.8) -- (0,3) -- (-.6,3) -- (-.6,0); 
\filldraw[primedregion=\QrColor] (1.2,0) -- (1.2,.6) .. controls ++(90:.4cm) and ++(-45:.2cm) .. (.6,1.2) .. controls ++(-135:.2cm) and ++(90:.4cm) .. (0,.6) -- (0,0);
\filldraw[boxregion=\PrColor] (0,3) -- (0,1.8) .. controls ++(270:.4cm) and ++(135:.2cm) .. (.6,1.2) .. controls ++(45:.2cm) and ++(270:.4cm) .. (1.2,1.8) -- (1.2,3);
\filldraw[boxregion=\QrColor] (1.2,0) -- (1.2,.6) .. controls ++(90:.4cm) and ++(-45:.2cm) .. (.6,1.2) .. controls ++(45:.2cm) and ++(270:.4cm) .. (1.2,1.8) -- (1.2,3) -- (1.8,3) -- (1.8,0);
\end{scope}
\draw[\XColor,thick] (0,0) -- (0,.6) .. controls ++(90:.6cm) and ++(270:.6cm) .. (1.2,1.8) -- (1.2,2.4);
\draw[\phiColor,thick] (1.2,0) -- (1.2,.6) .. controls ++(90:.6cm) and ++(270:.6cm) .. (0,1.8) -- (0,2.4);
\filldraw[white] (.6,1.2) circle (.1cm);
\draw[thick] (.6,1.2) circle (.1cm); 
\node at (0,-.2) {\scriptsize{$F(X)$}};
\node at (1.2,2.6) {\scriptsize{$G(X)$}};
\node at (1.2,-.2) {\scriptsize{$\QSys(\varphi)_Q$}};
\node at (0,2.6) {\scriptsize{$\QSys(\varphi)_P$}};
} 
:=
\tikzmath[scale=.7, transform shape]{
\begin{scope}
\clip[rounded corners = 5] (-3.6,-2.4) rectangle (3.6,2.4);
\filldraw[primedregion=\AColor] (0,0) -- (-1.2,1.2) .. controls ++(135:.2cm) and ++(270:.4cm) .. (-1.8,2.4) -- (-3,2.4) -- (-3,-.6) -- (-.6,-.6);
\filldraw[boxregion=\AColor] (0,0) -- (-1.2,1.2) .. controls ++(135:.2cm) and ++(270:.4cm) .. (-1.8,2.4) -- (-.6,2.4) -- (-.6,1.8) -- (.6,.6);
\filldraw[primedregion=\BColor] (1.8,-2.4) .. controls ++(90:.4cm) and ++(-45:.2cm) .. (1.2,-1.2) -- (0,0) -- (-.6,-.6) -- (.6,-1.8) -- (.6,-2.4);
\filldraw[boxregion=\BColor] (1.8,-2.4) .. controls ++(90:.4cm) and ++(-45:.2cm) .. (1.2,-1.2) -- (0,0) -- (.6,.6) -- (3,.6) -- (3,-2.4);
\filldraw[primedregion=\PrColor] (-.6,-2.4) -- (-.6,-.6) -- (-2.4,1.2) .. controls ++(135:.2cm) and ++(270:.4cm) .. (-3,2.4) -- (-3.6,2.4) -- (-3.6,-2.4);
\filldraw[boxregion=\PrColor] (-.6,2.4) -- (-.6,1.8) -- (.6,.6) -- (.6,2.4); 
\filldraw[primedregion=\QrColor] (.6,-2.4) -- (.6,-1.8) -- (-.6,-.6) -- (-.6,-2.4); 
\filldraw[boxregion=\QrColor] (3,-2.4) .. controls ++(90:.4cm) and ++(-45:.2cm) .. (2.4,-1.2) -- (.6,.6) -- (.6,2.4) -- (3.6,2.4) -- (3.6,-2.4);
\filldraw[\AColor] (-1.6,1.2) circle (.22cm);
\filldraw[\AColor] (-.4,0) circle (.22cm);
\filldraw[\BColor] (.8,-1.2) circle (.22cm);
\end{scope}
\draw[\phiColor,thick] (1.8,-2.4) .. controls ++(90:.4cm) and ++(-45:.2cm) .. (1.2,-1.2) -- (-1.2,1.2) .. controls ++(135:.2cm) and ++(270:.4cm) .. (-1.8,2.4);
\draw[\XColor,thick] (-.6,-2.4) -- (-.6,-.6) -- (.6,.6) -- (.6,2.4);
\draw[\PsColor,thick] (-.6,-.6) -- (-2.4,1.2) .. controls ++(135:.2cm) and ++(270:.4cm) .. (-3,2.4);
\draw[\PsColor,thick] (.6,.6) -- (-.6,1.8) -- (-.6,2.4);
\draw[\PsColor,thick] (-1.8,.6) -- (-.6,1.8);
\draw[\QsColor,thick] (-.6,-.6) -- (.6,-1.8) -- (.6,-2.4);
\draw[\QsColor,thick] (3,-2.4) .. controls ++(90:.4cm) and ++(-45:.2cm) .. (2.4,-1.2) -- (.6,.6);
\draw[\QsColor,thick] (1.8,-.6) -- (.6,-1.8);
\filldraw[\XColor] (-.6,-.6) circle (.05cm);
\filldraw[\XColor] (.6,.6) circle (.05cm);
\filldraw[\PsColor] (-1.8,.6) circle (.05cm);
\filldraw[\PsColor] (-.6,1.8) circle (.05cm);
\filldraw[\QsColor] (1.8,-.6) circle (.05cm);
\filldraw[\QsColor] (.6,-1.8) circle (.05cm);
\filldraw[white] (0,0) circle (.1cm);
\draw[thick] (0,0) circle (.1cm); 
\filldraw[white] (-1.2,1.2) circle (.1cm);
\draw[thick] (-1.2,1.2) circle (.1cm); 
\filldraw[white] (1.2,-1.2) circle (.1cm);
\draw[thick] (1.2,-1.2) circle (.1cm); 
\node at (-.6,-2.6) {\scriptsize{$F(X)$}};
\node at (.6,-2.6) {\scriptsize{$F(Q)$}};
\node at (1.8,-2.6) {\scriptsize{$\varphi_b$}};
\node at (3,-2.6) {\scriptsize{$G(Q)$}};
\node at (.6,2.6) {\scriptsize{$G(X)$}};
\node at (-.6,2.6) {\scriptsize{$G(P)$}};
\node at (-1.8,2.6) {\scriptsize{$\varphi_a$}};
\node at (-3,2.6) {\scriptsize{$F(P)$}};
\node at (-.4,0) {\scriptsize{$\varphi_X$}};
\node at (-1.6,1.2) {\scriptsize{$\varphi_P$}};
\node at (.8,-1.2) {\scriptsize{$\varphi_Q$}};
}
\]
To see that $\QSys(\varphi)_X$ is unitary, we observe
\begin{align*}
\QSys(\varphi)_X^\dag\xxt\QSys(\varphi)_X 
&=
\tikzmath[scale=.7, transform shape]{
\begin{scope}
\clip[rounded corners = 5] (-.6,0) rectangle (1.8,4.2);
\filldraw[primedregion=\QrColor] (-.6,0) rectangle (1.8,4.2);
\filldraw[boxregion=\PrColor] (0,3) -- (0,1.8) .. controls ++(270:.4cm) and ++(135:.2cm) .. (.6,1.2) .. controls ++(45:.2cm) and ++(270:.4cm) .. (1.2,1.8) -- (1.2,3);
\filldraw[primedregion=\PrColor] (0,0) -- (0,.6) .. controls ++(90:.4cm) and ++(-135:.2cm) .. (.6,1.2) .. controls ++(135:.2cm) and ++(270:.4cm) .. (0,1.8) -- (0,2.4) .. controls ++(90:.4cm) and ++(-135:.2cm) .. (.6,3) .. controls ++(135:.2cm) and ++(270:.4cm) .. (0,3.6) -- (0,4.2) -- (-.6,4.2) -- (-.6,0); 
\filldraw[boxregion=\QrColor] (1.2,0) -- (1.2,.6) .. controls ++(90:.4cm) and ++(-45:.2cm) .. (.6,1.2) .. controls ++(45:.2cm) and ++(270:.4cm) .. (1.2,1.8) -- (1.2,2.4) .. controls ++(90:.4cm) and ++(-45:.2cm) .. (.6,3) .. controls ++(45:.2cm) and ++(270:.4cm) .. (1.2,3.6) -- (1.2,4.2) -- (1.8,4.2) -- (1.8,0);
\end{scope}
\draw[\XColor,thick] (0,0) -- (0,.6) .. controls ++(90:.6cm) and ++(270:.6cm) .. (1.2,1.8) -- (1.2,2.4) .. controls ++(90:.6cm) and ++(270:.6cm) .. (0,3.6) -- (0,4.2);
\draw[\phiColor,thick] (1.2,0) -- (1.2,.6) .. controls ++(90:.6cm) and ++(270:.6cm) .. (0,1.8) -- (0,2.4) .. controls ++(90:.6cm) and ++(270:.6cm) .. (1.2,3.6) -- (1.2,4.2);
\filldraw[white] (.6,1.2) circle (.1cm);
\draw[thick] (.6,1.2) circle (.1cm); 
\filldraw[white] (.6,3) circle (.1cm);
\draw[thick] (.6,3) circle (.1cm); 
\node at (0,-.2) {\scriptsize{$F(X)$}};
\node at (0,4.4) {\scriptsize{$F(X)$}};
\node at (1.2,-.2) {\scriptsize{$\QSys(\varphi)_Q$}};
\node at (1.2,4.4) {\scriptsize{$\QSys(\varphi)_Q$}};
} 
=
\tikzmath[scale=.45, transform shape]{
\begin{scope}
\clip[rounded corners = 5] (-3.6,-4.8) rectangle (3.6,4.8);
\filldraw[primedregion=\AColor] (-3.6,-4.8) rectangle (3.6,4.8);
\filldraw[boxregion=\AColor] (.6,1.8) -- (.6,-1.8) -- (0,-2.4) -- (-1.2,-1.2) .. controls ++(135:.2cm) and ++(270:.4cm) .. (-1.8,-.2) -- (-1.8,.2) .. controls ++(90:.4cm) and ++(-135:.2cm) .. (-1.2,1.2) -- (0,2.4);
\filldraw[primedregion=\BColor] (1.8,-4.8) .. controls ++(90:.4cm) and ++(-45:.2cm) .. (1.2,-3.6) -- (0,-2.4) -- (-.6,-3) -- (-.6,-4.8);
\filldraw[primedregion=\BColor] (-.6,4.8) -- (-.6,3) -- (0,2.4) -- (1.2,3.6) .. controls ++(45:.2cm) and ++(270:.4cm) .. (1.8,4.8);
\filldraw[boxregion=\BColor] (1.8,-4.8) .. controls ++(90:.4cm) and ++(-45:.2cm) .. (1.2,-3.6) -- (0,-2.4) -- (.6,-1.8) -- (.6,1.8) -- (0,2.4) -- (1.2,3.6) .. controls ++(45:.2cm) and ++(270:.4cm) .. (1.8,4.8) -- (3.6,4.8) -- (3.6,-4.8);
\filldraw[boxregion=\PrColor] (.6,-1.8) -- (-.6,-.6) -- (-.6,.6) -- (.6,1.8);
\filldraw[primedregion=\QrColor] (-.6,3) -- (.6,4.2) -- (.6,4.8) -- (-.6,4.8);
\filldraw[primedregion=\QrColor] (-.6,-3) -- (.6,-4.2) -- (.6,-4.8) -- (-.6,-4.8);
\filldraw[primedregion=\PrColor] (-.6,-4.8) -- (-.6,-3) -- (-2.4,-1.2) .. controls ++(135:.2cm) and ++(270:.4cm) .. (-3,-.2) -- (-3,.2) .. controls ++(90:.4cm) and ++(-135:.2cm) .. (-2.4,1.2) -- (-.6,3) -- (-.6,4.8) -- (-3.6,4.8) -- (-3.6,-4.8);
\filldraw[boxregion=\QrColor] (3,-4.8) .. controls ++(90:.4cm) and ++(-45:.2cm) .. (2.4,-3.6) -- (.6,-1.8) -- (.6,1.8) -- (2.4,3.6) .. controls ++(45:.2cm) and ++(270:.4cm) .. (3,4.8) -- (3.6,4.8) -- (3.6,-4.8);
\end{scope}
\draw[\phiColor,thick] (1.8,-4.8) .. controls ++(90:.4cm) and ++(-45:.2cm) .. (1.2,-3.6) -- (-1.2,-1.2) .. controls ++(135:.2cm) and ++(270:.4cm) .. (-1.8,-.2) -- (-1.8,.2) .. controls ++(90:.4cm) and ++(-135:.2cm) .. (-1.2,1.2) -- (1.2,3.6) .. controls ++(45:.2cm) and ++(270:.4cm) .. (1.8,4.8);
\draw[\XColor,thick] (-.6,-4.8) -- (-.6,-3) -- (.6,-1.8) -- (.6,1.8) -- (-.6,3) -- (-.6,4.8);
\draw[\PsColor,thick] (-.6,-3) -- (-2.4,-1.2) .. controls ++(135:.2cm) and ++(270:.4cm) .. (-3,-.2) -- (-3,.2) .. controls ++(90:.4cm) and ++(-135:.2cm) .. (-2.4,1.2) -- (-.6,3);
\draw[\PsColor,thick] (.6,-1.8) -- (-.6,-.6) -- (-.6,.6) -- (.6,1.8);
\draw[\PsColor,thick] (-1.8,-1.8) -- (-.6,-.6);
\draw[\PsColor,thick] (-1.8,1.8) -- (-.6,.6);
\draw[\QsColor,thick] (-.6,-3) -- (.6,-4.2) -- (.6,-4.8);
\draw[\QsColor,thick] (3,-4.8) .. controls ++(90:.4cm) and ++(-45:.2cm) .. (2.4,-3.6) -- (.6,-1.8);
\draw[\QsColor,thick] (1.8,-3) -- (.6,-4.2);
\draw[\QsColor,thick] (-.6,3) -- (.6,4.2) -- (.6,4.8);
\draw[\QsColor,thick] (.6,1.8) -- (2.4,3.6) .. controls ++(45:.2cm) and ++(270:.4cm) .. (3,4.8);
\draw[\QsColor,thick] (1.8,3) -- (.6,4.2);
\filldraw[\XColor] (-.6,-3) circle (.07cm);
\filldraw[\XColor] (.6,-1.8) circle (.07cm);
\filldraw[\XColor] (-.6,3) circle (.07cm);
\filldraw[\XColor] (.6,1.8) circle (.07cm);
\filldraw[\PsColor] (-1.8,-1.8) circle (.07cm);
\filldraw[\PsColor] (-.6,-.6) circle (.07cm);
\filldraw[\PsColor] (-1.8,1.8) circle (.07cm);
\filldraw[\PsColor] (-.6,.6) circle (.07cm);
\filldraw[\QsColor] (1.8,3) circle (.07cm);
\filldraw[\QsColor] (.6,4.2) circle (.07cm);
\filldraw[\QsColor] (1.8,-3) circle (.07cm);
\filldraw[\QsColor] (.6,-4.2) circle (.07cm);
\filldraw[white] (0,-2.4) circle (.1cm);
\draw[thick] (0,-2.4) circle (.1cm); 
\filldraw[white] (0,2.4) circle (.1cm);
\draw[thick] (0,2.4) circle (.1cm); 
\filldraw[white] (-1.2,-1.2) circle (.1cm);
\draw[thick] (-1.2,-1.2) circle (.1cm); 
\filldraw[white] (-1.2,1.2) circle (.1cm);
\draw[thick] (-1.2,1.2) circle (.1cm); 
\filldraw[white] (1.2,-3.6) circle (.1cm);
\draw[thick] (1.2,-3.6) circle (.1cm); 
\filldraw[white] (1.2,3.6) circle (.1cm);
\draw[thick] (1.2,3.6) circle (.1cm); 
}
=
\tikzmath[scale=.45, transform shape]{
\begin{scope}
\clip[rounded corners = 5] (-3.6,-4.8) rectangle (3.6,4.8);
\filldraw[primedregion=\AColor] (-3.6,-4.8) rectangle (3.6,4.8);
\filldraw[boxregion=\AColor] (.6,1.8) -- (.6,-1.8) -- (0,-2.4) -- (-1.2,-1.2) .. controls ++(135:.2cm) and ++(270:.4cm) .. (-1.8,-.2) -- (-1.8,.2) .. controls ++(90:.4cm) and ++(-135:.2cm) .. (-1.2,1.2) -- (0,2.4);
\filldraw[primedregion=\BColor] (1.8,-4.8) .. controls ++(90:.4cm) and ++(-45:.2cm) .. (1.2,-3.6) -- (0,-2.4) -- (-.6,-3) -- (-.6,-4.8);
\filldraw[primedregion=\BColor] (-.6,4.8) -- (-.6,3) -- (0,2.4) -- (1.2,3.6) .. controls ++(45:.2cm) and ++(270:.4cm) .. (1.8,4.8);
\filldraw[boxregion=\BColor] (1.8,-4.8) .. controls ++(90:.4cm) and ++(-45:.2cm) .. (1.2,-3.6) -- (0,-2.4) -- (.6,-1.8) -- (.6,1.8) -- (0,2.4) -- (1.2,3.6) .. controls ++(45:.2cm) and ++(270:.4cm) .. (1.8,4.8) -- (3.6,4.8) -- (3.6,-4.8);
\filldraw[primedregion=\QrColor] (-.6,3) -- (.6,4.2) -- (.6,4.8) -- (-.6,4.8);
\filldraw[primedregion=\QrColor] (-.6,-3) -- (.6,-4.2) -- (.6,-4.8) -- (-.6,-4.8);
\filldraw[primedregion=\PrColor] (-.6,-4.8) -- (-.6,-3) -- (-2.4,-1.2) .. controls ++(135:.2cm) and ++(270:.4cm) .. (-3,-.2) -- (-3,.2) .. controls ++(90:.4cm) and ++(-135:.2cm) .. (-2.4,1.2) -- (-.6,3) -- (-.6,4.8) -- (-3.6,4.8) -- (-3.6,-4.8);
\filldraw[boxregion=\QrColor] (3,-4.8) .. controls ++(90:.4cm) and ++(-45:.2cm) .. (2.4,-3.6) -- (1.8,-3) -- (1.8,3) -- (2.4,3.6) .. controls ++(45:.2cm) and ++(270:.4cm) .. (3,4.8) -- (3.6,4.8) -- (3.6,-4.8);
\end{scope}
\draw[\phiColor,thick] (1.8,-4.8) .. controls ++(90:.4cm) and ++(-45:.2cm) .. (1.2,-3.6) -- (-1.2,-1.2) .. controls ++(135:.2cm) and ++(270:.4cm) .. (-1.8,-.2) -- (-1.8,.2) .. controls ++(90:.4cm) and ++(-135:.2cm) .. (-1.2,1.2) -- (1.2,3.6) .. controls ++(45:.2cm) and ++(270:.4cm) .. (1.8,4.8);
\draw[\XColor,thick] (-.6,-4.8) -- (-.6,-3) -- (.6,-1.8) -- (.6,1.8) -- (-.6,3) -- (-.6,4.8);
\draw[\PsColor,thick] (-.6,-3) -- (-2.4,-1.2) .. controls ++(135:.2cm) and ++(270:.4cm) .. (-3,-.2) -- (-3,.2) .. controls ++(90:.4cm) and ++(-135:.2cm) .. (-2.4,1.2) -- (-.6,3);
\draw[\QsColor,thick] (-.6,-3) -- (.6,-4.2) -- (.6,-4.8);
\draw[\QsColor,thick] (3,-4.8) .. controls ++(90:.4cm) and ++(-45:.2cm) .. (2.4,-3.6) -- (1.8,-3) -- (1.8,3);
\draw[\QsColor,thick] (1.8,-3) -- (.6,-4.2);
\draw[\QsColor,thick] (-.6,3) -- (.6,4.2) -- (.6,4.8);
\draw[\QsColor,thick] (1.8,3) -- (2.4,3.6) .. controls ++(45:.2cm) and ++(270:.4cm) .. (3,4.8);
\draw[\QsColor,thick] (1.8,3) -- (.6,4.2);
\filldraw[\XColor] (-.6,-3) circle (.07cm);
\filldraw[\XColor] (-.6,3) circle (.07cm);
\filldraw[\QsColor] (1.8,3) circle (.07cm);
\filldraw[\QsColor] (.6,4.2) circle (.07cm);
\filldraw[\QsColor] (1.8,-3) circle (.07cm);
\filldraw[\QsColor] (.6,-4.2) circle (.07cm);
\filldraw[white] (0,-2.4) circle (.1cm);
\draw[thick] (0,-2.4) circle (.1cm); 
\filldraw[white] (0,2.4) circle (.1cm);
\draw[thick] (0,2.4) circle (.1cm); 
\filldraw[white] (1.2,-3.6) circle (.1cm);
\draw[thick] (1.2,-3.6) circle (.1cm); 
\filldraw[white] (1.2,3.6) circle (.1cm);
\draw[thick] (1.2,3.6) circle (.1cm); 
}
\displaybreak[1]\\
&=
\tikzmath[scale=.45, transform shape]{
\begin{scope}
\clip[rounded corners = 5] (-3,-4.8) rectangle (3.6,4.8);
\filldraw[primedregion=\BColor] (1.8,-4.8) .. controls ++(90:.4cm) and ++(-45:.2cm) .. (1.2,-3.6) .. controls ++(135:.2cm) and ++(270:.4cm) .. (.6,-2.4) -- (.6,2.4) .. controls ++(90:.4cm) and ++(-135:.2cm) .. (1.2,3.6) .. controls ++(45:.2cm) and ++(270:.4cm) .. (1.8,4.8) -- (-.6,4.8) -- (-.6,-4.8);
\filldraw[boxregion=\BColor] (1.8,-4.8) .. controls ++(90:.4cm) and ++(-45:.2cm) .. (1.2,-3.6) .. controls ++(135:.2cm) and ++(270:.4cm) .. (.6,-2.4) -- (.6,2.4) .. controls ++(90:.4cm) and ++(-135:.2cm) .. (1.2,3.6) .. controls ++(45:.2cm) and ++(270:.4cm) .. (1.8,4.8) -- (3.6,4.8) -- (3.6,-4.8);
%
\filldraw[primedregion=\QrColor] (.6,-4.8) -- (.6,-4.2) -- (0,-3.6)  .. controls ++(135:.2cm) and ++(270:.4cm) .. (-.6,-2.4) -- (-.6,2.4) .. controls ++(90:.4cm) and ++(-135:.2cm) .. (0,3.6) -- (.6,4.2) -- (.6,4.8) -- (-1.8,4.8) -- (-1.8,-4.8);
\filldraw[primedregion=\PrColor] (-3,-4.8) rectangle (-1.8,4.8);
\filldraw[boxregion=\QrColor] (3,-4.8) .. controls ++(90:.4cm) and ++(-45:.2cm) .. (2.4,-3.6) -- (1.8,-3) -- (1.8,3) -- (2.4,3.6) .. controls ++(45:.2cm) and ++(270:.4cm) .. (3,4.8) -- (3.6,4.8) -- (3.6,-4.8);
\end{scope}
\draw[\phiColor,thick] (1.8,-4.8) .. controls ++(90:.4cm) and ++(-45:.2cm) .. (1.2,-3.6) .. controls ++(135:.2cm) and ++(270:.4cm) .. (.6,-2.4) -- (.6,2.4) .. controls ++(90:.4cm) and ++(-135:.2cm) .. (1.2,3.6) .. controls ++(45:.2cm) and ++(270:.4cm) .. (1.8,4.8);
\draw[\XColor,thick] (-1.8,-4.8) -- (-1.8,4.8);
\draw[\QsColor,thick] (.6,-4.8) -- (.6,-4.2) -- (0,-3.6)  .. controls ++(135:.2cm) and ++(270:.4cm) .. (-.6,-2.4) -- (-.6,2.4) .. controls ++(90:.4cm) and ++(-135:.2cm) .. (0,3.6) -- (.6,4.2) -- (.6,4.8);
\draw[\QsColor,thick] (3,-4.8) .. controls ++(90:.4cm) and ++(-45:.2cm) .. (2.4,-3.6) -- (1.8,-3) -- (1.8,3);
\draw[\QsColor,thick] (1.8,-3) -- (.6,-4.2);
\draw[\QsColor,thick] (1.8,3) -- (2.4,3.6) .. controls ++(45:.2cm) and ++(270:.4cm) .. (3,4.8);
\draw[\QsColor,thick] (1.8,3) -- (.6,4.2);
\filldraw[\QsColor] (1.8,3) circle (.07cm);
\filldraw[\QsColor] (.6,4.2) circle (.07cm);
\filldraw[\QsColor] (1.8,-3) circle (.07cm);
\filldraw[\QsColor] (.6,-4.2) circle (.07cm);
\filldraw[white] (1.2,-3.6) circle (.1cm);
\draw[thick] (1.2,-3.6) circle (.1cm); 
\filldraw[white] (1.2,3.6) circle (.1cm);
\draw[thick] (1.2,3.6) circle (.1cm); 
}
=
\tikzmath[scale=.5, transform shape]{
\begin{scope}
\clip[rounded corners = 5] (-3,-1.2) rectangle (2.4,1.2);
\filldraw[primedregion=\BColor] (.6,-1.2) .. controls ++(90:.4cm) and ++(-45:.2cm) .. (0,0) .. controls ++(135:.2cm) and ++(270:.4cm) .. (-.6,1.2) -- (-2.4,1.2) -- (-2.4,-1.2);
\filldraw[boxregion=\BColor] (.6,-1.2) .. controls ++(90:.4cm) and ++(-45:.2cm) .. (0,0) .. controls ++(135:.2cm) and ++(270:.4cm) .. (-.6,1.2) -- (2.4,1.2) -- (2.4,-1.2);
\filldraw[primedregion=\PrColor] (-3,-1.2) rectangle (-2.4,1.2);
\filldraw[primedregion=\QrColor] (-.6,-1.2) -- (-.6,-.6) -- (-1.2,0) .. controls ++(135:.2cm) and ++(270:.4cm) .. (-1.8,1.2) -- (-2.4,1.2) -- (-2.4,-1.2);
\filldraw[boxregion=\QrColor] (1.8,-1.2) .. controls ++(90:.4cm) and ++(-45:.2cm) .. (1.2,0) -- (.6,.6) -- (.6,1.2) -- (2.4,1.2) -- (2.4,-1.2);
\end{scope}
\draw[\XColor,thick] (-2.4,-1.2) -- (-2.4,1.2);
\draw[\phiColor,thick] (.6,-1.2) .. controls ++(90:.4cm) and ++(-45:.2cm) .. (0,0) .. controls ++(135:.2cm) and ++(270:.4cm) .. (-.6,1.2);
\draw[\QsColor,thick] (-.6,-1.2) -- (-.6,-.6) -- (-1.2,0) .. controls ++(135:.2cm) and ++(270:.4cm) .. (-1.8,1.2);
\draw[\QsColor,thick] (1.8,-1.2) .. controls ++(90:.4cm) and ++(-45:.2cm) .. (1.2,0) -- (.6,.6) -- (.6,1.2);
\draw[\QsColor,thick] (-.6,-.6) -- (.6,.6);
\filldraw[\QsColor] (.6,.6) circle (.07cm);
\filldraw[\QsColor] (-.6,-.6) circle (.07cm);
\filldraw[white] (0,0) circle (.1cm);
\draw[thick] (0,0) circle (.1cm); 
\node at (-.6,-1.4) {\small{$F(Q)$}};
\node at (.6,-1.4) {\small{$\varphi_b$}};
\node at (1.8,-1.4) {\small{$G(Q)$}};
}
=
\tikzmath{
\begin{scope}
\clip[rounded corners = 5] (-.5,-.5) rectangle (.5,.5);
\filldraw[primedregion=\PrColor] (-.5,-.5) rectangle (-.2,.5);
\filldraw[primedregion=\QrColor] (-.2,-.5) rectangle (.2,.5);
\filldraw[boxregion=\QrColor] (.2,-.5) rectangle (.5,.5);
\end{scope}
\draw[\XColor,thick] (-.2,-.5) -- (-.2,.5);
\draw[\phiColor,thick] (.2,-.5) -- (.2,.5);
}\,.
\displaybreak[1]
\end{align*}
Similarly, $\QSys(\varphi)_X\xxt \QSys(\varphi)_X^\dag=1_{\QSys(\varphi)_P\xxo_{G(P)}G(X)}$.

To see that $\QSys(\varphi):\QSys(F)\Rightarrow\QSys(G)$ is a 2-transformation, we observe
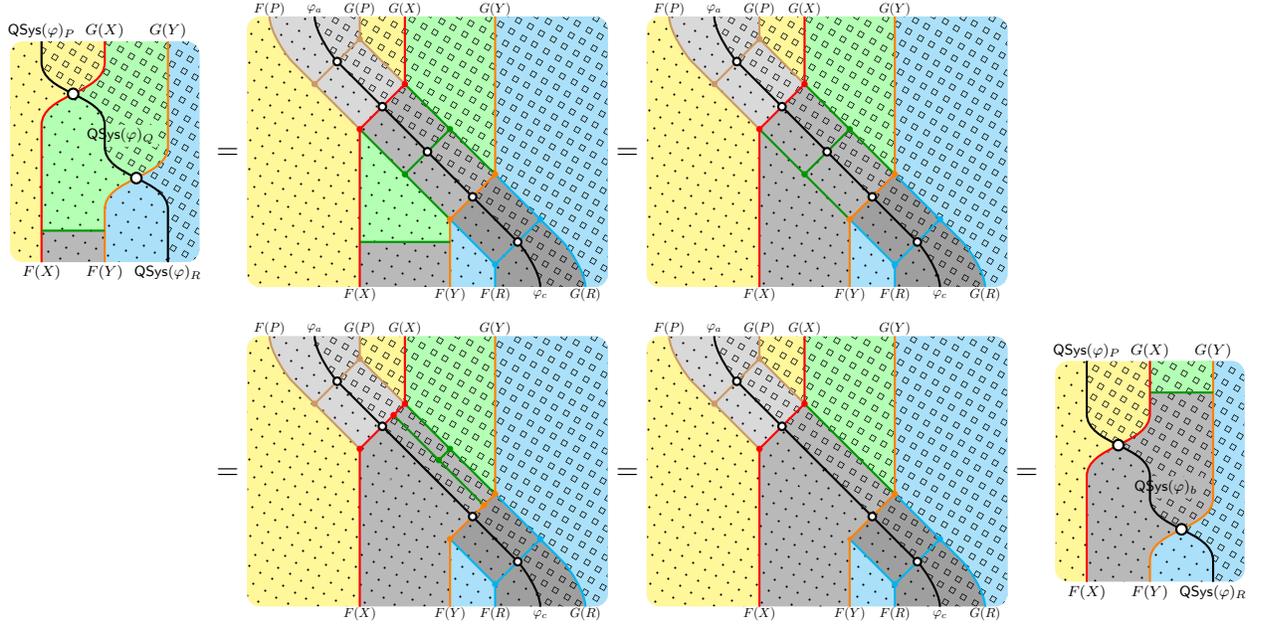
\begin{align*}
\tikzmath[scale=.7, transform shape]{
\begin{scope}
\clip[rounded corners = 5] (-.6,-.6) rectangle (3,3.6);
\filldraw[primedregion=\PrColor] (0,-.6) -- (0,2) .. controls ++(90:.4cm) and ++(-135:.2cm) .. (.6,2.6) .. controls ++(135:.2cm) and ++(270:.4cm) .. (0,3.2) -- (0,4.2) -- (-.6,4.2) -- (-.6,-.6); 
\filldraw[primedregion=\QrColor] (1.2,0) -- (1.2,.4) .. controls ++(90:.4cm) and ++(-135:.2cm) .. (1.8,1) .. controls ++(135:.2cm) and ++(270:.4cm) .. (1.2,1.6) -- (1.2,2) .. controls ++(90:.4cm) and ++(-45:.2cm) .. (.6,2.6) .. controls ++(-135:.2cm) and ++(90:.4cm) .. (0,2) -- (0,0);
\filldraw[primedregion=\RrColor] (2.4,-.6) -- (2.4,.4) .. controls ++(90:.4cm) and ++(-45:.2cm) .. (1.8,1) .. controls ++(-135:.2cm) and ++(90:.4cm) .. (1.2,.4) -- (1.2,-.6);
\filldraw[boxregion=\PrColor] (0,4.2) -- (0,3.2) .. controls ++(270:.4cm) and ++(135:.2cm) .. (.6,2.6) .. controls ++(45:.2cm) and ++(270:.4cm) .. (1.2,3.2) -- (1.2,4.2);
\filldraw[boxregion=\QrColor] (2.4,3.6) -- (2.4,1.6) .. controls ++(270:.4cm) and ++(45:.2cm) .. (1.8,1) .. controls ++(135:.2cm) and ++(270:.4cm) .. (1.2,1.6) -- (1.2,2) .. controls ++(90:.4cm) and ++(-45:.2cm) .. (.6,2.6) .. controls ++(45:.2cm) and ++(270:.4cm) .. (1.2,3.2) -- (1.2,3.6);
\filldraw[primedregion=\BColor] (0,-.6) rectangle (1.2,0);
\filldraw[boxregion=\RrColor] (2.4,-.6) -- (2.4,.4) .. controls ++(90:.4cm) and ++(-45:.2cm) .. (1.8,1) .. controls ++(45:.2cm) and ++(270:.4cm) .. (2.4,1.6) -- (2.4,4.2) -- (3,4.2) -- (3,-.6);
\filldraw[\QrColor] (.8,1.6) rectangle (2.2,2);
\end{scope}
\draw[\XColor,thick] (0,-.6) -- (0,.4) -- (0,2) .. controls ++(90:.6cm) and ++(270:.6cm) .. (1.2,3.2) -- (1.2,3.6);
\draw[\YColor,thick] (1.2,-.6) -- (1.2,.4) .. controls ++(90:.6cm) and ++(270:.6cm) .. (2.4,1.6) -- (2.4,3.6);
\draw[\phiColor,thick] (2.4,-.6) -- (2.4,.4) .. controls ++(90:.6cm) and ++(270:.6cm) .. (1.2,1.6) -- (1.2,2) .. controls ++(90:.6cm) and ++(270:.6cm) .. (0,3.2) -- (0,3.6);
\draw[\QsColor,thick] (0,0) -- (1.2,0);
\filldraw[white] (1.8,1) circle (.1cm);
\draw[thick] (1.8,1) circle (.1cm); 
\filldraw[white] (.6,2.6) circle (.1cm);
\draw[thick] (.6,2.6) circle (.1cm); 
\node at (0,-.8) {\scriptsize{$F(X)$}};
\node at (1.2,-.8) {\scriptsize{$F(Y)$}};
\node at (1.2,3.8) {\scriptsize{$G(X)$}};
\node at (2.4,3.8) {\scriptsize{$G(Y)$}};
\node at (2.4,-.8) {\scriptsize{$\QSys(\varphi)_R$}};
\node at (1.5,1.8) {\scriptsize{$\QSys(\varphi)_Q$}};
\node at (0,3.8) {\scriptsize{$\QSys(\varphi)_P$}};
}
&=
\tikzmath[scale=.5, transform shape]{
\begin{scope}
\clip[rounded corners = 5] (-4.8,-3.6) rectangle (4.8,3.6);
\filldraw[primedregion=\AColor] (-1.2,1.2) -- (-2.4,2.4) .. controls ++(135:.2cm) and ++(270:.4cm) .. (-3,3.6) -- (-4.2,3.6) -- (-4.2,-.6) -- (-1.8,.6);
\filldraw[boxregion=\AColor] (-1.2,1.2) -- (-2.4,2.4) .. controls ++(135:.2cm) and ++(270:.4cm) .. (-3,3.6) -- (-1.8,3.6) -- (-1.8,3) -- (-.6,1.8);
\filldraw[primedregion=\BColor] (.6,-1.8) -- (-1.8,.6) -- (-1.2,1.2) -- (1.2,-1.2);
\filldraw[boxregion=\BColor] (1.8,-.6) -- (-.6,1.8)-- (-1.2,1.2) -- (1.2,-1.2);
\filldraw[primedregion=\CColor] (3,-3.6) .. controls ++(90:.4cm) and ++(-45:.2cm) .. (2.4,-2.4) -- (1.2,-1.2) -- (.6,-1.8) -- (1.8,-3) -- (1.8,-3.6);
\filldraw[boxregion=\CColor] (3,-3.6) .. controls ++(90:.4cm) and ++(-45:.2cm) .. (2.4,-2.4) -- (1.2,-1.2) -- (1.8,-.6) -- (4.2,-.6) -- (4.2,-3.6);
\filldraw[primedregion=\PrColor] (-1.8,-3.6) -- (-1.8,.6) -- (-3.6,2.4) .. controls ++(135:.2cm) and ++(270:.4cm) .. (-4.2,3.6) -- (-4.8,3.6) -- (-4.8,-3.6);
\filldraw[boxregion=\PrColor] (-1.8,3.6) -- (-1.8,3) -- (-.6,1.8) -- (-.6,3.6);
\filldraw[primedregion=\QrColor] (-1.8,-3.6) -- (-1.8,.6) -- (.6,-1.8) -- (.6,-3.6);
\filldraw[boxregion=\QrColor] (-.6,1.8) -- (-.6,3.6) -- (1.8,3.6) -- (1.8,-.6);
\filldraw[primedregion=\RrColor] (1.8,-3.6) -- (1.8,-3) -- (.6,-1.8) -- (.6,-3.6);
\filldraw[boxregion=\RrColor] (4.2,-3.6) .. controls ++(90:.4cm) and ++(-45:.2cm) .. (3.6,-2.4) -- (1.8,-.6) -- (1.8,3.6) -- (4.8,3.6) -- (4.8,-3.6);
\filldraw[primedregion=\BColor] (.6,-3.6) -- (.6,-2.4) -- (-1.8,-2.4) -- (-1.8,-3.6);
\end{scope}
\draw[\phiColor,thick] (3,-3.6) .. controls ++(90:.4cm) and ++(-45:.2cm) .. (2.4,-2.4) -- (-2.4,2.4) .. controls ++(135:.2cm) and ++(270:.4cm) .. (-3,3.6);
\draw[\XColor,thick] (-1.8,-3.6) -- (-1.8,.6) -- (-.6,1.8) -- (-.6,3.6);
\draw[\YColor,thick] (1.8,3.6) -- (1.8,-.6) -- (.6,-1.8) -- (.6,-3.6);
\draw[\QsColor,thick] (.6,-1.8) -- (-1.8,.6);
\draw[\QsColor,thick] (1.8,-.6) -- (-.6,1.8);
\draw[\QsColor,thick] (-.6,-.6) -- (.6,.6);
\draw[\PsColor,thick] (-1.8,3.6) -- (-1.8,3) -- (-.6,1.8);
\draw[\PsColor,thick] (-1.8,.6) -- (-3.6,2.4) .. controls ++(135:.2cm) and ++(270:.4cm) .. (-4.2,3.6);
\draw[\PsColor,thick] (-3,1.8) -- (-1.8,3);
\draw[\RsColor,thick] (.6,-1.8) -- (1.8,-3) -- (1.8,-3.6);
\draw[\RsColor,thick] (4.2,-3.6) .. controls ++(90:.4cm) and ++(-45:.2cm) .. (3.6,-2.4) -- (1.8,-.6);
\draw[\RsColor,thick] (1.8,-3) -- (3,-1.8);
\draw[\QsColor,thick] (.6,-2.4) -- (-1.8,-2.4);
\filldraw[white] (-1.2,1.2) circle (.1cm);
\draw[thick] (-1.2,1.2) circle (.1cm); 
\filldraw[white] (1.2,-1.2) circle (.1cm);
\draw[thick] (1.2,-1.2) circle (.1cm); 
\filldraw[white] (-2.4,2.4) circle (.1cm);
\draw[thick] (-2.4,2.4) circle (.1cm); 
\filldraw[white] (0,0) circle (.1cm);
\draw[thick] (0,0) circle (.1cm); 
\filldraw[white] (2.4,-2.4) circle (.1cm);
\draw[thick] (2.4,-2.4) circle (.1cm); 
\filldraw[\XColor] (-1.8,.6) circle (.07cm);
\filldraw[\XColor] (-.6,1.8) circle (.07cm);
\filldraw[\YColor] (1.8,-.6) circle (.07cm);
\filldraw[\YColor] (.6,-1.8) circle (.07cm);
\filldraw[\QsColor] (-.6,-.6) circle (.07cm);
\filldraw[\QsColor] (.6,.6) circle (.07cm);
\filldraw[\PsColor] (-3,1.8) circle (.07cm);
\filldraw[\PsColor] (-1.8,3) circle (.07cm);
\filldraw[\RsColor] (3,-1.8) circle (.07cm);
\filldraw[\RsColor] (1.8,-3) circle (.07cm);
\node at (-1.8,-3.8) {\small{$F(X)$}};
\node at (.6,-3.8) {\small{$F(Y)$}};
\node at (1.8,-3.8) {\small{$F(R)$}};
\node at (3,-3.8) {\small{$\varphi_c$}};
\node at (4.2,-3.8) {\small{$G(R)$}};
\node at (1.8,3.8) {\small{$G(Y)$}};
\node at (-.6,3.8) {\small{$G(X)$}};
\node at (-1.8,3.8) {\small{$G(P)$}};
\node at (-3,3.8) {\small{$\varphi_a$}};
\node at (-4.2,3.8) {\small{$F(P)$}};
}
=
\tikzmath[scale=.5, transform shape]{
\begin{scope}
\clip[rounded corners = 5] (-4.8,-3.6) rectangle (4.8,3.6);
\filldraw[primedregion=\AColor] (-1.2,1.2) -- (-2.4,2.4) .. controls ++(135:.2cm) and ++(270:.4cm) .. (-3,3.6) -- (-4.2,3.6) -- (-4.2,-.6) -- (-1.8,.6);
\filldraw[boxregion=\AColor] (-1.2,1.2) -- (-2.4,2.4) .. controls ++(135:.2cm) and ++(270:.4cm) .. (-3,3.6) -- (-1.8,3.6) -- (-1.8,3) -- (-.6,1.8);
\filldraw[primedregion=\BColor] (.6,-1.8) -- (-1.8,.6) -- (-1.2,1.2) -- (1.2,-1.2);
\filldraw[boxregion=\BColor] (1.8,-.6) -- (-.6,1.8)-- (-1.2,1.2) -- (1.2,-1.2);
\filldraw[primedregion=\CColor] (3,-3.6) .. controls ++(90:.4cm) and ++(-45:.2cm) .. (2.4,-2.4) -- (1.2,-1.2) -- (.6,-1.8) -- (1.8,-3) -- (1.8,-3.6);
\filldraw[boxregion=\CColor] (3,-3.6) .. controls ++(90:.4cm) and ++(-45:.2cm) .. (2.4,-2.4) -- (1.2,-1.2) -- (1.8,-.6) -- (4.2,-.6) -- (4.2,-3.6);
\filldraw[primedregion=\PrColor] (-1.8,-3.6) -- (-1.8,.6) -- (-3.6,2.4) .. controls ++(135:.2cm) and ++(270:.4cm) .. (-4.2,3.6) -- (-4.8,3.6) -- (-4.8,-3.6);
\filldraw[boxregion=\PrColor] (-1.8,3.6) -- (-1.8,3) -- (-.6,1.8) -- (-.6,3.6);
\filldraw[primedregion=\BColor] (-1.8,-3.6) -- (-1.8,.6) -- (.6,-1.8) -- (.6,-3.6);
\filldraw[boxregion=\QrColor] (-.6,1.8) -- (-.6,3.6) -- (1.8,3.6) -- (1.8,-.6);
\filldraw[primedregion=\RrColor] (1.8,-3.6) -- (1.8,-3) -- (.6,-1.8) -- (.6,-3.6);
\filldraw[boxregion=\RrColor] (4.2,-3.6) .. controls ++(90:.4cm) and ++(-45:.2cm) .. (3.6,-2.4) -- (1.8,-.6) -- (1.8,3.6) -- (4.8,3.6) -- (4.8,-3.6);
\end{scope}
\draw[\phiColor,thick] (3,-3.6) .. controls ++(90:.4cm) and ++(-45:.2cm) .. (2.4,-2.4) -- (-2.4,2.4) .. controls ++(135:.2cm) and ++(270:.4cm) .. (-3,3.6);
\draw[\XColor,thick] (-1.8,-3.6) -- (-1.8,.6) -- (-.6,1.8) -- (-.6,3.6);
\draw[\YColor,thick] (1.8,3.6) -- (1.8,-.6) -- (.6,-1.8) -- (.6,-3.6);
\draw[\QsColor,thick] (.6,-1.8) -- (-1.8,.6);
\draw[\QsColor,thick] (1.8,-.6) -- (-.6,1.8);
\draw[\QsColor,thick] (-.6,-.6) -- (.6,.6);
\draw[\PsColor,thick] (-1.8,3.6) -- (-1.8,3) -- (-.6,1.8);
\draw[\PsColor,thick] (-1.8,.6) -- (-3.6,2.4) .. controls ++(135:.2cm) and ++(270:.4cm) .. (-4.2,3.6);
\draw[\PsColor,thick] (-3,1.8) -- (-1.8,3);
\draw[\RsColor,thick] (.6,-1.8) -- (1.8,-3) -- (1.8,-3.6);
\draw[\RsColor,thick] (4.2,-3.6) .. controls ++(90:.4cm) and ++(-45:.2cm) .. (3.6,-2.4) -- (1.8,-.6);
\draw[\RsColor,thick] (1.8,-3) -- (3,-1.8);
\filldraw[white] (-1.2,1.2) circle (.1cm);
\draw[thick] (-1.2,1.2) circle (.1cm); 
\filldraw[white] (1.2,-1.2) circle (.1cm);
\draw[thick] (1.2,-1.2) circle (.1cm); 
\filldraw[white] (-2.4,2.4) circle (.1cm);
\draw[thick] (-2.4,2.4) circle (.1cm); 
\filldraw[white] (0,0) circle (.1cm);
\draw[thick] (0,0) circle (.1cm); 
\filldraw[white] (2.4,-2.4) circle (.1cm);
\draw[thick] (2.4,-2.4) circle (.1cm); 
\filldraw[\XColor] (-1.8,.6) circle (.07cm);
\filldraw[\XColor] (-.6,1.8) circle (.07cm);
\filldraw[\YColor] (1.8,-.6) circle (.07cm);
\filldraw[\YColor] (.6,-1.8) circle (.07cm);
\filldraw[\QsColor] (-.6,-.6) circle (.07cm);
\filldraw[\QsColor] (.6,.6) circle (.07cm);
\filldraw[\PsColor] (-3,1.8) circle (.07cm);
\filldraw[\PsColor] (-1.8,3) circle (.07cm);
\filldraw[\RsColor] (3,-1.8) circle (.07cm);
\filldraw[\RsColor] (1.8,-3) circle (.07cm);
\node at (-1.8,-3.8) {\small{$F(X)$}};
\node at (.6,-3.8) {\small{$F(Y)$}};
\node at (1.8,-3.8) {\small{$F(R)$}};
\node at (3,-3.8) {\small{$\varphi_c$}};
\node at (4.2,-3.8) {\small{$G(R)$}};
\node at (1.8,3.8) {\small{$G(Y)$}};
\node at (-.6,3.8) {\small{$G(X)$}};
\node at (-1.8,3.8) {\small{$G(P)$}};
\node at (-3,3.8) {\small{$\varphi_a$}};
\node at (-4.2,3.8) {\small{$F(P)$}};
}
\\
&=
\tikzmath[scale=.5, transform shape]{
\begin{scope}
\clip[rounded corners = 5] (-4.8,-3.6) rectangle (4.8,3.6);
\filldraw[primedregion=\AColor] (-1.2,1.2) -- (-2.4,2.4) .. controls ++(135:.2cm) and ++(270:.4cm) .. (-3,3.6) -- (-4.2,3.6) -- (-4.2,-.6) -- (-1.8,.6);
\filldraw[boxregion=\AColor] (-1.2,1.2) -- (-2.4,2.4) .. controls ++(135:.2cm) and ++(270:.4cm) .. (-3,3.6) -- (-1.8,3.6) -- (-1.8,3) -- (-.6,1.8);
\filldraw[primedregion=\BColor] (.6,-1.8) -- (-1.8,.6) -- (-1.2,1.2) -- (1.2,-1.2);
\filldraw[boxregion=\BColor] (1.8,-.6) -- (-.6,1.8)-- (-1.2,1.2) -- (1.2,-1.2);
\filldraw[primedregion=\CColor] (3,-3.6) .. controls ++(90:.4cm) and ++(-45:.2cm) .. (2.4,-2.4) -- (1.2,-1.2) -- (.6,-1.8) -- (1.8,-3) -- (1.8,-3.6);
\filldraw[boxregion=\CColor] (3,-3.6) .. controls ++(90:.4cm) and ++(-45:.2cm) .. (2.4,-2.4) -- (1.2,-1.2) -- (1.8,-.6) -- (4.2,-.6) -- (4.2,-3.6);
\filldraw[primedregion=\PrColor] (-1.8,-3.6) -- (-1.8,.6) -- (-3.6,2.4) .. controls ++(135:.2cm) and ++(270:.4cm) .. (-4.2,3.6) -- (-4.8,3.6) -- (-4.8,-3.6);
\filldraw[boxregion=\PrColor] (-1.8,3.6) -- (-1.8,3) -- (-.6,1.8) -- (-.6,3.6);
\filldraw[primedregion=\BColor] (-1.8,-3.6) -- (-1.8,.6) -- (.6,-1.8) -- (.6,-3.6);
\filldraw[boxregion=\QrColor] (-.6,1.8) -- (-.6,3.6) -- (1.8,3.6) -- (1.8,-.6);
\filldraw[primedregion=\RrColor] (1.8,-3.6) -- (1.8,-3) -- (.6,-1.8) -- (.6,-3.6);
\filldraw[boxregion=\RrColor] (4.2,-3.6) .. controls ++(90:.4cm) and ++(-45:.2cm) .. (3.6,-2.4) -- (1.8,-.6) -- (1.8,3.6) -- (4.8,3.6) -- (4.8,-3.6);
\end{scope}
\draw[\phiColor,thick] (3,-3.6) .. controls ++(90:.4cm) and ++(-45:.2cm) .. (2.4,-2.4) -- (-2.4,2.4) .. controls ++(135:.2cm) and ++(270:.4cm) .. (-3,3.6);
\draw[\XColor,thick] (-1.8,-3.6) -- (-1.8,.6) -- (-.6,1.8) -- (-.6,3.6);
\draw[\YColor,thick] (1.8,3.6) -- (1.8,-.6) -- (.6,-1.8) -- (.6,-3.6);
\draw[\QsColor,thick] (1.5,-.9) -- (-.9,1.5); 
\draw[\QsColor,thick] (1.8,-.6) -- (-.6,1.8); 
\draw[\QsColor,thick] (.3,.3) -- (.6,.6); 
\draw[\PsColor,thick] (-1.8,3.6) -- (-1.8,3) -- (-.6,1.8);
\draw[\PsColor,thick] (-1.8,.6) -- (-3.6,2.4) .. controls ++(135:.2cm) and ++(270:.4cm) .. (-4.2,3.6);
\draw[\PsColor,thick] (-3,1.8) -- (-1.8,3);
\draw[\RsColor,thick] (.6,-1.8) -- (1.8,-3) -- (1.8,-3.6);
\draw[\RsColor,thick] (4.2,-3.6) .. controls ++(90:.4cm) and ++(-45:.2cm) .. (3.6,-2.4) -- (1.8,-.6);
\draw[\RsColor,thick] (1.8,-3) -- (3,-1.8);
\filldraw[white] (-1.2,1.2) circle (.1cm);
\draw[thick] (-1.2,1.2) circle (.1cm); 
\filldraw[white] (1.2,-1.2) circle (.1cm);
\draw[thick] (1.2,-1.2) circle (.1cm); 
\filldraw[white] (-2.4,2.4) circle (.1cm);
\draw[thick] (-2.4,2.4) circle (.1cm); 
\filldraw[white] (2.4,-2.4) circle (.1cm);
\draw[thick] (2.4,-2.4) circle (.1cm); 
\filldraw[\XColor] (-1.8,.6) circle (.07cm);
\filldraw[\XColor] (-.6,1.8) circle (.07cm);
\filldraw[\XColor] (-.9,1.5) circle (.07cm);
\filldraw[\YColor] (1.5,-.9) circle (.07cm);
\filldraw[\YColor] (1.8,-.6) circle (.07cm);
\filldraw[\YColor] (.6,-1.8) circle (.07cm);
\filldraw[\QsColor] (.3,.3) circle (.07cm); 
\filldraw[\QsColor] (.6,.6) circle (.07cm);
\filldraw[\PsColor] (-3,1.8) circle (.07cm);
\filldraw[\PsColor] (-1.8,3) circle (.07cm);
\filldraw[\RsColor] (3,-1.8) circle (.07cm);
\filldraw[\RsColor] (1.8,-3) circle (.07cm);
\node at (-1.8,-3.8) {\small{$F(X)$}};
\node at (.6,-3.8) {\small{$F(Y)$}};
\node at (1.8,-3.8) {\small{$F(R)$}};
\node at (3,-3.8) {\small{$\varphi_c$}};
\node at (4.2,-3.8) {\small{$G(R)$}};
\node at (1.8,3.8) {\small{$G(Y)$}};
\node at (-.6,3.8) {\small{$G(X)$}};
\node at (-1.8,3.8) {\small{$G(P)$}};
\node at (-3,3.8) {\small{$\varphi_a$}};
\node at (-4.2,3.8) {\small{$F(P)$}};
}
=
\tikzmath[scale=.5, transform shape]{
\begin{scope}
\clip[rounded corners = 5] (-4.8,-3.6) rectangle (4.8,3.6);
\filldraw[primedregion=\AColor] (-1.2,1.2) -- (-2.4,2.4) .. controls ++(135:.2cm) and ++(270:.4cm) .. (-3,3.6) -- (-4.2,3.6) -- (-4.2,-.6) -- (-1.8,.6);
\filldraw[boxregion=\AColor] (-1.2,1.2) -- (-2.4,2.4) .. controls ++(135:.2cm) and ++(270:.4cm) .. (-3,3.6) -- (-1.8,3.6) -- (-1.8,3) -- (-.6,1.8);
\filldraw[primedregion=\BColor] (.6,-1.8) -- (-1.8,.6) -- (-1.2,1.2) -- (1.2,-1.2);
\filldraw[boxregion=\BColor] (1.8,-.6) -- (-.6,1.8)-- (-1.2,1.2) -- (1.2,-1.2);
\filldraw[primedregion=\CColor] (3,-3.6) .. controls ++(90:.4cm) and ++(-45:.2cm) .. (2.4,-2.4) -- (1.2,-1.2) -- (.6,-1.8) -- (1.8,-3) -- (1.8,-3.6);
\filldraw[boxregion=\CColor] (3,-3.6) .. controls ++(90:.4cm) and ++(-45:.2cm) .. (2.4,-2.4) -- (1.2,-1.2) -- (1.8,-.6) -- (4.2,-.6) -- (4.2,-3.6);
\filldraw[primedregion=\PrColor] (-1.8,-3.6) -- (-1.8,.6) -- (-3.6,2.4) .. controls ++(135:.2cm) and ++(270:.4cm) .. (-4.2,3.6) -- (-4.8,3.6) -- (-4.8,-3.6);
\filldraw[boxregion=\PrColor] (-1.8,3.6) -- (-1.8,3) -- (-.6,1.8) -- (-.6,3.6);
\filldraw[primedregion=\BColor] (-1.8,-3.6) -- (-1.8,.6) -- (.6,-1.8) -- (.6,-3.6);
\filldraw[boxregion=\QrColor] (-.6,1.8) -- (-.6,3.6) -- (1.8,3.6) -- (1.8,-.6);
\filldraw[primedregion=\RrColor] (1.8,-3.6) -- (1.8,-3) -- (.6,-1.8) -- (.6,-3.6);
\filldraw[boxregion=\RrColor] (4.2,-3.6) .. controls ++(90:.4cm) and ++(-45:.2cm) .. (3.6,-2.4) -- (1.8,-.6) -- (1.8,3.6) -- (4.8,3.6) -- (4.8,-3.6);
\end{scope}
\draw[\phiColor,thick] (3,-3.6) .. controls ++(90:.4cm) and ++(-45:.2cm) .. (2.4,-2.4) -- (-2.4,2.4) .. controls ++(135:.2cm) and ++(270:.4cm) .. (-3,3.6);
\draw[\XColor,thick] (-1.8,-3.6) -- (-1.8,.6) -- (-.6,1.8) -- (-.6,3.6);
\draw[\YColor,thick] (1.8,3.6) -- (1.8,-.6) -- (.6,-1.8) -- (.6,-3.6);
\draw[\QsColor,thick] (1.8,-.6) -- (-.6,1.8); 
\draw[\PsColor,thick] (-1.8,3.6) -- (-1.8,3) -- (-.6,1.8);
\draw[\PsColor,thick] (-1.8,.6) -- (-3.6,2.4) .. controls ++(135:.2cm) and ++(270:.4cm) .. (-4.2,3.6);
\draw[\PsColor,thick] (-3,1.8) -- (-1.8,3);
\draw[\RsColor,thick] (.6,-1.8) -- (1.8,-3) -- (1.8,-3.6);
\draw[\RsColor,thick] (4.2,-3.6) .. controls ++(90:.4cm) and ++(-45:.2cm) .. (3.6,-2.4) -- (1.8,-.6);
\draw[\RsColor,thick] (1.8,-3) -- (3,-1.8);
\filldraw[white] (-1.2,1.2) circle (.1cm);
\draw[thick] (-1.2,1.2) circle (.1cm); 
\filldraw[white] (1.2,-1.2) circle (.1cm);
\draw[thick] (1.2,-1.2) circle (.1cm); 
\filldraw[white] (-2.4,2.4) circle (.1cm);
\draw[thick] (-2.4,2.4) circle (.1cm); 
\filldraw[white] (2.4,-2.4) circle (.1cm);
\draw[thick] (2.4,-2.4) circle (.1cm); 
\filldraw[\XColor] (-1.8,.6) circle (.07cm);
\filldraw[\XColor] (-.6,1.8) circle (.07cm);
\filldraw[\YColor] (1.8,-.6) circle (.07cm);
\filldraw[\YColor] (.6,-1.8) circle (.05cm);
\filldraw[\PsColor] (-3,1.8) circle (.07cm);
\filldraw[\PsColor] (-1.8,3) circle (.07cm);
\filldraw[\RsColor] (3,-1.8) circle (.07cm);
\filldraw[\RsColor] (1.8,-3) circle (.07cm);
\node at (-1.8,-3.8) {\small{$F(X)$}};
\node at (.6,-3.8) {\small{$F(Y)$}};
\node at (1.8,-3.8) {\small{$F(R)$}};
\node at (3,-3.8) {\small{$\varphi_c$}};
\node at (4.2,-3.8) {\small{$G(R)$}};
\node at (1.8,3.8) {\small{$G(Y)$}};
\node at (-.6,3.8) {\small{$G(X)$}};
\node at (-1.8,3.8) {\small{$G(P)$}};
\node at (-3,3.8) {\small{$\varphi_a$}};
\node at (-4.2,3.8) {\small{$F(P)$}};
}
=
\tikzmath[scale=.7, transform shape]{
\begin{scope}
\clip[rounded corners = 5] (-.6,0) rectangle (3,4.2);
\filldraw[primedregion=\PrColor] (0,0) -- (0,2) .. controls ++(90:.4cm) and ++(-135:.2cm) .. (.6,2.6) .. controls ++(135:.2cm) and ++(270:.4cm) .. (0,3.2) -- (0,4.2) -- (-.6,4.2) -- (-.6,0); 
\filldraw[primedregion=\BColor] (1.2,0) -- (1.2,.4) .. controls ++(90:.4cm) and ++(-135:.2cm) .. (1.8,1) .. controls ++(135:.2cm) and ++(270:.4cm) .. (1.2,1.6) -- (1.2,2) .. controls ++(90:.4cm) and ++(-45:.2cm) .. (.6,2.6) .. controls ++(-135:.2cm) and ++(90:.4cm) .. (0,2) -- (0,0);
\filldraw[primedregion=\RrColor] (2.4,0) -- (2.4,.4) .. controls ++(90:.4cm) and ++(-45:.2cm) .. (1.8,1) .. controls ++(-135:.2cm) and ++(90:.4cm) .. (1.2,.4) -- (1.2,0);
\filldraw[boxregion=\PrColor] (0,4.2) -- (0,3.2) .. controls ++(270:.4cm) and ++(135:.2cm) .. (.6,2.6) .. controls ++(45:.2cm) and ++(270:.4cm) .. (1.2,3.2) -- (1.2,4.2);
\filldraw[boxregion=\BColor] (2.4,3.6) -- (2.4,1.6) .. controls ++(270:.4cm) and ++(45:.2cm) .. (1.8,1) .. controls ++(135:.2cm) and ++(270:.4cm) .. (1.2,1.6) -- (1.2,2) .. controls ++(90:.4cm) and ++(-45:.2cm) .. (.6,2.6) .. controls ++(45:.2cm) and ++(270:.4cm) .. (1.2,3.2) -- (1.2,3.6);
\filldraw[boxregion=\QrColor] (1.2,3.6) rectangle (2.4,4.2);
\filldraw[boxregion=\RrColor] (2.4,0) -- (2.4,.4) .. controls ++(90:.4cm) and ++(-45:.2cm) .. (1.8,1) .. controls ++(45:.2cm) and ++(270:.4cm) .. (2.4,1.6) -- (2.4,4.2) -- (3,4.2) -- (3,0);
\filldraw[\BColor] (.8,1.6) rectangle (2.2,2);
\end{scope}
\draw[\XColor,thick] (0,0) -- (0,.4) -- (0,2) .. controls ++(90:.6cm) and ++(270:.6cm) .. (1.2,3.2) -- (1.2,4.2);
\draw[\YColor,thick] (1.2,0) -- (1.2,.4) .. controls ++(90:.6cm) and ++(270:.6cm) .. (2.4,1.6) -- (2.4,4.2);
\draw[\phiColor,thick] (2.4,0) -- (2.4,.4) .. controls ++(90:.6cm) and ++(270:.6cm) .. (1.2,1.6) -- (1.2,2) .. controls ++(90:.6cm) and ++(270:.6cm) .. (0,3.2) -- (0,4.2);
\draw[\QsColor,thick] (1.2,3.6) -- (2.4,3.6);
\filldraw[white] (1.8,1) circle (.1cm);
\draw[thick] (1.8,1) circle (.1cm); 
\filldraw[white] (.6,2.6) circle (.1cm);
\draw[thick] (.6,2.6) circle (.1cm); 
\node at (0,-.2) {\scriptsize{$F(X)$}};
\node at (1.2,-.2) {\scriptsize{$F(Y)$}};
\node at (1.2,4.4) {\scriptsize{$G(X)$}};
\node at (2.4,4.4) {\scriptsize{$G(Y)$}};
\node at (2.4,-.2) {\scriptsize{$\QSys(\varphi)_R$}};
\node at (1.5,1.8) {\scriptsize{$\QSys(\varphi)_b$}};
\node at (0,4.4) {\scriptsize{$\QSys(\varphi)_P$}};
}
\end{align*}
This relation implies the monoidality coherence condition:
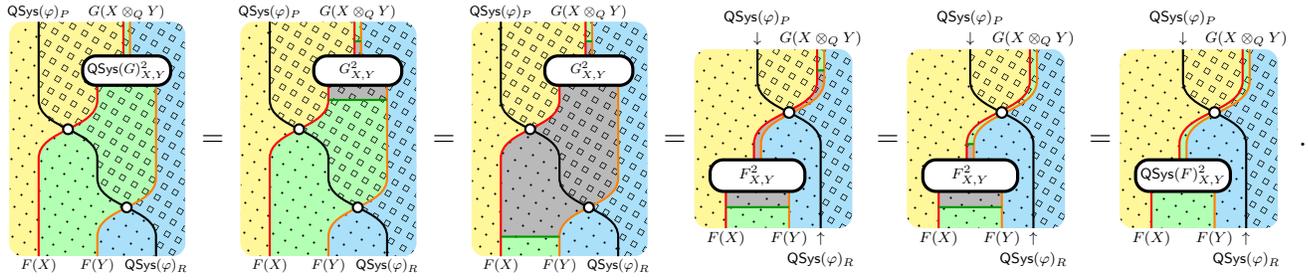
\[
\tikzmath[scale=.65, transform shape]{
\begin{scope}
\clip[rounded corners = 5] (-.6,0) rectangle (3,4.8);
\filldraw[primedregion=\PrColor] (0,0) -- (0,2) .. controls ++(90:.4cm) and ++(-135:.2cm) .. (.6,2.6) .. controls ++(135:.2cm) and ++(270:.4cm) .. (0,3.2) -- (0,4.8) -- (-.6,4.8) -- (-.6,0); 
\filldraw[primedregion=\QrColor] (1.2,0) -- (1.2,.4) .. controls ++(90:.4cm) and ++(-135:.2cm) .. (1.8,1) .. controls ++(135:.2cm) and ++(270:.4cm) .. (1.2,1.6) -- (1.2,2) .. controls ++(90:.4cm) and ++(-45:.2cm) .. (.6,2.6) .. controls ++(-135:.2cm) and ++(90:.4cm) .. (0,2) -- (0,0);
\filldraw[primedregion=\RrColor] (2.4,0) -- (2.4,.4) .. controls ++(90:.4cm) and ++(-45:.2cm) .. (1.8,1) .. controls ++(-135:.2cm) and ++(90:.4cm) .. (1.2,.4) -- (1.2,0);
\filldraw[boxregion=\PrColor] (0,4.8) -- (0,3.2) .. controls ++(270:.4cm) and ++(135:.2cm) .. (.6,2.6) .. controls ++(45:.2cm) and ++(270:.4cm) .. (1.2,3.2) -- (1.2,3.8) -- (1.733,3.8) -- (1.733,4.8);
\filldraw[boxregion=\RrColor] (2.4,0) -- (2.4,.4) .. controls ++(90:.4cm) and ++(-45:.2cm) .. (1.8,1) .. controls ++(45:.2cm) and ++(270:.4cm) .. (2.4,1.6) -- (2.4,3.2) -- (1.867,3.6) -- (1.867,4.8) -- (3,4.8) -- (3,0);
\filldraw[boxregion=\QrColor] (1.867,4.8) -- (1.867,3.8) -- (2.4,3.8) -- (2.4,1.6) .. controls ++(270:.4cm) and ++(45:.2cm) .. (1.8,1) .. controls ++(135:.2cm) and ++(270:.4cm) .. (1.2,1.6) -- (1.2,2) .. controls ++(90:.4cm) and ++(-45:.2cm) .. (.6,2.6) .. controls ++(45:.2cm) and ++(270:.4cm) .. (1.2,3.2) -- (1.2,3.8)-- (1.733,3.8) -- (1.733,4.8);
\end{scope}
\draw[\XColor,thick] (0,0) -- (0,.4) -- (0,2) .. controls ++(90:.6cm) and ++(270:.6cm) .. (1.2,3.2) -- (1.2,3.6);
\draw[\YColor,thick] (1.2,0) -- (1.2,.4) .. controls ++(90:.6cm) and ++(270:.6cm) .. (2.4,1.6) -- (2.4,3.6);
\draw[\phiColor,thick] (2.4,0) -- (2.4,.4) .. controls ++(90:.6cm) and ++(270:.6cm) .. (1.2,1.6) -- (1.2,2) .. controls ++(90:.6cm) and ++(270:.6cm) .. (0,3.2) -- (0,4.8);
\draw[\XColor,thick] (1.733,4) -- (1.733,4.8);
\draw[\YColor,thick] (1.867,4) -- (1.867,4.8);
\roundNbox{unshaded}{(1.8,3.8)}{.3}{.6}{.6}{\scriptsize{$\QSys(G)^2_{X,Y}$}}; 
\filldraw[white] (1.8,1) circle (.1cm);
\draw[thick] (1.8,1) circle (.1cm); 
\filldraw[white] (.6,2.6) circle (.1cm);
\draw[thick] (.6,2.6) circle (.1cm); 
\node at (0,-.2) {\scriptsize{$F(X)$}};
\node at (1.2,-.2) {\scriptsize{$F(Y)$}};
\node at (1.8,5) {\scriptsize{$G(X\xxo_Q Y)$}};
\node at (2.4,-.2) {\scriptsize{$\QSys(\varphi)_R$}};
\node at (0,5) {\scriptsize{$\QSys(\varphi)_P$}};
}
=
\tikzmath[scale=.65, transform shape]{
\begin{scope}
\clip[rounded corners = 5] (-.6,0) rectangle (3,4.8);
\filldraw[primedregion=\PrColor] (0,0) -- (0,2) .. controls ++(90:.4cm) and ++(-135:.2cm) .. (.6,2.6) .. controls ++(135:.2cm) and ++(270:.4cm) .. (0,3.2) -- (0,4.8) -- (-.6,4.8) -- (-.6,0); 
\filldraw[primedregion=\QrColor] (1.2,0) -- (1.2,.4) .. controls ++(90:.4cm) and ++(-135:.2cm) .. (1.8,1) .. controls ++(135:.2cm) and ++(270:.4cm) .. (1.2,1.6) -- (1.2,2) .. controls ++(90:.4cm) and ++(-45:.2cm) .. (.6,2.6) .. controls ++(-135:.2cm) and ++(90:.4cm) .. (0,2) -- (0,0);
\filldraw[primedregion=\RrColor] (2.4,0) -- (2.4,.4) .. controls ++(90:.4cm) and ++(-45:.2cm) .. (1.8,1) .. controls ++(-135:.2cm) and ++(90:.4cm) .. (1.2,.4) -- (1.2,0);
\filldraw[boxregion=\PrColor] (0,4.8) -- (0,3.2) .. controls ++(270:.4cm) and ++(135:.2cm) .. (.6,2.6) .. controls ++(45:.2cm) and ++(270:.4cm) .. (1.2,3.2) -- (1.2,3.8) -- (1.733,3.8) -- (1.733,4.8);
\filldraw[boxregion=\RrColor] (2.4,0) -- (2.4,.4) .. controls ++(90:.4cm) and ++(-45:.2cm) .. (1.8,1) .. controls ++(45:.2cm) and ++(270:.4cm) .. (2.4,1.6) -- (2.4,3.2) -- (1.867,3.6) -- (1.867,4.8) -- (3,4.8) -- (3,0);
\filldraw[boxregion=\QrColor] (1.867,4.8) -- (1.867,3.8) -- (2.4,3.8) -- (2.4,1.6) .. controls ++(270:.4cm) and ++(45:.2cm) .. (1.8,1) .. controls ++(135:.2cm) and ++(270:.4cm) .. (1.2,1.6) -- (1.2,2) .. controls ++(90:.4cm) and ++(-45:.2cm) .. (.6,2.6) .. controls ++(45:.2cm) and ++(270:.4cm) .. (1.2,3.2) -- (1.2,3.8)-- (1.733,3.8) -- (1.733,4.8);
\filldraw[boxregion=\BColor] (1.2,3.2) rectangle (2.4,3.8);
\filldraw[boxregion=\BColor] (1.733,4.4) rectangle (1.867,3.8);
\end{scope}
\draw[\XColor,thick] (0,0) -- (0,.4) -- (0,2) .. controls ++(90:.6cm) and ++(270:.6cm) .. (1.2,3.2) -- (1.2,3.6);
\draw[\YColor,thick] (1.2,0) -- (1.2,.4) .. controls ++(90:.6cm) and ++(270:.6cm) .. (2.4,1.6) -- (2.4,3.6);
\draw[\phiColor,thick] (2.4,0) -- (2.4,.4) .. controls ++(90:.6cm) and ++(270:.6cm) .. (1.2,1.6) -- (1.2,2) .. controls ++(90:.6cm) and ++(270:.6cm) .. (0,3.2) -- (0,4.8);
\draw[\XColor,thick] (1.733,4) -- (1.733,4.8);
\draw[\YColor,thick] (1.867,4) -- (1.867,4.8);
\draw[\QsColor,thick] (1.2,3.2) -- (2.4,3.2);
\draw[\QsColor,thick] (1.733,4.4) -- (1.867,4.4);
\roundNbox{unshaded}{(1.8,3.8)}{.3}{.6}{.6}{\scriptsize{$G^2_{X,Y}$}}; 
\filldraw[white] (1.8,1) circle (.1cm);
\draw[thick] (1.8,1) circle (.1cm); 
\filldraw[white] (.6,2.6) circle (.1cm);
\draw[thick] (.6,2.6) circle (.1cm); 
\node at (0,-.2) {\scriptsize{$F(X)$}};
\node at (1.2,-.2) {\scriptsize{$F(Y)$}};
\node at (1.8,5) {\scriptsize{$G(X\xxo_Q Y)$}};
\node at (2.4,-.2) {\scriptsize{$\QSys(\varphi)_R$}};
\node at (0,5) {\scriptsize{$\QSys(\varphi)_P$}};
}
=
\tikzmath[scale=.65, transform shape]{
\begin{scope}
\clip[rounded corners = 5] (-.6,0) rectangle (3,4.8);
\filldraw[primedregion=\PrColor] (0,0) -- (0,2) .. controls ++(90:.4cm) and ++(-135:.2cm) .. (.6,2.6) .. controls ++(135:.2cm) and ++(270:.4cm) .. (0,3.2) -- (0,4.8) -- (-.6,4.8) -- (-.6,0); 
\filldraw[primedregion=\BColor] (1.2,0) -- (1.2,.4) .. controls ++(90:.4cm) and ++(-135:.2cm) .. (1.8,1) .. controls ++(135:.2cm) and ++(270:.4cm) .. (1.2,1.6) -- (1.2,2) .. controls ++(90:.4cm) and ++(-45:.2cm) .. (.6,2.6) .. controls ++(-135:.2cm) and ++(90:.4cm) .. (0,2) -- (0,0);
\filldraw[primedregion=\RrColor] (2.4,0) -- (2.4,.4) .. controls ++(90:.4cm) and ++(-45:.2cm) .. (1.8,1) .. controls ++(-135:.2cm) and ++(90:.4cm) .. (1.2,.4) -- (1.2,0);
\filldraw[boxregion=\PrColor] (0,4.8) -- (0,3.2) .. controls ++(270:.4cm) and ++(135:.2cm) .. (.6,2.6) .. controls ++(45:.2cm) and ++(270:.4cm) .. (1.2,3.2) -- (1.2,3.8) -- (1.733,3.8) -- (1.733,4.8);
\filldraw[boxregion=\RrColor] (2.4,0) -- (2.4,.4) .. controls ++(90:.4cm) and ++(-45:.2cm) .. (1.8,1) .. controls ++(45:.2cm) and ++(270:.4cm) .. (2.4,1.6) -- (2.4,3.2) -- (1.867,3.6) -- (1.867,4.8) -- (3,4.8) -- (3,0);
\filldraw[boxregion=\BColor] (1.867,4.8) -- (1.867,3.8) -- (2.4,3.8) -- (2.4,1.6) .. controls ++(270:.4cm) and ++(45:.2cm) .. (1.8,1) .. controls ++(135:.2cm) and ++(270:.4cm) .. (1.2,1.6) -- (1.2,2) .. controls ++(90:.4cm) and ++(-45:.2cm) .. (.6,2.6) .. controls ++(45:.2cm) and ++(270:.4cm) .. (1.2,3.2) -- (1.2,3.8)-- (1.733,3.8) -- (1.733,4.8);
\filldraw[primedregion=\QrColor] (0,0) rectangle (1.2,.4);
\filldraw[boxregion=\QrColor] (1.733,4.4) rectangle (1.867,4.8);
\end{scope}
\draw[\XColor,thick] (0,0) -- (0,.4) -- (0,2) .. controls ++(90:.6cm) and ++(270:.6cm) .. (1.2,3.2) -- (1.2,3.6);
\draw[\YColor,thick] (1.2,0) -- (1.2,.4) .. controls ++(90:.6cm) and ++(270:.6cm) .. (2.4,1.6) -- (2.4,3.6);
\draw[\phiColor,thick] (2.4,0) -- (2.4,.4) .. controls ++(90:.6cm) and ++(270:.6cm) .. (1.2,1.6) -- (1.2,2) .. controls ++(90:.6cm) and ++(270:.6cm) .. (0,3.2) -- (0,4.8);
\draw[\XColor,thick] (1.733,4) -- (1.733,4.8);
\draw[\YColor,thick] (1.867,4) -- (1.867,4.8);
\draw[\QsColor,thick] (0,.4) -- (1.2,.4);
\draw[\QsColor,thick] (1.733,4.4) -- (1.867,4.4);
\roundNbox{unshaded}{(1.8,3.8)}{.3}{.6}{.6}{\scriptsize{$G^2_{X,Y}$}}; 
\filldraw[white] (1.8,1) circle (.1cm);
\draw[thick] (1.8,1) circle (.1cm); 
\filldraw[white] (.6,2.6) circle (.1cm);
\draw[thick] (.6,2.6) circle (.1cm); 
\node at (0,-.2) {\scriptsize{$F(X)$}};
\node at (1.2,-.2) {\scriptsize{$F(Y)$}};
\node at (1.8,5) {\scriptsize{$G(X\xxo_Q Y)$}};
\node at (2.4,-.2) {\scriptsize{$\QSys(\varphi)_R$}};
\node at (0,5) {\scriptsize{$\QSys(\varphi)_P$}};
}
=
\tikzmath[scale=.7, transform shape]{
\begin{scope}
\clip[rounded corners = 5] (-.6,-.4) rectangle (2.4,3);
\filldraw[primedregion=\PrColor] (0,-.4) -- (0,.8) -- (.533,.8) -- (.533,1.2) .. controls ++(90:.4cm) and ++(-135:.2cm) .. (1.2,1.8) .. controls ++(135:.2cm) and ++(270:.4cm) .. (.6,2.4) -- (.6,3) -- (-.6,3) -- (-.6,-.4); 
\filldraw[primedregion=\RrColor] (1.8,-.4) -- (1.8,1.2) .. controls ++(90:.4cm) and ++(-45:.2cm) .. (1.2,1.8) .. controls ++(-135:.2cm) and ++(90:.4cm) .. (.65,1.2) -- (.65,.8) -- (1.2,.8) -- (1.2,-.4);
\filldraw[primedregion=\BColor] (0,-.4) -- (0,.8) -- (.533,.8) -- (.533,1.2) .. controls ++(90:.4cm) and ++(-135:.2cm) .. (1.133,1.8) -- (1.267,1.8) .. controls ++(-135:.2cm) and ++(90:.4cm) .. (.667,1.2) -- (.667,.8) -- (1.2,.8) -- (1.2,-.4);
\filldraw[boxregion=\PrColor] (.6,3) -- (.6,2.4) .. controls ++(270:.4cm) and ++(135:.2cm) .. (1.2,1.8) .. controls ++(45:.2cm) and ++(270:.4cm) .. (1.75,2.4) -- (1.75,3);
\filldraw[boxregion=\RrColor] (1.8,-.4) -- (1.8,1.2) .. controls ++(90:.4cm) and ++(-45:.2cm) .. (1.2,1.8) .. controls ++(45:.2cm) and ++(270:.4cm) .. (1.85,2.4) -- (1.85,3) -- (2.4,3) -- (2.4,-.4);
\filldraw[boxregion=\BColor] (1.733,3) -- (1.733,2.4) .. controls ++(270:.4cm) and ++(45:.2cm) .. (1.133,1.8) -- (1.267,1.8) .. controls ++(45:.2cm) and ++(270:.4cm) .. (1.867,2.4) -- (1.867,3); 
\filldraw[primedregion=\QrColor] (0,0) rectangle (1.2,-.4);
\filldraw[boxregion=\QrColor] (1.733,2.6) rectangle (1.867,3);
\end{scope}
\draw[\XColor,thick] (0,-.4) -- (0,.6);
\draw[\YColor,thick] (1.2,-.4) -- (1.2,.6);
\draw[\XColor,thick] (.533,.6) -- (.533,1.2) .. controls ++(90:.6cm) and ++(270:.6cm) .. (1.733,2.4) -- (1.733,3);
\draw[\YColor,thick] (.667,.6) -- (.667,1.2) .. controls ++(90:.6cm) and ++(270:.6cm) .. (1.867,2.4) -- (1.867,3);
\draw[\phiColor,thick] (1.8,-.4) -- (1.8,1.2) .. controls ++(90:.6cm) and ++(270:.6cm) .. (.6,2.4) -- (.6,3);
\draw[\QsColor,thick] (0,0) -- (1.2,0);
\draw[\QsColor,thick] (1.733,2.6) -- (1.867,2.6);
\roundNbox{unshaded}{(.6,.6)}{.3}{.6}{.6}{\scriptsize{$F^2_{X,Y}$}}; 
\filldraw[white] (1.2,1.8) circle (.1cm);
\draw[thick] (1.2,1.8) circle (.1cm); 
\node at (0,-.6) {\scriptsize{$F(X)$}};
\node at (1.2,-.6) {\scriptsize{$F(Y)$}};
\node at (1.8,3.2) {\scriptsize{$G(X\xxo_Q Y)$}};
\node at (1.8,-.6) {\scriptsize{$\uparrow$}};
\node at (1.8,-1) {\scriptsize{$\QSys(\varphi)_R$}};
\node at (.6,3.2) {\scriptsize{$\downarrow$}};
\node at (.6,3.6) {\scriptsize{$\QSys(\varphi)_P$}};
}
=
\tikzmath[scale=.7, transform shape]{
\begin{scope}
\clip[rounded corners = 5] (-.6,-.4) rectangle (2.4,3);
\filldraw[primedregion=\PrColor] (0,-.4) -- (0,.8) -- (.533,.8) -- (.533,1.2) .. controls ++(90:.4cm) and ++(-135:.2cm) .. (1.2,1.8) .. controls ++(135:.2cm) and ++(270:.4cm) .. (.6,2.4) -- (.6,3) -- (-.6,3) -- (-.6,-.4); 
\filldraw[primedregion=\RrColor] (1.8,-.4) -- (1.8,1.2) .. controls ++(90:.4cm) and ++(-45:.2cm) .. (1.2,1.8) .. controls ++(-135:.2cm) and ++(90:.4cm) .. (.65,1.2) -- (.65,.8) -- (1.2,.8) -- (1.2,-.4);
\filldraw[primedregion=\QrColor] (0,-.4) -- (0,.8) -- (.533,.8) -- (.533,1.2) .. controls ++(90:.4cm) and ++(-135:.2cm) .. (1.133,1.8) -- (1.267,1.8) .. controls ++(-135:.2cm) and ++(90:.4cm) .. (.667,1.2) -- (.667,.8) -- (1.2,.8) -- (1.2,-.4);
\filldraw[boxregion=\PrColor] (.6,3) -- (.6,2.4) .. controls ++(270:.4cm) and ++(135:.2cm) .. (1.2,1.8) .. controls ++(45:.2cm) and ++(270:.4cm) .. (1.75,2.4) -- (1.75,3);
\filldraw[boxregion=\RrColor] (1.8,-.4) -- (1.8,1.2) .. controls ++(90:.4cm) and ++(-45:.2cm) .. (1.2,1.8) .. controls ++(45:.2cm) and ++(270:.4cm) .. (1.85,2.4) -- (1.85,3) -- (2.4,3) -- (2.4,-.4);
\filldraw[boxregion=\QrColor] (1.733,3) -- (1.733,2.4) .. controls ++(270:.4cm) and ++(45:.2cm) .. (1.133,1.8) -- (1.267,1.8) .. controls ++(45:.2cm) and ++(270:.4cm) .. (1.867,2.4) -- (1.867,3); 
\filldraw[primedregion=\BColor] (0,0) rectangle (1.2,.4);
\filldraw[boxregion=\BColor] (.533,.8) rectangle (.667,1.2);
\end{scope}
\draw[\XColor,thick] (0,-.4) -- (0,.6);
\draw[\YColor,thick] (1.2,-.4) -- (1.2,.6);
\draw[\XColor,thick] (.533,.6) -- (.533,1.2) .. controls ++(90:.6cm) and ++(270:.6cm) .. (1.733,2.4) -- (1.733,3);
\draw[\YColor,thick] (.667,.6) -- (.667,1.2) .. controls ++(90:.6cm) and ++(270:.6cm) .. (1.867,2.4) -- (1.867,3);
\draw[\phiColor,thick] (1.8,-.4) -- (1.8,1.2) .. controls ++(90:.6cm) and ++(270:.6cm) .. (.6,2.4) -- (.6,3);
\draw[\QsColor,thick] (0,0) -- (1.2,0);
\draw[\QsColor,thick] (.533,1.2) -- (.667,1.2);
\roundNbox{unshaded}{(.6,.6)}{.3}{.6}{.6}{\scriptsize{$F^2_{X,Y}$}}; 
\filldraw[white] (1.2,1.8) circle (.1cm);
\draw[thick] (1.2,1.8) circle (.1cm); 
\node at (0,-.6) {\scriptsize{$F(X)$}};
\node at (1.2,-.6) {\scriptsize{$F(Y)$}};
\node at (1.8,3.2) {\scriptsize{$G(X\xxo_Q Y)$}};
\node at (1.8,-.6) {\scriptsize{$\uparrow$}};
\node at (1.8,-1) {\scriptsize{$\QSys(\varphi)_R$}};
\node at (.6,3.2) {\scriptsize{$\downarrow$}};
\node at (.6,3.6) {\scriptsize{$\QSys(\varphi)_P$}};
}
=
\tikzmath[scale=.7, transform shape]{
\begin{scope}
\clip[rounded corners = 5] (-.6,-.4) rectangle (2.4,3);
\filldraw[primedregion=\PrColor] (0,-.4) -- (0,.8) -- (.533,.8) -- (.533,1.2) .. controls ++(90:.4cm) and ++(-135:.2cm) .. (1.2,1.8) .. controls ++(135:.2cm) and ++(270:.4cm) .. (.6,2.4) -- (.6,3) -- (-.6,3) -- (-.6,-.4); 
\filldraw[primedregion=\RrColor] (1.8,-.4) -- (1.8,1.2) .. controls ++(90:.4cm) and ++(-45:.2cm) .. (1.2,1.8) .. controls ++(-135:.2cm) and ++(90:.4cm) .. (.65,1.2) -- (.65,.8) -- (1.2,.8) -- (1.2,-.4);
\filldraw[primedregion=\QrColor] (0,-.4) -- (0,.8) -- (.533,.8) -- (.533,1.2) .. controls ++(90:.4cm) and ++(-135:.2cm) .. (1.133,1.8) -- (1.267,1.8) .. controls ++(-135:.2cm) and ++(90:.4cm) .. (.667,1.2) -- (.667,.8) -- (1.2,.8) -- (1.2,-.4);
\filldraw[boxregion=\PrColor] (.6,3) -- (.6,2.4) .. controls ++(270:.4cm) and ++(135:.2cm) .. (1.2,1.8) .. controls ++(45:.2cm) and ++(270:.4cm) .. (1.75,2.4) -- (1.75,3);
\filldraw[boxregion=\RrColor] (1.8,-.4) -- (1.8,1.2) .. controls ++(90:.4cm) and ++(-45:.2cm) .. (1.2,1.8) .. controls ++(45:.2cm) and ++(270:.4cm) .. (1.85,2.4) -- (1.85,3) -- (2.4,3) -- (2.4,-.4);
\filldraw[boxregion=\QrColor] (1.733,3) -- (1.733,2.4) .. controls ++(270:.4cm) and ++(45:.2cm) .. (1.133,1.8) -- (1.267,1.8) .. controls ++(45:.2cm) and ++(270:.4cm) .. (1.867,2.4) -- (1.867,3); 
\end{scope}
\draw[\XColor,thick] (0,-.4) -- (0,.6);
\draw[\YColor,thick] (1.2,-.4) -- (1.2,.6);
\draw[\XColor,thick] (.533,.6) -- (.533,1.2) .. controls ++(90:.6cm) and ++(270:.6cm) .. (1.733,2.4) -- (1.733,3);
\draw[\YColor,thick] (.667,.6) -- (.667,1.2) .. controls ++(90:.6cm) and ++(270:.6cm) .. (1.867,2.4) -- (1.867,3);
\draw[\phiColor,thick] (1.8,-.4) -- (1.8,1.2) .. controls ++(90:.6cm) and ++(270:.6cm) .. (.6,2.4) -- (.6,3);
\roundNbox{unshaded}{(.6,.6)}{.3}{.6}{.6}{\scriptsize{$\QSys(F)^2_{X,Y}$}}; 
\filldraw[white] (1.2,1.8) circle (.1cm);
\draw[thick] (1.2,1.8) circle (.1cm); 
\node at (0,-.6) {\scriptsize{$F(X)$}};
\node at (1.2,-.6) {\scriptsize{$F(Y)$}};
\node at (1.8,3.2) {\scriptsize{$G(X\xxo_Q Y)$}};
\node at (1.8,-.6) {\scriptsize{$\uparrow$}};
\node at (1.8,-1) {\scriptsize{$\QSys(\varphi)_R$}};
\node at (.6,3.2) {\scriptsize{$\downarrow$}};
\node at (.6,3.6) {\scriptsize{$\QSys(\varphi)_P$}};
}\,.
\]
Unitality is checked similarly.
Finally, to check naturality, for a 2-cell $f\in\cC({}_PX_Q\to {}_PZ_Q)$:
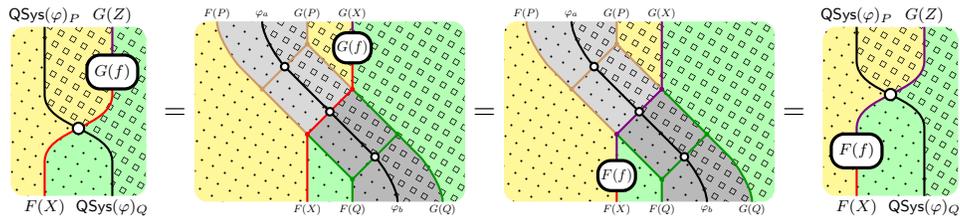
\[
\tikzmath[scale=.75, transform shape]{
\begin{scope}
\clip[rounded corners = 5] (-.6,0) rectangle (1.8,3);
\filldraw[primedregion=\PrColor] (0,0) -- (0,.6) .. controls ++(90:.4cm) and ++(-135:.2cm) .. (.6,1.2) .. controls ++(135:.2cm) and ++(270:.4cm) .. (0,1.8) -- (0,3) -- (-.6,3) -- (-.6,0); 
\filldraw[primedregion=\QrColor] (1.2,0) -- (1.2,.6) .. controls ++(90:.4cm) and ++(-45:.2cm) .. (.6,1.2) .. controls ++(-135:.2cm) and ++(90:.4cm) .. (0,.6) -- (0,0);
\filldraw[boxregion=\PrColor] (0,3) -- (0,1.8) .. controls ++(270:.4cm) and ++(135:.2cm) .. (.6,1.2) .. controls ++(45:.2cm) and ++(270:.4cm) .. (1.2,1.8) -- (1.2,3);
\filldraw[boxregion=\QrColor] (1.2,0) -- (1.2,.6) .. controls ++(90:.4cm) and ++(-45:.2cm) .. (.6,1.2) .. controls ++(45:.2cm) and ++(270:.4cm) .. (1.2,1.8) -- (1.2,3) -- (1.8,3) -- (1.8,0);
\end{scope}
\draw[\XColor,thick] (0,0) -- (0,.6) .. controls ++(90:.6cm) and ++(270:.6cm) .. (1.2,1.8) -- (1.2,2.2);
\draw[\ZColor,thick] (1.2,2.2) -- (1.2,3);
\draw[\phiColor,thick] (1.2,0) -- (1.2,.6) .. controls ++(90:.6cm) and ++(270:.6cm) .. (0,1.8) -- (0,3);
\roundNbox{unshaded}{(1.2,2.2)}{.3}{.15}{.15}{\scriptsize{$G(f)$}}; 
\filldraw[white] (.6,1.2) circle (.1cm);
\draw[\phiColor,thick] (.6,1.2) circle (.1cm); 
\node at (0,-.2) {\scriptsize{$F(X)$}};
\node at (1.2,3.2) {\scriptsize{$G(Z)$}};
\node at (1.2,-.2) {\scriptsize{$\QSys(\varphi)_Q$}};
\node at (0,3.2) {\scriptsize{$\QSys(\varphi)_P$}};
} 
=
\tikzmath[scale=.5, transform shape]{
\begin{scope}
\clip[rounded corners = 5] (-3.6,-2.4) rectangle (3.6,2.4);
\filldraw[primedregion=\AColor] (0,0) -- (-1.2,1.2) .. controls ++(135:.2cm) and ++(270:.4cm) .. (-1.8,2.4) -- (-3,2.4) -- (-3,-.6) -- (-.6,-.6);
\filldraw[boxregion=\AColor] (0,0) -- (-1.2,1.2) .. controls ++(135:.2cm) and ++(270:.4cm) .. (-1.8,2.4) -- (-.6,2.4) -- (-.6,1.8) -- (.6,.6);
\filldraw[primedregion=\BColor] (1.8,-2.4) .. controls ++(90:.4cm) and ++(-45:.2cm) .. (1.2,-1.2) -- (0,0) -- (-.6,-.6) -- (.6,-1.8) -- (.6,-2.4);
\filldraw[boxregion=\BColor] (1.8,-2.4) .. controls ++(90:.4cm) and ++(-45:.2cm) .. (1.2,-1.2) -- (0,0) -- (.6,.6) -- (3,.6) -- (3,-2.4);
\filldraw[primedregion=\PrColor] (-.6,-2.4) -- (-.6,-.6) -- (-2.4,1.2) .. controls ++(135:.2cm) and ++(270:.4cm) .. (-3,2.4) -- (-3.6,2.4) -- (-3.6,-2.4);
\filldraw[boxregion=\PrColor] (-.6,2.4) -- (-.6,1.8) -- (.6,.6) -- (.6,2.4); 
\filldraw[primedregion=\QrColor] (.6,-2.4) -- (.6,-1.8) -- (-.6,-.6) -- (-.6,-2.4); 
\filldraw[boxregion=\QrColor] (3,-2.4) .. controls ++(90:.4cm) and ++(-45:.2cm) .. (2.4,-1.2) -- (.6,.6) -- (.6,2.4) -- (3.6,2.4) -- (3.6,-2.4);
\end{scope}
\draw[\phiColor,thick] (1.8,-2.4) .. controls ++(90:.4cm) and ++(-45:.2cm) .. (1.2,-1.2) -- (-1.2,1.2) .. controls ++(135:.2cm) and ++(270:.4cm) .. (-1.8,2.4);
\draw[\XColor,thick] (-.6,-2.4) -- (-.6,-.6) -- (.6,.6) -- (.6,1.7);
\draw[\ZColor,thick] (.6,1.7) -- (.6,2.4);
\draw[\PsColor,thick] (-.6,-.6) -- (-2.4,1.2) .. controls ++(135:.2cm) and ++(270:.4cm) .. (-3,2.4);
\draw[\PsColor,thick] (.6,.6) -- (-.6,1.8) -- (-.6,2.4);
\draw[\PsColor,thick] (-1.8,.6) -- (-.6,1.8);
\draw[\QsColor,thick] (-.6,-.6) -- (.6,-1.8) -- (.6,-2.4);
\draw[\QsColor,thick] (3,-2.4) .. controls ++(90:.4cm) and ++(-45:.2cm) .. (2.4,-1.2) -- (.6,.6);
\draw[\QsColor,thick] (1.8,-.6) -- (.6,-1.8);
\filldraw[\XColor] (-.6,-.6) circle (.05cm);
\filldraw[\XColor] (.6,.6) circle (.05cm);
\filldraw[\PsColor] (-1.8,.6) circle (.05cm);
\filldraw[\PsColor] (-.6,1.8) circle (.05cm);
\filldraw[\QsColor] (1.8,-.6) circle (.05cm);
\filldraw[\QsColor] (.6,-1.8) circle (.05cm);
\filldraw[white] (0,0) circle (.1cm);
\draw[thick] (0,0) circle (.1cm); 
\filldraw[white] (-1.2,1.2) circle (.1cm);
\draw[thick] (-1.2,1.2) circle (.1cm); 
\filldraw[white] (1.2,-1.2) circle (.1cm);
\draw[thick] (1.2,-1.2) circle (.1cm); 
\roundNbox{unshaded}{(.6,1.7)}{.4}{.1}{.1}{\normalsize{$G(f)$}}; 
\node at (-.6,-2.6) {\scriptsize{$F(X)$}};
\node at (.6,-2.6) {\scriptsize{$F(Q)$}};
\node at (1.8,-2.6) {\scriptsize{$\varphi_b$}};
\node at (3,-2.6) {\scriptsize{$G(Q)$}};
\node at (.6,2.6) {\scriptsize{$G(X)$}};
\node at (-.6,2.6) {\scriptsize{$G(P)$}};
\node at (-1.8,2.6) {\scriptsize{$\varphi_a$}};
\node at (-3,2.6) {\scriptsize{$F(P)$}};
}
=
\tikzmath[scale=.5, transform shape]{
\begin{scope}
\clip[rounded corners = 5] (-3.6,-2.4) rectangle (3.6,2.4);
\filldraw[primedregion=\AColor] (0,0) -- (-1.2,1.2) .. controls ++(135:.2cm) and ++(270:.4cm) .. (-1.8,2.4) -- (-3,2.4) -- (-3,-.6) -- (-.6,-.6);
\filldraw[boxregion=\AColor] (0,0) -- (-1.2,1.2) .. controls ++(135:.2cm) and ++(270:.4cm) .. (-1.8,2.4) -- (-.6,2.4) -- (-.6,1.8) -- (.6,.6);
\filldraw[primedregion=\BColor] (1.8,-2.4) .. controls ++(90:.4cm) and ++(-45:.2cm) .. (1.2,-1.2) -- (0,0) -- (-.6,-.6) -- (.6,-1.8) -- (.6,-2.4);
\filldraw[boxregion=\BColor] (1.8,-2.4) .. controls ++(90:.4cm) and ++(-45:.2cm) .. (1.2,-1.2) -- (0,0) -- (.6,.6) -- (3,.6) -- (3,-2.4);
\filldraw[primedregion=\PrColor] (-.6,-2.4) -- (-.6,-.6) -- (-2.4,1.2) .. controls ++(135:.2cm) and ++(270:.4cm) .. (-3,2.4) -- (-3.6,2.4) -- (-3.6,-2.4);
\filldraw[boxregion=\PrColor] (-.6,2.4) -- (-.6,1.8) -- (.6,.6) -- (.6,2.4); 
\filldraw[primedregion=\QrColor] (.6,-2.4) -- (.6,-1.8) -- (-.6,-.6) -- (-.6,-2.4); 
\filldraw[boxregion=\QrColor] (3,-2.4) .. controls ++(90:.4cm) and ++(-45:.2cm) .. (2.4,-1.2) -- (.6,.6) -- (.6,2.4) -- (3.6,2.4) -- (3.6,-2.4);
\end{scope}
\draw[\phiColor,thick] (1.8,-2.4) .. controls ++(90:.4cm) and ++(-45:.2cm) .. (1.2,-1.2) -- (-1.2,1.2) .. controls ++(135:.2cm) and ++(270:.4cm) .. (-1.8,2.4);
\draw[\XColor,thick] (-.6,-2.4) -- (-.6,-1.7);
\draw[\ZColor,thick] (-.6,-1.7) -- (-.6,-.6) -- (.6,.6) -- (.6,2.4);
\draw[\PsColor,thick] (-.6,-.6) -- (-2.4,1.2) .. controls ++(135:.2cm) and ++(270:.4cm) .. (-3,2.4);
\draw[\PsColor,thick] (.6,.6) -- (-.6,1.8) -- (-.6,2.4);
\draw[\PsColor,thick] (-1.8,.6) -- (-.6,1.8);
\draw[\QsColor,thick] (-.6,-.6) -- (.6,-1.8) -- (.6,-2.4);
\draw[\QsColor,thick] (3,-2.4) .. controls ++(90:.4cm) and ++(-45:.2cm) .. (2.4,-1.2) -- (.6,.6);
\draw[\QsColor,thick] (1.8,-.6) -- (.6,-1.8);
\filldraw[\ZColor] (-.6,-.6) circle (.05cm);
\filldraw[\ZColor] (.6,.6) circle (.05cm);
\filldraw[\PsColor] (-1.8,.6) circle (.05cm);
\filldraw[\PsColor] (-.6,1.8) circle (.05cm);
\filldraw[\QsColor] (1.8,-.6) circle (.05cm);
\filldraw[\QsColor] (.6,-1.8) circle (.05cm);
\filldraw[white] (0,0) circle (.1cm);
\draw[thick] (0,0) circle (.1cm); 
\filldraw[white] (-1.2,1.2) circle (.1cm);
\draw[thick] (-1.2,1.2) circle (.1cm); 
\filldraw[white] (1.2,-1.2) circle (.1cm);
\draw[thick] (1.2,-1.2) circle (.1cm); 
\roundNbox{unshaded}{(-.6,-1.7)}{.4}{.1}{.1}{\normalsize{$F(f)$}}; 
\node at (-.6,-2.6) {\scriptsize{$F(X)$}};
\node at (.6,-2.6) {\scriptsize{$F(Q)$}};
\node at (1.8,-2.6) {\scriptsize{$\varphi_b$}};
\node at (3,-2.6) {\scriptsize{$G(Q)$}};
\node at (.6,2.6) {\scriptsize{$G(X)$}};
\node at (-.6,2.6) {\scriptsize{$G(P)$}};
\node at (-1.8,2.6) {\scriptsize{$\varphi_a$}};
\node at (-3,2.6) {\scriptsize{$F(P)$}};
}
=
\tikzmath[scale=.75, transform shape]{
\begin{scope}
\clip[rounded corners = 5] (-.6,0) rectangle (1.8,3);
\filldraw[primedregion=\PrColor] (0,0) -- (0,1.2) .. controls ++(90:.4cm) and ++(-135:.2cm) .. (.6,1.8) .. controls ++(135:.2cm) and ++(270:.4cm) .. (0,2.4) -- (0,3) -- (-.6,3) -- (-.6,0); 
\filldraw[primedregion=\QrColor] (1.2,0) -- (1.2,1.2) .. controls ++(90:.4cm) and ++(-45:.2cm) .. (.6,1.8) .. controls ++(-135:.2cm) and ++(90:.4cm) .. (0,1.2) -- (0,0);
\filldraw[boxregion=\PrColor] (0,3) -- (0,2.4) .. controls ++(270:.4cm) and ++(135:.2cm) .. (.6,1.8) .. controls ++(45:.2cm) and ++(270:.4cm) .. (1.2,2.4) -- (1.2,3);
\filldraw[boxregion=\QrColor] (1.2,0) -- (1.2,1.2) .. controls ++(90:.4cm) and ++(-45:.2cm) .. (.6,1.8) .. controls ++(45:.2cm) and ++(270:.4cm) .. (1.2,2.4) -- (1.2,3) -- (1.8,3) -- (1.8,0);
\end{scope}
\draw[\XColor,thick] (0,0) -- (0,.8);
\draw[\ZColor,thick] (0,.8) -- (0,1.2) .. controls ++(90:.6cm) and ++(270:.6cm) .. (1.2,2.4) -- (1.2,3);
\draw[\phiColor,thick] (1.2,0) -- (1.2,1.2) .. controls ++(90:.6cm) and ++(270:.6cm) .. (0,2.4) -- (0,3);
\roundNbox{unshaded}{(0,.8)}{.3}{.15}{.15}{\scriptsize{$F(f)$}}; 
\filldraw[white] (.6,1.8) circle (.1cm);
\draw[thick] (.6,1.8) circle (.1cm); 
\node at (0,-.2) {\scriptsize{$F(X)$}};
\node at (1.2,3.2) {\scriptsize{$G(Z)$}};
\node at (1.2,-.2) {\scriptsize{$\QSys(\varphi)_Q$}};
\node at (0,3.2) {\scriptsize{$\QSys(\varphi)_P$}};
}\,.
\]
\end{construction}

\begin{construction}
\label{construction:QSys(n)}
Suppose $n: \varphi \Rrightarrow \psi$ is a bounded modification between $\dag$-transformations.
We define a bounded modification $\QSys(n): \QSys(\varphi) \Rrightarrow \QSys(\psi)$ as follows.
Given a Q-system ${}_bQ_b\in\QSys(\cC)$, we define 
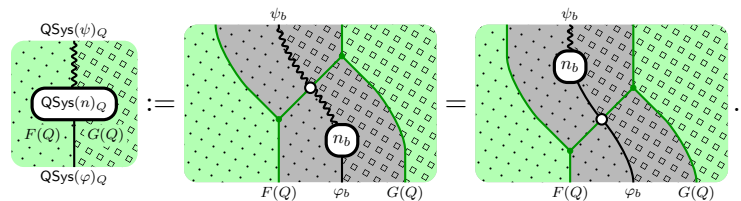
\[
\tikzmath[scale=.7, transform shape]{
\begin{scope}
\clip[rounded corners = 5] (-1.2,-1.2) rectangle (1.2,1.2);
\filldraw[primedregion=\QrColor] (-1.2,-1.2) rectangle (0,1.2);
\filldraw[boxregion=\QrColor] (1.2,1.2) rectangle (0,-1.2);
\filldraw[\QrColor] (-1,-.8) rectangle (-.2,-.4);
\filldraw[\QrColor] (.2,-.8) rectangle (1,-.4);
\end{scope}
\draw[\phiColor,thick] (0,-1.2) -- (0,0);
\draw[\psiColor,thick] (0,0) -- (0,1.2);
\roundNbox{unshaded}{(0,0)}{.3}{.5}{.5}{\scriptsize{$\QSys(n)_Q$}}; 
\node at (.6,-.6) {\scriptsize{$G(Q)$}};
\node at (-.6,-.6) {\scriptsize{$F(Q)$}};
\node at (0,-1.4) {\scriptsize{$\QSys(\varphi)_Q$}};
\node at (0,1.4) {\scriptsize{$\QSys(\psi)_Q$}};
}
:=
\tikzmath[scale=.7, transform shape]{
\begin{scope}
\clip[rounded corners = 5] (-2.4,-1.8) rectangle (2.4,1.2);
\filldraw[primedregion=\BColor] (.6,-1.8) -- (.6,-1.2) .. controls ++(90:.4cm) and ++(-45:.2cm) .. (0,0) .. controls ++(135:.2cm) and ++(270:.4cm) .. (-.6,1.2) -- (-2.4,1.2) -- (-2.4,-1.8);
\filldraw[boxregion=\BColor] (.6,-1.8) -- (.6,-1.2) .. controls ++(90:.4cm) and ++(-45:.2cm) .. (0,0) .. controls ++(135:.2cm) and ++(270:.4cm) .. (-.6,1.2) -- (2.4,1.2) -- (2.4,-1.8);
\filldraw[primedregion=\QrColor] (-.6,-1.8) -- (-.6,-.6) -- (-1.2,0) .. controls ++(135:.2cm) and ++(270:.4cm) .. (-1.8,1.2) -- (-2.4,1.2) -- (-2.4,-1.8);
\filldraw[boxregion=\QrColor] (1.8,-1.8) -- (1.8,-1.2) .. controls ++(90:.4cm) and ++(-45:.2cm) .. (1.2,0) -- (.6,.6) -- (.6,1.2) -- (2.4,1.2) -- (2.4,-1.8);
\end{scope}
\draw[\phiColor,thick] (.6,-1.8) -- (.6,-1.2);
\draw[\psiColor,thick] (.6,-1.2) .. controls ++(90:.4cm) and ++(-45:.2cm) .. (0,0) .. controls ++(135:.2cm) and ++(270:.4cm) .. (-.6,1.2);
\draw[\QsColor,thick] (-.6,-1.8) -- (-.6,-1.2) -- (-.6,-.6) -- (-1.2,0) .. controls ++(135:.2cm) and ++(270:.4cm) .. (-1.8,1.2);
\draw[\QsColor,thick] (1.8,-1.8) -- (1.8,-1.2) .. controls ++(90:.4cm) and ++(-45:.2cm) .. (1.2,0) -- (.6,.6) -- (.6,1.2);
\draw[\QsColor,thick] (-.6,-.6) -- (.6,.6);
\filldraw[\QsColor] (.6,.6) circle (.05cm);
\filldraw[\QsColor] (-.6,-.6) circle (.05cm);
\filldraw[white] (0,0) circle (.1cm);
\draw[thick] (0,0) circle (.1cm); 
\roundNbox{unshaded}{(.6,-1)}{.3}{0}{0}{\small{$n_b$}}; 
\node at (-.6,-2) {\scriptsize{$F(Q)$}};
\node at (.6,-2) {\scriptsize{$\varphi_b$}};
\node at (1.8,-2) {\scriptsize{$G(Q)$}};
\node at (-.6,1.4) {\scriptsize{$\psi_b$}};
}
=
\tikzmath[scale=.7, transform shape]{
\begin{scope}
\clip[rounded corners = 5] (-2.4,-1.2) rectangle (2.4,1.8);
\filldraw[primedregion=\BColor] (.6,-1.2) .. controls ++(90:.4cm) and ++(-45:.2cm) .. (0,0) .. controls ++(135:.2cm) and ++(270:.4cm) .. (-.6,1.8) -- (-2.4,1.8) -- (-2.4,-1.2);
\filldraw[boxregion=\BColor] (.6,-1.2) .. controls ++(90:.4cm) and ++(-45:.2cm) .. (0,0) .. controls ++(135:.2cm) and ++(270:.4cm) .. (-.6,1.8) -- (2.4,1.8) -- (2.4,-1.2);
\filldraw[primedregion=\QrColor] (-.6,-1.2) -- (-.6,-.6) -- (-1.2,0) .. controls ++(135:.2cm) and ++(270:.4cm) .. (-1.8,1.2) -- (-1.8,1.8) -- (-2.4,1.8) -- (-2.4,-1.2);
\filldraw[boxregion=\QrColor] (1.8,-1.2) .. controls ++(90:.4cm) and ++(-45:.2cm) .. (1.2,0) -- (.6,.6) -- (.6,1.8) -- (2.4,1.8) -- (2.4,-1.2);
\end{scope}
\draw[\phiColor,thick] (.6,-1.2) .. controls ++(90:.4cm) and ++(-45:.2cm) .. (0,0) .. controls ++(135:.2cm) and ++(270:.4cm) .. (-.6,1.2);
\draw[\psiColor,thick] (-.6,1.2) -- (-.6,1.8);
\draw[\QsColor,thick] (-.6,-1.2) -- (-.6,-.6) -- (-1.2,0) .. controls ++(135:.2cm) and ++(270:.4cm) .. (-1.8,1.2) -- (-1.8,1.8);
\draw[\QsColor,thick] (1.8,-1.2) .. controls ++(90:.4cm) and ++(-45:.2cm) .. (1.2,0) -- (.6,.6) -- (.6,1.2) -- (.6,1.8);
\draw[\QsColor,thick] (-.6,-.6) -- (.6,.6);
\filldraw[\QsColor] (.6,.6) circle (.05cm);
\filldraw[\QsColor] (-.6,-.6) circle (.05cm);
\filldraw[white] (0,0) circle (.1cm);
\draw[thick] (0,0) circle (.1cm); 
\roundNbox{unshaded}{(-.6,1)}{.3}{0}{0}{\small{$n_b$}}; 
\node at (-.6,-1.4) {\scriptsize{$F(Q)$}};
\node at (.6,-1.4) {\scriptsize{$\varphi_b$}};
\node at (1.8,-1.4) {\scriptsize{$G(Q)$}};
\node at (-.6,2) {\scriptsize{$\psi_b$}};
}\,.
\]
It is clear that $\QSys(n^\dag)=\QSys(n)^\dag$.
The modification coherence axiom is verified by
\[
\tikzmath[scale=.75, transform shape]{
\begin{scope}
\clip[rounded corners = 5] (-.6,0) rectangle (2.1,3);
\filldraw[primedregion=\PrColor] (0,0) -- (0,1.2) .. controls ++(90:.4cm) and ++(-135:.2cm) .. (.6,1.8) .. controls ++(135:.2cm) and ++(270:.4cm) .. (0,2.4) -- (0,3) -- (-.6,3) -- (-.6,0); 
\filldraw[primedregion=\QrColor] (1.2,0) -- (1.2,1.2) .. controls ++(90:.4cm) and ++(-45:.2cm) .. (.6,1.8) .. controls ++(-135:.2cm) and ++(90:.4cm) .. (0,1.2) -- (0,0);
\filldraw[boxregion=\PrColor] (0,3) -- (0,2.4) .. controls ++(270:.4cm) and ++(135:.2cm) .. (.6,1.8) .. controls ++(45:.2cm) and ++(270:.4cm) .. (1.2,2.4) -- (1.2,3);
\filldraw[boxregion=\QrColor] (1.2,0) -- (1.2,1.2) .. controls ++(90:.4cm) and ++(-45:.2cm) .. (.6,1.8) .. controls ++(45:.2cm) and ++(270:.4cm) .. (1.2,2.4) -- (1.2,3) -- (2.1,3) -- (2.1,0);
\end{scope}
\draw[\XColor,thick] (0,0) -- (0,1.2) .. controls ++(90:.6cm) and ++(270:.6cm) .. (1.2,2.4) -- (1.2,3);
\draw[\phiColor,thick] (1.2,0) -- (1.2,.8);
\draw[\psiColor,thick] (1.2,.8) -- (1.2,1.2) .. controls ++(90:.6cm) and ++(270:.6cm) .. (0,2.4) -- (0,3);
\roundNbox{unshaded}{(1.2,.8)}{.3}{.45}{.45}{\scriptsize{$\QSys(n)_Q$}}; 
\filldraw[white] (.6,1.8) circle (.1cm);
\draw[thick] (.6,1.8) circle (.1cm); 
\node at (0,-.2) {\scriptsize{$F(X)$}};
\node at (1.2,3.2) {\scriptsize{$G(X)$}};
\node at (1.2,-.2) {\scriptsize{$\QSys(\varphi)_Q$}};
\node at (0,3.2) {\scriptsize{$\QSys(\psi)_P$}};
%
}
=\tikzmath[scale=.7, transform shape]{
\begin{scope}
\clip[rounded corners = 5] (-3.6,-3) rectangle (3.6,2.4);
\filldraw[primedregion=\AColor] (0,0) -- (-1.2,1.2) .. controls ++(135:.2cm) and ++(270:.4cm) .. (-1.8,2.4) -- (-3,2.4) -- (-3,-.6) -- (-.6,-.6);
\filldraw[boxregion=\AColor] (0,0) -- (-1.2,1.2) .. controls ++(135:.2cm) and ++(270:.4cm) .. (-1.8,2.4) -- (-.6,2.4) -- (-.6,1.8) -- (.6,.6);
\filldraw[primedregion=\BColor] (1.8,-3) -- (1.8,-2.4) .. controls ++(90:.4cm) and ++(-45:.2cm) .. (1.2,-1.2) -- (0,0) -- (-.6,-.6) -- (.6,-1.8) -- (.6,-3);
\filldraw[boxregion=\BColor] (1.8,-3) -- (1.8,-2.4) .. controls ++(90:.4cm) and ++(-45:.2cm) .. (1.2,-1.2) -- (0,0) -- (.6,.6) -- (3,.6) -- (3,-3);
\filldraw[primedregion=\PrColor] (-.6,-3) -- (-.6,-.6) -- (-2.4,1.2) .. controls ++(135:.2cm) and ++(270:.4cm) .. (-3,2.4) -- (-3.6,2.4) -- (-3.6,-3);
\filldraw[boxregion=\PrColor] (-.6,2.4) -- (-.6,1.8) -- (.6,.6) -- (.6,2.4); 
\filldraw[primedregion=\QrColor] (.6,-3) -- (.6,-1.8) -- (-.6,-.6) -- (-.6,-3); 
\filldraw[boxregion=\QrColor] (3,-3) -- (3,-2.4) .. controls ++(90:.4cm) and ++(-45:.2cm) .. (2.4,-1.2) -- (.6,.6) -- (.6,2.4) -- (3.6,2.4) -- (3.6,-3);
\end{scope}
\draw[\phiColor,thick] (1.8,-3) -- (1.8,-2.4);
\draw[\psiColor,thick] (1.8,-2.4) .. controls ++(90:.4cm) and ++(-45:.2cm) .. (1.2,-1.2) -- (-1.2,1.2) .. controls ++(135:.2cm) and ++(270:.4cm) .. (-1.8,2.4);
\draw[\XColor,thick] (-.6,-3) -- (-.6,-.6) -- (.6,.6) -- (.6,2.4);
\draw[\PsColor,thick] (-.6,-.6) -- (-2.4,1.2) .. controls ++(135:.2cm) and ++(270:.4cm) .. (-3,2.4);
\draw[\PsColor,thick] (.6,.6) -- (-.6,1.8) -- (-.6,2.4);
\draw[\PsColor,thick] (-1.8,.6) -- (-.6,1.8);
\draw[\QsColor,thick] (-.6,-.6) -- (.6,-1.8) -- (.6,-3);
\draw[\QsColor,thick] (3,-3) -- (3,-2.4) .. controls ++(90:.4cm) and ++(-45:.2cm) .. (2.4,-1.2) -- (.6,.6);
\draw[\QsColor,thick] (1.8,-.6) -- (.6,-1.8);
\filldraw[\XColor] (-.6,-.6) circle (.05cm);
\filldraw[\XColor] (.6,.6) circle (.05cm);
\filldraw[\PsColor] (-1.8,.6) circle (.05cm);
\filldraw[\PsColor] (-.6,1.8) circle (.05cm);
\filldraw[\QsColor] (1.8,-.6) circle (.05cm);
\filldraw[\QsColor] (.6,-1.8) circle (.05cm);
\filldraw[white] (0,0) circle (.1cm);
\draw[thick] (0,0) circle (.1cm); 
\filldraw[white] (-1.2,1.2) circle (.1cm);
\draw[thick] (-1.2,1.2) circle (.1cm); 
\filldraw[white] (1.2,-1.2) circle (.1cm);
\draw[thick] (1.2,-1.2) circle (.1cm); 
\roundNbox{unshaded}{(1.8,-2.2)}{.3}{0}{0}{\small{$n_b$}}; 
\node at (-.6,-3.2) {\scriptsize{$F(X)$}};
\node at (.6,-3.2) {\scriptsize{$F(Q)$}};
\node at (1.8,-3.2) {\scriptsize{$\varphi_b$}};
\node at (3,-3.2) {\scriptsize{$G(Q)$}};
\node at (.6,2.6) {\scriptsize{$G(X)$}};
\node at (-.6,2.6) {\scriptsize{$G(P)$}};
\node at (-1.8,2.6) {\scriptsize{$\psi_a$}};
\node at (-3,2.6) {\scriptsize{$F(P)$}};
}
=
\tikzmath[scale=.7, transform shape]{
\begin{scope}
\clip[rounded corners = 5] (-3.6,-2.4) rectangle (3.6,3);
\filldraw[primedregion=\AColor] (0,0) -- (-1.2,1.2) .. controls ++(135:.2cm) and ++(270:.4cm) .. (-1.8,2.4) -- (-1.8,3) -- (-3,3) -- (-3,-.6) -- (-.6,-.6);
\filldraw[boxregion=\AColor] (0,0) -- (-1.2,1.2) .. controls ++(135:.2cm) and ++(270:.4cm) .. (-1.8,2.4) -- (-1.8,3) -- (-.6,3) -- (-.6,1.8) -- (.6,.6);
\filldraw[primedregion=\BColor] (1.8,-2.4) .. controls ++(90:.4cm) and ++(-45:.2cm) .. (1.2,-1.2) -- (0,0) -- (-.6,-.6) -- (.6,-1.8) -- (.6,-2.4);
\filldraw[boxregion=\BColor] (1.8,-2.4) .. controls ++(90:.4cm) and ++(-45:.2cm) .. (1.2,-1.2) -- (0,0) -- (.6,.6) -- (3,.6) -- (3,-2.4);
\filldraw[primedregion=\PrColor] (-.6,-2.4) -- (-.6,-.6) -- (-2.4,1.2) .. controls ++(135:.2cm) and ++(270:.4cm) .. (-3,2.4) -- (-3,3) -- (-3.6,3) -- (-3.6,-2.4);
\filldraw[boxregion=\PrColor] (-.6,3) -- (-.6,1.8) -- (.6,.6) -- (.6,3); 
\filldraw[primedregion=\QrColor] (.6,-2.4) -- (.6,-1.8) -- (-.6,-.6) -- (-.6,-2.4); 
\filldraw[boxregion=\QrColor] (3,-2.4) .. controls ++(90:.4cm) and ++(-45:.2cm) .. (2.4,-1.2) -- (.6,.6) -- (.6,3) -- (3.6,3) -- (3.6,-2.4);
\end{scope}
\draw[\phiColor,thick] (1.8,-2.4) .. controls ++(90:.4cm) and ++(-45:.2cm) .. (1.2,-1.2) -- (-1.2,1.2) .. controls ++(135:.2cm) and ++(270:.4cm) .. (-1.8,2.4);
\draw[\psiColor,thick] (-1.8,2.4) -- (-1.8,3);
\draw[\XColor,thick] (-.6,-2.4) -- (-.6,-.6) -- (.6,.6) -- (.6,3);
\draw[\PsColor,thick] (-.6,-.6) -- (-2.4,1.2) .. controls ++(135:.2cm) and ++(270:.4cm) .. (-3,2.4) -- (-3,3);
\draw[\PsColor,thick] (.6,.6) -- (-.6,1.8) -- (-.6,2.4) -- (-.6,3);
\draw[\PsColor,thick] (-1.8,.6) -- (-.6,1.8);
\draw[\QsColor,thick] (-.6,-.6) -- (.6,-1.8) -- (.6,-2.4);
\draw[\QsColor,thick] (3,-2.4) .. controls ++(90:.4cm) and ++(-45:.2cm) .. (2.4,-1.2) -- (.6,.6);
\draw[\QsColor,thick] (1.8,-.6) -- (.6,-1.8);
\filldraw[\XColor] (-.6,-.6) circle (.05cm);
\filldraw[\XColor] (.6,.6) circle (.05cm);
\filldraw[\PsColor] (-1.8,.6) circle (.05cm);
\filldraw[\PsColor] (-.6,1.8) circle (.05cm);
\filldraw[\QsColor] (1.8,-.6) circle (.05cm);
\filldraw[\QsColor] (.6,-1.8) circle (.05cm);
\filldraw[white] (0,0) circle (.1cm);
\draw[thick] (0,0) circle (.1cm); 
\filldraw[white] (-1.2,1.2) circle (.1cm);
\draw[thick] (-1.2,1.2) circle (.1cm); 
\filldraw[white] (1.2,-1.2) circle (.1cm);
\draw[thick] (1.2,-1.2) circle (.1cm); 
\roundNbox{unshaded}{(-1.8,2.2)}{.3}{0}{0}{\small{$n_a$}}; 
\node at (-.6,-2.6) {\scriptsize{$F(X)$}};
\node at (.6,-2.6) {\scriptsize{$F(Q)$}};
\node at (1.8,-2.6) {\scriptsize{$\varphi_b$}};
\node at (3,-2.6) {\scriptsize{$G(Q)$}};
\node at (.6,3.2) {\scriptsize{$G(X)$}};
\node at (-.6,3.2) {\scriptsize{$G(P)$}};
\node at (-1.8,3.2) {\scriptsize{$\psi_a$}};
\node at (-3,3.2) {\scriptsize{$F(P)$}};
}
=
\tikzmath[scale=.75, transform shape]{
\begin{scope}
\clip[rounded corners = 5] (-.9,0) rectangle (1.8,3);
\filldraw[primedregion=\PrColor] (0,0) -- (0,.6) .. controls ++(90:.4cm) and ++(-135:.2cm) .. (.6,1.2) .. controls ++(135:.2cm) and ++(270:.4cm) .. (0,1.8) -- (0,3) -- (-.9,3) -- (-.9,0); 
\filldraw[primedregion=\QrColor] (1.2,0) -- (1.2,.6) .. controls ++(90:.4cm) and ++(-45:.2cm) .. (.6,1.2) .. controls ++(-135:.2cm) and ++(90:.4cm) .. (0,.6) -- (0,0);
\filldraw[boxregion=\PrColor] (0,3) -- (0,1.8) .. controls ++(270:.4cm) and ++(135:.2cm) .. (.6,1.2) .. controls ++(45:.2cm) and ++(270:.4cm) .. (1.2,1.8) -- (1.2,3);
\filldraw[boxregion=\QrColor] (1.2,0) -- (1.2,.6) .. controls ++(90:.4cm) and ++(-45:.2cm) .. (.6,1.2) .. controls ++(45:.2cm) and ++(270:.4cm) .. (1.2,1.8) -- (1.2,3) -- (1.8,3) -- (1.8,0);
\end{scope}
\draw[\XColor,thick] (0,0) -- (0,.6) .. controls ++(90:.6cm) and ++(270:.6cm) .. (1.2,1.8) -- (1.2,3);
\draw[\phiColor,thick] (1.2,0) -- (1.2,.6) .. controls ++(90:.6cm) and ++(270:.6cm) .. (0,1.8) -- (0,2.2);
\draw[\psiColor,thick] (0,2.2) -- (0,3);
\roundNbox{unshaded}{(0,2.2)}{.3}{.45}{.45}{\scriptsize{$\QSys(n)_P$}}; 
\filldraw[white] (.6,1.2) circle (.1cm);
\draw[thick] (.6,1.2) circle (.1cm); 
\node at (0,-.2) {\scriptsize{$F(X)$}};
\node at (1.2,3.2) {\scriptsize{$G(X)$}};
\node at (1.2,-.2) {\scriptsize{$\QSys(\varphi)_Q$}};
\node at (0,3.2) {\scriptsize{$\QSys(\psi)_P$}};
}\,.
\]
By our construction, it is clear that when $n:\varphi\Rrightarrow \psi$ is invertible, $\QSys(n):\QSys(\varphi)\Rrightarrow\QSys(\psi)$ is also invertible.
\end{construction}

\begin{construction}
\label{construction:QSys-otimes}
Given $\cA,\cB\in 2\Cat$, $F,G,H:\cA\to\cB$, and $\varphi:F\Rightarrow G$, $\psi:G\Rightarrow H$, we construct $\QSys^{\xo}_{\varphi,\psi}:\QSys(\varphi)\xo\QSys(\psi)\Rrightarrow\QSys(\varphi\xo\psi)$
by
\[
\tikzmath{
\begin{scope}
\clip[rounded corners=5pt] (-.7,.3) rectangle (.7,1.7);
\filldraw[primedregion=white] (-.7,0) rectangle (.7,2);
\filldraw[boxregion=white] (.2,0) -- (.2,1) -- (.05,1) -- (.05,2) -- (-.05,2) -- (-.05,1) -- (-.2,1) -- (-.2,0);
\filldraw[plusregion=white] (.2,0) -- (.2,1) -- (.05,1) -- (.05,2) -- (.7,2) -- (.7,0);
\end{scope}
\draw[\phiColor,thick] (-.2,.3) -- (-.2,1);
\draw[\psiColor,thick] (.2,.3) -- (.2,1);
\draw[\phiColor,thick] (-.05,1) -- (-.05,1.7);
\draw[\psiColor,thick] (.05,1) -- (.05,1.7);
\roundNbox{unshaded}{(0,1)}{.3}{.25}{.25}{\tiny{$\QSys^\xo_{\varphi,\psi}$}};
\draw[thin, dotted, rounded corners = 5pt] (-.7,.3) rectangle (.7,1.7);
}
:
\qquad
\tikzmath[scale=.75, transform shape]{
\begin{scope}
\clip[rounded corners=5pt] (-1.2,0) rectangle (1.2,2);
\filldraw[primedregion=\QrColor] (-1.2,0) rectangle (1,2);
\filldraw[boxregion=\QrColor] (.6,0) -- (.6,1) -- (.05,1) -- (.05,2) -- (-.05,2) -- (-.05,1) -- (-.6,1) -- (-.6,0);
\filldraw[plusregion=\QrColor] (.6,0) -- (.6,1) -- (.05,1) -- (.05,2) -- (1.2,2) -- (1.2,0);
\end{scope}
\draw[\phiColor,thick] (-.6,0) -- (-.6,1);
\draw[\psiColor,thick] (.6,0) -- (.6,1);
\draw[\phiColor,thick] (-.05,1) -- (-.05,2);
\draw[\psiColor,thick] (.05,1) -- (.05,2);
\roundNbox{unshaded}{(0,1)}{.33}{.6}{.6}{\scriptsize{$\left(\QSys^\xo_{\varphi,\psi}\right)_Q$}};
}
:=
\tikzmath[scale=.5, transform shape]{
\begin{scope}
\clip[rounded corners = 5] (-4.8,-1.2) rectangle (2.4,1.2);
\filldraw[primedregion=\QrColor] (-4.8,-1.2) rectangle (2.4,1.2);
\filldraw[primedregion=\BColor] (-4.2,-1.2) .. controls ++(90:.4cm) and ++(-135:.2cm) .. (-3.6,0) -- (-3,.6) -- (-3,1.2) -- (-1.2,1.2) -- (-1.2,-1.2);
\filldraw[boxregion=\BColor] ((-3,-1.2) .. controls ++(90:.4cm) and ++(-135:.2cm) ..  (-2.4,0) .. controls ++(45:.2cm) and ++(270:.4cm) .. (-1.8,1.2) -- (.6,1.2) -- (.6,-1.2);
\filldraw[boxregion=\QrColor] (-1.8,-.6) rectangle (-.6,-1.2);
\filldraw[plusregion=\BColor] (.6,-1.2) .. controls ++(90:.4cm) and ++(-45:.2cm) .. (0,0) .. controls ++(135:.2cm) and ++(270:.4cm) .. (-.6,1.2) -- (2.4,1.2) -- (2.4,-1.2);
\filldraw[plusregion=\QrColor] (1.8,-1.2) .. controls ++(90:.4cm) and ++(-45:.2cm) .. (1.2,0) -- (.6,.6) -- (.6,1.2) -- (2.4,1.2) -- (2.4,-1.2);
\end{scope}
\draw[\phiColor,thick] (-3,-1.2) .. controls ++(90:.4cm) and ++(-135:.2cm) ..  (-2.4,0) .. controls ++(45:.2cm) and ++(270:.4cm) .. (-1.8,1.2);
\draw[\psiColor,thick] (.6,-1.2) .. controls ++(90:.4cm) and ++(45:.2cm) .. (0,0) .. controls ++(135:.2cm) and ++(270:.4cm) .. (-.6,1.2);
\draw[\QsColor,thick] (-4.2,-1.2) .. controls ++(90:.4cm) and ++(-135:.2cm) .. (-3.6,0) -- (-3,.6);
\draw[\QsColor,thick] (1.8,-1.2) .. controls ++(90:.4cm) and ++(-45:.2cm) .. (1.2,0) -- (.6,.6);
\draw[\QsColor,thick] (-.6,-1.2) -- (-.6,-.6) -- (.6,.6) -- (.6,1.2);
\draw[\QsColor,thick] (-1.8,-1.2) -- (-1.8,-.6) -- (-3,.6) -- (-3,1.2);
\draw[\QsColor,thick] (-1.8,-.6) -- (-.6,-.6);
\filldraw[\QsColor] (.6,.6) circle (.07cm);
\filldraw[\QsColor] (-.6,-.6) circle (.07cm);
\filldraw[\QsColor] (-3,.6) circle (.07cm);
\filldraw[\QsColor] (-1.8,-.6) circle (.07cm);
\filldraw[white] (0,0) circle (.1cm);
\draw[thick] (0,0) circle (.1cm); 
\filldraw[white] (-2.4,0) circle (.1cm);
\draw[thick] (-2.4,0) circle (.1cm); 
\node at (-4.2,-1.4) {\small{$F(Q)$}};
\node at (-3,-1.4) {\small{$\varphi_b$}};
\node at (-.6,-1.4) {\small{$G(Q)$}};
\node at (.6,-1.4) {\small{$\psi_b$}};
\node at (1.8,-1.4) {\small{$H(Q)$}};
}
=
\tikzmath[scale=.7, transform shape]{
\begin{scope}
\clip[rounded corners = 5] (-2.1,-1) rectangle (2.1,1);
\filldraw[primedregion=\QrColor] (-2.1,-1.2) rectangle (-1.5,1.2);
\filldraw[primedregion=\BColor] (-1.5,-1.2) rectangle (-.9,1.2);
\filldraw[boxregion=\BColor] (-.9,-1.2) rectangle (.9,1.2);
\filldraw[boxregion=\QrColor] (-.3,-1.2) rectangle (.3,0);
\filldraw[plusregion=\BColor] (.9,-1.2) rectangle (1.5,1.2);
\filldraw[plusregion=\QrColor] (1.5,-1.2) rectangle (2.1,1.2);
\end{scope}
\draw[\phiColor,thick] (-.9,-1) -- (-.9,1);
\draw[\psiColor,thick] (.9,-1) -- (.9,1);
\draw[\QsColor,thick] (-1.5,-1) -- (-1.5,1);
\draw[\QsColor,thick] (1.5,-1) -- (1.5,1);
\draw[\QsColor,thick] (.3,-1) -- (.3,0);
\draw[\QsColor,thick] (-.3,-1) -- (-.3,0);
\draw[\QsColor,thick] (-1.5,0) -- (1.5,0);
\filldraw[white] (-.9,0) circle (.08cm);
\draw[thick] (-.9,0) circle (.08cm); 
\filldraw[white] (.9,0) circle (.08cm);
\draw[thick] (.9,0) circle (.08cm); 
}\,.
\]
It is straightforward to verify $(\QSys^\xo_{\varphi,\psi})_Q$ is unitary. 
The following calculation shows $\QSys^\xo_{\varphi,\psi}$ is a modification:
\[
\tikzmath[scale=.7, transform shape]{
\begin{scope}
\clip[rounded corners = 5] (-1.2,-1.9) rectangle (1.8,1);
\filldraw[primedregion=\PrColor] (-.6,-1.9) -- (-.6,-.6) .. controls ++(90:.4cm) and ++(-135:.2cm) .. (0,0) -- (.05,0) .. controls ++(135:.2cm) and ++(270:.4cm) .. (-.5,.6) -- (-.5,1) -- (-1.2,1) -- (-1.2,-1.9);
\filldraw[primedregion=\QrColor] (-.6,-1.9) -- (-.6,-.6) .. controls ++(90:.4cm) and ++(-135:.2cm) .. (0,0) -- (.05,0) .. controls ++(-45:.2cm) and ++(90:.4cm) .. (.6,-.6) -- (.6,-1) -- (-.05,-1) -- (-.05,-1.9);
\filldraw[boxregion=\QrColor] (-.05,-1.9) rectangle (1.15,-1);
\filldraw[plusregion=\PrColor] (-.5,1) -- (-.5,.6) .. controls ++(270:.4cm) and ++(135:.2cm) .. (.05,0) -- (0,0) .. controls ++(45:.2cm) and ++(270:.4cm) .. (.6,.6) -- (.6,1);
\filldraw[plusregion=\QrColor] (1.15,-1.9) -- (1.15,-1) -- (.6,-1) -- (.6,-.6) .. controls ++(90:.4cm) and ++(-45:.2cm) .. (.05,0) -- (0,0) .. controls ++(45:.2cm) and ++(270:.4cm) .. (.6,.6) -- (.6,1) -- (1.8,1) -- (1.8,-1.9);
\end{scope}
\draw[\XColor,thick] (-.6,-1.9) -- (-.6,-.6) .. controls ++(90:.6cm) and ++(270:.6cm) .. (.6,.6) -- (.6,1);
\draw[\phiColor,thick] (.5,-1) -- (.5,-.6) .. controls ++(90:.6cm) and ++(270:.6cm) .. (-.6,.6) -- (-.6,1);
\draw[\psiColor,thick] (.6,-1) -- (.6,-.6) .. controls ++(90:.6cm) and ++(270:.6cm) .. (-.5,.6) -- (-.5,1);
\draw[\phiColor,thick] (-.05,-1) -- (-.05,-1.9);
\draw[\psiColor,thick] (1.15,-1) -- (1.15,-1.9);
\roundNbox{unshaded}{(.55,-1)}{.33}{.6}{.6}{\scriptsize{$\left(\QSys^\xo_{\varphi,\psi}\right)_Q$}};
\filldraw[white] (0,0) circle (.1cm);
\draw[thick] (0,0) circle (.1cm); 
}
=
\tikzmath[scale=.45, transform shape]{
\begin{scope}
\clip[rounded corners = 5] (-4.5,-5.1) rectangle (5.7,2.7);
\filldraw[plusregion=\QrColor] (5.1,-5.1) .. controls ++(90:.4cm) and ++(-45:.2cm) .. (4.5,-3.9) -- (3.9,-3.3) -- (3.9,-2.7) .. controls ++(90:.4cm) and ++(-45:.2cm) .. (3.3,-1.5) -- (.9,.9) -- (.9,2.7) -- (5.7,2.7) -- (5.7,-5.1);
\filldraw[plusregion=\BColor] (5.1,-5.1) .. controls ++(90:.4cm) and ++(-45:.2cm) .. (4.5,-3.9) -- (3.9,-3.3) -- (3.9,-2.7) .. controls ++(90:.4cm) and ++(-45:.2cm) .. (3.3,-1.5) -- (.9,.9) -- (-.9,-.9) -- (-.9,-5.1);
\filldraw[boxregion=\BColor] (3.9,-5.1) .. controls ++(90:.4cm) and ++(45:.2cm) .. (3.3,-3.9) .. controls ++(135:.2cm) and ++(270:.4cm) .. (2.7,-2.7) .. controls ++(90:.4cm) and ++(-45:.2cm) .. (2.1,-1.5) -- (.3,.3) -- (-.9,-.9) -- (-.9,-5.1);
\filldraw[boxregion=\QrColor] (1.5,-5.1) rectangle (2.7,-4.5);
\filldraw[primedregion=\BColor] (.3,-5.1) .. controls ++(90:.4cm) and ++(-135:.2cm) ..  (.9,-3.9) .. controls ++(45:.2cm) and ++(270:.4cm) .. (1.5,-2.7) .. controls ++(90:.4cm) and ++(-45:.2cm) .. (.9,-1.5) -- (-.3,-.3) -- (-.9,-.9) -- (-.9,-5.1);
\filldraw[primedregion=\QrColor] (-.9,-5.1) .. controls ++(90:.4cm) and ++(-135:.2cm) .. (-.3,-3.9) -- (.3,-3.3) -- (.3,-2.1) -- (-.9,-.9) -- (-2.1,-.9) -- (-2.1,-5.1);
\filldraw[primedregion=\PrColor] (-2.1,-5.1) -- (-2.1,-2.7)  .. controls ++(90:.4cm) and ++(-135:.2cm) .. (-1.5,-1.5) -- (-.9,-.9) -- (-3.3,1.5) .. controls ++(135:.2cm) and ++(270:.4cm) .. (-3.9,2.7) -- (-4.5,2.7) -- (-4.5,-5.1);
\filldraw[primedregion=\AColor] (-.9,-.9) -- (-3.3,1.5) .. controls ++(135:.2cm) and ++(270:.4cm) .. (-3.9,2.7) -- (.9,2.7) -- (.9,.9);
\filldraw[boxregion=\AColor] (-.3,-.3) -- (-2.1,1.5) .. controls ++(135:.2cm) and ++(270:.4cm) .. (-2.7,2.7) -- (.9,2.7) -- (.9,.9);
\filldraw[plusregion=\AColor] (.3,.3) -- (-.9,1.5) .. controls ++(135:.2cm) and ++(270:.4cm) .. (-1.5,2.7) -- (.9,2.7) -- (.9,.9);
\filldraw[plusregion=\PrColor] (.9,.9) -- (-.3,2.1) -- (-.3,2.7) -- (.9,2.7) -- (.9,.9);
\end{scope}
\draw[\phiColor,thick] (1.5,-2.7) .. controls ++(90:.4cm) and ++(-45:.2cm) .. (.9,-1.5) -- (-2.1,1.5) .. controls ++(135:.2cm) and ++(270:.4cm) .. (-2.7,2.7);
\draw[\psiColor,thick] (2.7,-2.7) .. controls ++(90:.4cm) and ++(-45:.2cm) .. (2.1,-1.5) -- (-.9,1.5) .. controls ++(135:.2cm) and ++(270:.4cm) .. (-1.5,2.7);
\draw[\XColor,thick] (-2.1,-5.1) -- (-2.1,-2.7)  .. controls ++(90:.4cm) and ++(-135:.2cm) .. (-1.5,-1.5) -- (-.9,-.9) -- (.9,.9) -- (.9,2.7);
\draw[\PsColor,thick] (-2.1,.3) -- (-.3,2.1);
\draw[\PsColor,thick] (.9,.9) -- (-.3,2.1) -- (-.3,2.7);
\draw[\PsColor,thick] (-.9,-.9) -- (-3.3,1.5) .. controls ++(135:.2cm) and ++(270:.4cm) .. (-3.9,2.7);
\draw[\QsColor,thick] (.3,-2.1) -- (2.1,-.3);
\draw[\QsColor,thick] (-.9,-.9) -- (.3,-2.1) -- (.3,-2.7);
\draw[\QsColor,thick] (3.9,-2.7) .. controls ++(90:.4cm) and ++(-45:.2cm) .. (3.3,-1.5) -- (.9,.9);
\filldraw[white] (-.3,-.3) circle (.1cm);
\draw[thick] (-.3,-.3) circle (.1cm); 
\filldraw[white] (.3,.3) circle (.1cm);
\draw[thick] (.3,.3) circle (.1cm); 
\filldraw[white] (-1.5,.9) circle (.1cm);
\draw[thick] (-1.5,.9) circle (.1cm); 
\filldraw[white] (.9,-1.5) circle (.1cm);
\draw[thick] (.9,-1.5) circle (.1cm); 
\filldraw[white] (1.5,-.9) circle (.1cm);
\draw[thick] (1.5,-.9) circle (.1cm); 
\filldraw[white] (-.9,1.5) circle (.1cm);
\draw[thick] (-.9,1.5) circle (.1cm); 
%
\draw[dashed,rounded corners = 5] (-1.2,-5.1) rectangle (5.4,-2.7);
\draw[\phiColor,thick] (.3,-5.1) .. controls ++(90:.4cm) and ++(-135:.2cm) ..  (.9,-3.9) .. controls ++(45:.2cm) and ++(270:.4cm) .. (1.5,-2.7);
\draw[\psiColor,thick] (3.9,-5.1) .. controls ++(90:.4cm) and ++(45:.2cm) .. (3.3,-3.9) .. controls ++(135:.2cm) and ++(270:.4cm) .. (2.7,-2.7) -- (2.7,-2.6);
\draw[\QsColor,thick] (-.9,-5.1) .. controls ++(90:.4cm) and ++(-135:.2cm) .. (-.3,-3.9) -- (.3,-3.3) -- (.3,-2.7);
\draw[\QsColor,thick] (5.1,-5.1) .. controls ++(90:.4cm) and ++(-45:.2cm) .. (4.5,-3.9) -- (3.9,-3.3) -- (3.9,-2.7);
\draw[\QsColor,thick] (.3,-3.3) -- (1.5,-4.5) -- (1.5,-5.1);
\draw[\QsColor,thick] (3.9,-3.3) -- (2.7,-4.5) -- (2.7,-5.1);
\draw[\QsColor,thick] (1.5,-4.5) -- (2.7,-4.5);
\filldraw[white] (.9,-3.9) circle (.1cm);
\draw[thick] (.9,-3.9) circle (.1cm); 
\filldraw[white] (3.3,-3.9) circle (.1cm);
\draw[thick] (3.3,-3.9) circle (.1cm); 
\filldraw[\XColor] (-.9,-.9) circle (.07cm);
\filldraw[\XColor] (.9,.9) circle (.07cm);
\filldraw[\PsColor] (-2.1,.3) circle (.07cm);
\filldraw[\PsColor] (-.3,2.1) circle (.07cm);
\filldraw[\QsColor] (.3,-2.1) circle (.07cm);
\filldraw[\QsColor] (2.1,-.3) circle (.07cm);
\filldraw[\QsColor] (.3,-3.3) circle (.07cm);
\filldraw[\QsColor] (1.5,-4.5) circle (.07cm);
\filldraw[\QsColor] (2.7,-4.5) circle (.07cm);
\filldraw[\QsColor] (3.9,-3.3) circle (.07cm);
}
=
\tikzmath[scale=.4, transform shape]{
\begin{scope}
\clip[rounded corners = 5] (-6.3,-3.3) rectangle (6.3,5.7);
\filldraw[primedregion=\AColor] (-1.5,-1.5) -- (-.9,-.9) --  (-3.9,2.1) .. controls ++(135:.2cm) and ++(270:.4cm) .. (-4.5,3.3) .. controls ++(90:.4cm) and ++(-135:.2cm) ..  (-3.9,4.5) .. controls ++(45:.2cm) and ++(270:.4cm) .. (-3.3,5.7) -- (-5.7,5.7) -- (-5.7,-1.5);
\filldraw[primedregion=\PrColor] (-1.5,-3.3) -- (-1.5,-1.5) -- (-5.1,2.1) .. controls ++(135:.2cm) and ++(270:.4cm) .. (-5.7,3.3) .. controls ++(90:.4cm) and ++(-135:.2cm) .. (-5.1,4.5) -- (-4.5,5.1) -- (-4.5,5.7) -- (-6.3,5.7) -- (-6.3,-3.3);
\filldraw[boxregion=\AColor] (-.9,-.9) -- (-3.9,2.1) .. controls ++(135:.2cm) and ++(270:.4cm) .. (-4.5,3.3).. controls ++(90:.4cm) and ++(-135:.2cm) ..  (-3.9,4.5) .. controls ++(45:.2cm) and ++(270:.4cm) .. (-3.3,5.7) -- (-.9,5.7) -- (1.5,5.7) -- (1.5,1.5);
\filldraw[boxregion=\PrColor] (-.3,-.3) -- (-2.7,2.1) .. controls ++(135:.2cm) and ++(270:.4cm) .. (-3.3,3.3) -- (-3.3,3.9) -- (-2.1,3.9) -- (.3,.3);
\filldraw[boxregion=\AColor] (.3,.3) -- (-1.5,2.1) .. controls ++(135:.2cm) and ++(270:.4cm) .. (-2.1,3.3) -- (-2.1,3.9) -- (-1.5,4.5) -- (1.5,1.5);
\filldraw[plusregion=\AColor] (.9,.9) -- (-.3,2.1) .. controls ++(135:.2cm) and ++(270:.4cm) .. (-.9,3.3)  .. controls ++(90:.4cm) and ++(45:.2cm) .. (-1.5,4.5) .. controls ++(135:.2cm) and ++(270:.4cm) .. (-2.1,5.7) -- (-.9,5.7) -- (1.5,5.7) -- (1.5,1.5);
\filldraw[plusregion=\PrColor] (.3,2.7) -- (.3,3.3) .. controls ++(90:.4cm) and ++(-45:.2cm) .. (-.3,4.5) -- (-.9,5.1) -- (-.9,5.7) -- (1.5,5.7) -- (1.5,1.5);
\filldraw[boxregion=\BColor] (4.5,-3.3) .. controls ++(90:.4cm) and ++(-45:.2cm) .. (3.9,-2.1) -- (.9,.9) -- (-1.5,-1.5) -- (-1.5,-3.3);
\filldraw[boxregion=\QrColor] (3.3,-3.3) .. controls ++(90:.4cm) and ++(-45:.2cm) .. (2.7,-2.1) -- (.3,.3) -- (-1.5,-1.5) -- (-1.5,-3.3);
\filldraw[boxregion=\BColor] (2.1,-3.3) .. controls ++(90:.4cm) and ++(-45:.2cm) .. (1.5,-2.1) -- (-.3,-.3) -- (-1.5,-1.5) -- (-1.5,-3.3);
\filldraw[primedregion=\BColor] (.9,-3.3) .. controls ++(90:.4cm) and ++(-45:.2cm) .. (.3,-2.1) -- (-.9,-.9) -- (-1.5,-1.5) -- (-1.5,-3.3);
\filldraw[primedregion=\QrColor] (-1.5,-3.3) -- (-1.5,-1.5) -- (-.3,-2.7) -- (-.3,-3.3);
\filldraw[plusregion=\BColor] (4.5,-3.3) .. controls ++(90:.4cm) and ++(-45:.2cm) .. (3.9,-2.1) -- (.9,.9) -- (1.5,1.5) -- (5.7,1.5) -- (5.7,-3.3);
\filldraw[plusregion=\QrColor] (5.7,-3.3) .. controls ++(90:.4cm) and ++(-45:.2cm) .. (5.1,-2.1) -- (1.5,1.5) -- (1.5,5.7) -- (6.3,5.7) -- (6.3,-3.3);
\end{scope}
\draw[\phiColor,thick] (.9,-3.3) .. controls ++(90:.4cm) and ++(-45:.2cm) .. (.3,-2.1) -- (-3.9,2.1) .. controls ++(135:.2cm) and ++(270:.4cm) .. (-4.5,3.3);
\draw[\psiColor,thick] (4.5,-3.3) .. controls ++(90:.4cm) and ++(-45:.2cm) .. (3.9,-2.1) -- (-.3,2.1) .. controls ++(135:.2cm) and ++(270:.4cm) .. (-.9,3.3);
\draw[\XColor,thick] (-1.5,-3.3) -- (-1.5,-1.5) -- (1.5,1.5) -- (1.5,5.7);
\draw[\PsColor,thick] (-2.7,-.3) -- (-1.5,.9);
\draw[\PsColor,thick] (-.9,1.5) -- (.3,2.7);
\draw[\PsColor,thick] (1.5,1.5) -- (.3,2.7) -- (.3,3.3);
\draw[\PsColor,thick] (-1.5,-1.5) -- (-5.1,2.1) .. controls ++(135:.2cm) and ++(270:.4cm) .. (-5.7,3.3);
\draw[\PsColor,thick] (-.3,-.3) -- (-2.7,2.1) .. controls ++(135:.2cm) and ++(270:.4cm) .. (-3.3,3.3);
\draw[\PsColor,thick] (.3,.3) -- (-1.5,2.1) .. controls ++(135:.2cm) and ++(270:.4cm) .. (-2.1,3.3);
\draw[\QsColor,thick] (-.3,-2.7) -- (.9,-1.5);
\draw[\QsColor,thick] (1.5,-.9) -- (2.7,.3);
\draw[\QsColor,thick] (-1.5,-1.5) -- (-.3,-2.7) -- (-.3,-3.3);
\draw[\QsColor,thick] (2.1,-3.3) .. controls ++(90:.4cm) and ++(-45:.2cm) .. (1.5,-2.1) -- (-.3,-.3);
\draw[\QsColor,thick] (3.3,-3.3) .. controls ++(90:.4cm) and ++(-45:.2cm) .. (2.7,-2.1) -- (.3,.3);
\draw[\QsColor,thick] (5.7,-3.3) .. controls ++(90:.4cm) and ++(-45:.2cm) .. (5.1,-2.1) -- (1.5,1.5);
\filldraw[white] (-.9,-.9) circle (.1cm);
\draw[thick] (-.9,-.9) circle (.1cm); 
\filldraw[white] (.9,.9) circle (.1cm);
\draw[thick] (.9,.9) circle (.1cm); 
\filldraw[white] (-2.1,.3) circle (.1cm);
\draw[thick] (-2.1,.3) circle (.1cm); 
\filldraw[white] (.3,-2.1) circle (.1cm);
\draw[thick] (.3,-2.1) circle (.1cm); 
\filldraw[white] (2.1,-.3) circle (.1cm);
\draw[thick] (2.1,-.3) circle (.1cm); 
\filldraw[white] (-.3,2.1) circle (.1cm);
\draw[thick] (-.3,2.1) circle (.1cm); 
%
\draw[dashed,rounded corners = 5] (-6,3.3) rectangle (.6,5.7);
\draw[\phiColor,thick] (-4.5,3.3) .. controls ++(90:.4cm) and ++(-135:.2cm) ..  (-3.9,4.5) .. controls ++(45:.2cm) and ++(270:.4cm) .. (-3.3,5.7);
\draw[\psiColor,thick] (-.9,3.3) .. controls ++(90:.4cm) and ++(45:.2cm) .. (-1.5,4.5) .. controls ++(135:.2cm) and ++(270:.4cm) .. (-2.1,5.7);
\draw[\PsColor,thick] (-5.7,3.3) .. controls ++(90:.4cm) and ++(-135:.2cm) .. (-5.1,4.5) -- (-4.5,5.1);
\draw[\PsColor,thick] (.3,3.3) .. controls ++(90:.4cm) and ++(-45:.2cm) .. (-.3,4.5) -- (-.9,5.1);
\draw[\PsColor,thick] (-2.1,3.3) -- (-2.1,3.9) -- (-.9,5.1) -- (-.9,5.7);
\draw[\PsColor,thick] (-3.3,3.3) -- (-3.3,3.9) -- (-4.5,5.1) -- (-4.5,5.7);
\draw[\PsColor,thick] (-2.1,3.9) -- (-3.3,3.9);
\filldraw[white] (-1.5,4.5) circle (.1cm);
\draw[thick] (-1.5,4.5) circle (.1cm); 
\filldraw[white] (-3.9,4.5) circle (.1cm);
\draw[thick] (-3.9,4.5) circle (.1cm); 
\filldraw[\XColor] (-1.5,-1.5) circle (.07cm);
\filldraw[\XColor] (-.3,-.3) circle (.07cm);
\filldraw[\XColor] (.3,.3) circle (.07cm);
\filldraw[\XColor] (1.5,1.5) circle (.07cm);
\filldraw[\PsColor] (-2.7,-.3) circle (.07cm);
\filldraw[\PsColor] (-1.5,.9) circle (.07cm);
\filldraw[\PsColor] (-.9,1.5) circle (.07cm);
\filldraw[\PsColor] (.3,2.7) circle (.07cm);
\filldraw[\PsColor] (-4.5,5.1) circle (.07cm);
\filldraw[\PsColor] (-3.3,3.9) circle (.07cm);
\filldraw[\PsColor] (-2.1,3.9) circle (.07cm);
\filldraw[\PsColor] (-.9,5.1) circle (.07cm);
\filldraw[\QsColor] (-.3,-2.7) circle (.07cm);
\filldraw[\QsColor] (.9,-1.5) circle (.07cm);
\filldraw[\QsColor] (1.5,-.9) circle (.07cm);
\filldraw[\QsColor] (2.7,.3) circle (.07cm);
}
=
\tikzmath[scale=.7, transform shape]{
\begin{scope}
\clip[rounded corners = 5] (-1.8,-1.8) rectangle (1.8,2.7);
\filldraw[primedregion=\PrColor] (-1.2,-1.8) -- (-1.2,-1.4) .. controls ++(90:.4cm) and ++(-135:.2cm) .. (-.6,-.8) .. controls ++(135:.2cm) and ++(270:.4cm) .. (-1.2,-.2) -- (-1.2,1.8) -- (-.55,1.8) -- (-.55,2.7) -- (-1.8,2.7) -- (-1.8,-1.8);
\filldraw[primedregion=\QrColor](0,-1.8) -- (0,-1.4) .. controls ++(90:.4cm) and ++(-45:.2cm) .. (-.6,-.8) .. controls ++(-135:.2cm) and ++(90:.4cm) .. (-1.2,-1.4) -- (-1.2,-1.8);
\filldraw[boxregion=\PrColor] (-1.2,1.8) -- (-1.2,-.2) .. controls ++(270:.4cm) and ++(135:.2cm) .. (-.6,-.8) .. controls ++(45:.2cm) and ++(270:.4cm) .. (0,-.2) -- (0,.2) .. controls ++(90:.4cm) and ++(-135:.2cm) .. (.6,.8) .. controls ++(135:.2cm) and ++(270:.4cm) .. (0,1.4) -- (0,1.8);
\filldraw[boxregion=\QrColor] (0,-1.8) -- (0,-1.4) .. controls ++(90:.4cm) and ++(-45:.2cm) .. (-.6,-.8) .. controls ++(45:.2cm) and ++(270:.4cm) .. (0,-.2) -- (0,.2) .. controls ++(90:.4cm) and ++(-135:.2cm) .. (.6,.8) .. controls ++(-45:.2cm) and ++(90:.4cm) .. (1.2,.2) -- (1.2,-1.8);
\filldraw[plusregion=\PrColor] (-.55,2.7) -- (-.55,1.8) -- (0,1.8) -- (0,1.4) .. controls ++(270:.4cm) and ++(135:.2cm) .. (.6,.8) .. controls ++(45:.2cm) and ++(270:.4cm) .. (1.2,1.4) -- (1.2,2.7);
\filldraw[plusregion=\QrColor] (1.2,-1.8) -- (1.2,.2) .. controls ++(90:.4cm) and ++(-45:.2cm) .. (.6,.8) .. controls ++(45:.2cm) and ++(270:.4cm) .. (1.2,1.4) -- (1.2,2.7) -- (1.8,2.7) -- (1.8,-1.8);
\end{scope}
\draw[\XColor,thick] (-1.2,-1.8) -- (-1.2,-1.4) .. controls ++(90:.6cm) and ++(270:.6cm) .. (0,-.2) -- (0,.2) .. controls ++(90:.6cm) and ++(270:.6cm) .. (1.2,1.4) -- (1.2,2.7);
\draw[\phiColor,thick] (0,-1.8) -- (0,-1.4) .. controls ++(90:.6cm) and ++(270:.6cm) .. (-1.2,-.2) -- (-1.2,1.8);
\draw[\psiColor,thick] (1.2,-1.8) -- (1.2,.2) .. controls ++(90:.6cm) and ++(270:.6cm) .. (0,1.4) -- (0,1.8);
\draw[\phiColor,thick] (-.65,1.8) -- (-.65,2.7);
\draw[\psiColor,thick] (-.55,1.8) -- (-.55,2.7);
\roundNbox{unshaded}{(-.6,1.8)}{.33}{.6}{.6}{\scriptsize{$\left(\QSys^\xo_{\varphi,\psi}\right)_P$}};
\filldraw[white] (-.6,-.8) circle (.1cm);
\draw[thick] (-.6,-.8) circle (.1cm); 
\filldraw[white] (.6,.8) circle (.1cm);
\draw[thick] (.6,.8) circle (.1cm); 
}\,.
\]
Finally, we check the monoidality coherence axiom for $\QSys^\otimes_{\bullet,\bullet}$, and we leave $\QSys^\otimes_{\bullet}$ to the reader:
\[
\tikzmath[scale=.5, transform shape]{
\begin{scope}
\clip[rounded corners=5pt] (-1.9,.6) rectangle (1.7,7.4);
\filldraw[primedregion=\QrColor] (-1.9,0) rectangle (1.7,8);
\filldraw[plusregion=\QrColor] (1,0) rectangle (-.4,4);
\filldraw[boxregion=\QrColor] (-1,0) rectangle (0,2);
\filldraw[boxregion=\QrColor] (-.6,2) rectangle (-.4,4);
\filldraw[boxregion=\QrColor] (-.2,4) rectangle (0,8);
\filldraw[plusregion=\QrColor] (0,4) rectangle (.2,8);
\filldraw[starregion=\QrColor] (1,0) -- (1,4) -- (.2,4) -- (.2,8) -- (1.7,8) -- (1.7,0);
\end{scope}
\draw[\phiColor,thick] (-1,.6) -- (-1,2);
\draw[\psiColor,thick] (0,.6) -- (0,2);
\draw[\phiColor,thick] (-.6,2) -- (-.6,4);
\draw[\psiColor,thick] (-.4,2) -- (-.4,4);
\draw[\gammaColor,thick] (1,.6) -- (1,4);
\draw[\phiColor,thick] (-.2,4) -- (-.2,7.4);
\draw[\psiColor,thick] (0,4) -- (0,7.4);
\draw[\gammaColor,thick] (.2,4) -- (.2,7.4);
\roundNbox{unshaded}{(0,6)}{.6}{.95}{.75}{\normalsize{$\left(\QSys(\alpha^\xo_{\varphi,\psi,\gamma})\right)_Q$}};
\roundNbox{unshaded}{(0,4)}{.6}{.95}{.75}{\normalsize{$\left(\QSys^\xo_{\varphi\xo\psi,\gamma}\right)_Q$}};
\roundNbox{unshaded}{(-.5,2)}{.6}{.45}{.45}{\normalsize{$\left(\QSys^\xo_{\varphi,\psi}\right)_Q$}};
\draw[thin, dotted, rounded corners = 5pt] (-1.9,.6) rectangle (1.7,7.4);
}
=
\tikzmath[scale=.7, transform shape]{
\begin{scope}
\clip[rounded corners = 5] (-2.2,-2.9) rectangle (2.2,1.2);
\filldraw[primedregion=\QrColor] (-2.2,-2.9) rectangle (-1.6,1.2);
\filldraw[primedregion=\BColor] (-1.6,-2.9) rectangle (-1.2,1.2);
\filldraw[boxregion=\BColor] (-1.2,-2.9) rectangle (0,1.2);
\filldraw[boxregion=\QrColor] (-.8,-2.9) rectangle (-.4,-2.3);
\filldraw[plusregion=\BColor] (0,-2.9) rectangle (1.2,1.2);
\filldraw[plusregion=\QrColor] (.4,-2.9) rectangle (.8,-1.5);
\filldraw[starregion=\BColor] (1.2,-2.9) rectangle (1.6,1.2);
\filldraw[starregion=\QrColor] (1.6,-2.9) rectangle (2.2,1.2);
\end{scope}
\draw[\phiColor,thick] (-1.2,-2.9) -- (-1.2,1.2); 
\draw[\psiColor,thick] (0,-2.9) -- (0,1.2);
\draw[\gammaColor,thick] (1.2,-2.9) -- (1.2,1.2); 
\draw[\QsColor,thick] (-1.6,-2.9) -- (-1.6,1.2);
\draw[\QsColor,thick] (1.6,-2.9) -- (1.6,1.2);
\draw[\QsColor,thick] (.4,-2.9) -- (.4,-1.5);
\draw[\QsColor,thick] (.8,-2.9) -- (.8,-1.5);
\draw[\QsColor,thick] (-.4,-2.9) -- (-.4,-2.3);
\draw[\QsColor,thick] (-.8,-2.9) -- (-.8,-2.3);
\draw[\QsColor,thick] (-1.6,.7) -- (1.6,.7);
\draw[\QsColor,thick] (-1.6,-.7) -- (1.6,-.7);
\draw[\QsColor,thick] (-1.6,-1.5) -- (1.6,-1.5);
\draw[\QsColor,thick] (-1.6,-2.3) -- (.4,-2.3);
\filldraw[white] (-1.2,.7) circle (.1cm);
\draw[thick] (-1.2,.7) circle (.1cm);
\filldraw[white] (-1.2,-.7) circle (.1cm);
\draw[thick] (-1.2,-.7) circle (.1cm);
\filldraw[white] (-1.2,-1.5) circle (.1cm);
\draw[thick] (-1.2,-1.5) circle (.1cm);
\filldraw[white] (-1.2,-2.3) circle (.1cm);
\draw[thick] (-1.2,-2.3) circle (.1cm);
\filldraw[white] (0,.7) circle (.1cm);
\draw[thick] (0,.7) circle (.1cm);
\filldraw[white] (0,-.7) circle (.1cm);
\draw[thick] (0,-.7) circle (.1cm);
\filldraw[white] (0,-1.5) circle (.1cm);
\draw[thick] (0,-1.5) circle (.1cm);
\filldraw[white] (0,-2.3) circle (.1cm);
\draw[thick] (0,-2.3) circle (.1cm);
\filldraw[white] (1.2,.7) circle (.1cm);
\draw[thick] (1.2,.7) circle (.1cm);
\filldraw[white] (1.2,-.7) circle (.1cm);
\draw[thick] (1.2,-.7) circle (.1cm);
\filldraw[white] (1.2,-1.5) circle (.1cm);
\draw[thick] (1.2,-1.5) circle (.1cm);
\roundNbox{unshaded}{(0,0)}{.33}{1.1}{1.1}{\scriptsize{$\left(\alpha^\xo_{\varphi,\psi,\gamma}\right)_b$}};
}
=
\tikzmath[scale=.7, transform shape]{
\begin{scope}
\clip[rounded corners = 5] (-2.2,-1.2) rectangle (2.2,1.2);
\filldraw[primedregion=\QrColor] (-2.2,-1.2) rectangle (-1.6,1.2);
\filldraw[primedregion=\BColor] (-1.6,-1.2) rectangle (-1.2,1.2);
\filldraw[boxregion=\BColor] (-1.2,-1.2) rectangle (0,1.2);
\filldraw[boxregion=\QrColor] (-.8,-1.2) rectangle (-.4,-.7);
\filldraw[plusregion=\BColor] (0,-1.2) rectangle (1.2,1.2);
\filldraw[plusregion=\QrColor] (.4,-1.2) rectangle (.8,-.7);
\filldraw[starregion=\BColor] (1.2,-1.2) rectangle (1.6,1.2);
\filldraw[starregion=\QrColor] (1.6,-1.2) rectangle (2.2,1.2);
\end{scope}
\draw[\phiColor,thick] (-1.2,-1.2) -- (-1.2,1.2); 
\draw[\psiColor,thick] (0,-1.2) -- (0,1.2);
\draw[\gammaColor,thick] (1.2,-1.2) -- (1.2,1.2); 
\draw[\QsColor,thick] (-1.6,-1.2) -- (-1.6,1.2);
\draw[\QsColor,thick] (1.6,-1.2) -- (1.6,1.2);
\draw[\QsColor,thick] (-.4,-1.2) -- (-.4,-.7);
\draw[\QsColor,thick] (-.8,-1.2) -- (-.8,-.7);
\draw[\QsColor,thick] (.4,-1.2) -- (.4,-.7);
\draw[\QsColor,thick] (.8,-1.2) -- (.8,-.7);
\draw[\QsColor,thick] (-1.6,.7) -- (1.6,.7);
\draw[\QsColor,thick] (-1.6,-.7) -- (1.6,-.7);
\filldraw[white] (-1.2,.7) circle (.1cm);
\draw[thick] (-1.2,.7) circle (.1cm);
\filldraw[white] (-1.2,-.7) circle (.1cm);
\draw[thick] (-1.2,-.7) circle (.1cm);
\filldraw[white] (0,.7) circle (.1cm);
\draw[thick] (0,.7) circle (.1cm);
\filldraw[white] (0,-.7) circle (.1cm);
\draw[thick] (0,-.7) circle (.1cm);
\filldraw[white] (1.2,.7) circle (.1cm);
\draw[thick] (1.2,.7) circle (.1cm);
\filldraw[white] (1.2,-.7) circle (.1cm);
\draw[thick] (1.2,-.7) circle (.1cm);
\roundNbox{unshaded}{(0,0)}{.33}{1.1}{1.1}{\scriptsize{$\left(\alpha^\xo_{\varphi,\psi,\gamma}\right)_b$}};
}
=
\tikzmath[scale=.7, transform shape]{
\begin{scope}
\clip[rounded corners = 5] (-2.2,-1.2) rectangle (2.2,2.9);
\filldraw[primedregion=\QrColor] (-2.2,-1.2) rectangle (-1.6,2.9);
\filldraw[primedregion=\BColor] (-1.6,-1.2) rectangle (-1.2,2.9);
\filldraw[boxregion=\BColor] (-1.2,-1.2) rectangle (0,2.9);
\filldraw[boxregion=\QrColor] (-.8,-1.2) rectangle (-.4,-.7);
\filldraw[boxregion=\QrColor] (-.8,.7) rectangle (-.4,2.3);
\filldraw[plusregion=\BColor] (0,-1.2) rectangle (1.2,2.9);
\filldraw[plusregion=\QrColor] (.4,-1.2) rectangle (.8,-.7);
\filldraw[plusregion=\QrColor] (.4,.7) rectangle (.8,1.5);
\filldraw[starregion=\BColor] (1.2,-1.2) rectangle (1.6,2.9);
\filldraw[starregion=\QrColor] (1.6,-1.2) rectangle (2.2,2.9);
\end{scope}
\draw[\phiColor,thick] (-1.2,-1.2) -- (-1.2,2.9); 
\draw[\psiColor,thick] (0,-1.2) -- (0,2.9);
\draw[\gammaColor,thick] (1.2,-1.2) -- (1.2,2.9); 
\draw[\QsColor,thick] (-1.6,-1.2) -- (-1.6,2.9);
\draw[\QsColor,thick] (1.6,-1.2) -- (1.6,2.9);
\draw[\QsColor,thick] (-.4,-1.2) -- (-.4,-.7);
\draw[\QsColor,thick] (-.8,-1.2) -- (-.8,-.7);
\draw[\QsColor,thick] (.4,-1.2) -- (.4,-.7);
\draw[\QsColor,thick] (.8,-1.2) -- (.8,-.7);
\draw[\QsColor,thick] (-.4,.7) -- (-.4,2.3);
\draw[\QsColor,thick] (-.8,.7) -- (-.8,2.3);
\draw[\QsColor,thick] (.4,.7) -- (.4,1.5);
\draw[\QsColor,thick] (.8,.7) -- (.8,1.5);
\draw[\QsColor,thick] (-1.6,.7) -- (1.6,.7);
\draw[\QsColor,thick] (-1.6,-.7) -- (1.6,-.7);
\draw[\QsColor,thick] (-1.6,2.3) -- (1.6,2.3);
\draw[\QsColor,thick] (-.4,1.5) -- (1.6,1.5);
\filldraw[white] (-1.2,.7) circle (.1cm);
\draw[thick] (-1.2,.7) circle (.1cm);
\filldraw[white] (-1.2,-.7) circle (.1cm);
\draw[thick] (-1.2,-.7) circle (.1cm);
\filldraw[white] (-1.2,2.3) circle (.1cm);
\draw[thick] (-1.2,2.3) circle (.1cm);
\filldraw[white] (0,.7) circle (.1cm);
\draw[thick] (0,.7) circle (.1cm);
\filldraw[white] (0,-.7) circle (.1cm);
\draw[thick] (0,-.7) circle (.1cm);
\filldraw[white] (0,1.5) circle (.1cm);
\draw[thick] (0,1.5) circle (.1cm);
\filldraw[white] (0,2.3) circle (.1cm);
\draw[thick] (0,2.3) circle (.1cm);
\filldraw[white] (1.2,.7) circle (.1cm);
\draw[thick] (1.2,.7) circle (.1cm);
\filldraw[white] (1.2,-.7) circle (.1cm);
\draw[thick] (1.2,-.7) circle (.1cm);
\filldraw[white] (1.2,1.5) circle (.1cm);
\draw[thick] (1.2,1.5) circle (.1cm);
\filldraw[white] (1.2,2.3) circle (.1cm);
\draw[thick] (1.2,2.3) circle (.1cm); 
\roundNbox{unshaded}{(0,0)}{.33}{1.1}{1.1}{\scriptsize{$\left(\alpha^\xo_{\varphi,\psi,\gamma}\right)_b$}};
}
=
\tikzmath[scale=.5, transform shape]{
\begin{scope}
\clip[rounded corners = 5] (-2.3,.6) rectangle (2.5,7.4);
\filldraw[primedregion=\QrColor] (-2.4,0) rectangle (2.6,8);
\filldraw[boxregion=\QrColor] (-1,0) rectangle (.4,6);
\filldraw[plusregion=\QrColor] (0,0) rectangle (1,4);
\filldraw[plusregion=\QrColor] (.6,4) rectangle (.4,6);
\filldraw[boxregion=\QrColor] (-.2,6) rectangle (0,8);
\filldraw[plusregion=\QrColor] (0,6) rectangle (.2,8);
\filldraw[starregion=\QrColor] (1,0) -- (1,4) -- (.6,4) -- (.6,6) -- (.2,6) -- (.2,8) -- (2.6,8) -- (2.6,.6);
\end{scope}
\draw[\phiColor,thick] (-1,.6) -- (-1,6);
\draw[\psiColor,thick] (0,.6) -- (0,4);
\draw[\gammaColor,thick] (1,.6) -- (1,4);
\draw[\psiColor,thick] (.4,4) -- (.4,6);
\draw[\gammaColor,thick] (.6,4) -- (.6,6);
\draw[\phiColor,thick] (-.2,6) -- (-.2,7.4);
\draw[\psiColor,thick] (0,6) -- (0,7.4);
\draw[\gammaColor,thick] (.2,6) -- (.2,7.4);
\roundNbox{unshaded}{(0,6)}{.6}{.75}{.95}{\normalsize{$\left(\QSys^\xo_{\varphi,\psi\xo\gamma}\right)_Q$}};
\roundNbox{unshaded}{(.5,4)}{.6}{.45}{.45}{\normalsize{$\left(\QSys^\xo_{\psi,\gamma}\right)_Q$}};
\roundNbox{unshaded}{(0,2)}{.6}{1.3}{1.5}{\small{$\left(\alpha^\xo_{\QSys(\varphi),\QSys(\psi),\QSys(\gamma)}\right)_Q$}};
\draw[thin, dotted, rounded corners = 5pt] (-2.3,.6) rectangle (2.5,7.4);
}
\,.
\]
\end{construction}

Constructions \ref{construction:Qsys(F)}, \ref{construction:QSys(phi)}, \ref{construction:QSys(n)}, and \ref{construction:QSys-otimes} immediately imply the following proposition.

\begin{prop}
$\QSys$ as defined above is a 2-functor on every hom 2-category $\Fun^\dag(\cA\to \cB)$.
\end{prop}

\begin{lem}
\label{lem:Strict1Associator}
For $F\in 2\Cat(\cA\to\cB)$ and $G\in  2\Cat(\cB\to \cC)$, 
$\QSys(G)\xz\QSys(F)=\QSys(G\xz F)$.
\end{lem}
\begin{proof}
By Constructions \ref{const:1CompositionIn2Cat} and \ref{construction:Qsys(F)},
for a 0-cell $Q\in\QSys(\cA)$, 
$$
\QSys(G\xz F)(Q)=G(F(Q))=\QSys(G)(\QSys(F)(Q))=[\QSys(G)\xz \QSys(F)](Q),
$$
for a 1-cell $X\in\QSys(\cA)(P\to Q)$, 
$$
\QSys(G\xz F)(X)=G(F(X))=[\QSys(G)\xz \QSys(F)](X),
$$ 
and for a 2-cell $f\in\QSys(\cA)(X\Rightarrow Y)$, 
$$
\QSys(G\xz F)(f)=G(F(f))=[\QSys(G)\xz \QSys(F)](f).
$$
For a 0-cell $Q\in\QSys(\cA)$, $\QSys(F)^1_Q=\id$, so $\QSys(G\xz F)^1_Q=\id=(\QSys(G)\xz\QSys(F))^1_Q$.
For 1-cells $X\in\QSys(\cA)(P\to Q)$ and $Y\in\QSys(\cA)(Q\to R)$, we have 
\begin{small}
\begin{align*}
(\QSys(G)&\xz\QSys(F))^2_{X,Y}
\\ &= 
\QSys(G)(\QSys(F)^2_{X,Y})\xxt \QSys(G)^2_{\QSys(F)(X),\QSys(F)(Y)}
\\
& = \QSys(G)(F(u^Q_{X,Y})\xxt F^2_{X,Y}\xxt (u^{F(Q)}_{F(X),F(Y)})^\dag)\xxt \QSys(G)^2_{F(X),F(Y)}
\\
& = G(F(u^Q_{X,Y}))\xxt G(F^2_{X,Y})\xxt G({(u^{F(Q)}_{F(X),F(Y)})}^\dag)\xxt G(u^{F(Q)}_{F(X),F(Y)})
\xxt G^2_{F(X),F(Y)}\xxt {(u^{G(F(Q))}_{G(F(X)),G(F(Y))})}^\dag
\\
& = G(F(u^Q_{X,Y}))\xxt G(F^2_{X,Y})\xxt G^2_{F(X),F(Y)}\xxt {(u^{G(F(Q))}_{G(F(X)),G(F(Y))})}^\dag 
\\
& = (G\xz F)(u^Q_{X,Y})\xxt (G\xz F)^2_{X,Y}\xxt {(u^{(G\xz F)(Q)}_{(G\xz F)(X),(G\xz F)(Y)})}^\dag
\\
& = \QSys(G\xz F)^2_{X,Y}.
\end{align*}
\end{small}
Hence $\QSys(G)\xz\QSys(F)=\QSys(G\xz F)$ as claimed.
\end{proof}

By Lemma \ref{lem:Strict1Associator},
we may define each $\QSys^\xz_{G,F}:\QSys(G)\xz\QSys(F)\Rightarrow \QSys(G\xz F)$ to be the identity transformation,
and we may define each 1-associator modification $\omega^\xz_{H,G,F}$ to be the identity modification, as well as each unitor modification $\ell^\xz_F$ and $r^\xz_G$.
Theorem \ref{thm:QSys3Functor} follows immediately, i.e., $\QSys$ is a $\dag$ 3-functor.
\qed

\begin{rem}
Since 1-composition is strict in $2\Cat$, 2-categories and 2-functors form a 1-category where we forget all transformations and modifications.
(Observe we have \emph{not} truncated, as this would identify equivalent 2-functors.)
Lemma \ref{lem:Strict1Associator} shows that $\QSys$ is a functor on this 1-category.
\end{rem}

\begin{rem}
It was pointed out to us by 
Thibault D\'ecoppet 
and
David Reutter
that our $\dag$ 3-endofunctor $\QSys$ on $\rm C^*/W^*$ $2\Cat$ should be left 3-adjoint to the inclusion of the full 3-subcategory on the Q-system complete $\rm C^*/W^*$ 2-categories.
We will not prove this here as it would take us too far afield.
We note, however, that this would endow $\QSys$ with the structure of a symmetric lax monoidal $\dag$ 3-functor.
\end{rem}

\section{Universal property of Q-system completion}
\label{sec:UniversalProperty}

In this section, we give the strongest possible universal property which is satisfied by Q-system completion.
Namely, we prove Theorem \ref{thm:UniqueLift}, which states that the \emph{lift 2-category} of a $\dag$ 2-functor $F: \cC\to \cD$ from a $\rm C^*/W^*$ 2-category $\cC$ into a Q-system complete $\rm C^*/W^*$ 2-category $\cD$ is \emph{$(-2)$-truncated}, i.e., equivalent to a point.
We now define the necessary terms to prove this theorem, and we explain the proof strategy from \cite[\S3.1]{1910.03178}.

\subsection{Lift categories and homotopy fibers}

\begin{defn}
\label{defn:Lift2Category}
Suppose $\cC,\cD,\cE$ are $\rm C^*/W^*$ 2-categories and $F: \cC\to \cD$ and  $G: \cC \to \cE$ are $\dag$ 2-functors.
The \emph{lift 2-category} of $F$ along $G$ is the \emph{homotopy fiber} 2-category of the pre-composition 2-functor $-\xz G : \Fun^\dag(\cE\to \cD) \to \Fun^\dag(\cC\to \cD)$ at $F\in \Fun^\dag(\cC\to \cD)$.
We remind the reader that the definition of $-\xz-$ in $2\Cat$ is detailed in Construction \ref{const:1CompositionIn2Cat} above.
\end{defn}

\begin{rem}
We now further unpack Defintion \ref{defn:Lift2Category}.
The lift 2-category of $F$ along $G$ has:
\begin{itemize}
\item \underline{objects:}
pairs $(A,\alpha)$, where $A:\cE\to \cD$ is a $\dag$ 2-functor and $\alpha:F\Rightarrow A\xz G$ is a unitary 2-transformation.
\[
    \begin{tikzcd}[row sep=1.2em]
    \cC\arrow[rr, "G"]
    \arrow[ddrr, swap, "F"]
    &&
    \cE\arrow[dd,"A"]
    \arrow[dl,Leftarrow,shorten <= 1em, shorten >= 1em, "\alpha"]
    \\
    &\mbox{}&
    \\
    &&
    \cD.
    \end{tikzcd}
\]

\item \underline{1-morphisms:}
pairs $(\varphi,m):(A,\alpha)\to(B,\beta)$, where $\varphi:A\Rightarrow B$ is a $\dag$ 2-transformation and $m:\beta\Rrightarrow \alpha\xo (\varphi\xz G)$ is a unitary 2-modification: 
\[
    \begin{tikzcd}[row sep=1.2em]
    \cC\arrow[rr, "G"]
    \arrow[ddrr, swap, "F"]
    &&
    \cE\arrow[dd,"B"]
    \arrow[dl,Leftarrow,shorten <= 1em, shorten >= 1em, "\beta"]
    \\
    &\mbox{}&
    \\
    &&
    \cD.
    \end{tikzcd}
    \qquad
    \overset{m}{\Rrightarrow}
    \qquad
    \begin{tikzcd}[row sep=1.2em]
    \cC\arrow[rr, "G"]
    \arrow[ddrr, swap, "F"]
    &&
    \cE\arrow[dd,swap,"A"]
    \arrow[dl,Leftarrow,shorten <= 1em, shorten >= 1em, "\alpha"]
    \arrow[dd,bend left = 90, "B"]
    \\
    &\mbox{}&
    \arrow[r,Rightarrow, shorten >= 1.5em, "\!\!\!\!\!\!\!\!\!\!\!\varphi"]
    &\mbox{}
    \\
    &&
    \cD
    \end{tikzcd}
\]

\item \underline{2-morphisms:}
$p:(\varphi,m) \Rightarrow (\psi,n)$, where $p:\varphi\Rrightarrow\psi$ is a $\dag$ 2-modification such that
\[
\begin{tikzcd}[row sep=1.2em]
F\arrow[rr,Rightarrow, "\alpha"]
\arrow[ddrr, Rightarrow,swap, "\beta"]
&&
A\xz G\arrow[dd,Rightarrow, "\psi\xz G"]
\\
&
\mbox{}
\arrow[ur,triplecd,shorten <= 1em, shorten >= 1em, "n"]
&
\\
&&
B\xz G.
\end{tikzcd}
\qquad
=
\qquad
\begin{tikzcd}[row sep=1.2em]
F\arrow[rr,Rightarrow, "\alpha"]
\arrow[ddrr, Rightarrow,swap, "\beta"]
&&
A\xz G\arrow[dd,swap,Rightarrow, "\varphi\xz G"]
\arrow[dd,Rightarrow, bend left = 90, "\psi \xz G"]
\\
&
\arrow[ur,triplecd,shorten <= 1em, shorten >= 1em, "m"]
\mbox{}
&
\arrow[r,triplecd,shorten >=1.5em, "\!\!\!\!\!\!\!\!\!\!\!p \xz G"]
&\mbox{}
\\
&&
B\xz G
\end{tikzcd}
\]
\end{itemize}
\end{rem}

Recall that for a 2-category $\cC$, its \emph{core} is the 2-subcategory $\core(\cC)$ with only invertible 1-cells and invertible 2-cells.
When $\cC$ is $\rm C^*/W^*$, its \emph{unitary core} $\core^\dag(\cC)$ is the 2-subcategory of $\core(\cC)$ with only unitary 2-cells.
In a $\rm C^*/W^*$ 2-category, by polar decomposition for invertible 2-cells, there exists an invertible 2-cell $\cC({}_aX_b\Rightarrow {}_aY_b)$ if and only if there exists a unitary 2-cell, so the connectivity of $\core(\cC)$ and $\core^\dag(\cC)$ agree.

We pass to cores in order to utilize the notion of $k$-truncated 2-functor between 2-groupoids from \cite[\S3.1]{1910.03178}.

\begin{defn}[{cf.~\cite[Def.~3.3]{1910.03178}}]
Suppose $\cC, \cD$ are $2$-groupoids and $U: \cC \to \cD$ is a 2-functor.
We call $U$ $k$-\emph{truncated} or $(k+1)$-\emph{monic} \cite[\S5.5]{MR2664619} if $U$ is:
\begin{itemize}
\item
$k=2$: (no condition)
\item
$k=1$: faithful on $2$-cells
\item
$k=0$: fully faithful on $2$-cells
\item
$k=-1$: an equivalence on hom-categories
\item
$k=-2$: an equivalence of 2-categories.
\end{itemize}
\end{defn}

The following proposition connects the notions of a $k$-truncated 2-functor between 2-groupoids and its homotopy fibers.

\begin{prop}[{cf.~\cite[Prop.~3.4]{1910.03178}}]
\label{prop:TruncatedGroupoids}
Suppose $\cC,\cD$ are $2$-groupoids, and $U: \cC \to \cD$ is a 2-functor.
For every $-2\leq k\leq 2$,
$U$ is $k$-truncated if and only if 
at each object $d\in \cD$, the homotopy fiber $\mathsf{hoFib}_d(U)$ is $k$-truncated as an $2$-groupoid, i.e., a $k$-groupoid.\footnote{We use `negative categorical thinking' \cite{MR2664619} when $k=-2,-1,0$.
That is, a $0$-groupoid is a set, a $(-1)$-groupoid is either a point or the empty set, and a $(-2)$-groupoid is a point.}
\end{prop}

\subsection{Dominance and truncation}
\label{sec:Dominance}

Observe that $-\xz G: \Fun^\dag(\cE \to \cD) \to \Fun^\dag(\cC \to \cD)$ restricts to a $\dag$ 2-functor $-\xz G: \core^\dag(\Fun^\dag(\cE \to \cD)) \to \core^\dag(\Fun^\dag(\cC \to \cD))$.
Hence, in order to apply Proposition \ref{prop:TruncatedGroupoids} to the $\dag$ 2-functor $-\xz G$, we need (essential) surjectivity conditions on $-\xz G$.
(Being faithful on 2-morphisms is being surjective on equalities between 2-morphisms.)
A suitable notion of (essential) surjectivity for a linear 2-functor is \emph{dominance}, which we define via the notion of \emph{condensation} in a 2-category \cite{1905.09566}.

\begin{defn}
\label{defn:Condensation}
Suppose $\cC$ is a 2-category and $a,b\in \cC$ are 0-cells.
A \emph{condensation} $X: a\condense b$ consists of 1-cells ${}_aX_b, {}_bX_a^{\bullet}$ and 2-cells $\varepsilon_X : {}_bX^\bullet\xxo_a X_b\to 1_b$ and $\delta_X: 1_b \to {}_bX^\bullet\xxo_a X_b$
such that $\varepsilon_X\xxt \delta_X = 1_{1_b}$.
Graphically, we denote $X: a\condense b$ by
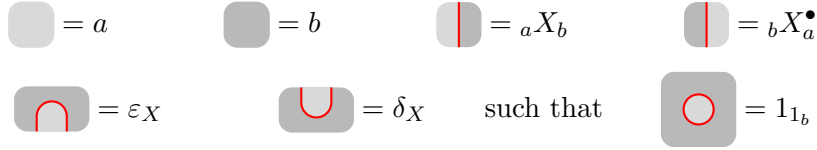
\[
\tikzmath{\filldraw[\AColor, rounded corners=5, very thin, baseline=1cm] (0,0) rectangle (.6,.6);}=a
\qquad\qquad
\tikzmath{\filldraw[\BColor, rounded corners=5, very thin, baseline=1cm] (0,0) rectangle (.6,.6);}=b
\qquad\qquad
\tikzmath{
\begin{scope}
\clip[rounded corners=5pt] (-.3,0) rectangle (.3,.6);
\fill[\AColor] (0,0) rectangle (-.3,.6);
\fill[\BColor] (0,0) rectangle (.3,.6);
\end{scope}
\draw[thick, \XColor] (0,0) -- (0,.6);
}
={}_{a}X_b
\qquad\qquad
\tikzmath{
\begin{scope}
\clip[rounded corners=5pt] (-.3,0) rectangle (.3,.6);
\fill[\BColor] (0,0) rectangle (-.3,.6);
\fill[\AColor] (0,0) rectangle (.3,.6);
\end{scope}
\draw[thick, \XColor] (0,0) -- (0,.6);
}
={}_{b}X^\bullet_a
\]
\[
\tikzmath{
\begin{scope}
\clip[rounded corners=5pt] (-.5,0) rectangle (.5,.6);
\filldraw[\BColor] (-.5,0) rectangle (.5,.8);
\filldraw[\AColor] (-.2,0) -- (-.2,.2) arc (180:0:.2cm) -- (.2,0);
\end{scope}
\draw[thick,\XColor] (-.2,0) -- (-.2,.2) arc (180:0:.2cm) -- (.2,0);
}
=\varepsilon_X
\qquad\qquad
\tikzmath{
\begin{scope}
\clip[rounded corners=5pt] (-.5,0) rectangle (.5,.6);
\filldraw[\BColor] (-.5,0) rectangle (.5,.8);
\filldraw[\AColor] (-.2,.6) -- (-.2,.4) arc (-180:0:.2cm) -- (.2,.6);
\end{scope}
\draw[thick,\XColor] (-.2,.6) -- (-.2,.4) arc (-180:0:.2cm) -- (.2,.6);
}
=\delta_X
\qquad\text{such that}\qquad
\tikzmath{
\begin{scope}
\clip[rounded corners=5pt] (-.5,-.5) rectangle (.5,.5);
\filldraw[\BColor] (-.5,-.5) rectangle (.5,.5);
\filldraw[\AColor] (0,0) circle (.2cm);
\end{scope}
\draw[thick,\XColor] (0,0) circle (.2cm);
}
=1_{1_b}
\]
When $\cC$ is $\rm C^*/W^*$, a condensation $X: a\condense b$ is called a \emph{dagger condensation} if $\delta_X = \varepsilon_X^\dag$.
\end{defn}

\begin{defn}
A 2-functor $G: \cC \to \cE$ is called:
\begin{itemize}
    \item \emph{0-dominant} if for all $e\in\cE$, there is a condensation $G(c)\condense e$ for some $c\in \cC$,
    \item \emph{locally dominant} if every hom functor $G_{a\to b}: \cC(a\to b) \to \cE(G(a)\to G(b))$ is dominant as a linear functor, and
    \item \emph{dominant} if $G$ is both 0-dominant and locally dominant.
\end{itemize}
When $G$ is a $\dag$ 2-functor between $\rm C^*/W^*$ 2-categories,
we call $G$ 
\begin{itemize}
\item 
\emph{orthogonally} 0-dominant if for all $e\in \cE$, there is a dagger condensation $G(c) \condense e$ for some $c\in \cC$,
\item
locally \emph{orthogonally} dominant if every hom functor $G_{a\to b}$ is orthogonally dominant as a linear $\dag$-functor, i.e., every 1-cell ${}_{G(a)}Y_{G(b)}\in \cE$ is unitarily isomorphic to an orthogonal direct summand of some ${}_{G(a)}G(X)_{G(b)}$, and
\item
\emph{orthogonally} dominant if $G$ is both orthogonally 0-dominant and locally dominant.
\end{itemize}
\end{defn}

\begin{rem}
There is an analogous notion of
$k$-dominance for an $n$-functor $G$ between $n$-categories for $0\leq k\leq n-1$:
every $k$-morphism in the target should admit a 
condensation from a source in the image of $G$.
\end{rem}

For the propositions in this section, we work with algebraic 2-categories and 2-functors, and we make particular comments about the $\rm C^*/W^*$ setting.

\begin{prop}
If a 2-functor $G:\cC\to \cE$ is 0-dominant, then $-\xz G:\Fun(\cC\to \cE)\to \Fun(\cC\to \cD)$ is faithful on 2-morphisms.
In the $\rm C^*/W^*$ setting,
if $G: \cC\to \cE$ is orthogonally 0-dominant, then $-\xz G: \Fun^\dag(\cC\to \cE) \to \Fun^\dag(\cC\to \cD)$ is faithful on 2-morphisms.
\end{prop}
\begin{proof}
Let $A,B\in\Fun(\cE\to\cD)$ and $\varphi,\psi:A\Rightarrow B$.
Suppose $m,n:\varphi\Rrightarrow\psi$ and $m\xz G=n\xz G$.
We show $m=n$.
For each $e\in\cE$, there exists a 0-cell $c\in\cC$ and a condensation $X:G(c)\condense e$.
We denote $G(c),e\in \cD$, ${}_{G(c)}X_e\in \cD(G(c)\to e)$, and the functors $A,B$ graphically by
\begin{equation}
\label{eq:GraphicalCalculusFaithfulOn2Mor}
\tikzmath{\filldraw[\AColor, rounded corners=5, very thin, baseline=1cm] (0,0) rectangle (.6,.6);}=G(c)
\qquad\qquad
\tikzmath{\filldraw[\BColor, rounded corners=5, very thin, baseline=1cm] (0,0) rectangle (.6,.6);}=e
\qquad\qquad
\tikzmath{
\begin{scope}
\clip[rounded corners=5pt] (-.3,0) rectangle (.3,.6);
\fill[\AColor] (0,0) rectangle (-.3,.6);
\fill[\BColor] (0,0) rectangle (.3,.6);
\end{scope}
\draw[thick, \XColor] (0,0) -- (0,.6);
}
={}_{G(c)}X_e
\qquad\qquad
\tikzmath{
\filldraw[primedregion=white, rounded corners = 5pt] (0,0) rectangle (.6,.6);
\draw[thin, dotted, rounded corners = 5pt] (0,0) rectangle (.6,.6);
}
=
A
\qquad\qquad
\tikzmath{
\filldraw[boxregion=white, rounded corners = 5pt] (0,0) rectangle (.6,.6);
\draw[thin, dotted, rounded corners = 5pt] (0,0) rectangle (.6,.6);
}
=
B.
\end{equation}
The modification axiom implies the following:
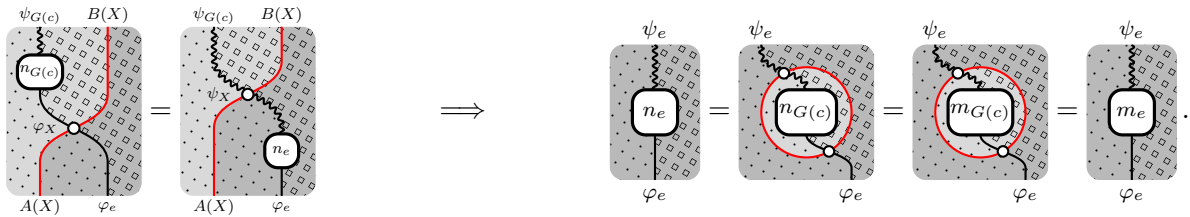
\[
\tikzmath[scale=.75, transform shape]{
\begin{scope}
\clip[rounded corners = 5] (-.6,0) rectangle (1.8,3);
\filldraw[primedregion=\AColor] (0,0) -- (0,.6) .. controls ++(90:.4cm) and ++(-135:.2cm) .. (.6,1.2) .. controls ++(135:.2cm) and ++(270:.4cm) .. (0,1.8) -- (0,3) -- (-.6,3) -- (-.6,0); 
\filldraw[primedregion=\BColor] (1.2,0) -- (1.2,.6) .. controls ++(90:.4cm) and ++(-45:.2cm) .. (.6,1.2) .. controls ++(-135:.2cm) and ++(90:.4cm) .. (0,.6) -- (0,0);
\filldraw[boxregion=\AColor] (0,3) -- (0,1.8) .. controls ++(270:.4cm) and ++(135:.2cm) .. (.6,1.2) .. controls ++(45:.2cm) and ++(270:.4cm) .. (1.2,1.8) -- (1.2,3);
\filldraw[boxregion=\BColor] (1.2,0) -- (1.2,.6) .. controls ++(90:.4cm) and ++(-45:.2cm) .. (.6,1.2) .. controls ++(45:.2cm) and ++(270:.4cm) .. (1.2,1.8) -- (1.2,3) -- (1.8,3) -- (1.8,0);
\filldraw[\AColor] (-.1,1) rectangle (.29,1.4);
\end{scope}
\draw[\XColor,thick] (0,0) -- (0,.6) .. controls ++(90:.6cm) and ++(270:.6cm) .. (1.2,1.8) -- (1.2,3);
\draw[\phiColor,thick] (1.2,0) -- (1.2,.6) .. controls ++(90:.6cm) and ++(270:.6cm) .. (0,1.8) -- (0,2.2);
\draw[\psiColor,thick] (0,2.2) -- (0,3);
\roundNbox{unshaded}{(0,2.2)}{.3}{.1}{.1}{\scriptsize{$n_{G(c)}$}}; 
\filldraw[white] (.6,1.2) circle (.1cm);
\draw[thick] (.6,1.2) circle (.1cm); 
\node at (0,-.2) {\scriptsize{$A(X)$}};
\node at (1.2,3.2) {\scriptsize{$B(X)$}};
\node at (1.2,-.2) {\scriptsize{$\varphi_e$}};
\node at (0,3.2) {\scriptsize{$\psi_{G(c)}$}};
\node at (.1,1.2) {\scriptsize{$\varphi_{X}$}};
} 
=
\tikzmath[scale=.75, transform shape]{
\begin{scope}
\clip[rounded corners = 5] (-.6,0) rectangle (1.8,3);
\filldraw[primedregion=\AColor] (0,0) -- (0,1.2) .. controls ++(90:.4cm) and ++(-135:.2cm) .. (.6,1.8) .. controls ++(135:.2cm) and ++(270:.4cm) .. (0,2.4) -- (0,3) -- (-.6,3) -- (-.6,0); 
\filldraw[primedregion=\BColor] (1.2,0) -- (1.2,1.2) .. controls ++(90:.4cm) and ++(-45:.2cm) .. (.6,1.8) .. controls ++(-135:.2cm) and ++(90:.4cm) .. (0,1.2) -- (0,0);
\filldraw[boxregion=\AColor] (0,3) -- (0,2.4) .. controls ++(270:.4cm) and ++(135:.2cm) .. (.6,1.8) .. controls ++(45:.2cm) and ++(270:.4cm) .. (1.2,2.4) -- (1.2,3);
\filldraw[boxregion=\BColor] (1.2,0) -- (1.2,1.2) .. controls ++(90:.4cm) and ++(-45:.2cm) .. (.6,1.8) .. controls ++(45:.2cm) and ++(270:.4cm) .. (1.2,2.4) -- (1.2,3) -- (1.8,3) -- (1.8,0);
\filldraw[\AColor] (-.1,1.6) rectangle (.29,2);
\end{scope}
\draw[\XColor,thick] (0,0) -- (0,1.2) .. controls ++(90:.6cm) and ++(270:.6cm) .. (1.2,2.4) -- (1.2,3);
\draw[\phiColor,thick] (1.2,0) -- (1.2,.8);
\draw[\psiColor,thick] (1.2,.8) -- (1.2,1.2) .. controls ++(90:.6cm) and ++(270:.6cm) .. (0,2.4) -- (0,3);
\roundNbox{unshaded}{(1.2,.8)}{.3}{0}{0}{\scriptsize{$n_e$}}; 
\filldraw[white] (.6,1.8) circle (.1cm);
\draw[thick] (.6,1.8) circle (.1cm); 
\node at (0,-.2) {\scriptsize{$A(X)$}};
\node at (1.2,3.2) {\scriptsize{$B(X)$}};
\node at (1.2,-.2) {\scriptsize{$\varphi_e$}};
\node at (0,3.2) {\scriptsize{$\psi_{G(c)}$}};
\node at (.1,1.8) {\scriptsize{$\psi_X$}};
}
\qquad\qquad
\Longrightarrow
\qquad\qquad
\tikzmath{
\begin{scope}
\clip[rounded corners = 5] (-.6,-.9) rectangle (.6,.9);
\filldraw[primedregion=\BColor] (-.6,-.9) rectangle (0,.9);
\filldraw[boxregion=\BColor] (0,-.9) rectangle (.6,.9);
\end{scope}
\draw[\phiColor,thick] (0,-.9) -- (0,0);
\draw[\psiColor,thick] (0,.9) -- (0,0);
\roundNbox{unshaded}{(0,0)}{.3}{0}{0}{\scriptsize{$n_e$}}; 
\node at (0,1.1) {\scriptsize{$\psi_{e}$}};
\node at (0,-1.1) {\scriptsize{$\varphi_e$}};
}
=
\tikzmath{
\begin{scope}
\clip[rounded corners = 5] (-.9,-.9) rectangle (.9,.9);
\filldraw[primedregion=\BColor] (-.9,-.9) rectangle (.9,.9);
\filldraw[boxregion=\BColor] (.6,-.9) -- (.6,-.74) .. controls ++(90:.2cm) and ++(270:.2cm) .. (0,-.3) -- (0,.3) .. controls ++(90:.2cm) and ++(270:.2cm) .. (-.6,.74) -- (-.6,.9) -- (.9,.9) -- (.9,-.9);
\filldraw[primedregion=\AColor] (0,0) circle (.6cm);
\filldraw[boxregion=\AColor] (.3,-.52)  .. controls ++(135:.0667cm) and ++(270:.1333cm) .. (0,-.3) -- (0,.3)  .. controls ++(90:.1333cm) and ++(-45:.0667cm) .. (-.3,.52) arc (120:-60:.6cm);
\end{scope}
\draw[\XColor,thick] (0,0) circle (.6cm);
\draw[\phiColor,thick] (0,-.3) .. controls ++(270:.2cm) and ++(90:.2cm) .. (.6,-.74) -- (.6,-.9);
\draw[\psiColor,thick] (0,.3) .. controls ++(90:.2cm) and ++(270:.2cm) .. (-.6,.74) -- (-.6,.9);
\roundNbox{unshaded}{(0,0)}{.3}{.1}{.1}{\scriptsize{$n_{G(c)}$}}; %
\filldraw[white] (.3,-.52) circle (.07cm);
\draw[thick] (.3,-.52) circle (.07cm); 
\filldraw[white] (-.3,.52) circle (.07cm);
\draw[thick] (-.3,.52) circle (.07cm); 
\node at (.6,-1.1) {\scriptsize{$\varphi_e$}};
\node at (-.6,1.1) {\scriptsize{$\psi_{e}$}};
}
=
\tikzmath{
\begin{scope}
\clip[rounded corners = 5] (-.9,-.9) rectangle (.9,.9);
\filldraw[primedregion=\BColor] (-.9,-.9) rectangle (.9,.9);
\filldraw[boxregion=\BColor] (.6,-.9) -- (.6,-.74) .. controls ++(90:.2cm) and ++(270:.2cm) .. (0,-.3) -- (0,.3) .. controls ++(90:.2cm) and ++(270:.2cm) .. (-.6,.74) -- (-.6,.9) -- (.9,.9) -- (.9,-.9);
\filldraw[primedregion=\AColor] (0,0) circle (.6cm);
\filldraw[boxregion=\AColor] (.3,-.52)  .. controls ++(135:.0667cm) and ++(270:.1333cm) .. (0,-.3) -- (0,.3)  .. controls ++(90:.1333cm) and ++(-45:.0667cm) .. (-.3,.52) arc (120:-60:.6cm);
\end{scope}
\draw[\XColor,thick] (0,0) circle (.6cm);
\draw[\phiColor,thick] (0,-.3) .. controls ++(270:.2cm) and ++(90:.2cm) .. (.6,-.74) -- (.6,-.9);
\draw[\psiColor,thick] (0,.3) .. controls ++(90:.2cm) and ++(270:.2cm) .. (-.6,.74) -- (-.6,.9);
\roundNbox{unshaded}{(0,0)}{.3}{.13}{.13}{\scriptsize{$m_{G(c)}$}}; %
\filldraw[white] (.3,-.52) circle (.07cm);
\draw[thick] (.3,-.52) circle (.07cm); 
\filldraw[white] (-.3,.52) circle (.07cm);
\draw[thick] (-.3,.52) circle (.07cm); 
\node at (.6,-1.1) {\scriptsize{$\varphi_e$}};
\node at (-.6,1.1) {\scriptsize{$\psi_{e}$}};
}
=
\tikzmath{
\begin{scope}
\clip[rounded corners = 5] (-.6,-.9) rectangle (.6,.9);
\filldraw[primedregion=\BColor] (-.6,-.9) rectangle (0,.9);
\filldraw[boxregion=\BColor] (0,-.9) rectangle (.6,.9);
\end{scope}
\draw[\phiColor,thick] (0,-.9) -- (0,0);
\draw[\psiColor,thick] (0,.9) -- (0,0);
\roundNbox{unshaded}{(0,0)}{.3}{0}{0}{\scriptsize{$m_e$}}; 
\node at (0,1.1) {\scriptsize{$\psi_{e}$}};
\node at (0,-1.1) {\scriptsize{$\varphi_e$}};
}\,.
\]
Hence $m=n$, as claimed.
\end{proof}

\begin{prop}
\label{prop:DominantToFullyFaithful}
If a 2-functor $G:\cC\to \cE$ is dominant, then $-\xz G:\Fun(\cC\to \cE)\to \Fun(\cC\to \cD)$ is fully faithful on 2-morphisms.
An analogous statement holds in the $\rm C^*/W^*$ setting.
\end{prop}
\begin{proof}
It suffices to show $-\xz G$ is full on 2-morphisms.
Suppose $A,B\in\Fun(\cE\to\cD)$, $\varphi,\psi:A\Rightarrow B$, and $p:\varphi\xz G\Rrightarrow \psi\xz G$.
We show there exists $n:\varphi\Rrightarrow\psi$ such that $p=n\xz G$.

First, for each 1-cell $X\in\cE(G(c)\to G(c'))$, there exists a 1-cell $Y\in\cC(c\to c')$ such that $G(Y)\mathrel{\mathop{\rightleftarrows}\limits^{r}_{s}}X$ is a retract, i.e., $rs=1_X$.
Since $p:\varphi\xz G\Rrightarrow\psi\xz G$ is a 2-modification,
building on our graphical conventions \eqref{eq:GraphicalCalculusFaithfulOn2Mor},
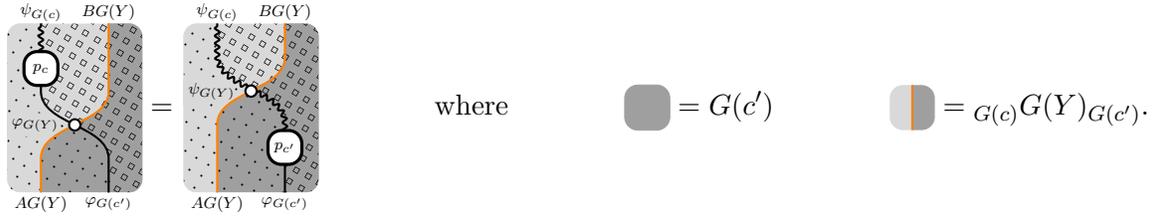
\[
\tikzmath[scale=.75, transform shape]{
\begin{scope}
\clip[rounded corners = 5] (-.6,0) rectangle (1.8,3);
\filldraw[primedregion=\AColor] (0,0) -- (0,.6) .. controls ++(90:.4cm) and ++(-135:.2cm) .. (.6,1.2) .. controls ++(135:.2cm) and ++(270:.4cm) .. (0,1.8) -- (0,3) -- (-.6,3) -- (-.6,0); 
\filldraw[primedregion=\CColor] (1.2,0) -- (1.2,.6) .. controls ++(90:.4cm) and ++(-45:.2cm) .. (.6,1.2) .. controls ++(-135:.2cm) and ++(90:.4cm) .. (0,.6) -- (0,0);
\filldraw[boxregion=\AColor] (0,3) -- (0,1.8) .. controls ++(270:.4cm) and ++(135:.2cm) .. (.6,1.2) .. controls ++(45:.2cm) and ++(270:.4cm) .. (1.2,1.8) -- (1.2,3);
\filldraw[boxregion=\CColor] (1.2,0) -- (1.2,.6) .. controls ++(90:.4cm) and ++(-45:.2cm) .. (.6,1.2) .. controls ++(45:.2cm) and ++(270:.4cm) .. (1.2,1.8) -- (1.2,3) -- (1.8,3) -- (1.8,0);
\filldraw[\AColor] (-.6,1) rectangle (.25,1.4);
\end{scope}
\draw[\YColor,thick] (0,0) -- (0,.6) .. controls ++(90:.6cm) and ++(270:.6cm) .. (1.2,1.8) -- (1.2,3);
\draw[\phiColor,thick] (1.2,0) -- (1.2,.6) .. controls ++(90:.6cm) and ++(270:.6cm) .. (0,1.8) -- (0,2.2);
\draw[\psiColor,thick] (0,2.2) -- (0,3);
\roundNbox{unshaded}{(0,2.2)}{.3}{0}{0}{\scriptsize{$p_c$}}; 
\filldraw[white] (.6,1.2) circle (.1cm);
\draw[thick] (.6,1.2) circle (.1cm); 
\node at (0,-.2) {\scriptsize{$AG(Y)$}};
\node at (1.2,3.2) {\scriptsize{$BG(Y)$}};
\node at (1.2,-.2) {\scriptsize{$\varphi_{G(c')}$}};
\node at (0,3.2) {\scriptsize{$\psi_{G(c)}$}};
\node at (-.1,1.2) {\scriptsize{$\varphi_{G(Y)}$}};
} 
=
\tikzmath[scale=.75, transform shape]{
\begin{scope}
\clip[rounded corners = 5] (-.6,0) rectangle (1.8,3);
\filldraw[primedregion=\AColor] (0,0) -- (0,1.2) .. controls ++(90:.4cm) and ++(-135:.2cm) .. (.6,1.8) .. controls ++(135:.2cm) and ++(270:.4cm) .. (0,2.4) -- (0,3) -- (-.6,3) -- (-.6,0); 
\filldraw[primedregion=\CColor] (1.2,0) -- (1.2,1.2) .. controls ++(90:.4cm) and ++(-45:.2cm) .. (.6,1.8) .. controls ++(-135:.2cm) and ++(90:.4cm) .. (0,1.2) -- (0,0);
\filldraw[boxregion=\AColor] (0,3) -- (0,2.4) .. controls ++(270:.4cm) and ++(135:.2cm) .. (.6,1.8) .. controls ++(45:.2cm) and ++(270:.4cm) .. (1.2,2.4) -- (1.2,3);
\filldraw[boxregion=\CColor] (1.2,0) -- (1.2,1.2) .. controls ++(90:.4cm) and ++(-45:.2cm) .. (.6,1.8) .. controls ++(45:.2cm) and ++(270:.4cm) .. (1.2,2.4) -- (1.2,3) -- (1.8,3) -- (1.8,0);
\filldraw[\AColor] (-.6,1.6) rectangle (.25,2);
\end{scope}
\draw[\YColor,thick] (0,0) -- (0,1.2) .. controls ++(90:.6cm) and ++(270:.6cm) .. (1.2,2.4) -- (1.2,3);
\draw[\phiColor,thick] (1.2,0) -- (1.2,.8);
\draw[\psiColor,thick] (1.2,.8) -- (1.2,1.2) .. controls ++(90:.6cm) and ++(270:.6cm) .. (0,2.4) -- (0,3);
\roundNbox{unshaded}{(1.2,.8)}{.3}{0}{0}{\scriptsize{$p_{c'}$}}; 
\filldraw[white] (.6,1.8) circle (.1cm);
\draw[thick] (.6,1.8) circle (.1cm); 
\node at (0,-.2) {\scriptsize{$AG(Y)$}};
\node at (1.2,3.2) {\scriptsize{$BG(Y)$}};
\node at (1.2,-.2) {\scriptsize{$\varphi_{G(c')}$}};
\node at (0,3.2) {\scriptsize{$\psi_{G(c)}$}};
\node at (-.1,1.8) {\scriptsize{$\psi_{G(Y)}$}};
}
\qquad\qquad
\text{where}
\qquad\qquad
\tikzmath{\filldraw[\CColor, rounded corners=5, very thin, baseline=1cm] (0,0) rectangle (.6,.6);}=G(c')
\qquad\qquad
\tikzmath{
\begin{scope}
\clip[rounded corners=5pt] (-.3,0) rectangle (.3,.6);
\fill[\AColor] (0,0) rectangle (-.3,.6);
\fill[\CColor] (0,0) rectangle (.3,.6);
\end{scope}
\draw[thick, \YColor] (0,0) -- (0,.6);
}
={}_{G(c)}G(Y)_{G(c')}.
\]
This implies that for \emph{any} $X\in \cE(G(c)\to G(c'))$ (and not just 1-cells in the image of $G$!),
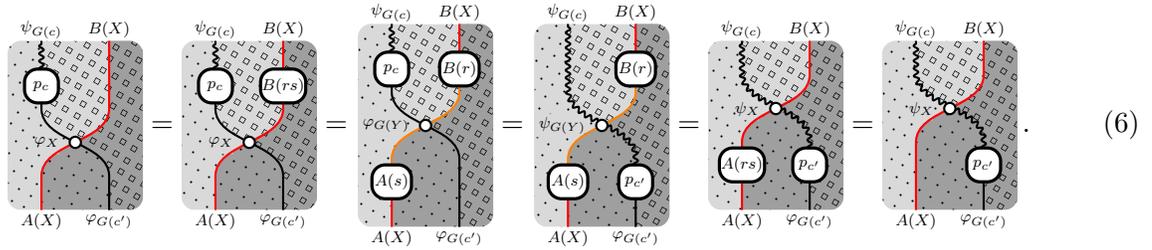
\begin{equation}
\label{eq:ModificationForAll1Cells}
\tikzmath[scale=.75, transform shape]{
\begin{scope}
\clip[rounded corners = 5] (-.6,0) rectangle (1.8,3);
\filldraw[primedregion=\AColor] (0,0) -- (0,.6) .. controls ++(90:.4cm) and ++(-135:.2cm) .. (.6,1.2) .. controls ++(135:.2cm) and ++(270:.4cm) .. (0,1.8) -- (0,3) -- (-.6,3) -- (-.6,0); 
\filldraw[primedregion=\CColor] (1.2,0) -- (1.2,.6) .. controls ++(90:.4cm) and ++(-45:.2cm) .. (.6,1.2) .. controls ++(-135:.2cm) and ++(90:.4cm) .. (0,.6) -- (0,0);
\filldraw[boxregion=\AColor] (0,3) -- (0,1.8) .. controls ++(270:.4cm) and ++(135:.2cm) .. (.6,1.2) .. controls ++(45:.2cm) and ++(270:.4cm) .. (1.2,1.8) -- (1.2,3);
\filldraw[boxregion=\CColor] (1.2,0) -- (1.2,.6) .. controls ++(90:.4cm) and ++(-45:.2cm) .. (.6,1.2) .. controls ++(45:.2cm) and ++(270:.4cm) .. (1.2,1.8) -- (1.2,3) -- (1.8,3) -- (1.8,0);
\filldraw[\AColor] (-.1,1) rectangle (.3,1.4);
\end{scope}
\draw[\XColor,thick] (0,0) -- (0,.6) .. controls ++(90:.6cm) and ++(270:.6cm) .. (1.2,1.8) -- (1.2,3);
\draw[\phiColor,thick] (1.2,0) -- (1.2,.6) .. controls ++(90:.6cm) and ++(270:.6cm) .. (0,1.8) -- (0,2.2);
\draw[\psiColor,thick] (0,2.2) -- (0,3);
\roundNbox{unshaded}{(0,2.2)}{.3}{0}{0}{\scriptsize{$p_c$}}; 
\filldraw[white] (.6,1.2) circle (.1cm);
\draw[thick] (.6,1.2) circle (.1cm); 
\node at (0,-.2) {\scriptsize{$A(X)$}};
\node at (1.2,3.2) {\scriptsize{$B(X)$}};
\node at (1.2,-.2) {\scriptsize{$\varphi_{G(c')}$}};
\node at (0,3.2) {\scriptsize{$\psi_{G(c)}$}};
\node at (.1,1.2) {\scriptsize{$\varphi_{X}$}};
} 
=
\tikzmath[scale=.75, transform shape]{
\begin{scope}
\clip[rounded corners = 5] (-.6,0) rectangle (1.8,3);
\filldraw[primedregion=\AColor] (0,0) -- (0,.6) .. controls ++(90:.4cm) and ++(-135:.2cm) .. (.6,1.2) .. controls ++(135:.2cm) and ++(270:.4cm) .. (0,1.8) -- (0,3) -- (-.6,3) -- (-.6,0); 
\filldraw[primedregion=\CColor] (1.2,0) -- (1.2,.6) .. controls ++(90:.4cm) and ++(-45:.2cm) .. (.6,1.2) .. controls ++(-135:.2cm) and ++(90:.4cm) .. (0,.6) -- (0,0);
\filldraw[boxregion=\AColor] (0,3) -- (0,1.8) .. controls ++(270:.4cm) and ++(135:.2cm) .. (.6,1.2) .. controls ++(45:.2cm) and ++(270:.4cm) .. (1.2,1.8) -- (1.2,3);
\filldraw[boxregion=\CColor] (1.2,0) -- (1.2,.6) .. controls ++(90:.4cm) and ++(-45:.2cm) .. (.6,1.2) .. controls ++(45:.2cm) and ++(270:.4cm) .. (1.2,1.8) -- (1.2,3) -- (1.8,3) -- (1.8,0);
\filldraw[\AColor] (-.1,1) rectangle (.3,1.4);
\end{scope}
\draw[\XColor,thick] (0,0) -- (0,.6) .. controls ++(90:.6cm) and ++(270:.6cm) .. (1.2,1.8) -- (1.2,3);
\draw[\phiColor,thick] (1.2,0) -- (1.2,.6) .. controls ++(90:.6cm) and ++(270:.6cm) .. (0,1.8) -- (0,2.2);
\draw[\psiColor,thick] (0,2.2) -- (0,3);
\roundNbox{unshaded}{(0,2.2)}{.3}{0}{0}{\scriptsize{$p_c$}}; 
\roundNbox{unshaded}{(1.2,2.2)}{.3}{.1}{.1}{\scriptsize{$B(rs)$}}; 
\filldraw[white] (.6,1.2) circle (.1cm);
\draw[thick] (.6,1.2) circle (.1cm); 
\node at (0,-.2) {\scriptsize{$A(X)$}};
\node at (1.2,3.2) {\scriptsize{$B(X)$}};
\node at (1.2,-.2) {\scriptsize{$\varphi_{G(c')}$}};
\node at (0,3.2) {\scriptsize{$\psi_{G(c)}$}};
\node at (.1,1.2) {\scriptsize{$\varphi_{X}$}};
} 
=
\tikzmath[scale=.75, transform shape]{
\begin{scope}
\clip[rounded corners = 5] (-.6,-.6) rectangle (1.8,3);
\filldraw[primedregion=\AColor] (0,-.6) -- (0,.6) .. controls ++(90:.4cm) and ++(-135:.2cm) .. (.6,1.2) .. controls ++(135:.2cm) and ++(270:.4cm) .. (0,1.8) -- (0,3) -- (-.6,3) -- (-.6,-.6); 
\filldraw[primedregion=\CColor] (1.2,-.6) -- (1.2,.6) .. controls ++(90:.4cm) and ++(-45:.2cm) .. (.6,1.2) .. controls ++(-135:.2cm) and ++(90:.4cm) .. (0,.6) -- (0,-.6);
\filldraw[boxregion=\AColor] (0,3) -- (0,1.8) .. controls ++(270:.4cm) and ++(135:.2cm) .. (.6,1.2) .. controls ++(45:.2cm) and ++(270:.4cm) .. (1.2,1.8) -- (1.2,3);
\filldraw[boxregion=\CColor] (1.2,-.6) -- (1.2,.6) .. controls ++(90:.4cm) and ++(-45:.2cm) .. (.6,1.2) .. controls ++(45:.2cm) and ++(270:.4cm) .. (1.2,1.8) -- (1.2,3) -- (1.8,3) -- (1.8,-.6);
\filldraw[\AColor] (-.6,1) rectangle (.25,1.4);
\end{scope}
\draw[\XColor,thick] (0,-.6) -- (0,.2);
\draw[\XColor,thick] (1.2,2.2) -- (1.2,3);
\draw[\YColor,thick] (0,.2) -- (0,.6) .. controls ++(90:.6cm) and ++(270:.6cm) .. (1.2,1.8) -- (1.2,2.2);
\draw[\phiColor,thick] (1.2,-.6) -- (1.2,.6) .. controls ++(90:.6cm) and ++(270:.6cm) .. (0,1.8) -- (0,2.2);
\draw[\psiColor,thick] (0,2.2) -- (0,3);
\roundNbox{unshaded}{(0,2.2)}{.3}{0}{0}{\scriptsize{$p_c$}}; 
\roundNbox{unshaded}{(1.2,2.2)}{.3}{.05}{.05}{\scriptsize{$B(r)$}}; %
\roundNbox{unshaded}{(0,.2)}{.3}{.05}{.05}{\scriptsize{$A(s)$}}; %
\filldraw[white] (.6,1.2) circle (.1cm);
\draw[thick] (.6,1.2) circle (.1cm); 
\node at (0,-.8) {\scriptsize{$A(X)$}};
\node at (1.2,3.2) {\scriptsize{$B(X)$}};
\node at (1.2,-.8) {\scriptsize{$\varphi_{G(c')}$}};
\node at (0,3.2) {\scriptsize{$\psi_{G(c)}$}};
\node at (-.1,1.2) {\scriptsize{$\varphi_{G(Y)}$}};
} 
=
\tikzmath[scale=.75, transform shape]{
\begin{scope}
\clip[rounded corners = 5] (-.6,0) rectangle (1.8,3.6);
\filldraw[primedregion=\AColor] (0,0) -- (0,1.2) .. controls ++(90:.4cm) and ++(-135:.2cm) .. (.6,1.8) .. controls ++(135:.2cm) and ++(270:.4cm) .. (0,2.4) -- (0,3.6) -- (-.6,3.6) -- (-.6,0); 
\filldraw[primedregion=\CColor] (1.2,0) -- (1.2,1.2) .. controls ++(90:.4cm) and ++(-45:.2cm) .. (.6,1.8) .. controls ++(-135:.2cm) and ++(90:.4cm) .. (0,1.2) -- (0,0);
\filldraw[boxregion=\AColor] (0,3.6) -- (0,2.4) .. controls ++(270:.4cm) and ++(135:.2cm) .. (.6,1.8) .. controls ++(45:.2cm) and ++(270:.4cm) .. (1.2,2.4) -- (1.2,3.6);
\filldraw[boxregion=\CColor] (1.2,0) -- (1.2,1.2) .. controls ++(90:.4cm) and ++(-45:.2cm) .. (.6,1.8) .. controls ++(45:.2cm) and ++(270:.4cm) .. (1.2,2.4) -- (1.2,3.6) -- (1.8,3.6) -- (1.8,0);
\filldraw[\AColor] (-.6,1.6) rectangle (.25,2);
\end{scope}
\draw[\XColor,thick] (0,0) -- (0,.8);
\draw[\XColor,thick] (1.2,2.8) -- (1.2,3.6);
\draw[\YColor,thick] (0,.8) -- (0,1.2) .. controls ++(90:.6cm) and ++(270:.6cm) .. (1.2,2.4) -- (1.2,2.8);
\draw[\phiColor,thick] (1.2,0) -- (1.2,.8);
\draw[\psiColor,thick] (1.2,.8) -- (1.2,1.2) .. controls ++(90:.6cm) and ++(270:.6cm) .. (0,2.4) -- (0,3.6);
\roundNbox{unshaded}{(1.2,.8)}{.3}{0}{0}{\scriptsize{$p_{c'}$}}; 
\roundNbox{unshaded}{(1.2,2.8)}{.3}{.05}{.05}{\scriptsize{$B(r)$}}; %
\roundNbox{unshaded}{(0,.8)}{.3}{.05}{.05}{\scriptsize{$A(s)$}}; %
\filldraw[white] (.6,1.8) circle (.1cm);
\draw[thick] (.6,1.8) circle (.1cm); 
\node at (0,-.2) {\scriptsize{$A(X)$}};
\node at (1.2,3.8) {\scriptsize{$B(X)$}};
\node at (1.2,-.2) {\scriptsize{$\varphi_{G(c')}$}};
\node at (0,3.8) {\scriptsize{$\psi_{G(c)}$}};
\node at (-.1,1.8) {\scriptsize{$\psi_{G(Y)}$}};
}
=
\tikzmath[scale=.75, transform shape]{
\begin{scope}
\clip[rounded corners = 5] (-.6,0) rectangle (1.8,3);
\filldraw[primedregion=\AColor] (0,0) -- (0,1.2) .. controls ++(90:.4cm) and ++(-135:.2cm) .. (.6,1.8) .. controls ++(135:.2cm) and ++(270:.4cm) .. (0,2.4) -- (0,3) -- (-.6,3) -- (-.6,0); 
\filldraw[primedregion=\CColor] (1.2,0) -- (1.2,1.2) .. controls ++(90:.4cm) and ++(-45:.2cm) .. (.6,1.8) .. controls ++(-135:.2cm) and ++(90:.4cm) .. (0,1.2) -- (0,0);
\filldraw[boxregion=\AColor] (0,3) -- (0,2.4) .. controls ++(270:.4cm) and ++(135:.2cm) .. (.6,1.8) .. controls ++(45:.2cm) and ++(270:.4cm) .. (1.2,2.4) -- (1.2,3);
\filldraw[boxregion=\CColor] (1.2,0) -- (1.2,1.2) .. controls ++(90:.4cm) and ++(-45:.2cm) .. (.6,1.8) .. controls ++(45:.2cm) and ++(270:.4cm) .. (1.2,2.4) -- (1.2,3) -- (1.8,3) -- (1.8,0);
\filldraw[\AColor] (-.1,1.6) rectangle (.29,2);
\end{scope}
\draw[\XColor,thick] (0,0) -- (0,1.2) .. controls ++(90:.6cm) and ++(270:.6cm) .. (1.2,2.4) -- (1.2,3);
\draw[\phiColor,thick] (1.2,0) -- (1.2,.8);
\draw[\psiColor,thick] (1.2,.8) -- (1.2,1.2) .. controls ++(90:.6cm) and ++(270:.6cm) .. (0,2.4) -- (0,3);
\roundNbox{unshaded}{(1.2,.8)}{.3}{0}{0}{\scriptsize{$p_{c'}$}}; 
\roundNbox{unshaded}{(0,.8)}{.3}{.1}{.1}{\scriptsize{$A(rs)$}}; %
\filldraw[white] (.6,1.8) circle (.1cm);
\draw[thick] (.6,1.8) circle (.1cm); 
\node at (0,-.2) {\scriptsize{$A(X)$}};
\node at (1.2,3.2) {\scriptsize{$B(X)$}};
\node at (1.2,-.2) {\scriptsize{$\varphi_{G(c')}$}};
\node at (0,3.2) {\scriptsize{$\psi_{G(c)}$}};
\node at (.1,1.8) {\scriptsize{$\psi_{X}$}};
}
=
\tikzmath[scale=.75, transform shape]{
\begin{scope}
\clip[rounded corners = 5] (-.6,0) rectangle (1.8,3);
\filldraw[primedregion=\AColor] (0,0) -- (0,1.2) .. controls ++(90:.4cm) and ++(-135:.2cm) .. (.6,1.8) .. controls ++(135:.2cm) and ++(270:.4cm) .. (0,2.4) -- (0,3) -- (-.6,3) -- (-.6,0); 
\filldraw[primedregion=\CColor] (1.2,0) -- (1.2,1.2) .. controls ++(90:.4cm) and ++(-45:.2cm) .. (.6,1.8) .. controls ++(-135:.2cm) and ++(90:.4cm) .. (0,1.2) -- (0,0);
\filldraw[boxregion=\AColor] (0,3) -- (0,2.4) .. controls ++(270:.4cm) and ++(135:.2cm) .. (.6,1.8) .. controls ++(45:.2cm) and ++(270:.4cm) .. (1.2,2.4) -- (1.2,3);
\filldraw[boxregion=\CColor] (1.2,0) -- (1.2,1.2) .. controls ++(90:.4cm) and ++(-45:.2cm) .. (.6,1.8) .. controls ++(45:.2cm) and ++(270:.4cm) .. (1.2,2.4) -- (1.2,3) -- (1.8,3) -- (1.8,0);
\filldraw[\AColor] (-.1,1.6) rectangle (.29,2);
\end{scope}
\draw[\XColor,thick] (0,0) -- (0,1.2) .. controls ++(90:.6cm) and ++(270:.6cm) .. (1.2,2.4) -- (1.2,3);
\draw[\phiColor,thick] (1.2,0) -- (1.2,.8);
\draw[\psiColor,thick] (1.2,.8) -- (1.2,1.2) .. controls ++(90:.6cm) and ++(270:.6cm) .. (0,2.4) -- (0,3);
\roundNbox{unshaded}{(1.2,.8)}{.3}{0}{0}{\scriptsize{$p_{c'}$}}; 
\filldraw[white] (.6,1.8) circle (.1cm);
\draw[thick] (.6,1.8) circle (.1cm); 
\node at (0,-.2) {\scriptsize{$A(X)$}};
\node at (1.2,3.2) {\scriptsize{$B(X)$}};
\node at (1.2,-.2) {\scriptsize{$\varphi_{G(c')}$}};
\node at (0,3.2) {\scriptsize{$\psi_{G(c)}$}};
\node at (.1,1.8) {\scriptsize{$\psi_{X}$}};
}\,.
\end{equation}

Next we construct $n:\varphi\Rrightarrow \psi$ such that $p=n\xz G$.
For each $c\in \cC$, we define $n_{G(c)}:=p_c$ so that $p_c=(n\xz G)_c$, and $p=n\xz G$, provided we can extend $n$ to a modification.
For each $e\in \cE$, there exists a 0-cell $c\in\cC$ and a condensation $X:G(c)\condense e$.
We define $n_e$ as follows
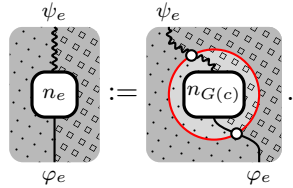
\[
\tikzmath{
\begin{scope}
\clip[rounded corners = 5] (-.6,-.9) rectangle (.6,.9);
\filldraw[primedregion=\BColor] (-.6,-.9) rectangle (0,.9);
\filldraw[boxregion=\BColor] (0,-.9) rectangle (.6,.9);
\end{scope}
\draw[\phiColor,thick] (0,-.9) -- (0,0);
\draw[\psiColor,thick] (0,.9) -- (0,0);
\roundNbox{unshaded}{(0,0)}{.3}{0}{0}{\scriptsize{$n_e$}}; 
\node at (0,1.1) {\scriptsize{$\psi_{e}$}};
\node at (0,-1.1) {\scriptsize{$\varphi_e$}};
}
:=
\tikzmath{
\begin{scope}
\clip[rounded corners = 5] (-.9,-.9) rectangle (.9,.9);
\filldraw[primedregion=\BColor] (-.9,-.9) rectangle (.9,.9);
\filldraw[boxregion=\BColor] (.6,-.9) -- (.6,-.74) .. controls ++(90:.2cm) and ++(270:.2cm) .. (0,-.3) -- (0,.3) .. controls ++(90:.2cm) and ++(270:.2cm) .. (-.6,.74) -- (-.6,.9) -- (.9,.9) -- (.9,-.9);
\filldraw[primedregion=\AColor] (0,0) circle (.6cm);
\filldraw[boxregion=\AColor] (.3,-.52)  .. controls ++(135:.0667cm) and ++(270:.1333cm) .. (0,-.3) -- (0,.3)  .. controls ++(90:.1333cm) and ++(-45:.0667cm) .. (-.3,.52) arc (120:-60:.6cm);
\end{scope}
\draw[\XColor,thick] (0,0) circle (.6cm);
\draw[\phiColor,thick] (0,-.3) .. controls ++(270:.2cm) and ++(90:.2cm) .. (.6,-.74) -- (.6,-.9);
\draw[\psiColor,thick] (0,.3) .. controls ++(90:.2cm) and ++(270:.2cm) .. (-.6,.74) -- (-.6,.9);
\roundNbox{unshaded}{(0,0)}{.3}{.1}{.1}{\scriptsize{$n_{G(c)}$}}; %
\filldraw[white] (.3,-.52) circle (.07cm);
\draw[thick] (.3,-.52) circle (.07cm); 
\filldraw[white] (-.3,.52) circle (.07cm);
\draw[thick] (-.3,.52) circle (.07cm); 
\node at (.6,-1.1) {\scriptsize{$\varphi_e$}};
\node at (-.6,1.1) {\scriptsize{$\psi_{e}$}};
}\,.
\]
We prove $n$ is a 2-modification $\varphi\Rrightarrow\psi$. 
Suppose $e'\in\cE$ is a 0-cell and $Z\in\cE(e\to e')$ is a 1-cell. 
Let $X':G(c')\condense e'$ be a condensation for some 0-cell $c'\in\cC$.
Using the graphical conventions
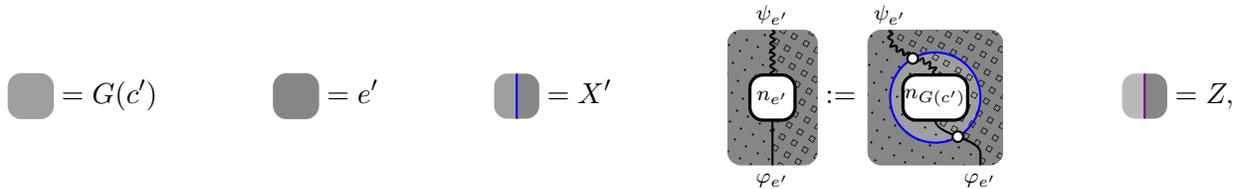
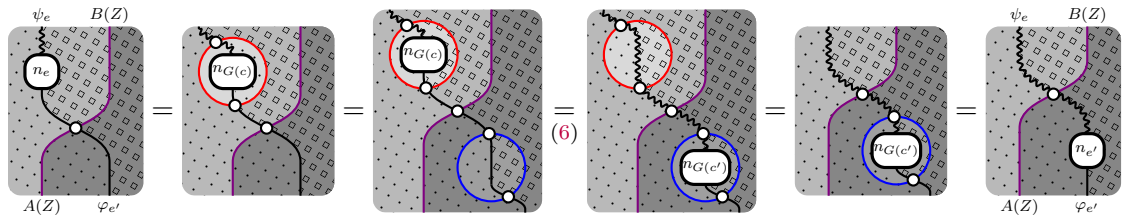
\[
\tikzmath{\filldraw[\CColor, rounded corners=5, very thin, baseline=1cm] (0,0) rectangle (.6,.6);}=G(c')
\qquad\qquad
\tikzmath{\filldraw[\DColor, rounded corners=5, very thin, baseline=1cm] (0,0) rectangle (.6,.6);}=e'
\qquad\qquad
\tikzmath{
\begin{scope}
\clip[rounded corners=5pt] (-.3,0) rectangle (.3,.6);
\fill[\CColor] (0,0) rectangle (-.3,.6);
\fill[\DColor] (0,0) rectangle (.3,.6);
\end{scope}
\draw[thick, \WColor] (0,0) -- (0,.6);
}=X'
\qquad\qquad
\tikzmath{
\begin{scope}
\clip[rounded corners = 5] (-.6,-.9) rectangle (.6,.9);
\filldraw[primedregion=\DColor] (-.6,-.9) rectangle (0,.9);
\filldraw[boxregion=\DColor] (0,-.9) rectangle (.6,.9);
\end{scope}
\draw[\phiColor,thick] (0,-.9) -- (0,0);
\draw[\psiColor,thick] (0,.9) -- (0,0);
\roundNbox{unshaded}{(0,0)}{.3}{0}{0}{\scriptsize{$n_{e'}$}}; 
\node at (0,1.1) {\scriptsize{$\psi_{e'}$}};
\node at (0,-1.1) {\scriptsize{$\varphi_{e'}$}};
}
:=
\tikzmath{
\begin{scope}
\clip[rounded corners = 5] (-.9,-.9) rectangle (.9,.9);
\filldraw[primedregion=\DColor] (-.9,-.9) rectangle (.9,.9);
\filldraw[boxregion=\DColor] (.6,-.9) -- (.6,-.74) .. controls ++(90:.2cm) and ++(270:.2cm) .. (0,-.3) -- (0,.3) .. controls ++(90:.2cm) and ++(270:.2cm) .. (-.6,.74) -- (-.6,.9) -- (.9,.9) -- (.9,-.9);
\filldraw[primedregion=\CColor] (0,0) circle (.6cm);
\filldraw[boxregion=\CColor] (.3,-.52) .. controls ++(135:.0667cm) and ++(270:.1333cm) .. (0,-.3) -- (0,.3) .. controls ++(90:.1333cm) and ++(-45:.0667cm) .. (-.3,.52) arc (120:-60:.6cm);
\end{scope}
\draw[\WColor,thick] (0,0) circle (.6cm);
\draw[\phiColor,thick] (0,-.3) .. controls ++(270:.2cm) and ++(90:.2cm) .. (.6,-.74) -- (.6,-.9);
\draw[\psiColor,thick] (0,.3) .. controls ++(90:.2cm) and ++(270:.2cm) .. (-.6,.74) -- (-.6,.9);
\roundNbox{unshaded}{(0,0)}{.3}{.13}{.13}{\scriptsize{$n_{G(c')}$}}; %
\filldraw[white] (.3,-.52) circle (.07cm);
\draw[thick] (.3,-.52) circle (.07cm); 
\filldraw[white] (-.3,.52) circle (.07cm);
\draw[thick] (-.3,.52) circle (.07cm); 
\node at (.6,-1.1) {\scriptsize{$\varphi_{e'}$}};
\node at (-.6,1.1) {\scriptsize{$\psi_{e'}$}};
}
\qquad\qquad
\tikzmath{
\begin{scope}
\clip[rounded corners=5pt] (-.3,0) rectangle (.3,.6);
\fill[\BColor] (0,0) rectangle (-.3,.6);
\fill[\DColor] (0,0) rectangle (.3,.6);
\end{scope}
\draw[thick, \ZColor] (0,0) -- (0,.6);
}=Z,
\]
we see that
\[
\tikzmath[scale=.75, transform shape]{
\begin{scope}
\clip[rounded corners = 5] (-.6,0) rectangle (1.8,3);
\filldraw[primedregion=\BColor] (0,0) -- (0,.6) .. controls ++(90:.4cm) and ++(-135:.2cm) .. (.6,1.2) .. controls ++(135:.2cm) and ++(270:.4cm) .. (0,1.8) -- (0,3) -- (-.6,3) -- (-.6,0); 
\filldraw[primedregion=\DColor] (1.2,0) -- (1.2,.6) .. controls ++(90:.4cm) and ++(-45:.2cm) .. (.6,1.2) .. controls ++(-135:.2cm) and ++(90:.4cm) .. (0,.6) -- (0,0);
\filldraw[boxregion=\BColor] (0,3) -- (0,1.8) .. controls ++(270:.4cm) and ++(135:.2cm) .. (.6,1.2) .. controls ++(45:.2cm) and ++(270:.4cm) .. (1.2,1.8) -- (1.2,3);
\filldraw[boxregion=\DColor] (1.2,0) -- (1.2,.6) .. controls ++(90:.4cm) and ++(-45:.2cm) .. (.6,1.2) .. controls ++(45:.2cm) and ++(270:.4cm) .. (1.2,1.8) -- (1.2,3) -- (1.8,3) -- (1.8,0);
\end{scope}
\draw[\ZColor,thick] (0,0) -- (0,.6) .. controls ++(90:.6cm) and ++(270:.6cm) .. (1.2,1.8) -- (1.2,3);
\draw[\phiColor,thick] (1.2,0) -- (1.2,.6) .. controls ++(90:.6cm) and ++(270:.6cm) .. (0,1.8) -- (0,2.2);
\draw[\psiColor,thick] (0,2.2) -- (0,3);
\roundNbox{unshaded}{(0,2.2)}{.3}{0}{0}{\scriptsize{$n_e$}}; 
\filldraw[white] (.6,1.2) circle (.1cm);
\draw[thick] (.6,1.2) circle (.1cm); 
\node at (0,-.2) {\scriptsize{$A(Z)$}};
\node at (1.2,3.2) {\scriptsize{$B(Z)$}};
\node at (1.2,-.2) {\scriptsize{$\varphi_{e'}$}};
\node at (0,3.2) {\scriptsize{$\psi_{e}$}};
} 
=
\tikzmath[scale=.75, transform shape]{
\begin{scope}
\clip[rounded corners = 5] (-.9,0) rectangle (1.8,3);
\filldraw[primedregion=\BColor] (0,0) -- (0,.6) .. controls ++(90:.6cm) and ++(270:.6cm) .. (1.2,1.8) -- (1.2,3) -- (-.9,3) -- (-.9,0); 
\filldraw[primedregion=\DColor] (0,0) -- (0,.6) .. controls ++(90:.6cm) and ++(270:.6cm) .. (1.2,1.8) -- (1.2,3) -- (1.8,3) -- (1.8,0);
\filldraw[boxregion=\BColor] (-.6,3) -- (-.6,2.94) .. controls ++(270:.2cm) and ++(90:.2cm) .. (0,2.5) -- (0,1.8) .. controls ++(270:.4cm) and ++(135:.2cm) .. (.6,1.2) .. controls ++(45:.2cm) and ++(270:.4cm) .. (1.2,1.8) -- (1.2,3);
\filldraw[boxregion=\DColor] (1.2,0) -- (1.2,.6) .. controls ++(90:.4cm) and ++(-45:.2cm) .. (.6,1.2) .. controls ++(45:.2cm) and ++(270:.4cm) .. (1.2,1.8) -- (1.2,3) -- (1.8,3) -- (1.8,0);
\filldraw[\AColor] (0,2.2) circle (.6cm);
\filldraw[primedregion=\AColor] (.04,1.6) -- (0,2.5) .. controls ++(90:.1333cm) and ++(-45:.0667cm) .. (-.3,2.72) arc (120:-86:.6cm);
\end{scope}
\draw[\ZColor,thick] (0,0) -- (0,.6) .. controls ++(90:.6cm) and ++(270:.6cm) .. (1.2,1.8) -- (1.2,3);
\draw[\phiColor,thick] (1.2,0) -- (1.2,.6) .. controls ++(90:.6cm) and ++(270:.6cm) .. (0,1.8) -- (0,2.2);
\draw[\psiColor,thick] (0,2.5) .. controls ++(90:.2cm) and ++(270:.2cm) .. (-.6,2.94) -- (-.6,3);
\draw[thick,\XColor] (0,2.2) circle (.6cm);
\roundNbox{unshaded}{(0,2.2)}{.3}{.1}{.1}{\scriptsize{$n_{G(c)}$}}; %
\filldraw[white] (.6,1.2) circle (.1cm);
\draw[thick] (.6,1.2) circle (.1cm); 
\filldraw[white] (.04,1.6) circle (.1cm);
\draw[thick] (.04,1.6) circle (.1cm); 
\filldraw[white] (-.3,2.72) circle (.1cm);
\draw[thick] (-.3,2.72) circle (.1cm); 
} 
=
\tikzmath[scale=.75, transform shape]{
\begin{scope}
\clip[rounded corners = 5] (-.9,-.6) rectangle (2.1,3);
\filldraw[primedregion=\BColor] (0,-.6) -- (0,.6) .. controls ++(90:.6cm) and ++(270:.6cm) .. (1.2,1.8) -- (1.2,3) -- (-.9,3) -- (-.9,-.6); 
\filldraw[primedregion=\DColor] (0,-.6) -- (0,.6) .. controls ++(90:.6cm) and ++(270:.6cm) .. (1.2,1.8) -- (1.2,3) -- (2.1,3) -- (2.1,-.6);
\filldraw[boxregion=\BColor] (-.6,3) -- (-.6,2.94) .. controls ++(270:.2cm) and ++(90:.2cm) .. (0,2.5) -- (0,1.8) .. controls ++(270:.4cm) and ++(135:.2cm) .. (.6,1.2) .. controls ++(45:.2cm) and ++(270:.4cm) .. (1.2,1.8) -- (1.2,3);
\filldraw[boxregion=\DColor] (1.8,-.6) -- (1.8,-.54) .. controls ++(90:.2cm) and ++(270:.2cm) .. (1.2,-.1) -- (1.2,.6) .. controls ++(90:.4cm) and ++(-45:.2cm) .. (.6,1.2) .. controls ++(45:.2cm) and ++(270:.4cm) .. (1.2,1.8) -- (1.2,3) -- (2.1,3) -- (2.1,0);
\filldraw[primedregion=\AColor] (0,2.2) circle (.6cm);
\filldraw[boxregion=\AColor] (.04,1.6) -- (0,2.5) .. controls ++(90:.1333cm) and ++(-45:.0667cm) .. (-.3,2.72) arc (120:-86:.6cm);
\filldraw[primedregion=\CColor] (1.2,.2) circle (.6cm);
\filldraw[boxregion=\CColor] (1.5,-.32) .. controls ++(135:.0667cm) and ++(270:.1333cm) .. (1.2,-.1) -- (1.16,.8) arc (94:-60:.6cm);
\end{scope}
\draw[\ZColor,thick] (0,-.6) -- (0,.6) .. controls ++(90:.6cm) and ++(270:.6cm) .. (1.2,1.8) -- (1.2,3);
\draw[\phiColor,thick] (1.2,-.1) -- (1.2,.6) .. controls ++(90:.6cm) and ++(270:.6cm) .. (0,1.8) -- (0,2.2);
\draw[\phiColor,thick] (1.2,-.1) .. controls ++(270:.2cm) and ++(90:.2cm) .. (1.8,-.54) -- (1.8,-.6);
\draw[\psiColor,thick] (0,2.5) .. controls ++(90:.2cm) and ++(270:.2cm) .. (-.6,2.94) -- (-.6,3);
\draw[thick,\XColor] (0,2.2) circle (.6cm);
\draw[thick,\WColor] (1.2,.2) circle (.6cm);
\roundNbox{unshaded}{(0,2.2)}{.3}{.1}{.1}{\scriptsize{$n_{G(c)}$}}; %
\filldraw[white] (.6,1.2) circle (.1cm);
\draw[thick] (.6,1.2) circle (.1cm); 
\filldraw[white] (.04,1.61) circle (.1cm);
\draw[thick] (.04,1.61) circle (.1cm); 
\filldraw[white] (-.3,2.72) circle (.1cm);
\draw[thick] (-.3,2.72) circle (.1cm); 
\filldraw[white] (1.16,.8) circle (.1cm);
\draw[thick] (1.16,.8) circle (.1cm); 
\filldraw[white] (1.5,-.32) circle (.1cm);
\draw[thick] (1.5,-.32) circle (.1cm); 
} 
\underset{\text{\eqref{eq:ModificationForAll1Cells}}}{=}
\tikzmath[scale=.75, transform shape]{
\begin{scope}
\clip[rounded corners = 5] (-.9,0) rectangle (2.1,3.6);
\filldraw[primedregion=\BColor] (0,0) -- (0,1.2) .. controls ++(90:.6cm) and ++(270:.6cm) .. (1.2,2.4) -- (1.2,3.6) -- (-.9,3.6) -- (-.9,0); 
\filldraw[primedregion=\DColor] (0,0) -- (0,1.2) .. controls ++(90:.6cm) and ++(270:.6cm) .. (1.2,2.4) -- (1.2,3.6) -- (2.1,3.6) -- (2.1,0);
\filldraw[boxregion=\BColor]  (-.6,3.6) -- (-.6,3.54) .. controls ++(270:.2cm) and ++(90:.2cm) .. (0,3.1) -- (0,2.4) .. controls ++(270:.4cm) and ++(135:.2cm) .. (.6,1.8) .. controls ++(45:.2cm) and ++(270:.4cm) .. (1.2,2.4) -- (1.2,3.6);
\filldraw[boxregion=\DColor] (1.8,0) -- (1.8,.06) .. controls ++(90:.2cm) and ++(270:.2cm) .. (1.2,.5) -- (1.2,1.2) .. controls ++(90:.4cm) and ++(-45:.2cm) .. (.6,1.8) .. controls ++(45:.2cm) and ++(270:.4cm) .. (1.2,2.4) -- (1.2,3.6) -- (2.1,3.6) -- (2.1,0);
\filldraw[primedregion=\AColor] (0,2.8) circle (.6cm);
\filldraw[boxregion=\AColor] (.04,2.2) -- (0,3.1) .. controls ++(90:.1333cm) and ++(-45:.0667cm) .. (-.3,3.32) arc (120:-86:.6cm);
\filldraw[primedregion=\CColor] (1.2,.8) circle (.6cm);
\filldraw[boxregion=\CColor] (1.5,.28) .. controls ++(135:.0667cm) and ++(270:.1333cm) .. (1.2,.5) -- (1.16,1.4) arc (94:-60:.6cm);
\end{scope}
\draw[\ZColor,thick] (0,0) -- (0,1.2) .. controls ++(90:.6cm) and ++(270:.6cm) .. (1.2,2.4) -- (1.2,3.6);
\draw[\phiColor,thick] (1.2,.5) .. controls ++(270:.2cm) and ++(90:.2cm) .. (1.8,.06) -- (1.8,0);
\draw[\psiColor,thick] (1.2,.8) -- (1.2,1.2) .. controls ++(90:.6cm) and ++(270:.6cm) .. (0,2.4) -- (0,3.1) .. controls ++(90:.2cm) and ++(270:.2cm) .. (-.6,3.54) -- (-.6,3.6);
\draw[thick,\XColor] (0,2.8) circle (.6cm);
\draw[thick,\WColor] (1.2,.8) circle (.6cm);
\roundNbox{unshaded}{(1.2,.8)}{.3}{.13}{.13}{\scriptsize{$n_{G(c')}$}}; 
\filldraw[white] (1.5,.28) circle (.1cm);
\draw[thick] (1.5,.28) circle (.1cm); 
\filldraw[white] (1.16,1.4) circle (.1cm);
\draw[thick] (1.16,1.4) circle (.1cm); 
\filldraw[white] (.6,1.8) circle (.1cm);
\draw[thick] (.6,1.8) circle (.1cm); 
\filldraw[white] (.04,2.21) circle (.1cm);
\draw[thick] (.04,2.21) circle (.1cm); 
\filldraw[white] (-.3,3.32) circle (.1cm);
\draw[thick] (-.3,3.32) circle (.1cm); 
}
=
\tikzmath[scale=.75, transform shape]{
\begin{scope}
\clip[rounded corners = 5] (-.6,0) rectangle (2.1,3);
\filldraw[primedregion=\BColor] (0,0) -- (0,1.2) .. controls ++(90:.6cm) and ++(270:.6cm) .. (1.2,2.4) -- (1.2,3) -- (-.6,3) -- (-.6,0); 
\filldraw[primedregion=\DColor] (0,0) -- (0,1.2) .. controls ++(90:.6cm) and ++(270:.6cm) .. (1.2,2.4) -- (1.2,3) -- (2.1,3) -- (2.1,0);
\filldraw[boxregion=\BColor] (0,3) -- (0,2.4) .. controls ++(270:.4cm) and ++(135:.2cm) .. (.6,1.8) .. controls ++(45:.2cm) and ++(270:.4cm) .. (1.2,2.4) -- (1.2,3);
\filldraw[boxregion=\DColor] (1.8,0) -- (1.8,.06) .. controls ++(90:.2cm) and ++(270:.2cm) .. (1.2,.5) -- (1.2,1.2) .. controls ++(90:.4cm) and ++(-45:.2cm) .. (.6,1.8) .. controls ++(45:.2cm) and ++(270:.4cm) .. (1.2,2.4) -- (1.2,3) -- (2.1,3) -- (2.1,0);
\filldraw[primedregion=\CColor] (1.2,.8) circle (.6cm);
\filldraw[boxregion=\CColor] (1.5,.28) .. controls ++(135:.0667cm) and ++(270:.1333cm) .. (1.2,.5) -- (1.16,1.4) arc (94:-60:.6cm);
\end{scope}
\draw[\ZColor,thick] (0,0) -- (0,1.2) .. controls ++(90:.6cm) and ++(270:.6cm) .. (1.2,2.4) -- (1.2,3);
\draw[\phiColor,thick] (1.2,.5) .. controls ++(270:.2cm) and ++(90:.2cm) .. (1.8,.06) -- (1.8,0);
\draw[\psiColor,thick] (1.2,.8) -- (1.2,1.2) .. controls ++(90:.6cm) and ++(270:.6cm) .. (0,2.4) -- (0,3);
\draw[thick,\WColor] (1.2,.8) circle (.6cm);
\roundNbox{unshaded}{(1.2,.8)}{.3}{.13}{.13}{\scriptsize{$n_{G(c')}$}}; 
\filldraw[white] (1.5,.28) circle (.1cm);
\draw[thick] (1.5,.28) circle (.1cm); 
\filldraw[white] (1.16,1.4) circle (.1cm);
\draw[thick] (1.16,1.4) circle (.1cm); 
\filldraw[white] (.6,1.8) circle (.1cm);
\draw[thick] (.6,1.8) circle (.1cm); 
}
=
\tikzmath[scale=.75, transform shape]{
\begin{scope}
\clip[rounded corners = 5] (-.6,0) rectangle (1.8,3);
\filldraw[primedregion=\BColor] (0,0) -- (0,1.2) .. controls ++(90:.4cm) and ++(-135:.2cm) .. (.6,1.8) .. controls ++(135:.2cm) and ++(270:.4cm) .. (0,2.4) -- (0,3) -- (-.6,3) -- (-.6,0); 
\filldraw[primedregion=\DColor] (1.2,0) -- (1.2,1.2) .. controls ++(90:.4cm) and ++(-45:.2cm) .. (.6,1.8) .. controls ++(-135:.2cm) and ++(90:.4cm) .. (0,1.2) -- (0,0);
\filldraw[boxregion=\BColor] (0,3) -- (0,2.4) .. controls ++(270:.4cm) and ++(135:.2cm) .. (.6,1.8) .. controls ++(45:.2cm) and ++(270:.4cm) .. (1.2,2.4) -- (1.2,3);
\filldraw[boxregion=\DColor] (1.2,0) -- (1.2,1.2) .. controls ++(90:.4cm) and ++(-45:.2cm) .. (.6,1.8) .. controls ++(45:.2cm) and ++(270:.4cm) .. (1.2,2.4) -- (1.2,3) -- (1.8,3) -- (1.8,0);
\end{scope}
\draw[\ZColor,thick] (0,0) -- (0,1.2) .. controls ++(90:.6cm) and ++(270:.6cm) .. (1.2,2.4) -- (1.2,3);
\draw[\phiColor,thick] (1.2,0) -- (1.2,.8);
\draw[\psiColor,thick] (1.2,.8) -- (1.2,1.2) .. controls ++(90:.6cm) and ++(270:.6cm) .. (0,2.4) -- (0,3);
\roundNbox{unshaded}{(1.2,.8)}{.3}{0}{0}{\scriptsize{$n_{e'}$}}; 
\filldraw[white] (.6,1.8) circle (.1cm);
\draw[thick] (.6,1.8) circle (.1cm); 
\node at (0,-.2) {\scriptsize{$A(Z)$}};
\node at (1.2,3.2) {\scriptsize{$B(Z)$}};
\node at (1.2,-.2) {\scriptsize{$\varphi_{e'}$}};
\node at (0,3.2) {\scriptsize{$\psi_{e}$}};
}\,.
\]
In the third equality above, we used the fact that $X\xxo_e Z\xxo_{e'}(X')^\bullet \in \cE(G(c)\to G(c'))$ to apply \eqref{eq:ModificationForAll1Cells}.
This completes the proof.
\end{proof}

\subsection{Proof of Theorem \ref{thm:UniqueLift}}

In this section, we prove Theorem \ref{thm:UniqueLift}.
We begin by recalling the construction of the canonical inclusion $\iota_\cC : \cC \hookrightarrow \QSys(\cC)$.

\begin{construction}[{\cite[Const.~3.24]{2105.12010}}]
\label{const:iota}
For each $\cA\in 2\Cat$, there is a canonical inclusion strict $\dag$ 2-functor $\iota_\cA: \cA \to \QSys(\cA)$ defined as follows:
\begin{itemize}
\item
For $a\in \cA$, $a\mapsto 1_a$, the trivial Q-system.
\item
For ${}_aX_b\in \cA(a\to b)$, $X$ is a separable $1_a-1_b$ bimodule, so $X$ maps to itself.
\item
For $f\in \cA(X\Rightarrow Y)$, $f$ is automatically $1_a-1_b$ bimodular, so $f$ maps to itself.
\end{itemize}
\end{construction}

\begin{construction}
\label{construction:LiftExists}
Suppose $F\in \Fun^\dag(\cA\to \cB)$.
We construct an invertible transformation $\psi^F:\iota_\cB \xz F \Rightarrow \QSys(F)\xz \iota_\cA$.

By Constructions \ref{construction:Qsys(F)} and \ref{const:iota},
for a 0-cell $b\in\cA$, 
we have
$$
(\iota_\cB\xz F)(b) = \iota_\cB(F(b)) = 1_{F(b)}
\qquad\text{and}\qquad 
(\QSys(F)\xz \iota_\cA)(b) = \QSys(F)(1_b)=F(1_b).
$$
For a 1-cell $X\in\QSys(\cA)(P\to Q)$, we have an equality
$$
(\iota_\cB\xz F)(X) = \iota_\cB(F(X)) = F(X) = \QSys(F)(X) = (\QSys(F)\xz \iota_\cA)(X),
$$ 
as well as for a 2-cell $f\in\QSys(\cA)(X\Rightarrow X')$:
$$
(\iota_\cB\xz F)(f) = \iota_\cB(F(f)) = F(f) = \QSys(F)(f) = (\QSys(F)\xz \iota_\cA)(f).
$$
Now $F(1_b)$ is equivalent to the trivial Q-system $1_{F(b)}$, and thus
for every $X\in\cA(a\to b)$,
$$
u^{F(1_b)}_{F(X),1_{F(b)}}: F(X)\otimes_{1_{F(b)}}1_{F(b)}\Rightarrow F(X)\otimes_{F(1_b)} 1_{F(b)}
$$ 
from \eqref{nota:QSys(C)andOther} is
unitary;
similarly, $u^{F(1_a)}_{1_{F(a)},F(X)}$ is a unitary.
We define:
\begin{itemize}
\item 
For 0-cell $a,b\in\cA$ and 1-cell $X\in\cA(a\to b)$, 
we define $\psi^F_b: = 1_{F(b)}$ as an $F(1_b)-1_{F(b)}$ bimodule, which is clearly invertible.
\item
For ${}_aX_b\in \cA$, we define
\[
\psi^F_X : = 
\tikzmath[scale=.7, transform shape]{
\begin{scope}
\clip[rounded corners = 5] (-.6,0) rectangle (1.8,2.4);
\filldraw[primedregion=\AColor] (0,0) -- (0,.6) .. controls ++(90:.4cm) and ++(-135:.2cm) .. (.6,1.2) .. controls ++(135:.2cm) and ++(270:.4cm) .. (0,1.8) -- (0,3) -- (-.6,3) -- (-.6,0); 
\filldraw[primedregion=\BColor] (1.2,0) -- (1.2,.6) .. controls ++(90:.4cm) and ++(-45:.2cm) .. (.6,1.2) .. controls ++(-135:.2cm) and ++(90:.4cm) .. (0,.6) -- (0,0);
\filldraw[\AColor] (0,3) -- (0,1.8) .. controls ++(270:.4cm) and ++(135:.2cm) .. (.6,1.2) .. controls ++(45:.2cm) and ++(270:.4cm) .. (1.2,1.8) -- (1.2,3);
\filldraw[\BColor] (1.2,0) -- (1.2,.6) .. controls ++(90:.4cm) and ++(-45:.2cm) .. (.6,1.2) .. controls ++(45:.2cm) and ++(270:.4cm) .. (1.2,1.8) -- (1.2,3) -- (1.8,3) -- (1.8,0);
\end{scope}
\draw[\XColor,thick] (0,0) -- (0,.6) .. controls ++(90:.6cm) and ++(270:.6cm) .. (1.2,1.8) -- (1.2,2.4);
\draw[\phiColor,thick,dotted] (1.2,0) -- (1.2,.6) .. controls ++(90:.6cm) and ++(270:.6cm) .. (0,1.8) -- (0,2.4);
\filldraw[white] (.6,1.2) circle (.1cm);
\draw[thick] (.6,1.2) circle (.1cm); 
\node at (0,-.2) {\scriptsize{$F(X)$}};
\node at (1.2,2.6) {\scriptsize{$F(X)$}};
\node at (1.2,-.2) {\scriptsize{$1_{F(b)}$}};
\node at (0,2.6) {\scriptsize{$1_{F(a)}$}};
} 
:= 
\tikzmath{
\begin{scope}
\clip[rounded corners = 5] (-.9,-.8) rectangle (.9,1);
\filldraw[primedregion=\AColor] (-.9,-.8) rectangle (0,1);
\filldraw[\AColor] (-.4,.6) rectangle (0,1);
\filldraw[\BColor] (0,-.8) rectangle (.9,1);
\filldraw[primedregion=\BColor] (.4,-.8) -- (.4,-.5) arc (0:90:.4cm) -- (0,-.8);
\end{scope}
\draw[\XColor,thick] (0,-.8) -- (0,1);
\draw[thick,dotted] (.4,-.8) -- (.4,-.5) arc (0:90:.4cm);
\draw[thick,dotted] (-.4,1) -- (-.4,.6) arc (180:270:.4cm);
\draw[thick] (-.4,.6) -- (0,.6);
\filldraw[\XColor] (0,-.1) circle (.05cm);
\filldraw[\XColor] (0,.2) circle (.05cm);
}
\qquad\qquad
\begin{aligned}
\tikzmath{
\filldraw[primedregion=\AColor, rounded corners = 5pt] (0,0) rectangle (.6,.6);
}
&=
F(1_a)
&
\tikzmath{
\filldraw[primedregion=\BColor, rounded corners = 5pt] (0,0) rectangle (.6,.6);
}
&=
F(1_b)
\\
\tikzmath{
\filldraw[\AColor, rounded corners = 5pt] (0,0) rectangle (.6,.6);
}
&=
1_{F(a)}
&
\tikzmath{
\filldraw[\BColor, rounded corners = 5pt] (0,0) rectangle (.6,.6);
}
&=
1_{F(b)}
\end{aligned}
\qquad
\tikzmath{
\begin{scope}
\clip[rounded corners = 5pt] (-.5,-.5) rectangle (.5,.5);
\filldraw[\AColor] (-.2,0) rectangle (.2,.5);
\filldraw[primedregion=\AColor] (-.2,-.5) rectangle (.2,0);
\end{scope}
\draw[thick, dotted] (-.2,-.5) -- (-.2,.5);
\draw[thick] (-.2,0) -- (.2,0);
\draw[thick, \XColor] (.2,-.5) -- (.2,.5);
}
=\left(u^{F(1_a)}_{1_{F(a)},F(X)}\right)^\dag
\]
Clearly $\psi^F_X$ is unitary. 
\end{itemize}
We leave the verification that $\psi^F$ is a 2-transformation to the reader.
\end{construction}

\begin{rem}
We expect that $\QSys$ is actually a symmetric lax monoidal 3-functor on the symmetric monoidal 3-category $2\Cat$.
Indeed, for each $\cA,\cB\in 2\Cat$, there is a canonical 2-functor
$\QSys(\cA) \boxtimes \QSys(\cB)\Rightarrow \QSys(\cA\boxtimes \cB)$ which satisfies various coherences.

At this time, we are unaware of a definition of a symmetric monoidal structure on an algebraic tricategory, as well as a definition of symmetric lax monoidal 3-functor on an algebraic tricategory.
We leave this exploration to the interested reader. 
A more tractable goal would be to produce a symmetric lax monoidal functor on the 1-category of 2-categories and equivalence classes of 2-functors, which is equivalent to the localization of the 1-category $\Gray$ of strict 2-categories and strict 2-functors with the $\Gray$ tensor product at the weak equivalences in the sense of a model category structure \cite{MR2138540}.
\end{rem}

Suppose now $\cC,\cD$ are $\rm C^*/W^*$ 2-categories with $\cD$ Q-system complete.
We apply the propositions from \S\ref{sec:Dominance} in the case that $\cE=\QSys(\cC)$ and $G=\iota_\cC$. 

\begin{lem}
\label{lem:IotaDominant}
$\iota_\cC$ is dominant. 
\end{lem}
\begin{proof}
For each 0-cell/Q-system ${}_bQ_b\in\QSys(\cC)$ where $b\in\cC$,
$Q:\iota_\cC(b)=1_b\condense Q$ is a dagger condensation when equipped with the 1-cells 
${}_bQ_Q 
=
\tikzmath{
\begin{scope}
\clip[rounded corners=5pt] (-.3,0) rectangle (.3,.6);
\fill[\BColor] (0,0) rectangle (-.3,.6);
\fill[\QrColor] (0,0) rectangle (.3,.6);
\end{scope}
\draw[thick, \QsColor] (0,0) -- (0,.6);
}$\,,
${}_QQ_b^\bullet:={}_QQ_b
=
\tikzmath{
\begin{scope}
\clip[rounded corners=5pt] (-.3,0) rectangle (.3,.6);
\fill[\QrColor] (0,0) rectangle (-.3,.6);
\fill[\BColor] (0,0) rectangle (.3,.6);
\end{scope}
\draw[thick, \QsColor] (0,0) -- (0,.6);
}
$\,,
and the 2-cells
\[
\tikzmath{
\fill[\QrColor, rounded corners=5pt] (-.3,0) rectangle (.9,.6);
\filldraw[\BColor] (0,0) arc (180:0:.3cm);
\draw[\QsColor,thick] (0,0) arc (180:0:.3cm);
\draw[\QsColor,thick] (.3,.3) -- (.3,.6);
\filldraw[\QsColor] (.3,.3) circle (.05cm);
}=\varepsilon_Q
\qquad\qquad
\tikzmath{
\fill[\QrColor, rounded corners=5pt] (-.3,0) rectangle (.9,-.6);
\filldraw[\BColor] (0,0) arc (-180:0:.3cm);
\draw[\QsColor,thick] (0,0) arc (-180:0:.3cm);
\draw[\QsColor,thick] (.3,-.3) -- (.3,-.6);
\filldraw[\QsColor] (.3,-.3) circle (.05cm);
}=\delta_Q=\varepsilon_Q^\dag
\qquad\quad\underset{\text{\ref{Q:separable}}}{\Longrightarrow}
\qquad\quad
\varepsilon_Q\varepsilon_Q^\dag = 
\tikzmath{
\fill[\QrColor, rounded corners=5pt] (-.3,0) rectangle (.9,1.2);
\filldraw[\BColor] (0,.6) arc (180:-180:.3cm);
\draw[\QsColor,thick] (0,.6) arc (180:-180:.3cm);
\draw[\QsColor,thick] (.3,1.2) -- (.3,.9);
\draw[\QsColor,thick] (.3,0) -- (.3,.3);
\filldraw[\QsColor] (.3,.3) circle (.05cm);
\filldraw[\QsColor] (.3,.9) circle (.05cm);
}
=
\tikzmath{
\fill[\QrColor, rounded corners=5pt ] (0,0) rectangle (.6,1.2);
\draw[\QsColor,thick] (.3,0) -- (.3,1.2);
}=\id_{{}_QQ_Q}\,.
\] 
The result now follows as $\iota_\cC$ is a local equivalence on hom categories by definition.
\end{proof}

\begin{prop}
\label{prop:PrecompositionEquivalenceOnHomCats}
$-\xz\iota_\cC: \Fun^\dag(\QSys(\cC)\to \cD)\to \Fun^\dag(\cC\to \cD)$ is a dagger equivalence on hom categories.
\end{prop}
\begin{proof}
By Lemma \ref{lem:IotaDominant}, $\iota_\cC$ is dominant, so by Proposition \ref{prop:DominantToFullyFaithful},
$-\xz \iota_\cC$ is fully faithful on 2-morphism. 
To prove $-\xz \iota_\cC$ is a dagger equivalence on hom categories, it remains to prove $-\xz \iota_\cC$ is unitarily essentially surjective on 1-morphisms, i.e., 
for all $A,B\in \Fun^\dag(\QSys(\cC) \to \cD)$
and
each 1-morphism $\gamma:A\xz \iota_\cC\Rightarrow B\xz \iota_\cC$, there exists $\varphi:A\Rightarrow B$ such that $\gamma\cong\varphi\xz \iota_\cC$. 

For 0-cells/Q-systems ${}_aP_a,{}_bQ_b\in \QSys(\cC)$ and a 1-cell ${}_PX_Q\in\QSys(\cC)(P\to Q)$, we define $\varphi_Q\in \cD(A(Q)\to B(Q))$ and $\varphi_X\in \cD({}_{A(P)}{A(X)\xxo_{A(Q)}\varphi_Q}_{B(Q)}\Rightarrow {}_{A(P)}{\varphi_P\xxo_{B(P)}B(X)}_{B(Q)})$ by
\[
\varphi_Q :=
\tikzmath[scale=.7, transform shape]{
\begin{scope}
\clip[rounded corners = 5] (-2.4,-1.2) rectangle (2.4,1.2);
\filldraw[primedregion=\BColor] (.6,-1.2) .. controls ++(90:.4cm) and ++(-45:.2cm) .. (0,0) .. controls ++(135:.2cm) and ++(270:.4cm) .. (-.6,1.2) -- (-2.4,1.2) -- (-2.4,-1.2);
\filldraw[boxregion=\BColor] (.6,-1.2) .. controls ++(90:.4cm) and ++(-45:.2cm) .. (0,0) .. controls ++(135:.2cm) and ++(270:.4cm) .. (-.6,1.2) -- (2.4,1.2) -- (2.4,-1.2);
\filldraw[primedregion=\QrColor] (-.6,-1.2) -- (-.6,-.6) -- (-1.2,0) .. controls ++(135:.2cm) and ++(270:.4cm) .. (-1.8,1.2) -- (-2.4,1.2) -- (-2.4,-1.2);
\filldraw[boxregion=\QrColor] (1.8,-1.2) .. controls ++(90:.4cm) and ++(-45:.2cm) .. (1.2,0) -- (.6,.6) -- (.6,1.2) -- (2.4,1.2) -- (2.4,-1.2);
\filldraw[\BColor] (-.6,-.2) rectangle (-.2,.2);
\end{scope}
\draw[\phiColor,thick] (.6,-1.2) .. controls ++(90:.4cm) and ++(-45:.2cm) .. (0,0) .. controls ++(135:.2cm) and ++(270:.4cm) .. (-.6,1.2);
\draw[\QsColor,thick] (-.6,-1.2) -- (-.6,-.6) -- (-1.2,0) .. controls ++(135:.2cm) and ++(270:.4cm) .. (-1.8,1.2);
\draw[\QsColor,thick] (1.8,-1.2) .. controls ++(90:.4cm) and ++(-45:.2cm) .. (1.2,0) -- (.6,.6) -- (.6,1.2);
\draw[\QsColor,thick] (-.6,-.6) -- (.6,.6);
\filldraw[\QsColor] (.6,.6) circle (.05cm);
\filldraw[\QsColor] (-.6,-.6) circle (.05cm);
\filldraw[white] (0,0) circle (.1cm);
\draw[thick] (0,0) circle (.1cm); 
\node at (-.6,-1.4) {\scriptsize{$A(Q)$}};
\node at (.6,-1.4) {\scriptsize{$\gamma_b$}};
\node at (1.8,-1.4) {\scriptsize{$B(Q)$}};
\node at (-.4,0) {\scriptsize{$\gamma_Q$}};
}
\qquad\qquad
\varphi_X :=
\tikzmath[scale=.5, transform shape]{
\begin{scope}
\clip[rounded corners = 5] (-3.6,-2.4) rectangle (3.6,2.4);
\filldraw[primedregion=\AColor] (0,0) -- (-1.2,1.2) .. controls ++(135:.2cm) and ++(270:.4cm) .. (-1.8,2.4) -- (-3,2.4) -- (-3,-.6) -- (-.6,-.6);
\filldraw[boxregion=\AColor] (0,0) -- (-1.2,1.2) .. controls ++(135:.2cm) and ++(270:.4cm) .. (-1.8,2.4) -- (-.6,2.4) -- (-.6,1.8) -- (.6,.6);
\filldraw[primedregion=\BColor] (1.8,-2.4) .. controls ++(90:.4cm) and ++(-45:.2cm) .. (1.2,-1.2) -- (0,0) -- (-.6,-.6) -- (.6,-1.8) -- (.6,-2.4);
\filldraw[boxregion=\BColor] (1.8,-2.4) .. controls ++(90:.4cm) and ++(-45:.2cm) .. (1.2,-1.2) -- (0,0) -- (.6,.6) -- (3,.6) -- (3,-2.4);
\filldraw[primedregion=\PrColor] (-.6,-2.4) -- (-.6,-.6) -- (-2.4,1.2) .. controls ++(135:.2cm) and ++(270:.4cm) .. (-3,2.4) -- (-3.6,2.4) -- (-3.6,-2.4);
\filldraw[boxregion=\PrColor] (-.6,2.4) -- (-.6,1.8) -- (.6,.6) -- (.6,2.4); 
\filldraw[primedregion=\QrColor] (.6,-2.4) -- (.6,-1.8) -- (-.6,-.6) -- (-.6,-2.4); 
\filldraw[boxregion=\QrColor] (3,-2.4) .. controls ++(90:.4cm) and ++(-45:.2cm) .. (2.4,-1.2) -- (.6,.6) -- (.6,2.4) -- (3.6,2.4) -- (3.6,-2.4);
\filldraw[\AColor] (-1.6,1.2) circle (.22cm);
\filldraw[\AColor] (-.4,0) circle (.22cm);
\filldraw[\BColor] (.8,-1.2) circle (.22cm);
\end{scope}
\draw[\phiColor,thick] (1.8,-2.4) .. controls ++(90:.4cm) and ++(-45:.2cm) .. (1.2,-1.2) -- (-1.2,1.2) .. controls ++(135:.2cm) and ++(270:.4cm) .. (-1.8,2.4);
\draw[\XColor,thick] (-.6,-2.4) -- (-.6,-.6) -- (.6,.6) -- (.6,2.4);
\draw[\PsColor,thick] (-.6,-.6) -- (-2.4,1.2) .. controls ++(135:.2cm) and ++(270:.4cm) .. (-3,2.4);
\draw[\PsColor,thick] (.6,.6) -- (-.6,1.8) -- (-.6,2.4);
\draw[\PsColor,thick] (-1.8,.6) -- (-.6,1.8);
\draw[\QsColor,thick] (-.6,-.6) -- (.6,-1.8) -- (.6,-2.4);
\draw[\QsColor,thick] (3,-2.4) .. controls ++(90:.4cm) and ++(-45:.2cm) .. (2.4,-1.2) -- (.6,.6);
\draw[\QsColor,thick] (1.8,-.6) -- (.6,-1.8);
\filldraw[\XColor] (-.6,-.6) circle (.05cm);
\filldraw[\XColor] (.6,.6) circle (.05cm);
\filldraw[\PsColor] (-1.8,.6) circle (.05cm);
\filldraw[\PsColor] (-.6,1.8) circle (.05cm);
\filldraw[\QsColor] (1.8,-.6) circle (.05cm);
\filldraw[\QsColor] (.6,-1.8) circle (.05cm);
\filldraw[white] (0,0) circle (.1cm);
\draw[thick] (0,0) circle (.1cm); 
\filldraw[white] (-1.2,1.2) circle (.1cm);
\draw[thick] (-1.2,1.2) circle (.1cm); 
\filldraw[white] (1.2,-1.2) circle (.1cm);
\draw[thick] (1.2,-1.2) circle (.1cm); 
\node at (-.6,-2.6) {\normalsize{$A(X)$}};
\node at (.6,-2.6) {\normalsize{$A(Q)$}};
\node at (1.8,-2.6) {\normalsize{$\gamma_b$}};
\node at (3,-2.6) {\normalsize{$B(Q)$}};
\node at (.6,2.6) {\normalsize{$B(X)$}};
\node at (-.6,2.6) {\normalsize{$B(P)$}};
\node at (-1.8,2.6) {\normalsize{$\gamma_a$}};
\node at (-3,2.6) {\normalsize{$A(P)$}};
\node at (-.4,0) {\normalsize{$\gamma_X$}};
\node at (-1.6,1.2) {\normalsize{$\gamma_P$}};
\node at (.8,-1.2) {\normalsize{$\gamma_Q$}};
}
\qquad\qquad
\begin{aligned}
\tikzmath{
\filldraw[\AColor, rounded corners = 5pt] (0,0) rectangle (.6,.6);
}
&=
1_a
&
\tikzmath{
\filldraw[\BColor, rounded corners = 5pt] (0,0) rectangle (.6,.6);
}
&=
1_b
\\
\tikzmath{
\filldraw[\PrColor, rounded corners = 5pt] (0,0) rectangle (.6,.6);
}
&=
P
&
\tikzmath{
\filldraw[\QrColor, rounded corners = 5pt] (0,0) rectangle (.6,.6);
}
&=
Q
\\
\tikzmath{
\filldraw[primedregion=white, rounded corners = 5pt] (0,0) rectangle (.6,.6);
\draw[thin, dotted, rounded corners = 5pt] (0,0) rectangle (.6,.6);
}
&=
A
&
\tikzmath{
\filldraw[boxregion=white, rounded corners = 5pt] (0,0) rectangle (.6,.6);
\draw[thin, dotted, rounded corners = 5pt] (0,0) rectangle (.6,.6);
}
&=
B
\end{aligned}
\]
Then for each 0-cell $b\in \cC$, and 1-cell ${}_aX_b\in\cC(a\to b)$,
by Construction \ref{const:1CompositionIn2Cat},
$$
(\varphi\xz \iota_\cC)_b 
= \varphi_{\iota_\cC(b)} 
= \varphi_{1_b}
= 
\gamma_b
\qquad
\text{and}
\qquad
(\varphi\xz \iota_\cC)_X 
= 
\varphi_{\iota_\cC(X)} 
= 
\varphi_{X} 
= 
\gamma_X
$$
where the latter is viewed as $1_a-1_b$ bimodular.
Therefore $\varphi\xz \iota_\cC\cong \gamma$ as desired, so $-\circ \iota_\cC$ is gives a dagger functor on hom 1-categories whose underlying functor is an equivalence.
Since $\Fun^\dag(\cC\to \cD)$ is $\rm C^*$, $-\circ \iota_\cC$ is a dagger equivalence on hom 1-categories by polar decomposition as discussed in Remark \ref{rem:Underlying2FunctorEquivalence}.
\end{proof}

By Propositions \ref{prop:TruncatedGroupoids} and \ref{prop:PrecompositionEquivalenceOnHomCats}, $-\xz \iota_\cC$ is $(-1)$-truncated when restricted to unitary cores, i.e., the homotopy fiber at each $F\in \core^\dag(\Fun^\dag(\cC\to \cD))$ is either empty or equivalent to a point.
By Constructions \ref{construction:Qsys(F)} and \ref{construction:LiftExists}, the homotopy fiber of $-\xz \iota_\cC$ is non-empty at each $F\in \Fun^\dag(\cC\to \cD)$.
Indeed, since $\cD$ is Q-system complete, $\iota_\cD$ is invertible, so 
there exists a $\dag$ 2-functor $\iota^{-1}_\cD:\QSys(\cD)\to\cD$ together with an invertible $\dag$ 2-transformation $\theta_\cD:1_\cD\Rightarrow \iota_\cD^{-1}\xz\iota_\cD$. 
Thus $\iota_\cD^{-1}\circ \QSys(F)$ provides the desired lift together with the composite invertible transformation
\begin{equation}
\label{eq:LiftExists}
\begin{tikzcd}
\QSys(\cC)
\arrow[rr, "\QSys(F)"]
&&
\QSys(\cD)
\arrow[dr, "\iota_\cD^{-1}"]
\\
\cC
\arrow[r,swap, "F"]
\arrow[u, "\iota_\cC"]
&
\cD
\arrow[rr,swap,"1_\cD"]
\arrow[ur, "\iota_\cD"]
\arrow[ul,Rightarrow,shorten <= 1em, shorten >= 1em,"\cong","\psi^F"']
&
\arrow[u,Rightarrow,shorten <= .5em, 
shorten >= .5em,
"\cong","\theta"']
&
\cD.
\end{tikzcd}
\end{equation}
Thus the homotopy fiber of $-\xz \iota_\cC$ at each $F\in \core^\dag(\Fun^\dag(\cC\to \cD))$ is equivalent to a point.
By Proposition \ref{prop:TruncatedGroupoids}, $-\xz\iota_\cC$ is (-2)-truncated when restricted to unitary cores.
This implies $-\xz\iota_\cC : \Fun^\dag(\cE\to \cD) \to \Fun^\dag(\cC\to \cD)$ is essentially surjective on objects.
Again by 
Remark \ref{rem:Underlying2FunctorEquivalence},
Proposition \ref{prop:PrecompositionEquivalenceOnHomCats}, and \cite[Thm.~7.4.1]{2002.06055}, $-\xz\iota_\cC : \Fun^\dag(\cE\to \cD) \to \Fun^\dag(\cC\to \cD)$ is a $\dag$-equivalence of $\rm C^*/W^*$ 2-categories.
\qed

\begin{rem}
Observe that we did not really need to pass to (unitary) cores, nor use Proposition \ref{prop:TruncatedGroupoids}.
Indeed, $-\circ \iota_\cC$ is an equivalence on hom categories by Proposition \ref{prop:PrecompositionEquivalenceOnHomCats} and essentially surjective on objects by \eqref{eq:LiftExists}, and thus an equivalence by \cite[Thm.~7.4.1]{2002.06055} and Remark \ref{rem:Underlying2FunctorEquivalence}.
\end{rem}

\bibliographystyle{alpha}
{\footnotesize{
\bibliography{bibliography}
}}
\end{document}